\documentclass[12pt]{article}
\usepackage[T2A]{fontenc}
\usepackage[utf8]{inputenc}
\usepackage{amssymb,amsmath,amsfonts,amsthm,amscd,latexsym,verbatim,graphics,epsfig,indentfirst,graphicx,caption}
\usepackage[bookmarks=true,pdfborder=0 0 0,colorlinks=true,linkcolor=blue,dvips=true,final=true]{hyperref}
\usepackage{geometry} 
\def\evv{\mathrm{ev}\,}

\def\tr{\mathrm{tr}\,}
\def\thr{\mathrm{thr}\,}
\def\Aut{\mathrm{Aut}}
\def\Ass{\mathrm{As}}

\def\Lie{\mathrm{Lie}\,}
\def\Supp{\mathrm{Supp}\,}

\begin{document}

\begin{center}
{\Large $PC$-polynomial of graph}

Vsevolod Gubarev
\end{center}

\begin{abstract}
We define $PC$-polynomial of graph
which is related to clique, (in)dependence and matching polynomials. 
The growth rate of partially commutative monoid is equal to
the largest root $\beta(G)$ of $PC$-polynomial of the corresponding graph.

The random algebra is defined in such way that its growth rate 
equals the largest root of $PC$-polynomial of random graph. 
We prove that for almost all graphs 
all sufficiently large real roots of $PC$-polynomial 
lie in neighbourhoods of roots of $PC$-polynomial of random graph.
We show how to calculate the series expansions of the latter roots.
The average value of $\beta(G)$ over all graphs 
with the same number of vertices is computed.

We found the graphs on which the maximal value of $\beta(G)$
with fixed numbers of vertices and edges is reached.
From this, we derive the upper bound of $\beta(G)$.
Modulo one assumption, we do the same for minimal value of $\beta(G)$.

We study the Nordhaus---Gaddum bounds of
$\beta(G)+\beta(\bar{G})$ and $\beta(G)\beta(\bar{G})$.

\end{abstract}

\newpage
\tableofcontents

\eject

\newcounter{Pic}
\setcounter{section}{-1}

\newpage
\section*{Introduction}
\addcontentsline{toc}{section}{Introduction}

In 1969~\cite{Cartier69}, P. Cartier and D. Foata defined a partially commutative monoid 
by some combinatorial reasons (reproving MacMahon Master theorem).
Given a simple finite graph $G(V,E)$, a partially commutative monoid
$M(G)$ is in algebraic language 
$M\langle V\mid v_1v_2 = v_2v_1,\,(v_1,v_2)\in E\rangle$, i.e., 
a quotient of the free monoid 
$M\langle V\rangle$ by the congruence relation generated by 
$\{v_1v_2 = v_2v_1,\,(v_1,v_2)\in E\}$.
Informally, $M(G)$ is a set of words in alphabet $V$ with
the operation of concatenation. Moreover, given a word,
we can interchange neighbor letters if they are connected in $G$.
Two words $u,v$ are equal if $u$ could be obtained from $v$
by a finite number of such interchanges of neighbor connected letters.

Further, partially commutative groups, algebras and Lie algebras appeared. 
Partially commutative structures are very natural object for study.
On the one hand, they are defined very simply, their properties
are close to the properties of free ones. On the other hand,
p.c. structures are enough reach to formulate a lot of problems
with nontrivial solutions. Such objects have been investigated 
in combinatorics, formal languages, auto\-mata, computer 
science~\cite{Cartier69,Diekert,Diekert2,Diekert3,Mazurkiewicz,TracesNew};
in algebra~\cite{Baudisch,DuchampKrob1992,DuchampKrob1993,Duncan,Gub2016,Kim,Poroshenko};
in topology~\cite{Agol,Charney,KimKoberda,Margolis}; 
in logics~\cite{Kazachkov,Gupta,Poroshenko2}; 
in robotics~\cite{Charney}.

We should clarify that elements of partially commutative monoids sometimes are called {\it traces},
a partially commutative group is called {\it right-angled Artin group} (RAAG).
There is a big area of mathematics devoted to more general groups: Artin and Coxeter groups.
In the first ones, we have the relations 
$[x_i,x_j]^{m_{ij}} = 1$ for all pairs of connected generators $x_i,x_j$ and some $m_{ij}\geq 2$.
A Coxeter group is an Artin group with additional relations $x_i^2 = 1$,
so, it is a generalization of Weil groups which play the significant role in study
of simple Lie groups and simple Lie algebras.

The natural question about the numbers of different words of given length
and, more specifically, about the growth rate of partially commutative monoid $M(G)$ 
leads to graph polynomials\footnote{About growth rates of partially 
commutative groups see~\cite{GeodesicGrowth,Athreya,Baik}.}.
For associative/Lie algebras we should ask about the dimensions respectively.

Dependence polynomial actually appeared in the work of P. Cartier and D. Foata~\cite{Cartier69}.
This polynomial was wrtitten down by D.C. Fisher and A.E. Solow in 1990~\cite{Fisher1990}.
Depen\-dence polynomial is defined as $D(G,x) = 1+\sum\limits_{k=1}^{\omega(G)}(-1)^k c_k(G)x^k$,
where $c_i(G)$ denotes a number of distinct cliques in $G$ of the size $i$
and $\omega(G)$ is the {\it clique number} of graph, the size of a maximal clique in~$G$.
The reciprocal of dependence polynomial is the generating function of partially commutative monoid $M(G)$. 

By analogy, in 1994 C. Hoede and X. Li defined~\cite{Clique1994} 
the clique polynomial as 
$C(G,x) = 1+\sum\limits_{k=1}^{\omega(G)}c_k(G)x^k$.

Let $s_i(G)$ denote the number of anticliques of size~$i$ in a graph $G$, i.e., 
$s_i(G) = c_i(\bar{G})$, where $\bar{G}$ is a complement graph.
The polynomial
$I(G,x) = 1+\sum\limits_{k=1}^{\alpha(G)}s_k(G)x^k$
is called the independence polynomial of a graph $G$.
Here $\alpha(G)$ is the size of the maximal independent set in $G$.
Independence polynomial firstly appeared in the work of 
A.~Motoyama and H.~Hosoya~\cite{Motoyama} in 1977 for lattice graphs. 
A systematic study of the subject began in 1983, with the work of 
I. Gutman and F. Harary~\cite{Gutman1983}.
In 2005, V. Levit and E.~Mandrescu wrote the interesting survey 
about independence polynomial~\cite{Levit2005}.

In 1971, H. Hosoya defined~\cite{Hosoya} matching polynomial,
actually it is independence polynomial of the line graph $L(G)$ of $G$
(see also the work of H. Kunz~\cite{Kunz}).
Let us refer to the important works~\cite{Farrell,Godsil,Heilman72}
and also monographs~\cite{Godsil93,LovaszPlummer} devoted to matching polynomial.
Note that independence and matching polynomials have connections 
with chemistry and physics~\cite{Balasubramanian,Heilman72,Motoyama,Scott}.

In the work, we define $PC$-polynomial 
(short for partially commutative polynomial) 
of a graph~$G$ as 
$PC(G,x) = \sum\limits_{k=0}^{\omega(G)}(-1)^k c_k(G)x^{\omega(G)-k}$, where $c_0(G):=1$. 
Actually $PC$-poly\-nomial have already appeared in the works~\cite{BencsAdjont,BrownNowak,Fisher1990} 
without any special name. The growth rate of partially commutative monoid/associative algebra/Lie algebra 
equals the largest real root $\beta(G)$ of $PC(G,x)$~\cite{Csikvari,DuchampKrob1992,Fisher1990,Santini2000}.
That is why $\beta(G)$ in \cite{Fisher1990} was called as the {\it growth factor}. 

Let us call (in)dependence polynomials, clique polynomial 
and $PC$-polynomial as clique-type polynomials.
They are connected in the following way:
\begin{gather*}
I(\bar{G},x) = C(G,x),\\
C(G,-x) = D(G,x) = x^{t_0}PC(G,1/x).
\end{gather*}

Clique-type polynomials have different applications in counting of 
(anti)cliques, mat\-chings, perfect matchings, homo\-morphisms 
and colorings in graphs with different 
constra\-ints~\cite{Fisher1989,Galvin,Galvin2,Zhao2018,Zhao2,Zhao2017}, 
in the Ramsey theory and sphere packings~\cite{Davies2017,Davies,Perkins}.
See also~\cite{Barvinok,ApplIndep,Widom-Rowlinson,CsikvariInfin}.

In 2005, A. Scott and A. Sokal showed~\cite{Scott} the deep connection 
between $\beta(G)$ and Lov\'{a}sz local lemma.

The main {\bf goals} of the current paper are the following:

0) give a survey on clique-type polynomials,

1) present new results about $PC$-polynomials of random and planar graphs,

2) state asymptotically tight lower and upper bounds on $\beta(G)$
in terms of $n = |V(G)|$ and $k = |E(G)|$,

3) find graphs on which $\beta(G)$ reaches the minimum or maximum among 
all graphs with fixed values of $n$ and $k$,

4) study Nordhaus---Gadddum~\cite{NordhausGadddumSurvey,NordhausGaddum}
inequalities for the expressions $\beta(G)+\beta(\bar{G})$ and $\beta(G)\beta(\bar{G})$,

5) find the average value of $\beta(G)$ among graphs of the same size.

Let us give an {\bf exposition} of the work.

In $\S0$, required preliminaries on polynomials and sequences are written down. 

The main goal of $\S1$ is to prove 

{\bf Theorem 1.1} \cite{Cartier69,Fisher1989-1}.
The numbers $m_n$, $n\geq1$, satisfy the reccurence relation
$$
m_n = c_1(G)m_{n-1}-c_2(G)m_{n-2} + \ldots + (-1)^{t_0+1}c_{t_0}(G)m_{n-t_0}
$$
with initial data $m_0 = 1$, $m_{-1} = \ldots = m_{-t_0+1} = 0$.

We state this result in different from the proofs of P. Cartier with D. Foata and D.C.~Fisher way.
As a corollary, we get

{\bf Lemma 1.4} \cite{Cartier69,Fisher1990}.
The function $1/D(G,x)$ is a generating function for the sequence $\{m_n\}$, i.e.,
$$
\frac{1}{D(G,x)} = \sum\limits_{n\geq0}m_n x^n.
$$

At the end of the section, we state the very useful equalities of clique-type polynomials
related to deletion of a vertex or an edge, to graphs join and union and to derivative. 

{\bf Lemma 1.5} \cite{Gutman1990,Clique1994}.
Given a graph $G$, we have 

a) $D(G,x) = D(G \setminus v,x) - xD(G[N(v)],x)$, $v\in V(G)$;

b) $D(G,x) = D(G\setminus uv,x) + x^2D(G[N(u)\cap N(v)],x)$, $(u,v)\in E(G)$;

c) $D'(G,x) = -\sum\limits_{v\in V(G)}D(G[N(v)],x)$;

d) $D(G_1\cup G_2) = D(G_1) + D(G_2) - 1$;

e) $D(G_1 + G_2) = D(G_1)D(G_2)$.

Here by $G[A]$ for $A\subset V(G)$ we denote the subgraph of $G$ induced 
by the set of vertices~$A$. By $N(v)$ we mean an open neighbourhood of $v\in V(G)$.

Exposition of~$\S2$ in general follows the work of P. Csikv\'{a}ri of 2013~\cite{Csikvari}.
Let $z_0$ denote a root of $PC(G,x)$ with the largest modulus and $\beta(G) = |z_0|$. 
D.C. Fisher and A.E. Solow in 1990 stated 

{\bf Lemma 2.1} \cite{Fisher1990}.
The number $\beta(G)$ is a root of $PC(G,x)$.

D.C. Fisher and A.E. Solow also pretended that they had proved the following theorem 

{\bf Theorem 2.1} \cite{Csikvari,Santini2000}.
The number $\beta(G)$ is the only complex root of $PC(G,x)$ 
with modulus greater or equal to $\beta(G)$.

In 2000, M. Goldwurm and M. Santini finally proved Theorem~2.1~\cite{Santini2000} 
via the Perron---Frobenius theory. 
In $\S2$, we consider the excellent proof of P. Csikv\'{a}ri~\cite{Csikvari}.

Now we are ready to connect the growth rate of p.c. structures 
with $PC$-polynomial:

{\bf Corollary 2.1} \cite{Santini2000}.
a) The growth rate of partially commutative monoid $M(X,G)$ equals~$\beta(G)$.

b) The growth rate of partially commutative associative algebra $\Ass(X,G)$ equals~$\beta(G)$\!.

In Lemmas~2.3 and 2.4 we prove the inequalities
$\beta(H)\leq \beta(G)\leq \beta(F)$~\cite{Csikvari,Hajiabolhassan}
for any induced subgraph $H$ of $G$ and any spanning subgraph~$F$ of $G$.

In~$\S3$, some applications of clique-type polynomials in graph theory
are considered. 

In 2009, D. Galvin~\cite{Galvin} proved the result based on the work of V.~Alekseev~\cite{Alekseev}.

{\bf Theorem 3.1} \cite{Alekseev,Davies,Galvin}.
Given a graph $G$ with $w = \omega(G)$, $n = |V(G)|$,
we have $C(G,x)\leq \big(1+\frac{nx}{w}\big)^w$ for all $x>0$
with equality if and only if $G$ is a complete multipartite graph with equal parts.

Theorem 3.1 allows to prove the following corollaries:

{\bf Corollary 3.1} \cite{Alekseev,Erdos}.
Let $G$ be a graph with $w = \omega(G)$, $|V(G)| = n$
and $c(G)$ denotes the number of all cliques in $G$. Then 
$c(G)\leq \big(1+\frac{n}{w}\big)^w$.
We have equality if and only if $G$ is a complete multipartite 
graph with equal parts.

{\bf Corollary 3.3 (Tur\'{a}n's Theorem)}.
Given a graph $G$ with $w = \omega(G)$, $|V(G)| = n$, 
we have $k = |E(G)|\leq \frac{n^2}{2}\big(1-\frac{1}{w}\big)$.

We write down the proof of Y. Zhao of the next result 
(modulo the case of bipartite graphs stated by J. Kahn in 2001~\cite{Kahn}).

{\bf Theorem 3.2} \cite{Zhao2}.
Given a $d$-regular graph $G$ with $n = |V(G)|$, for all $x\geq 0$ we have 
$I(G,x) \leq (2(1+x)^d-1)^{n/(2d)}$.

It implies the solution of the conjecture of N. Alon \cite{Alon}:

{\bf Corollary 3.4} \cite{Zhao2}.
For any $n$-vertex $d$-regular graph $G$, 
$i(G)\leq (2^{d+1} - 1)^{\frac{n}{2d}}$.

Let us introduce clique-type polynomials from the point of 
view of statistical physics, see, e.g.,~\cite{Perkins}.
By the {\it occupancy fraction} we mean the expected fraction of vertices 
that appear in the random independent set
$$
\alpha(G,x) 
 = \frac{E(|I|)}{|V(G)|}
 = \frac{1}{|V(G)|}\sum\limits_{I\in\mathcal{I}(G)}|I|\cdot \mathrm{Pr}[I]
 = \frac{xI'(G,x)}{|V(G)|I(G,x)}.
$$
Here $\mathrm{Pr}[I] = \dfrac{x^{|I|}}{\sum\limits_{J\in\mathcal{I}(G)}x^{|J|}}$
is so called {\it hard-core distribution} which is simply the uniform
distribution over all independent sets of $G$ at fugacity~$x$.
The independence polynomial is interpreted as the {\it partition function} 
of the hard-core model on $G$ at fugacity $x$.

{\bf Statement 3.7}~\cite{Davies}.
a) For any graph $G$, $\alpha (G,x)$ is monotone increasing in~$x$.

b) Let $G$ be a triangle-free graph on $n$ vertices with maximum degree $d$,
we have $\alpha(G,x)\geq (1+o(1))\frac{\ln d}{d}$ for any $x\geq 1/\ln d$.

E. Davies et al~\cite{Davies} applied Statement 3.7 to reprove 
the best known upper bound on the Ramsey numbers $R(3,k)$.

{\bf Corollary 3.8}~\cite{Davies,Shearer2}.
For the Ramsey numbers $R(3,k)$, we have the upper bound
$R(3,k)\leq (1+o(1))\frac{k^2}{\ln k}$.

In 2002, V. Nikiforov proved~\cite{NikiforovIneq} for 
the spectral radius $\rho = \rho(G)$ the inequality 
$$
\rho^w \leq c_2(G)\rho^{w-2}+2c_3(G)\rho^{w-2}+\ldots+(i-1)c_i(G)\rho^{w-i}+\ldots+(w-1)c_w(G),
$$
where $w = \omega(G)$. We show that this inequality is equivalent to 

{\bf Statement 3.8}.
Let $G$ be not empty graph with $n$ vertices and the spectral radius~$\rho$. 
Then $\alpha(\bar{G},x)\geq \frac{1}{n}$ for any $x\geq 1/\rho$.

The section~$\S4$ is devoted to partially commutative Lie algebras. 
By the definition, a partially commutative Lie algebra $L(X,G)$ equals
$\Lie\langle X|[a,b] = 0,(a,b)\in E(G)\rangle$.
Denote the dimension of the homogeneous space of all products
of length~$n$ in the alphabet~$X$ in $L(X,G)$ as $l_n$.

Given roots $x_1,\ldots,x_{\omega(G)}$ of $PC(G,x)$, define the numbers 
$p_n = x_1^n+x_2^n+\ldots +x_{\omega(G)}^n$.

In 1992, G. Duchamp and D. Krob actually proved~\cite{DuchampKrob1992} 
the following result (but not in the most comfortable form).

{\bf Theorem 4.1} \cite{DuchampKrob1992}.
We have
$$
l_n = \frac{1}{n}\sum\limits_{d|n}\mu(d)p_{n/d},
$$
where $\mu$ is the M\"{o}bius function.

{\bf Corollary 4.1}.
If $G = K_n$, then the growth rate of partially commutative Lie algebra $L(X,G)$ equals~0.
Otherwise, it equals~$\beta(G)$.

\eject

Also, we show that par\-tially commutative Lie algebras are connected
with Yang---Mills algebras~\cite{YM2} (Corollary~4.3) and the Lie algebras 
$L^{(l)}$ of primitive elements in the connected cocommutative Hopf algebra 
$\mathrm{Sym}^{(l)}$~\cite{Thibon} (Remark~4.4).

In~$\S5$, we study clique and $PC$-polynomial of the random graph
$G_{n,p}$ with $n$ vertices and edge probability $p$, 
$$
PC(G_{n,p},x) = x^n {-} \binom{n}{1}x^{n-1} {+} \binom{n}{2} px^{n-2} 
 {-} \ldots {+} (-1)^k \binom{n}{k} p^{\frac{k(k-1)}{2}}x^k
 {+} \ldots {+} (-1)^n p^{\frac{n(n-1)}{2}},
$$
and $C(G_{n,p},x) = (-x)^n PC(G,-1/x)$.
The main result of $\S5.1$ is 

{\bf Theorem 5.1} \cite{Brown2}.
Let $p\in(0,1)$. 

a) All roots of $C(G_{n,p},x)$ are real and simple.

b) Write roots of $C(G_{n,p},x)$ in ascending order
$r_{n}<\ldots<r_1<0$. Then $pr_{i+1}<r_{i}$ for all $i=1,\ldots,n-1$.

Theorem 5.1 was stated by J. Brown and R. Nowakowski in 2005~\cite{BrownNowak} for $p=1/2$
and by J. Brown et al in 2012~\cite{Brown2} for any $p$.
We show that Theorem 5.1a immediately follows from the result of E. Laguerre~\cite{Laguerre}
(Remark 5.1).

Let $q = \sqrt{p}$, $y = {q}^{n-1}x$. 
Define the polynomial
$$
\widetilde{C}(G_{n,q},y)
 = C(G_{n,p},x)
 = \sum\limits_{k=0}^n \binom{n}{k} \frac{y^k}{{q}^{k(n-k)}}.
$$

A polynomial $F(x) = \sum\limits_{i=0}^n a_i x^i$ of degree~$n$
is called symmetric if $a_i = a_{n-i}$ for all~$i$. 

{\bf Statement 5.1}.
a) The polynomial $\widetilde{C}(G_{n,q},y)$ is symmetric.

b) For odd $n$, $p^{\frac{n-1}{2}}$ is a middle root of $PC(G_{n,p},x)$.
All other roots of $PC(G_{n,p},x)$ for odd $n$
and all roots for even $n$ could be gathered in pairs
with the roots product equal~$p^{n-1}$.

Denote by $\beta(G_{n,p})$ the largest root of $PC(G_{n,p},x)$.
Statement~5.1 allows us to find $\beta(G_{n,p})$ for all $n\leq5$ (Corollary~5.1).

In~$\S5.2$, we define the {\it random algebra}. 
Let $X =\{x_1,\ldots,x_n\}$ be a finite set. Fix an order on $X$ such that $x_i>x_j$ if $i<j$.
Consider a word $w = w_1\ldots w_m\in X^*$ of length $m = |w|$.
Let a letter $x_{i_j}$ occurs in $w$ exactly $m_i\geq1$ times, $i=1,\ldots,k$.
We suppose that $x_{i_1}>\ldots>x_{i_k}$. Consider a new alphabet
$$
X' = X'(w) = \{x_{i_1}^1,\ldots,x_{i_1}^{m_1},x_{i_2}^1,\ldots,x_{i_2}^{m_2},\ldots,x_{i_k}^1,\ldots,x_{i_k}^{m_k}\}.
$$
Define an order on the set $X'$: $x_{i_a}^s>x_{i_b}^t$ if $a<b$ or $a=b$ and $s<t$.

Given a word $w\in X^*$, let us construct a word $w'\in (X')^*$ of the same length as follows.
If $w_j$ is the $t$-th occurrence (counting from the left) of a letter $x_s$ in~$w$,
then the $j$-th letter of $w'$ equals $x_s^t$.
Denote the set of all multipartite graphs with parts 
$\{x_{i_1}^1,\ldots,x_{i_1}^{m_1}\}$, \ldots, $\{x_{i_k}^1,\ldots,x_{i_k}^{m_k}\}$
as $MP(w)$. Define $M = M(w)$ be equal to the product $m_1m_2\ldots m_k$.

Let $p\in[0;1]$. Define a weight $s_p(w)$ of a word $w$ as 
$$
s_p(w) = \sum\limits_{G\in MP(w)}p^{|E(G)|}(1-p)^{M-|E(G)|}
 I(w'\ \mbox{is in n.f. in }\,M(X',G)),
$$
where n.f. means ``normal form'' 
(the maximal word among all words in a partially commutative monoid equal to it),
$M = \prod\limits_{i=1}^k m_i$,
$I(A) = \begin{cases}
1, & A\ \mbox{is true}, \\
0, & \mbox{otherwise}.
\end{cases}$

Actually $s_p(w)$ equals a probability of the event that 
$w'$ is in the normal form in hypothetical partially commutative monoid
with commutativity graph $G_{p}(w)$, 
the random multipartite graph with fixed parts 
with $m_1,\ldots,m_k$ vertices and edge probabi\-lity~$p$.

Define on the free associative algebra $\Ass\langle X\rangle$
as on the vector space a new product $\cdot$.
Let $X_n$ denote the set of all words of length $n$ in the alphabet~$X$. 
For $w_i\in X_{n_i}$, $i=1,2$, 
\begin{multline*}
w_1\cdot w_2
 = \frac{1}{s_p(w_1)s_p(w_2)}
\sum\limits_{G\in MP(w_1w_2)}p^{|E(G)|}(1-p)^{M-|E(G)|}
I(w_1',w_2'\ \mbox{in n.f.})[(w_1w_2)'] \\
 = \sum\limits_{u\in X_{n_1+n_2}}P(u = [(w_1w_2)'] \mid w_1',w_2'\ \mbox{in n.f.})u,
\end{multline*}
where $[w]$ denotes the normal form of~$w$,
$P(A)$ denotes the probability of an event $A$ 
in the probability theory model constructed 
by the random multipartite graph~$G_p(w_1w_2)$
and $P(A\mid B)$ denotes the conditional probability of~$A$ given~$B$. 

Let us call the space $\Ass\langle X\rangle$ under the product $\cdot$
as random algebra, notation:~$\Ass(X,p)$.

{\bf Lemma 5.3}.
a) For $0\leq p<1$, the algebra $\Ass(X,p)$ is isomorphic to 
the free asso\-ciative algebra $\Ass\langle X\rangle$.
For $p = 1$, $\Ass(X,p)$ is isomorphic to the polynomial algebra $\mathbb{R}[X]$.

b) The map $s_p\colon \Ass(X,p)\to \langle\mathbb{R},\cdot\rangle$
is a semigroup homomorphism.

Extend a weight $s_p$ on $\Ass(X,p)$ by linearity and put $m_t(p) = s_p(X_t)$.

{\bf Theorem 5.2}.
The numbers $m_t(p)$, $t\geq1$, satisfy the reccurence relation
$$
m_t(p)
= \binom{n}{1} m_{t-1}(p) - \binom{n}{2}p m_{t-2}(p) 
 + \ldots  + (-1)^{k+1}\binom{n}{k} p^{\binom{k}{2}} 
 + \ldots + (-1)^{n+1}p^{\binom{n}{2}}m_{t-n}(p)
$$
with initial data $m_0(p) = 1$, $m_{-1}(p) = \ldots = m_{-n+1}(p) = 0$.

{\bf Corollary 5.3}.
The polynomial $PC(G_{n,p},x)$ is a characteristic polynomial
for the sequence $\{m_t(p)\}$ and $\beta(G_{n,p})$ equals its growth rate.

{\bf Lemma 5.4}.
a) The following inequalities hold
$$
1+(n-1)(1-p)\leq \beta(G_{n,p})\leq 1+(n-1)\sqrt{1-p}.
$$

b) The number $\beta(G_{n,p})$ for fixed~$n$
is strictly monotonic function on $p\in[0;1]$ decreasing from $n$ to~1.

\begin{figure}[h]
\centering
\includegraphics[height = 4cm]{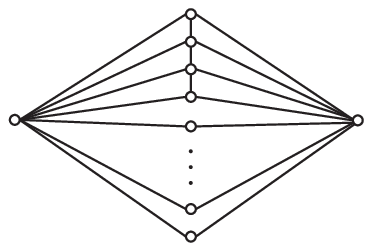}
\caption{Graph $G_-$ for the case $n\geq5$, $2n-4 < k < 3n-6$.} \label{Pic1}
\end{figure}

In~$\S6$, we solve the following problem.
Denote the set of all planar graphs with $n$ vertices and $k$ edges as $Pl(n,k)$.
We find minimal and maximal values of $\beta(G)$ for $G\in Pl(n,k)$
and the graphs for these extremal values (Theorem~6.1).
Let $n\geq6$, then the graph $G_-\in Pl(n,k)$ with the minimal $\beta$
is a~triangle-free graph for $k\leq 2n-4$.
For $2n-4 < k < 3n-6$, we construct $G_-$ as a~supergraph 
of $K_{2,n-2}$ in which $k-2n+4$ edges in the big part form a tree (see Picture~\ref{Pic1}). 
For $k = 3n-6$, put $G_- = \bar{K}_2\cup C_{n-2}$.

Let $n\geq4$, $k\geq3$, then the graph $G_+\in Pl(n,k)$ with the maximal $\beta$
is constructed as follows. We start with $K_3$.
On each step, we add one new vertex inside of some 
(triangle) face of the graph and connect it with each vertex of the face.
We proceed on while we have edges (see Picture~\ref{Pic3}). 
Sometimes, such graphs are called Apollonian networks.

Let $Pl(n)$ denote the set of all planar graphs with $n$ vertices.
In 2007, O. Gim\'{e}nez and M. Noy stated~\cite{Noy} that
$|Pl(n)|\sim C_0n^{-7/2}\gamma^n n!$
for $\gamma \approx 27.227$ and 
a planar graph in average contains $\kappa \approx 2.213 n$ edges. 

{\bf Theorem 6.2}.
a) $PC(G,x)$ of almost all planar graphs is a polynomial of 4-th degree,
has two complex and real roots. Complex roots lie in the right half-plane,
real roots are simple.

b) Let $\varepsilon>0$, then for almost all planar graphs
we have $|\rho|>\frac{1-\varepsilon}{6\gamma^3\kappa}$
for any root $\rho$ of $PC$-polynomial.

{\bf Statement 6.1}.
The average value of the growth rate of partially commutative monoid
with planar commutativity graph equals
$$
\beta_{\mathrm{ev},Pl}(n)
=\frac{1}{|Pl(n)|}\sum\limits_{G\in Pl(n)}\beta(G) = n-\kappa + O(1/n).
$$

In~$\S7$, we formulate the main problems of the paper.
By analogy with the Nordhaus---Gadddum inequalities
for chromatic number of a graph~$G$ and its complement, we post 

{\bf Problem~1}. To find the tight bounds for 
the expressions $\beta(G)+\beta(\bar{G})$ and $\beta(G)\beta(\bar{G})$.

Denote by $G(n,k)$ the set of all graphs with $n$ vertices and $k$ edges. Introduce 
$$
\beta_-(n,k) = \min\limits_{G\in G(n,k)}\beta(G),\quad
\beta_+(n,k) = \max\limits_{G\in G(n,k)}\beta(G).
$$

{\bf Problem~2}.
To find the values $\beta_{\pm}(n,k)$ and the graphs on which they are reached.

The following lower bound for $\beta(G)$ (part a)) was proved by D.C. Fisher in 1989.

{\bf Theorem 7.1}.
a) \cite{Fisher1989} Given a graph $G$ with $n$ vertices and $k$ edges, 
$\beta(G)\geq n-\frac{2k}{n}$. 

b) The bound from a) is reached if and only if $G$ is an empty graph or 
$G$ is a complete multipartite graph with equal parts.

For a graph $G\in G(n,k)$, define the edge $PC$-density 
$e(G) = \frac{n-\beta(G)}{k}$, $e(\bar{K}_n) := \frac{1}{n}$.

{\bf Corollary 7.1}.
For any graph $G$ with $n$ vertices, $e(G)\leq\frac{2}{n}$.

{\bf Corollary 7.3}.
For any graph $G$ with $n$ vertices, 

a) $n+1\leq \beta(G)+\beta(\bar{G})$,

b) $n\leq \beta(G)\beta(\bar{G})$.

\noindent Moreover, the bounds are reached only if and only if 
$\{G,\bar{G}\} = \{K_n,\bar{K}_n\}$.

{\bf Corollary 7.5} \cite{Fisher1989}.
For any graph $G$, we have $\beta(G)\geq 1 +\rho(\bar{G})$.

{\bf Statement 7.1}.
Let $k\leq\frac{n^2}{4}$, then 
$\beta_-(n,k) = \frac{n+\sqrt{n^2-4k}}{2}$.
This value is reached for a~graph $G\in G(n,k)$
if and only if $G$ is a triangle-free graph.

In 1990, D.C. Fisher and J.M. Nonis proved the strong lower bound.

{\bf Theorem 7.2}~\cite{FisherUpper}.
Given a graph $G$ with $n$ vertices and $k$ edges, let $w$ be such a~natural number that
$\big(1-\frac{1}{w-1}\big)\frac{n^2}{2}<k\leq \big(1-\frac{1}{w}\big)\frac{n^2}{2}$. Then 
$$
\beta(G)\geq \frac{n}{w}\left(1+\sqrt{1 -\frac{2kw}{n^2(w-1)}}\right).
$$

The section~$\S8$ could be considered as the central section of the work.
In~$\S8.1$, we find the graph $G\in G(n,k)$ with $\beta(G) = \beta_+(n,k)$.

Let us define a relation $\geq$ on simple graphs as follows.
We write that $G\geq H$ if $D(G,x)\leq D(H,x)$ on the line segment $[0;1/\beta(G)]$. 
From $G\geq H$, we have $\beta(G)\geq \beta(H)$.

{\bf Lemma 8.1} \cite{CsikvariKelman}.
a) Given an induced subgraph $H$ of $G$, we have $G\geq H$; \\
b) Given a spanning subgraph $H$ of $G$, we have $H\geq G$.

Let $G$ be a graph, $u,v\in V(G)$, $u\neq v$.
In 1981, A. Kelmans defined \cite{Kelmans} 
so called {\it Kelmans transformation} which transfers a graph $G$ into a graph $G' = KT(G,u,v)$.
To get $G'$, we erase all edges between $v$ and $N(v)\setminus (N(u)\cup\{u\})$
and add all edges between $u$ and $N(v)\setminus (N(u)\cup\{u\})$
(see Picture~\ref{Pic4}). Note that $|E(G')| = |E(G)|$.

In~2011, P. Csikv\'{a}ri stated~\cite{CsikvariKelman} 
some important properties of Kelmans transformations.

{\bf Lemma 8.2} \cite{CsikvariKelman,CsikvariThesis}.
Let $G$ be a graph and $G'$ be a graph obtained from~$G$ by a Kelmans transformation. Then 

a) $G'\geq G$ and so, $\beta(G')\geq\beta(G)$,

b) $c_k(G')\geq c_k(G)$ for all $k$.

It is easy to show that (see \cite{CsikvariNikif})
any graph by a series of Kelmans transformations
can be transformed to a threshold graph.

Let $G$ be a graph such that $V(G) = V(G_1)\cdot \hspace{-9.5pt}\cup V(G_2)$.
Moreover, let any vertex of $G_2$ be either connected or
disconnected with all vertices of $G_1$ and
$u\in V(G_1)$ be such hanging vertex in $G_1$ that $G_1\setminus u$ is not complete.
Define the {\it isolating transformation} which transforms~$G$ to a graph 
$G'=I(G,u)$ as follows.
We obtain a $G'$ by arising the only edge in $G_1$
incident to~$u$ and adding an edge in $G_1\setminus u$.

{\bf Lemma 8.3}.
Let $G$ be a graph such that $V(G) = V(G_1)\cdot \hspace{-8pt}\cup V(G_2)$.
Moreover, let any vertex of $G_2$ be either connected or
disconnected with all vertices of $G_1$ and
$u\in V(G_1)$ be such hanging vertex in $G_1$ that $G_1\setminus u$ is not complete.
There exists an isolating transformation $G' = I(G,u)$
such that $\beta(G')\geq\beta(G)$.

In the next statement, we prove Conjecture~1 from~\cite{FisherUpper}.

{\bf Theorem 8.1}.
Let $n,k$ be natural numbers, $k = \binom{d}{2} + e\leq \binom{n}{2}$ for $0\leq e<d$.
Construct a graph $G$ with $n$ vertices and $k$ edges as follows.
We add a vertex of degree $e$ to the complete graph $K_d$
and leave all other vertices to be isolated. Then $\beta_+(n,k) = \beta(G)$.

{\bf Corollary 8.1} \cite{CutlerRadcliffe,Wood2007}.
The constructed graph $G$ from Theorem~8.1 maximizes
all numbers $c_k$ among graphs from $G(n,k)$.
In particular, $C(H,x)\leq C(G,x)$ for any graph $H\in G(n,k)$
and for all $x\geq 0$.

In Remark 8.1, we discuss how the strategy of the proof of Theorem~8.1
could be applied to reprove the analogous result for the spectral radius.
This problem was initially posed by Brualdi and Hoffman in 1976~\cite{BrualdiHoffman}
and solved by P. Rowlinson in 1988 \cite{Rowlinson}.

In~$\S8.2$, we want to derive the upper bound on $\beta(G)$ applying Theorem~8.1.
Before this, we easily prove the necessary condition for real-rootedness of $PC(G,x)$.

{\bf Statement 8.1}.
Let $G$ be a graph with $n$ vertices such that 
all roots of $PC(G,x)$ are real. Then $\beta\leq n-\frac{k}{n}$
and $e(G)\geq \frac{1}{n}$.

{\bf Lemma 8.4}.
Let $G$ be a graph with $n\gg1$ vertices and 
$k\geq \frac{n^2}{2}\big(1-\frac{1}{pe^{p+2}}\big)^2$ edges,
then $\beta(G)<\frac{n}{p}$ and $e(G)>\frac{2}{n}\big(1-\frac{1}{p}\big)$.

{\bf Theorem 8.2}.
Let $n\gg1$.
For any graph $G$ with $n$ vertices and $k$ edges, we have

a) $\beta(G)\leq n-\frac{\alpha k}{n}$,

b) $e(G)\geq \frac{\alpha}{n}$, \\
where $\alpha\approx 0.9408008$.

Another important fact which follows from the proof of Theorem~8.2 is the following:
we have the asymptotically tight upper bound provided $k - n^2/2 = O(n^2)$
$$
\beta(G) \lesssim n\cdot\frac{\sqrt{x}}{W(\frac{\sqrt{x}}{1-\sqrt{x}})},\quad x = \frac{2k}{n^2},
$$
where $W(x)$ is the Lambert $W$-function, the inverse function to $f(y) = ye^y$.

{\bf Corollary 8.3}.
Let $G$ be a graph with $n\gg1$ vertices and $k = O(n^{2\theta})$, $\theta<1$, edges.
Then $\beta(G)\sim n -\frac{k}{n}$. 

{\bf Corollary 8.4}.
Let $G$ be a graph with $n\gg1$ vertices and $k$ edges.

a) For $k\geq 0.256736n^2$, we have 
$e(G)\geq\frac{1}{n}$ and $\beta(G)\leq n-\frac{k}{n} < \frac{3n}{4}$.

b) For $k\geq n^2/4$, we have 
$e(G)\geq\frac{0.996}{n}$
and $\beta(G)\leq n-\frac{0.996k}{n} < 0.751n$.

{\bf Corollary 8.5}.
For any graph $G$ with $n\gg1$ vertices the following inequalities hold

a) $\beta(G)+\beta(\bar{G})<1.50197n$,

b) $\beta(G)\beta(\bar{G})<0.56398n^2$.

We conjecture that the maximal values of the expressions 
$\beta(G)+\beta(\bar{G})$ and $\beta(G)\beta(\bar{G})$
are reached on graphs from the class 
$A(n) = \{K_s \cup \bar{K}_{n-s},\bar{K}_s + K_{n-s}\}$.
By Example 8.2, the following values are maximal among graphs from $A(n)$ 
\begin{gather*}
\beta(G)+\beta(\bar{G})\approx 1.46594n, \\
\beta(G)\beta(\bar{G})\approx 0.535919n^2.
\end{gather*}

In~$\S9$, we study the minimum value $\beta_-(n,k)$.
At first, we state some results devoted to transformations of graphs.
Given a graph $G$ and two distinct vertices $u,v\in V(G)$, let us call 
a Kelmans transformation $KT(G,u,v)$ as a {\it nontrivial} one, if $c(G')-c(G)>0$,
where $c(H)$ denotes the number of all cliques in a graph $H$.

{\bf Lemma 9.1}.
Let $G$ be a graph with connected complement.
If a graph $G'$ is a result of a nontrivial Kelmans transformation of~$G$,
then $\beta(G')>\beta(G)$.

{\bf Lemma 9.2}.
Let $G$ be a graph with connected complement.
Let $e\in E(G)$ be an edge lying in a clique of size~$t\geq3$
and $a,b\in V(G)$ such vertices that $(a,b)\not\in E(G)$
and $N_G(a)\cap N_G(b) = \emptyset$.
Consider the graph $G'$ obtained by removing an edge $e$
and adding an edge $(a,b)$. Then $G'\leq G$ and $\beta(G')<\beta(G)$.

{\bf Corollary 9.1}.
The graph $G$ constructed in Theorem~8.1 is a unique graph with the maximal~$\beta(G)$
among all graphs with $n$ vertices and $k$ edges with one exception:
when $k = \binom{d}{2}+1$ for some~$d$. In this case, 
the set $\{H\in G(n,k)\mid \beta(H) = \beta_+(n,k)\}$
consists of all graphs obtained from $K_d\cup \bar{K}_{n-d}$ by adding one edge.

{\bf Corollary 9.2}.
Let $k>[n^2/4]$ and $G$ be a graph such that $\beta(G) = \beta_-(n,k)$.
Then $G$ is connected graph having diameter~2.

We will find the exact value of $\beta_-(n,k)$ modulo the following conjecture.

{\bf Conjecture 9.1}.
Let $k>[n^2/4]$ and $G$ be a graph such that $\beta(G) = \beta_-(n,k)$.
Then $\bar{G}$ is disconnected.

Let $n^2/4<k$ and 
$\big(1-\frac{1}{w-1}\big)\frac{n^2}{2}<k<\big(1-\frac{1}{w}\big)\frac{n^2}{2}$
for a natural number $w$. If $k$ is not enough large to construct $K_{a,a,\ldots,a,b}$
(with $w-1$ parts of $a$ and $a\geq b\geq0$), then 
we construct the graph $G_-\in G(n,k)$ as a supergraph of 
$K_{l+1,l+1,\ldots,l+1,l,\ldots,l}$ 
with $p$ parts with $l+1$ edges and $q$ parts with $l$ edges,
where $l = \big[\frac{n}{w-1}\big]$, $p = n-l(w-1)$ and $q = w-1-p$.
Introduce 
$$
k' = k - \left(\binom{p}{2}(l+1)^2+\binom{q}{2}l^2+pql(l+1)\right).
$$
Further, we draw a triangle-free graphs with $[k'/p]$ vertices inside all $p$
parts with $l+1$ vertices and draw remaining edges anywhere. So, 
$$
\beta(G_-) = \frac{l+1}{2}\left(1+\sqrt{1 - \frac{4[k'/p]}{(l+1)^2}}\right).
$$

Let $k$ be enough large to construct a supergraph of $K_{a,a,\ldots,a,b}$ 
with prescribed conditions on parts.
Find a natural number $n_1$ such that
$$
(w-1)n_1\left(n-\frac{wn_1}{2}\right)\leq k 
< (w-1)(n_1-1)\left(n-\frac{w(n_1-1)}{2}\right).
$$

Denote $k' = k - \big((w-1)n_1n-\binom{w}{2}n_1^2\big)$ and 
$0\leq k'< (w-1)(wn_1-n-w/2)$. 
We construct the graph $G_-$ as a supergraph of 
$K_{n_1,\ldots,n_1,n-(w-1)n_1}$
in which the $[k'/(w-1)]$ edges form a triangle-free graph
in each part with $n_1$ vertices and we put remaining 
$k' - (w-1)[k'/(w-1)]$ edges anywhere. Hence, 
$$
\beta(G_-) = \frac{n_1}{2}\left(1 + \sqrt{1-\frac{4k'}{(w-1)n_1^2}+\frac{4\varepsilon}{(w-1)n_1^2}}\right),
$$
where $\varepsilon = k'/(w-1) - [k'/(w-1)]\in\{0,1\}$.

{\bf Theorem 9.2}.
Let $\big(1-\frac{1}{w-1}\big)\frac{n^2}{2}<k<\big(1-\frac{1}{w}\big)\frac{n^2}{2}$. 
If Conjecture~9.1 holds, then $\beta_-(n,k) = \beta(G_-)$
for the constructed graph $G_-$. 

{\bf Corollary 9.4}.
Let 
$\big(1-\frac{1}{w-1}\big)\frac{n^2}{2}<k\leq\big(1-\frac{1}{w}\big)\frac{n^2}{2}$, $w\geq2$, then 
$$
\frac{n}{w}+\frac{1}{w}\sqrt{n^2-\frac{2kw}{w-1}}\leq 
 \beta_{-}(n,k) < \frac{n}{w}+\frac{1}{w}\sqrt{n^2-\frac{2kw}{w-1}}+1.
$$

In~$\S10.1$, we are interested on the asymptotics of the average 
growth rate of partially commutative monoid with $n$-vertex commutativity graph: 
$$
\beta_{\evv}(n)
 =\frac{1}{2^{\binom{n}{2}}}\sum\limits_{G\colon |V(G)| = n}\beta(G).
$$

{\bf Lemma 10.1}.
For all $p\in[0;1]$, there exists the limit 
$\lim\limits_{n\to\infty}\frac{\beta(G_{n,p})}{n}$.

{\bf Theorem 10.1}.
The average value of $\beta(G)$ on graphs with $n\gg1$ vertices
asymptoti\-cally equals
$$
\beta_{\mathrm{ev}}(n)
 \sim n\lim\limits_{n\to\infty}\frac{\beta(G_{n,1/2})}{n} = \beta_0 n \approx 0.672008n.
$$

The constant $\beta_0$ firstly appeared 
in the article of R. Stanley~\cite{Stanley} in 1973:
the number of all acyclic orientations of a digraph was counted as 
$A\cdot n!2^{n(n-1)/2}\beta_0^n$ for $A\approx 1.741$.

{\bf Statement 10.1}.
For almost all graphs with $n$ vertices,
$\beta(G)$ lies in a neighbourhood of $\beta_0 n$
and is the unique root $PC(G,x)$ which modulus is greater than $n/2$.

Define the number $\beta_{\evv}(n,k)$
as the average value of $\beta(G)$ for the set $G(n,k)$
of all graphs with $n$ vertices and $k$ edges:
$$
\beta_{\evv}(n,k)
 =\frac{1}{\binom{\binom{n}{2} }{k}}\sum\limits_{G\in G(n,k)}\beta(G).
$$

{\bf Statement 10.2}.
Let $k(n)$ be such integer-valued function that $0\leq k(n)\leq \frac{n(n-1)}{2}$
and there exists $\lim\limits_{n\to\infty}\frac{2k(n)}{n^2} = k_0<1$. Then
$\beta_{\evv}(n,k) \sim n\lim\limits_{n\to\infty}\frac{\beta(G_{n,k_0})}{n}$.

In~$\S10.2$, we are interested on all roots of $PC(G_{n,p},x)$. 
Let $p>0$ and $\beta_r = \beta_r(G_{n,p})$ denote the $r$-th
largest root of $PC(G_{n,p},x)$.

{\bf Theorem 10.2}.
For all $p\in(0;1]$ and $r\geq1$, there exists the limit 
$\lim\limits_{n\to\infty}\frac{\beta_r(G_{n,p})}{n}$.

{\bf Corollary 10.2}.
Let $r>0$. For almost all graphs with $n\gg1$ vertices, 
the real roots of $PC$-polynomial which moduli is not less than $n/r$
lie in neighbourhoods of the roots of $PC(G_{n,1/2},x)$.

One can approximately compute six largest roots of $PC(G_{n,1/2},x)$: 
\begin{gather*}
0.672008n,\quad
0.204871n,\quad
0.073744n,\quad
0.028756n,\quad
0.011768n,\quad
0.004975n.
\end{gather*}

Further, we show how all coefficients of the series expansion of 
$\lim\limits_{n\to\infty}\frac{\beta_r(G_{n,p})}{n}$ on $p$
could be calculated. In particular, we find 
\begin{multline*}
\frac{\beta_1(G_{n,p})}{n}\sim 
1-\frac{p}{2}-\frac{p^2}{4}-\frac{p^3}{12}-\frac{p^4}{16}
 -\frac{p^5}{48}-\frac{7p^6}{288}-\frac{p^7}{96} -\frac{7p^8}{768} \\
 - \frac{49p^9}{6912}-\frac{113p^{10}}{23040}-\frac{17p^{11}}{4608}
 -\frac{293p^{12}}{92160}-\frac{737p^{13}}{276480}-\frac{3107p^{14}}{1658880}+O(p^{15}),
\end{multline*}
$$
\frac{\beta_2(G_{n,p})}{n}\sim 
\frac{p}{2}\left(1-\frac{p}{6}-\frac{5p^2}{18} -\frac{29p^3}{216} 
 - \frac{85p^4}{648}-\frac{163p^5}{3888}-\frac{1387p^6}{19440}\right) + O(p^8).
$$

{\bf Conjecture 10.1}.
For any~$r\geq1$, we have 
$$
\lim\limits_{n\to\infty}\frac{\beta_r(G_{n,p})}{n}
 = \frac{p^{r-1}}{r}\left(1-\frac{p}{r(r+1)}+O(p^2)\right).
$$ 

There are some articles devoted to the study of the behavior 
and properties of $C(G_{n,p},1)$. 
In 2014, W. Gawronski and T. Neuschel stated (see also \cite{Paris})

{\bf Theorem 10.3}~\cite{Gawronski}.
For a fixed $p\in(0;1)$, as $n\to\infty$, we have 
$$
C(G_{n,p},1)
 = \frac{1}{\sqrt{r(n)}}
 \left(\theta_3 \left(\frac{\pi r(n)}{\ln(1/p)}, e^{-2\pi^2/{\ln(1/p)}}\right)+o(1)\right)
 \exp\left(\frac{r(n)^2 + 2r(n)}{2\ln(1/p)}\right),
$$
where $\theta_3(z,q)$ is the Jacobi's third Theta function,
$r(n)$ is defined as the positive solution $t$ of the equation
$t\left(e^t+\sqrt{1/p}\right) = n\sqrt{1/p}\ln(1/p)$.

In 2012--13 \cite{BrownAsy,BrownFn}, 
J. Brown, K. Dilcher and V. Manna initiated to study the poly\-nomials 
$f_n(z) = \sum\limits_{k=0}^n\binom{n}{k}z^{\binom{k}{2}}$
for complex variable $z$.
For real $z\in[0;1]$, $f_n(z) = C(G_{n,z},1)$.

{\bf Statement 10.4} \cite{BrownFn}.
a) For each $n\geq 0$, $f_{2n+1}(z)$ is divisible by $z^n + 1$.

b) For each $n\geq 3$ the roots of $f_n(z)$ lie inside a circle of radius $1 + \frac{3\ln n}{n}$.

c) For each $n\geq 3$ the roots of $f_n(z)$ lie outside a circle of radius $2/n$.

d) For each $n\geq 4$ there is a negative real root of $f_n(z)$ in the interval
$-\frac{2+\frac{2}{n}}{n}< z < -\frac{2}{n}$.

Finally, in Picture~\ref{Pic8}, the borders of the possible
values of $\beta(G)/n$ are drawn (asymp\-totically). 

Weighted clique-type polynomials are defined in~$\S11$.
Given a simple graph $G$ and the set $Cl(G)$ of all cliques of $G$,
define the weighted dependence polynomial of $G$ as 
$$
D_w(G,x) = \sum\limits_{B\in Cl(G)}(-1)^{|B|}w(B),\quad 
w(B) = \prod\limits_{v\in B}\alpha_v x^{d_v},\
\alpha_v,d_v\in \mathbb{R}_{+}.
$$
If $\alpha_v = d_v = 1$ for all $v\in V(G)$, 
then $D_w(G,x)$ coincides with $D(G,x)$.

Define the weighted $PC$-polynomial of $G$ as 
$PC_w(G,x) = x^{w_0}D_w(G,1/x)$ for $w_0 = \deg(D_w(G,x))$.
S. Lavall\'{e}e and C. Reutenauer in 2009 and 
S. Lavall\'{e}e in 2010 proved the following results
which we gathered in one statement.

{\bf Theorem 11.1}.
a) \cite{Reutenauer} 
The set of all weighted dependence polynomials 
with $\alpha_v,d_v\in \mathbb{N}_{>0}$
coincide with the set of all polynomials of the 
form $\det(E-xM)$, where $M$ is a square matrix with natural entries.

b) \cite{Lavallee}
The set of all weighted dependence polynomials 
with $d_v\in \mathbb{N}_{>0}$
coincide with the set of all polynomials of the 
form $\det(E-xM)$, where $M$ is a square matrix with nonnegative entries.

In the proof, we partly follow the original proofs and partly suggest new steps.
In particular, we state the next result about the largest real root $\beta_w(G)$ of $PC_w(G,x)$.

{\bf Lemma 11.1}.
Let $G$ be a graph such that $\bar{G}$ is connected.
Let $D_w(G,x)$ be a weighted dependence polynomial
such that the set $\{d_v\mid v\in V(G)\}$ is coprime.
Then the number $1/\beta_w(G)$ is the only complex root 
of $D_w(G,x)$ with modulus less or equal to $1/\beta_w(G)$.

{\bf Corollary 11.1}.
a) \cite{Lavallee}
The set of all weighted $PC$-polynomials with $d_v\in \mathbb{N}_{>0}$
multiplied by $x^k$, $k\in\mathbb{N}$,
coincide with the set of all characteristic polynomials
of square matrices with nonnegative entries.

b) \cite{Reutenauer} 
The set of all weighted $PC$-polyno\-mials with $\alpha_v,d_v\in \mathbb{N}_{>0}$
multiplied by $x^k$, $k\in\mathbb{N}$,
coincide with the set of all characteristic polynomials of
square matrices with natural entries.

In 2015, C. McMullen considered~\cite{McMullen} weighted dependence polynomial 
$D_w(G,x)$ with $\alpha_v=1$ for all $v\in V(G)$, see observation of his results directly in~$\S11$.

In 2005, A. Scott and A. Sokal found~\cite{Scott} the deep connection
between weighted (in)dependence polynomials and Lov\'{a}sz local lemma.

For a graph~$G$, we define the weighted independence polynomial 
$I_w(G,x)$ as $D_w(\bar{G},x)$.

{\bf Theorem 11.4} \cite{CsikvariThesis,Scott}.
Assume that given a graph $G$ and there is an event $A_i$
assigned to each vertex~$i$. 
Assume that $A_i$ is totally independent of the events
$\{A_k\mid (i,k)\not\in E(G)\}$.
Set $P(A_i) = p_i$ and define the weighted independence polynomial
$I_w(G,x)$ with $\alpha_i = p_i$ and $d_i = 1$ for all $i=1,\ldots,n$.

a) Assume that $I_w(G,t)>0$ for $t\in[0;1]$, then 
$P\big(\bigcap\limits_{i=1}^n \overline{A_i}\big)\geq I_w(G,1)>0$.

b) Assume that $I_w(G,t) = 0$ for some $t\in [0;1]$.
Then there exist a probability space and a~family of events $B_i$, $i=1,\ldots,n$,
with $P(B_i)\leq p_i$ and with dependency graph~$G$ such that
$P\big(\bigcap\limits_{i=1}^n \overline{B_i}\big) = 0$.

{\bf Corollary 11.2} \cite{CsikvariThesis,Scott}.
Let $A_i$, $i=1,\ldots,n$, be events with dependence graph $G$
such that $P(A_i)\leq t$ for all $i=1,\ldots,n$. Then 
$P\big(\bigcap\limits_{i=1}^n \overline{A_i}\big)>0$
if and only if $t\leq 1/\beta(\bar{G})$.

{\bf Corollary 11.3}.
a) \cite{Scott}
For any graph $G$, 
$\beta(G)\leq \frac{d^{d}}{(d-1)^{d-1}}<ed$,
where $d = \Delta(\bar{G}){\geq}2$.
For $d = 1$, we have $\beta(G) = 2$.

b) Let $G$ be a $(n-d-1)$-regular graph, then 
$1+d\leq\beta(G) < ed$.

{\bf Corollary 11.4}.
Let $G$ be a graph with $n\gg1$ vertices and $k$ edges,
$A_i$, $i=1,\ldots,n$, are events with dependence graph~$G$.
Then $P\big(\bigcap\limits_{i=1}^n \overline{A_i}\big)>0$, if 

a) $k\leq 0.24326n^2$ and $P(A_i)\leq 4/(3n)$, $i=1,\ldots,n$, or

b) $k\leq n^2/4$ and $P(A_i)\leq 1.3316/n$, $i=1,\ldots,n$.

Note that in analogous way one can produce another global versions
of Lov\'{a}sz local lemma in case $k\leq \alpha n^2$.
It is an open question: 
Wether such version of Lov\'{a}sz local lemma like in Corollary~11.4 is useful?

The main goal of~$\S12$ is to prove the Chudnovsky---Seymour theorem (2004).
It was conjectured by Y. Hamidoune in 1990 \cite{Hamidoune} and R. Stanley in 1998 \cite{Stanley98}.

{\bf Theorem 12.1} \cite{Bencs,Chudnovsk}.
If $G$ is a claw-free graph then all roots of $I(G,x)$ are real.

We follow the proof of F. Bencs~\cite{Bencs} (maybe, 	the simplest one).

{\bf Corollary 12.1}.
If $\bar{G}$ is a claw-free graph then all roots of $PC(G,x)$ are real.

{\bf Corollary 12.2}.
Let $G$ be a graph with $n$ vertices and $k$ edges.

a) If $\bar{G}$ is claw-free then 
$n-\frac{2k}{n}\leq \beta(G)\leq n-\frac{k}{n}$
and $\frac{1}{n}\leq e(G)\leq \frac{2}{n}$.

b) If $\beta(G)>n-\frac{k}{n}$ then $\bar{G}$ has a claw.

c) Let $G$, $\bar{G}$ be claw-free~\cite{Pouzet}, then
$\beta(G)+\beta(\bar{G})\leq \frac{3n+1}{2}$,
$\beta(G)\beta(\bar{G})\leq \big(\frac{3n+1}{4}\big)^2$.

In~$\S12.2$, we study the matching polynomial
$\mu(G,x)=\sum\limits_{k=0}^{\nu(G)}(-1)^k m_k x^{n-2k}$,
where $n = |V(G)|$, $m_k$ denotes the number of matchings with $k$ edges in $G$
and $\nu(G)$ is the maximum of sizes of matchings in $G$.
The matching-generating polynomial 
$M(G,x)=\sum\limits_{k=0}^{\nu(G)}m_kx^k$
and matching polynomial are are connected in the following way:
$$
\mu(G,x) = x^n M(G,-x^{-2}),\quad M(G,x) = I(L(G),x),
$$
where $L(G)$ denotes the line graph of $G$.

From Theorem~12.1 follows the result of O. Heilmann and E. Lieb from 1972:

{\bf Corollary 12.3} \cite{Heilman72}.
All roots of $\mu(G,x)$ and $M(G,x)$ are real.

Denote the largest root of $\mu(G,x)$ as $t(G)$. 
By the definition, $t^2(G) = \beta(\overline{L(G)})$. 

{\bf Statement 12.1} \cite{Fisher1992-Match}.
For a graph~$G$ with $n$ vertices and $k$ edges, 
$t^2(G)\geq \frac{4k}{n}-1$.

{\bf Statement 12.2} \cite{Fisher1992-Match,Heilman72}.
Let $\Delta>1$ denote the maximal degree of graph~$G$, then 
the following inequalities hold
$\sqrt{\Delta}\leq t(G)\leq 2\sqrt{\Delta-1}$.

{\bf Statement 12.3} \cite{LeakeRyder}.
If $\bar{G}$ is a claw-free graph then 
$\alpha(G)\leq \beta(G)\leq 4\max\{\alpha(G)-1,\Delta(G)\}$.

At the end of~$\S12$, we list some results about matchings. 
Denote by $m(G)$ the number of matchings in a graph~$G$
and by $pm(G) = m_{|V(G)|/2}(G)$ the number of perfect matchings. 

{\bf Statement 12.6} \cite{Davies2017}.
For all $d$-regular graphs $G$ on $n$ vertices (where $2d$ divides $n$), 
$m_k(G) \le 2 \sqrt{n} \cdot m_k(H_{d,n})$,
where $H_{d,n}$ denotes the $d$-regular, $n$-vertex graph that 
is the disjoint union of $n/(2d)$ copies of $K_{d,d}$.

By the edge occupancy fraction \cite{Davies2017} we mean 
$\alpha^M(G,x) = \dfrac{xM'(G,x)}{|E(G)|M(G,x)}$.

{\bf Statement 12.7} \cite{Davies2017}.
For any $d$-regular graph $G$, the following inequalities hold 

a) $\alpha^M(G,x)\leq\alpha^M(K_{d,d},x)$,

b) $M(G,x)\leq M(K_{d,d},x)^{n/(2d)}$.

Moreover, the maximum is achieved only by unions of copies of $K_{d,d}$.

{\bf Corollary 12.4} \cite{Davies2017}.
For any $d$-regular graph~$G$, 
$m(G)\leq m(K_{d,d})^{n/(2d)}$.

{\bf Corollary 12.5} \cite{Bregman,Davies2017}.
For any $d$-regular graph~$G$, 
$pm(G)\leq (d!)^{n/(2d)}$.

{\bf Theorem 12.2}~\cite{CsikvariMatch}.
Let $G$ be a $d$-regular bipartite graph on $2n$ vertices,
let $p = \frac{k}{n}$, and $p_{\mu}$ be the probability 
that a random variable with distribution Binomial$(n,p)$ 
takes its mean value $\mu=k$. Then 
$$
m_k(G)\geq p_{\mu}{n \choose k}^2\left(\frac{d-p}{d}\right)^{n(d-p)}(dp)^{np}.
$$

{\bf Corollary 12.6} \cite{CsikvariMatch,Schrijver}. 
Let $G$ be a $d$--regular bipartite graph on $2n$ vertices, then
$pm(G)\geq \left(\frac{(d-1)^{d-1}}{d^{d-2}}\right)^n$.

The~$\S13$ is actually survey section. In~$\S13.1$, we discuss 
real-rootedness of $PC$-polyno\-mials. 
In the paper of J. Brown and R. Nowakowski of 2005~\cite{BrownNowak},
it was stated that for almost all graphs $PC$-polynomial has a complex non-real root.
Unfortunately, their proof seems to be not complete. Thus, the next problem is still open.

{\bf Problem 13.1}.
To state if for almost all graphs $PC$-polynomial has a non-real root.

The subsection~$\S13.2$ is devoted to weaker conditions than real-rootedness.
A polyno\-mial $F(x) = \sum\limits_{i=0}^{n}a_i x^i$, $a_i\in\mathbb{R}$, 
is called {\it unimodal} if 
$a_0\leq \ldots \leq a_{k-1}\leq a_k\geq a_{k+1}\geq \ldots \geq a_n$
for some $k\in\{0,1,\ldots,n\}$ and {\it log-concave} if 
$a_i^2\geq a_{i-1}a_{i+1}$ for all $i=1,\ldots,n-1$.

In 1987, Y. Alavi, J. Malde, A. Schwenk, P. Erd\H{o}s proved that

{\bf Theorem 13.2} \cite{Schwenk}.
For every permutation $\pi\in S_{n}$ there exists a graph~$G$
with $\omega(G) = n$ such that
$c_{\pi(1)}<c_{\pi(2)}<\ldots <c_{\pi(n)}$.

The central problem of this area of research is

{\bf Conjecture 13.1} \cite{Schwenk}.
Let $\bar{G}$ be a tree. Then $C(G,x)$ is unimodal.

A graph $G$ is said to be {\it well-covered} if every maximal independent set of $G$ 
is also a maximum independent set.
In 2003, T. Michael and W. Traves stated that 

{\bf Theorem 13.4} \cite{Michael}.
Let $\bar{G}$ be a well-covered graph and $w = \omega(G)$. Then 
$c_1(G)\leq c_2(G)\leq \ldots\leq c_{\lceil w/2\rceil}(G)$.

In 2014, J. Cutler and L. Pebody proved

{\bf Theorem 13.5} \cite{Cutler}.
For every permutation $\pi$ of the set 
$\{\lceil w/2\rceil,\lceil w/2\rceil+1,\ldots,w\}$
there exists a well-covered graph~$G$ such that
$c_{\pi(\lceil w/2\rceil)}(G)
 <c_{\pi(\lceil w/2\rceil+1)}(G)
 {<}\ldots <c_{\pi(w)}(G)$
and $w = \omega(G)$.

{\bf Corollary 13.1} \cite{Cutler,Levit2006}.
Let $k\geq4$. There exists a well-covered 
graph $\bar{G}$ with $\omega(G) = k$ such that $C(G,x)$ is not unimodal.

A well-covered graph~$G$ is called a {\it very well-covered graph}~\cite{Favaron},
if $G$ contains no isolated vertices and $|V(G)| = 2\alpha(G)$.

In 2006, V. Levit and E. Mandrescu stated the following theorem.

{\bf Theorem 13.6} \cite{Levit2006very}.
Let $\bar{G}$ be a very well-covered graph,
$|V(G)|\geq2$, $w = \omega(G)$. Then 
a) $c_1(G)\leq c_2(G)\leq \ldots\leq c_{\lceil w/2\rceil}(G)$ and
$c_{\lceil (2w-1)/3\rceil}(G)\geq \ldots\geq c_{w-1}(G)\geq c_w(G)$;

b) $C(G,x)$ is unimodal for $w\leq 9$ and log-concave for $w\leq 5$.

In~$\S13.3$, the following notion is considered.
A class of graphs $\mathcal{K}$ is called $PC$-recogni\-zable
if for any two graphs $G,H\in\mathcal{K}$
the equality $PC(G,x) = PC(H,x)$ implies $G\cong H$.
In 1994, C. Hoede and X. Li formulated~\cite{Clique1994}
the question: Which classes of graphs are $PC$-recognizable?
D. Stevanovic in 1997 proved that 

{\bf Theorem 13.8} \cite{Stevanovic97}.
The class of threshold graphs is $PC$-recognizable.

In 2008, V.~Levit and E. Mandrescu stated 

{\bf Conjecture 13.2} \cite{Levit2008}.
Let $\bar{G}$ be a connected graph, $\bar{T}$ be a well-covered tree. If 
$PC(\bar{G},x) = PC(\bar{T},x)$, then $\bar{G}$ is a well-covered tree.

In~$\S13.4$, the bounds on roots of $PC$-polynomials (including of random graph)
are discussed. Let us list only few results.

{\bf Theorem 13.10} \cite{BrownNowak01}.
For any graph $G$ with $n$ vertices and clique number 
$w\geq2$, modulus of any root of $PC(G,x)$ is not less than 
$\big(\frac{w-1}{n}\big)^{w-1} + O(n^{-w})$.
This bound is tight. 

{\bf Theorem 13.11} \cite{BrownNowak00}.
Given a well-covered graph~$G$, modulus of any root of $PC(G,x)$ 
is not less than $1/\omega(G)$.

The subsection~$\S13.5$ involves a new graph polynomial. 
In 1987, R.-Y. Liu defined~\cite{Adjont} for a graph $G$ 
the {\it adjoint polynomial} $h(G,x)=\sum\limits_{k=1}^{n}(-1)^{n-k}a_k(G)x^k$,
where $n = |V(G)|$ and $a_k(G)$ denotes the number of ways one can cover 
all vertices of~$G$ by exactly $k$ disjoint cliques of $G$.
In 2017, F. Bencs stated the following result. 

{\bf Theorem 13.13} \cite{BencsAdjont}.
For any graph $G$ there exists a graph $\hat{G}$ such that
$h(G,x) = x^n I(\hat{G},1/x) = x^{n-w}PC(\overline{\hat{G}},x)$,
where $n = |V(G)|$ and $w = \omega(G)$.

This result allows to apply the theory of clique-type polynomials
and results of $\beta(G)$ for the adjoint polynomial.

In~$\S13.6$, we investigate the sum 
$I(G,-1) = 1 - s_1(G)+s_2(G)-\ldots +(-1)^{\alpha(G)}s_{\alpha(G)}(G)$ 
which was called in~\cite{Bousquet} as the {\it alternating number of independent sets}.

Given a graph~$G$, the {\it decycling number} $\varphi(G)$~\cite{Beineke}
is the minimum number of vertices that need to be removed in order 
to eliminate all cycles in~$G$.

The following result was initially proved by A. Engstr\"{o}m in 2009,
V.~Levit and E.~Mandrescu found its elementary proof in 2011.

{\bf Theorem 13.13} \cite{Engstroem,Levit2011}.
For any graph~$G$, $|I(G,-1)|\leq 2^{\varphi(G)}$.

J. Cutler and N. Kahl in 2016 proved that 

{\bf Theorem 13.14} \cite{CutlerKahl}.
Given a positive integer $k$ and an integer $q$ with $|q| \leq 2^k$,
there is a connected graph~$G$ with $\varphi(G) = k$ and $I(G,-1) = q$.

In~$\S13.7$, we discuss different problems about the extremal 
values of $\beta_{\pm}(n,k)$ and possible strategies to struggle with them.

\newpage
\section{Preliminaries on polynomials and sequences}

{\bf Pringsheim's Theorem} \cite[p.~240--242]{Flajolet}.
If $f(z)$ is representable at the origin by a~series 
expansion that has non-negative coefficients and 
radius of convergence $R$, then the point $z = R$ is a singularity of $f(z)$.

Let $f(z) = \sum\limits_{n=0}^\infty f_n z^n$ and $\Supp(f) = \{k \mid f_k\neq0\}$.
The sequence $(f_n)$, as well as $f(z)$, is said to admit a span $d$ 
if for some $r$, there holds
$\Supp(f) \subset r + d\mathbb{Z}_{\geq0} = \{r,r+d,r+2d,\ldots\}$.
The largest span, $p$, is the period, all other spans being divisors of $p$. 
If the period is equal to~1, then the sequence $(f_n)$ and the function $f(z)$ are said to be aperiodic.

Applying triangle inequality, it is easy to prove

{\bf Daffodil Lemma} \cite[p.~266--267]{Flajolet}.
Let $f(z)\in\mathbb{C}[z]$ be analytic in $|z| < R$ and have nonnegative
coefficients at~0. Assume that $f$ does not reduce to a monomial and that for
some nonzero nonpositive $s$ satisfying $|s| < R$, one has $|f(s)| = f (|s|)$. 
Then, the following hold: 

(i) $s = |s|e^{i\theta}$ with $\theta/2\pi = r/p \in \mathbb{Q}$ (an irreducible fraction) and $0 < r < p$;

(ii) $f$ admits $p$ as a span.

{\bf Fekete's Lemma \cite[Lemma~1.2.2]{Madras}}.
Let $\{a_n\}$, $n\geq1$, be a sequence of real numbers such that 
$a_{s+t}\leq a_s + a_t$ for all $s,t\in\mathbb{N}_{>0}$.
Then there exists a limit $\lim\limits_{n\to\infty}\frac{a_n}{n}$ 
and the limit equals $\inf\limits_{n\geq1}\frac{a_n}{n}$.

{\bf Enestr\"{o}m---Kakeya Theorem \cite{Kakeya}}.
Let $f(x) = \sum\limits_{j=0}^n a_j x^j\in\mathbb{R}[x]$
be any polynomial with $a_j>0$ for all $j=0,\ldots,n$.
Then for any zero $\rho$ of $f(x)$, we have 
$$
\min\limits_{j=0,\ldots,n-1}\{|a_j/a_{j+1}|\}
\leq|\rho|\leq\max\limits_{j=0,\ldots,n-1}\{|a_j/a_{j+1}|\}.
$$

\textbf{Samuelson's Inequality \cite{Samuelson}}.
Let all roots of a polynomial 
$f(x) = \sum\limits_{j=0}^n a_j x^j\in\mathbb{R}[x]$
be real. Then all roots of $f(x)$ lie in the segment $[x_-;x_+]$, where
$$
x_{\pm} = -\frac{a_{n-1}}{na_n}\pm\frac{n-1}{na_n}\sqrt{a_{n-1}^2-\frac{2n}{n-1}a_n a_{n-2}}.
$$

{\bf Ostrowski's Theorem \cite{Ostrowski-type}}.
Let $f(z),g(z)$ be two monic polynomials of degree $n$ with complex coefficients:
$$
f(z) = z^n + a_1z^{n-1}+\ldots + a_n, \quad g(z) = z^n + b_1z^{n-1}+\ldots + b_n.
$$
Then the roots of $f$ and $g$ can be enumerated as 
$\alpha_1,\ldots,\alpha_n$, $\beta_1,\ldots,\beta_n$
respectively in such a way that 
$$
\max\limits_{k=1,\ldots,n}\{|\alpha_k-\beta_k|\}
\leq 4\cdot 2^{-1/n}\left(\sum\limits_{i=1}^n|a_i-b_i|\gamma^{n-i}\right)^{1/n},
$$
where $\gamma = 2\max\limits_{k=1,\ldots,n}\{|a_k|^{1/k},|b_k|^{1/k}\}$.

A matrix~$M\in M_n(\mathbb{R})$ is primitive if all entries of $M$ 
are nonnegative and for some~$t$, all entries of $M^t$ are strictly positive.

{\bf Boyle---Handelman Theorem \cite{Boyle,KOR}}.
Let $\Lambda = (\lambda_1,\lambda_2,\ldots,\lambda_d)$
be a tuple of nonzero complex numbers and let $S = \mathbb{Z}$ or $S = \mathbb{R}$. Also,
$$
s(\Lambda,n) = \sum\limits_{i=1}^d \lambda_i^n,\quad
t(\Lambda,n) = \sum\limits_{k|n}\mu(n/k)s(\Lambda,k).
$$
There is a primitive matrix~$M$ with entries from~$S$ 
and with characteristic polynomial 
$\chi_M(x) = x^m\prod\limits_{i=1}^d (x - \lambda_i)$
for some $m\geq 0$ if and only if 

(1) the coefficients of the polynomial $\chi_M(x)$ belong to~$S$,

(2)  there exists $\lambda_j\in\Lambda$ such that $\lambda_j > |\lambda_i|$ for all $i\neq j$,

(3) if $S = \mathbb{Z}$, then $t(\Lambda,n)\geq0$ for all $n \geq 1$,

(3$'$) if $S = \mathbb{R}$, then for all $n \geq 1$,
$s(\Lambda,n)\geq0$ and for all $k\geq1$, 
$s(\Lambda,n)>0$ implies $s(\Lambda,nk)>0$.

\newpage
\section{$PC$-polynomial of graph}

Let $X$ be a finite alphabet, $G = G(V,E)$ be a simple graph whose vertices
are in one-to-one correspondence with the elements of $X$.

Denote by $M(X,G)$ a partially commutative monoid
$M\langle X\mid ab = ba,\,(a,b)\in E\rangle$.
We may regard elements of $M(X,G)$ as words in the alphabet $X$
and two words $x$ and $y$ are equal if there exists a sequence
$z_1,z_2,\ldots,z_k$ of words such that $x = z_1$, $y = z_k$, 
and for all~$i$, $1\leq i < k$, there exist words 
$z_i'$, $z_i''$ and letters $a_i$, $b_i$ satisfying:
$$
z_i = z_i' a_ib_i z_i'',\quad  
z_{i+1} = z_i' b_ia_i z_i'',\quad (a_i,b_i)\in E.
$$

{\bf Lemma 1.1} \cite{Diekert}. Let $u,v\in X^*$. 

a) If $au = av$ in $M(X,G)$ for $a\in X$, then $u = v$.

b) If $su = sv$ in $M(X,G)$ for $s\in X^*$, then $u = v$.

{\sc Proof}.
a) Let $av$ to be obtained from $au$ by a finite sequence 
of arrangements in a pair of neighbour letters.
If the initial letter~$a$ does not participate in all such arrangements, then 
we have the equality $u = v$ by the same sequence of arrangements.
Otherwise, call all arrangements in which initial~$a$ participates
as 0-arrangements and all other ones as 1-arrangements.
By the sequence of all 1-arrangements with the preserved order
we also have $u = v$. 

b) It follows from a).

We assume that $X$ is totally ordered, and we consider the corresponding
lexicographi\-cal ordering $<$ on $X^*$.
Given a word $w\in M(X,G)$, its normal form $[w]$ is defined as a~maximal 
word among the set of all words in $M(X,G)$ which are equal to $w$.
As the lexicographic ordering is a total order and $|X|$ is finite, 
each word has a unique normal form.

{\bf Lemma 1.2} \cite{Diekert}.
A word $w\in M(X,G)$ is in normal form if and only if 
for any equality $w = xaybz$, $x,y,z\in X^*$, $a,b\in X$, $a < b$,
the inclusion $(a,b)\in E$ implies that there exists a letter~$t$ of $y$ such that $(t,b)\not\in E$.

{\sc Proof}.
It is easy to check that any word in normal form satisfies the conditions of Lemma.

Let a word $w$ satisfy the conditions of Lemma.
To the contrary, suppose that there exists a word $u$ equal to $w$ such that $u > w$.
Denote $w = xay$, $u = xbz$ for some $x,y,z\in X^*$, $a,b \in X$, $a < b$. 
By Lemma~1.1, we have $ay = bz$.
Let us find the first occurrences of the letter~$b$ in $y$: $y = y'by''$.
Then $ay'by'' = bz$. 
There exists a sequence of arrangements converting $ay'by''$ to $bz$, 
thus, $b$ is adjacent with all letters of $ay'$ in~$G$.
We arrive at a contradiction to the condition holding for $w = xay'by''$.

{\bf Corollary 1.1} \cite{Diekert}.
Let $A\subset X^*$ be a set consisting of all normal forms of words in $M(X,G)$.
Then $A$ is a regular language which is closed under taking subwords.

Let $M_n = M_n(X,G)$ denote the set of all normal forms of words 
of length $n$ in $M(X,G)$ and $m_n = |M_n|$. Assume that $M_1 = \{1\}$.

A clique is a subset of vertices of a graph such that every 
two distinct vertices in the clique are adjacent.
The number of vertices of a clique is called a size of the clique.
A~maximum clique of a graph~$H$ is a clique of maximum possible size for $H$.
The size of the maximum clique of~$H$ is called a clique number of~$H$
and is denoted by $\omega(H)$. 

Let $t_0 = \omega(G)$. Denote by $c_i(G)$ the number of all cliques of size $i$ in $G$. 
In particular, $c_1(G) = |V|$, $c_2(G) = |E|$. Put $c_0(G) = 1$.

The following Theorem~1.1, which is the main result of $\S1$, 
was actually stated by P.~Cartier and D.~Foata in 1969~\cite{Cartier69}.
The precise form of the result was written down by D.C. Fisher in 1989~\cite{Fisher1989-1}.
Below the new proof of this result is proposed. 
Of course, the proofs of P. Cartier, D. Foata and D.C. Fisher as well as 
the one presented below are very similar in their ideas.

{\bf Theorem 1.1} \cite{Cartier69,Fisher1989-1}.
The numbers $m_n$, $n\geq1$, satisfy the reccurence relation
\begin{equation}\label{rekurrent}
m_n = c_1(G)m_{n-1}-c_2(G)m_{n-2} + \ldots + (-1)^{t_0+1}c_{t_0}(G)m_{n-t_0}
\end{equation}
with initial data $m_0 = 1$, $m_{-1} = \ldots = m_{-t_0+1} = 0$.

{\sc Proof}.
Define sets 
\begin{gather*}
N_1 = \{(u, x) \mid u\in M_{n-1}, x \in X\}, \\
N_2 = \{(u, x_1, x_2) \mid u \in M_{n-2}, x_1,x_2\in X, x_1 > x_2,(x_1,x_2)\in E\}, \\
\ldots
\end{gather*}
\vspace{-1cm}
\begin{multline*}
N_{t_0} = \{(u, x_1,\ldots,x_{t_0}) \mid u \in M_{n-t_0}, \\
x_1,\ldots, x_{t_0}\in X, x_1 > x_2 > \ldots > x_{t_0},(x_i,x_j)\in E, i,j=1,\ldots,t_0\}
\end{multline*}
and maps 
$$
f_i\colon N_i \to M_n,\quad (u, x_1,\ldots,x_i) \to [ux_1\ldots x_i],\ i=1,\ldots,t_0.
$$

Given a word $w \in M_n$, define the numbers
$k_i(w) = |f^{-1}_i[\{w\}]|$, $i=1,\ldots,t_0$, i.e., 
the cardinality of the preimage of $\{w\}$ under $f_i$.
Also, define $k_0(w) = \sum\limits_{i=1}^{t_0}(-1)^{i+1}k_i(w)$.

Applying the equalities 
$c_i(G)m_{n-i} = |N_i| = \sum\limits_{w\in M_{n}}k_i(w)$,
conclude
\begin{multline}\label{rek:right-part}
c_1(G)m_{n-1}- c_2(G)m_{n-2} + \ldots + (-1)^{t_0+1}c_{t_0}(G)m_{n-t_0} \\
= |N_1| - |N_2| + \ldots + (-1)^{t_0+1}|N_{t_0}| \\
= \sum\limits_{i=1}^{t_0}(-1)^{i+1}\sum\limits_{w\in M_n}k_i(w) = \sum\limits_{w\in M_n}k_0(w).
\end{multline}

If $k_0(w) = 1$ for all $w \in M_n$, then the RHS of~\eqref{rekurrent}
equals $\sum\limits_{w\in M_n}k_0(w) = \sum\limits_{w\in M_n}1 = m_n$
by~\eqref{rek:right-part}, we are done. 
Thus, it remains to prove 

{\bf Lemma 1.3}.
For $w \in M_n$, $n\geq1$, we have $k_0(w) = 1$.

{\sc Proof}.
Let $(u, x_1, \ldots, x_i)\in f^{-1}_i[\{w\}]$, $i\in \{1,\ldots,t_0\}$, $u \in M_{n-i}$.
By Lemma~1.2, the normal form $w$ of the word $ux_1 \ldots x_i$
could be obtained by inserting $x_1,\ldots,x_i$ somewhere in $u$
and preserving in $w$ the order of all letters from~$u$.
Moreover, by Lemma~1.2, we have the next remark: 
if a letter $x_s$, $s\in \{1,\ldots,i\}$,
stays in the $k$-th position in the word $w = w_1w_2\ldots w_n$, $w_i\in X$,
then $(x_s,w_l)\in E$ for all $l = k+1,\ldots,n$.
Therefore, $k_i(w)$ equals a number of tuples of $i$ letters 
$w_{p_1}$, \ldots, $w_{p_i}$, $p_1<p_2<\ldots<p_i$, 
of the word $w$ such that they are pairwise adjacent 
and they are adjacent to all letters $w_s$ with $p_i<s$.

Let $k\in\{1,\ldots,n\}$ be a minimal number such that 
$(w_k,w_l)\in E$ for $l = k + 1, \ldots,n$. Define the set
$B_i = \{(u,y_1,y_2,\ldots,y_i)\in N_i\mid w_k\in\{y_1,\ldots,y_i\}\}$
and the numbers 
\begin{equation}\label{ka-itye}
k_i'(w)  = \big|f^{-1}_i[\{w\}]\cap B_i|,\quad
k_i''(w) = \big|f^{-1}_i[\{w\}]\cap (N_i\setminus B_i)|.
\end{equation}
Note that 
\begin{equation}\label{ka-itye2}
k_i = k_i'+k_i'',\ i=\overline{1,t_0},\quad k_1'(w) = 1,\quad k_{t_0}''(w)=0, \quad
k_i'(w) = k_{i-1}''(w),\ i=\overline{2,t_0}.
\end{equation}

Due to \eqref{ka-itye}, \eqref{ka-itye2}, 
we finish the proof of Lemma~1.3 by the following calculations
\begin{multline}
k_0(w) =
\sum\limits_{i=1}^{t_0}(-1)^{i+1}k_i(w)
= \sum\limits_{i=1}^{t_0}(-1)^{i+1}k_i'(w)+
\sum\limits_{i=1}^{t_0}(-1)^{i+1}k_i''(w) \\
= k_1'(w) + \sum\limits_{i=2}^{t_0}(-1)^{i+1}k_i'(w)
+\sum\limits_{i=1}^{t_0-1}(-1)^{i+1}k_i''(w) \\
= 1 + \sum\limits_{i=1}^{t_0-1}(-1)^{i+1}(
k_i''(w)-k_i''(w)) = 1.
\end{multline}

Define $PC$-polynomial on graph $G$ as a characteristic polynomial
of the recurrence relation~\eqref{rekurrent}:
\begin{equation}\label{PC-polynomial}
PC(G,x) = x^{t_0}-c_1(G)x^{t_0-1}+c_2(G)x^{t_0-2}
+\ldots+(-1)^{t_0-1}c_{t_0-1}(G)x+(-1)^{t_0}c_{t_0}(G),
\end{equation}
where $t_0 = \omega(G)$.

{\bf Example 1.1} \cite{Fisher1990}.
We have 

a) $PC(\bar{K}_n,x) = x-n$ for empty graph $\bar{K}_n$;

b) $PC(K_n,x) = (x-1)^n$ for complete graph $K_n$;

c) $PC(K_{q_1,\ldots,q_s}) = (x-q_1)\ldots(x-q_s)$
for complete multipartite graph with parts $q_1,\ldots,q_s$,
in particular, $PC(K_{n/p,\ldots,n/p}) = \big(x-\frac{n}p\big)^p$
for complete multipartite graph with $p$ equal parts;

d) $PC(T_n,x) = x^2-nx+n-1 = (x-1)(x-n+1)$ for a tree $T_n$ with $n$ vertices.

Note that 
\begin{equation}\label{PC-Dep}
x^{t_0}PC(G,1/x) = D(G,x),
\end{equation}
where 
$D(G,x) = \sum\limits_{k=0}^{\omega(G)}(-1)^k c_k(G)x^k$ 
is the dependence polynomial of~$G$.
Hence, the sets of roots of polynomials $PC(G,x)$ and $D(G,x)$
respectively are in one-to-one correspondence by the rule 
$x\leftrightarrow 1/x$,
as zero is not root of neither $PC(G,x)$ nor $D(G,x)$.

{\bf Lemma 1.4} \cite{Cartier69,Fisher1990}.
The function $1/D(G,x)$ is a generating function for the sequence $\{m_n\}$, i.e.,
\begin{equation}\label{GF-PC}
\frac{1}{D(G,x)} = \sum\limits_{n\geq0}m_n x^n.
\end{equation}

{\sc Proof}.
The statement of Lemma is equivalent to the equality
$$
1 = (1-c_1(G)x+c_2(G)x^2-\ldots+(-1)^{t_0}c_{t_0}(G)x^{t_0})
 \bigg(\sum\limits_{n\geq0}m_nx^n\bigg)
$$
holding by Theorem~1.1.
\eject

By $G[A]$ for $A\subset V(G)$
we denote the subgraph of $G$ induced by the set of vertices~$A$.
Let $N(v)$ be a neighbourhood of a vertex $v\in V(G)$, i.e.,
$N(v) = \{u\in V\mid (v,u)\in E\}$.
Given two graphs $G_1,G_2$ with $V(G_1)\cap V(G_2) = \emptyset$,
the graphs $G_1\cup G_2$ and $G_1 + G_2$ are defined as follows:
$V(G_1\cup G_2) = V(G_1+G_2) = V(G_1)\cup V(G_2)$,
$E(G_1\cup G_2) = E(G_1)\cup E(G_2)$,
$E(G_1+G_2) = E(G_1)\cup E(G_2)\cup \{(u,v)\mid u\in E(G_1),v\in E(G_2)\}$.

{\bf Lemma 1.5} \cite{Gutman1990,Clique1994}.
Given a graph $G$, we have 

a) $D(G,x) = D(G \setminus v,x) - xD(G[N(v)],x)$, $v\in V(G)$;

b) $D(G,x) = D(G\setminus uv,x) + x^2D(G[N(u)\cap N(v)],x)$, $(u,v)\in E(G)$;

c) $D'(G,x) = -\sum\limits_{v\in V(G)}D(G[N(v)],x)$;

d) $D(G_1\cup G_2) = D(G_1) + D(G_2) - 1$;

e) $D(G_1 + G_2) = D(G_1)D(G_2)$.

{\sc Proof}.
a) Let $v\in V(G)$. All cliques in $G$ could be divided into 
two types, the first ones contain $v$ but not the second ones.
So, all cliques of the first type are counted in a~summand $xD(G[N(v)],x)$
and all cliques of the second type --- in $D(G \setminus v,x)$.

b) All cliques in $G$ are of two types, 
the first ones contain the edge $(u,v)$ but not the second ones.
So, all cliques of the first type are counted in $x^2D(G[N(u)\cap N(v)],x)$
and all cliques of the second type --- in $D(G\setminus uv,x)$.

c) Prove the statement by induction on $n = |V(G)|$.
For $G = \{v\}$, we have 
$D(G,x) = 1 - x$ and $D'(G,x) = -1 = - D(\emptyset,x)$.

Differentiating the equality from a) and further applying 
the same equality and the induction hypothesis, we have 
\begin{multline*}\allowdisplaybreaks
D'(G, x) = D'(G \setminus v, x) - D(G[N(v)],x) - xD'(G[N(v)], x) \\
= - D(G[N(v)],x) + \sum\limits_{u\in V(G\setminus v))}D(G\setminus v[N(u)],x)
+ x\sum\limits_{u\in N(v)}D(G[N(v)\cap N(u)],x) \\
= - D(G[N(v)],x)
- \sum\limits_{u\in V(G\setminus N(v))}D(G\setminus v[N(u)],x) \\
-\sum\limits_{u\in N(v)}
(D(G \setminus v[N(u)],x) - xD(G[N(v)\cap N(u)],x)) \\
= - D(G[N(v)],x)
- \sum\limits_{u\in V(G\setminus N(v))}D(G[N(u)],x)
- \sum\limits_{u\in N(v)}D(G[N(u)],x) \\
= \sum\limits_{w\in V(G)}D(G[N(w)],x).
\end{multline*}

d) It follows from the fact that any clique in $G_1\cup G_2$
is a clique in either $G_1$ or $G_2$.

e) We conclude that $D(G_1 + G_2) = D(G_1)D(G_2)$ 
by the following equality 
$$
c_p(G_1+G_2)
 = \sum\limits_{k=0}^p c_k(G_1)c_{p-k}(G_2).
$$

\newpage
\section{Largest root of $PC$-polynomial}

Let $z_0(G)$ be a largest on modulus (complex) root of $PC(G,x)$,
$\beta(G) = |z_0(G)|$.

{\bf Lemma 2.1} \cite{Fisher1990}.
The number $\beta(G)$ is a root of $PC(G,x)$.

{\sc Proof}.
Let $s = 1/\beta(G)$.
By~\eqref{GF-PC} and minimality of modulus of the root $1/z_0(G)$ of $D(G,x)$,
we get that the convergence radius of $1/D(G,z)$ equals~$s$.
By Pringsheim's Theorem, the function $1/D(G,z)$ is not analytic in~$s$.
Thus, we have $D(G,s) = 0$ and $PC(G,\beta(G)) = 0$.

{\bf Lemma 2.2} \cite{Csikvari}.
Given a graph~$G$ and proper induced subgraph~$H$ of~$G$,
the rational function $D(H,x)/D(G,x)$ is representable 
at the origin by a series expansion that 
has positive integer coefficients.

{\sc Proof}.
Let us prove the statement by induction on $n = |V(G)|$.
For $n = 1$, the empty graph if the only proper induced subgraph. So, 
$$
\frac{D(H,x)}{D(G,x)} = \frac{1}{1-x} = \sum\limits_{k\geq0}x^k.
$$

It is enough to prove the statement for $H = G\setminus v$.
Indeed, let $H$ be a proper induced subgraph of~$G$, i.e.,
$H = G\setminus\{v_1,\ldots,v_l\}$. Then
$$
\frac{D(H,x)}{D(G,x)}
 = \frac{D(G\setminus v_1,x)}{D(G,x)} 
    \frac{D(G\setminus \{v_1,v_2\},x)}{D(G\setminus v_1,x)} \ldots 
    \frac{D(G\setminus\{v_1,\ldots,v_l\},x)}{D(G\setminus\{v_1,\ldots,v_{l-1}\},x)},
$$
where by induction all factors except the first one 
have a series expansion with positive coefficients. 

By Lemma~1.5a, 
\begin{multline*}
\frac{D(G \setminus v,x)}{D(G,x)}
 = \frac{D(G \setminus v,x)}{D(G \setminus v,x) - xD(G[N(v)],x)}  \\
 = \frac{1}{1 - x\frac{D(G[N(v)],x)}{D(G \setminus v,x)}}
 = \sum\limits_{k\geq0}\left(x\frac{D(G[N(v)],x)}{D(G \setminus v,x)}\right)^k.
\end{multline*}
If $G[N(v)]$ is a proper subgraph of~$G\setminus v$, then
by induction $\frac{D(G[N(v)],x)}{D(G \setminus v,x)}$
is representable by a series expansion with positive coefficients.
Otherwise, $G[N(v)] = G\setminus v$ and 
$\frac{D(G \setminus v,x)}{D(G,x)} = \frac{1}{1-x}$.
Lemma is proved.

{\bf Remark 2.1}.
In \cite{Csikvari}, P. Csikv\'{a}ri wondered if some interpretation 
of the series expansion of $D(H,x)/D(G,x)$, where $H$ is a subgraph of $G$, could be found. 
In 2006, C.~Krattenthaler proved the following result~\cite{Athreya,Krattenthaler}.
Let $w = w_1w_2\ldots w_n$ be a word from $M_n(G,X)$ (in the normal form).
Call the set 
$\{w_k\mid (w_k,w_l)\in E(G),\,l = k + 1,\ldots,n\}$ as the type $\tau(w)$ of~$w$.
So, the generating function of all words from $M(X,G)$
which type is a subset of $V(G)\setminus V(H)$
is equal to the ratio~$D(H,x)/D(G,x)$ \cite{Athreya,Krattenthaler}. 
Also, this result implies Lemma~2.2.

Let $G$ be a supergraph of a complete multipartite graph with parts 
$H_1,H_2,\ldots,H_k$, $k\geq2$, 
(i.e., $G = H_1 + H_2 + \ldots + H_k$)
such that $\bar{H}_1,\ldots,\bar{H}_k$, are connected.
By Lemma~5e, $PC(G,x) = \prod\limits_{i=1}^k PC(H_i,x)$.
A part $H_i$ is called a maximal one if $\beta(G)$ is a root of $PC(H_i,x)$.
A maximal part can be not unique, e.g., for $G = H + H$, where $\bar{H}$ is connected,
both parts are maximal.

{\bf Lemma 2.3} \cite{Csikvari,Hajiabolhassan}.
Let $G$ be a graph.

a) For any induced subgraph~$H$ of~$G$, we have $\beta(H)\leq \beta(G)$.

b) For any proper induced subgraph~$H$ of~$G$, we have $\beta(H) = \beta(G)$
if and only if $G$ is a supergraph of a complete multipartite graph 
with parts having connected complements and $H$~contains a maximal part of~$G$.

{\sc Proof}.
a) Prove the statement by induction on $n = |V(G)|$. 
For $n = 1,2$ it is true.
For inductive step, it is enough to show that 
$\beta(G\setminus v)\leq \beta(G)$ for all $v\in V(G)$.
Applying Lemma~1.5a for $\alpha = 1/\beta(G\setminus v,x)$
and the induction hypothesis, calculate
$$
D(G,\alpha) = D(G \setminus v,\alpha) - \alpha D(G[N(v)],\alpha)
 = -\alpha D(G[N(v)],\alpha)\leq0,
$$
so, $D(G,x)$ has a root in $(0,\alpha]$ as $D(G,0) = 1$
and $\beta(G)\geq \beta(G\setminus v)$.

b) It is easy to see that if $G$ and $H$ satisfy the conditions, then $\beta(H) = \beta(G)$.

Suppose that $\beta(H) = \beta(G)$. We prove the statement by induction on $n = |V(G)|$.
For $n = 2$ it is true. Consider the case when $G$ is such a graph that $\bar{G}$ is connected.
If $G$ is a supergraph of a complete multipartite graph 
with parts having connected complements, then by Lemma~1.5e 
and induction, we reduce such case to the one stated above.

Let us prove that $\beta(G)>\beta(G\setminus v)$ for all $v\in V(G)$,
then the statement for any proper induced subgraph~$H$ follows from a).
Note that the connectedness of $\bar{G}$ provides that 
$v$ is not connected to all vertices of $G$, i.e., 
$G[N(v)]$ is a proper induced subgraph in $G\setminus v$.

Assume that $\beta(G) = \beta(G\setminus v)$.
Then by Lemma~1.5a, $\beta(G\setminus v) = \beta(G[N(v)])$.
If the complement of $G\setminus v$ is connected,
then by the induction hypothesis
$\beta(G\setminus v)>\beta(G[N(v)])$, a contradiction.

Assume that the complement of $G\setminus v$ is disconnected,
then from the equality $\beta(G\setminus v) = \beta(G[N(v)])$
by the induction hypothesis we conclude that 
$G\setminus v$ is a supergraph of a complete multipartite graph
with parts $H_1,\ldots,H_k$, $\bar{H}_i$, is connected, 
and $G[N(v)]$ is an induced subgraph of $G\setminus v$ 
containing a maximal part $H_{i_0}$. 
Therefore, $\bar{G}$ is disconnected as the part $H_{i_0}$ 
has no any edges to all other vertices of $\bar{G}$.
Lemma is proved.

The following Theorem~2.1, the main result of the current paragraph, 
was initially formulated but only partly proved by D.C. Fisher and A.E. Solow in 1990~\cite{Fisher1990}.
Actually they proved Lemma~2.1.
In 2000, M. Goldwurm and M. Santini proved Theorem~2.1~\cite{Santini2000} 
with the help of Perron---Frobenius theory. 
We consider the proof of P. Csikv\'{a}ri stated in 2013~\cite{Csikvari}.

{\bf Theorem 2.1} \cite{Csikvari,Santini2000}.
The number $\beta(G)$ is the only complex root of $PC(G,x)$ 
with modulus greater or equal to $\beta(G)$.

{\sc Proof}.
Let us prove that $1/\beta(G)$ is the only complex root of $D(G,x)$ equal to $1/\beta(G)$.

If $G$ is a supergraph of a complete multipartite graph with parts $H_1,\ldots,H_k$,
$\bar{H}_i$ is connected, then by Lemma~1.5e $D(G,x) = \prod\limits_{i=1}^k D(H_i,x)$.
Thus, we may assume that $\bar{G}$ is connected.
Also, $|V(G)|\geq2$, as for $|V(G)| = 1$ the statement is trivial.

For any $v\in V(G)$, $G[N(v)]$ is a proper induced subgraph of $G\setminus v$.
Otherwise $v$ is connected to all vertices of $V(G)\setminus v$
and so $\bar{G}$ is disconnected.

Consider the functions
$$
g(x)
 = \frac{D(G\setminus v,x)}{D(G,x)}
 = \frac{1}{1 - x\frac{D(G[N(v)],x)}{D(G\setminus v,x)}},\quad
f(x) = x\frac{D(G[N(v)],x)}{D(G\setminus v,x)}.
$$

By Lemma~2.3a, $\beta(G[N(v)])\leq\beta(G\setminus v)$
and the convergence radius of the series expansion of $f(x)$ 
is greater or equal to $1/\beta(G\setminus v)$.
Let $\rho$ be a complex root of $D(G,x)$
with modulus $1/\beta(G)$, $\rho\neq 1/\beta(G)$.
Then by Lemma~2.3b, $\rho$ is a pole of $g(x)$.
Hence, $f(\rho) = 1$ and $|f(\rho)| = f(\rho)$.
By Daffodil Lemma, $\rho = (1/\beta(G))e^{2\pi ir/p}$, where $p$ is a span of $f(x)$.
By Lemma~2.2, all coefficients of the series expansion of $f(x)$
are positive, thus, $f(x)$ is aperiodic. We get a contradiction. Theorem is proved.

{\bf Example 2.1}~\cite{Fisher1990}.
a) $\beta(K_n) = 1$,

b) $\beta(\bar{K}_n) = n$,

c) $\beta(K_{n_1,\ldots,n_p}) = \max\{n_1,\ldots,n_p\}$,

d) $\beta(T_n) = n-1$ for a tree $T_n$ with $n$ vertices.

Let $\{a_n\}_{n\geq0}$ be a sequence of nonnegative real numbers.
If the limit $a = \lim\limits_{n\to\infty}\sqrt[n]{a_n}$ exists, 
then $a$~is called the (exponential) growth rate of the sequence $\{a_n\}$.

Given a graph $G(V,E)$, an associative algebra 
$\Ass(X,G) = \Ass\langle X\mid ab = ba,(a,b)\in E\rangle$ for $X = V$
is called a partially commutative algebra.
Denote by $a_n$ the dimension of homogeneous space of all words 
of length~$n$ in the alphabet~$X$ in $\Ass(X,G)$.

The growth rate of $M(X,G)$ and $\Ass(X,G)$ 
is called the growth rate of the sequence $\{m_n\}$ and $\{a_n\}$ respectively.

{\bf Corollary 2.1} \cite{Santini2000}.
a) The growth rate of partially commutative monoid $M(X,G)$ equals~$\beta(G)$.

b) The growth rate of partially commutative associative algebra $\Ass(X,G)$ equals~$\beta(G)$\!.

{\sc Proof}.
a) Prove an additional statement. 
Let $\alpha_1,\ldots,\alpha_k$ be roots of $PC(G,x)$ 
with multiplicities $s_1,\ldots,s_k$ respectively. Then 
\begin{equation}\label{Growth}
m_n = \sum\limits_{j=1}^k P_j(n)\alpha_j^n,
\end{equation}
where $P_j(x)$ is a polynomial of degree $s_j-1$.

Indeed, consider a full partial fraction expansion for the rational function $1/D(G,x)$:
\begin{equation}\label{frac-decomp}
\frac{1}{D(G,x)} = \sum\limits_{j=1}^k\sum\limits_{r=1}^{s_j}\frac{c_{j,r}}{(1-\alpha_j x)^r},\quad c_{j,r}\in\mathbb{R},
\end{equation}
where we have $c_{j,s_j}\neq0$.
Continuing on dealing with \eqref{frac-decomp}, we get
\begin{multline*}
\frac{1}{D(G,x)}
 = \sum\limits_{j=1}^k\sum\limits_{r=1}^{s_j}c_{j,r}
 \sum\limits_{n\geq0} \binom{n+r-1}{r-1}\alpha_j^n \\
 = \sum\limits_{n\geq0}\sum\limits_{j=1}^k \left(c_{j,1}\binom{n}{0} 
 + \ldots + c_{j,s_j}\binom{n+s_j-1}{s_j-1}\right)\alpha_j^n,
\end{multline*}
the expressions in brackets are exactly the polynomials $P_j(n)$ of degree $s_j-1$.

Let $\beta(G) = |\alpha_1|>|\alpha_2|\geq\ldots\geq|\alpha_k|$,
then we have an asymptotics 
$m_n\sim C_0 n^{s_1-1}\alpha_1^n$.
Moreover, the real constant $C_0$ is positive as $m_n$ are nonnegative. Therefore,
$$
\lim\limits_{n\to\infty}\sqrt[n]{m_n}
 = \lim\limits_{n\to\infty}\sqrt[n]{C_0 n^{s_1-1}\alpha_1^n}
 = \alpha_1\lim\limits_{n\to\infty}\sqrt[n]{C_0 n^{s_1-1}}
  = \alpha_1 = \beta(G).
$$

b) It follows from a) as $a_n = m_n$ for $n\geq1$.

{\bf Remark 2.2}.
It is easy to show the existence of the limit $\lim\limits_{n\to\infty}\sqrt[n]{m_n}$ directly.
Indeed, by Corollary~1.1 the set of all normal forms of words in $M(X,G)$ 
is closed under taking subwords. Thus, we have the inequality $m_{s+t}\leq m_s m_t$ for all $s,t\in\mathbb{N}$.
It remains to apply Fekete's Lemma for the sequence $\{\ln m_n\}$.

{\bf Lemma 2.4} \cite{Csikvari,Hajiabolhassan}.
Let $G$ be a graph.

a) For any spanning subgraph~$H$ of~$G$, we have $\beta(G)\leq \beta(H)$.

b) For any proper spanning subgraph~$H$ of~$G$, we have $\beta(H) < \beta(G)$
provided that $\bar{G}$ is connected.

c) For any spanning subgraph~$H$ of~$G$, we have $\beta(H) = \beta(G)$
if and only if $G$ is a supergraph of a complete multipartite graph 
with parts $H_1,\ldots,H_k$ having connected complements, moreover, 
if $H_i$ is a maximal part of~$G$, then $G[V(H_i)] = H[V(H_i)]$,
$\bar{H}_i$ is a connected component in $\bar{H}$ and $H_i$ is a maximal part of $H$
(the decomposition of $\bar{H}$ into connected components could be not exactly
$\bar{H}_1,\ldots,\bar{H}_k$, some of them could be joined).

{\sc Proof}.
a) By the condition, we have inequalities $m_n(G)\leq m_n(H)$, $n\geq1$. Hence,
$$
\beta(G) = \lim\limits_{n\to\infty}\sqrt[n]{m_n(G)}
\leq \lim\limits_{n\to\infty}\sqrt[n]{m_n(H)} = \beta(H).
$$

b) Let us prove that $\beta(G\setminus uv) > \beta(G)$ for any pair of $u,v\in V(G)$.
To the contrary, there exists a pair of vertices $u,v\in V(G)$ such that
$\beta(G\setminus uv) = \beta(G)$. Then by Lemma~1.5b
we have $\beta(G[N(u)\cap N(v)])\geq\beta(G)$.
As $\bar{G}$ is connected, the graph $G[N(u)\cap N(v)]$ is a proper induced subgraph of~$G$.
By Lemma~2.3b, $\beta(G[N(u)\cap N(v)])<\beta(G)$, a contradiction.

c) If $G,H$ satisfy the conditions, then $\beta(H) {=} \beta(G)$.
Conversely, it follows from a), b).

{\bf Lemma 2.5} \cite{Csikvari}.
a) If $G$ is not a supergraph of a complete multipartite graph, 
then $\beta(G)$ is a simple root of $PC(G,x)$.

b) If $G$ is a supergraph of a complete multipartite graph 
with parts having connected complements, then 
the multiplicity $k$ of the root $\beta(G)$ of $PC(G,x)$
equals the number of maximal parts of $G$.

{\sc Proof}.
a) Apply Lemma~1.5c: $D'(G,x) = -\sum\limits_{v\in V(G)}D(G[N(v)],x)$.
By Lemma~2.3b, we have $D(G[N(v)],1/\beta(G)) < 0$
for any $v\in V(G)$. Thus, $D'(G,1/\beta(G))<0$ and 
$1/\beta(G)$ is a simple root of $D(G,x)$.

b) It follows from Lemma~1.5e and a).

{\bf Corollary 2.2}.
a) For almost all graphs, the largest root of $PC$-polynomial is simple.

b) The largest root of either $PC(G,x)$ or $PC(\bar{G},x)$ is simple.

{\sc Proof}.
a) It follows from Lemma~2.5 and connectedness of almost all graphs~\cite{Bolobas}.

b) At least one of $G$ and $\bar{G}$ is connected and we are done.

{\bf Lemma 2.6} \cite{Fisher1990,McMullen}.
a) Given a graph~$G$ with $n = |V(G)|$, we have 
\begin{equation}\label{MaxRootBounds}
\frac{n}{\omega(G)}\leq \beta(G)\leq n.
\end{equation}
Moreover, the lower bound is reached if and only if 
$n/\omega(G)\in\mathbb{N}$ and $G$ is the complete multipartite graph 
with all parts having $n/\omega(G)$ vertices. 
The upper bound is reached if and only if $G$ is the empty graph.

b) If $G$ is not complete graph and $n>1$, then $2\leq\beta(G)$.

{\sc Proof}.
a) Assume that $\beta(G) < n/\omega(G)$, then all real roots and real parts 
of all complex non-real roots of $PC(G,x)$ are less than $n/\omega(G)$. 
Hence, the sum of all roots of $PC(G,x)$ is less than $n$, a contradiction. 
Suppose that $\beta(G) = n/\omega(G)$, it means that all roots of 
$PC(G,x)$ are real and equal $\beta(G)$.
Therefore, $PC(G,x) = (x-n/\omega(G))^{\omega(G)}$.
As all coefficients of $PC(G,x)$ are integers, $n/\omega(G)\in\mathbb{N}$.
So, $G$ is $(\omega(G)+1)$-clique-free graph with 
$\big(1-\frac{1}{\omega(G)}\big)\frac{n^2}{2}$ edges.
It's the Tur\'{a}n graph, i.e., $G$ is the complete multipartite graph 
with $\omega(G)$ parts of size $n/\omega(G)$.
Conversely, we apply Example~2.1.

Consider the empty spanning subgraph $H$ of $G$.
Then $PC(H,x) = x-n$ and $\beta(H) {=} n$. By Lemma~2.4a, we have $\beta(G)\leq n$.
If $G$ is a not empty graph such that $\beta(G) = n$, then we arrive at a contradiction by Lemma~2.4b.

b) Let $u$ and $v$ be disadjacent vertices in~$G$.
Consider in~$G$ a~subgraph~$H$ induced by the set of vertices $\{u,v\}$.
By Lemma~2.3a, $2 = \beta(H)\leq \beta(G)$.

{\bf Lemma 2.7} \cite{Fisher1989}.
Let $G$ be a graph with clique number $w$, $j\leq w$. Then 

a) $\beta^j-c_1\beta^{j-1}+c_2\beta^{j-2}-\ldots -c_j\leq 0$ for odd~$j$,

b) $\beta^j-c_1\beta^{j-1}+c_2\beta^{j-2}-\ldots +c_j\geq 0$ for even~$j$.

{\sc Proof}.
Following notations from Theorem~1.1 and Lemma~1.3, it is easy to show that 
$$
m_t-c_1m_{t-1}+c_2m_{t-2}-\ldots +(-1)^j c_j m_{t-j}
 = (-1)^j\sum\limits_{w\in M_t(X,G)} k_{j+1}'(w).
$$
Thus, 
$m_t-c_1m_{t-1}+c_2m_{t-2}-\ldots -c_j m_{t-j}\leq 0$ for odd~$j$ and
$m_t-c_1m_{t-1}+c_2m_{t-2}-\ldots +c_j m_{t-j}\geq 0$ for even~$j$.
It remains to divide the inequalities by $m_{t-j}$ and consider the limit $t\to\infty$.

\newpage
\section{Applications in graph theory}

{\bf Statement 3.1 (Mantel's Theorem)} \cite{Hajiabolhassan}.
The maximum number of edges in an $n$-vertex triangle-free 
graph is $[n^2/4]$.

{\sc Proof}.
Write down $PC$-polynomial of $G$ as $PC(G,x) = x^2-nx+k$, where $k = |E(G)|$. 
By Lemma~2.1, the equation $x^2 - nx + k = 0$ has real roots, i.e., $k\leq n^2/4$.

The complete bipartite graph $K_{[n/2],n-[n/2]}$ is the example of 
a triangle-free graph with $[n^2/4]$ edges.

{\bf Statement 3.2} \cite{Hajiabolhassan}.
Given a graph~$G$ with $n$ vertices, we have

a) $\alpha(G)\leq\beta(G)$, where $\alpha(G) = \omega(\bar{G})$ is the independence number of $G$,

b) $n/\beta(G)\leq \chi(G)$, where $\chi(G)$ is the chromatic number of $G$,

c) $g(G)\leq\dfrac{\beta^2(G)}{\beta(G)-1}$ provided that $g(G)<\infty$, where $g(G)$ is the girth of $G$.

{\sc Proof}.
a) Consider the subgraph~$H$ in~$G$ induced by the set of $\alpha(G)$ vertices 
forming the maximal independent set. 
Then $\beta(G) \geq \beta (H) = \alpha(G)$ by Example 2.1b and Lemma~2.3a.

b) It follows from a) and the well-known inequality $\alpha(G)\chi(G)\geq n$.

c) Consider the subgraph~$H$ in~$G$ induced by a cycle of length $g(G)$.
For $g(G) = 3$, we have $PC(H,x) = x^3 - g(G)x^2 + g(G)x-1 = (x-1)^3$.
So, $\beta(H) = 1$ and the statement follows from the inequalities
$$
\beta(G)^2 -3\beta(G)+3
 =\frac{(\beta(G)-1)^3+1}{\beta(G)}\geq\frac{1}{\beta(G)}>0
$$
holding by Lemma~2.6.

For $g(G)\geq4$, we have $PC(H,x) = x^2 - g(G)x + g(G)$. Thus, 
$$
\beta(H) = \frac{g(G)+\sqrt{g(G)^2-4g(G)}}{2}\leq\beta(G)
$$
or $(2\beta(G)-g(G))^2\geq g^2(G)-4g(G)$
which gives the required inequality.

{\bf Statement 3.3} \cite{Faal}.
Let $f\colon G\to H$ be a surjective homomorphism from~$G$ to~$H$,
i.e., a surjective map $f\colon V(G)\to V(H)$ such that 
$(f(u),f(v))\in E(H)$ for all $(u,v)\in E(G)$. Then $\beta(H)\leq\beta(G)$.

{\sc Proof}.
Let $V(H) = \{v_1,\ldots,v_m\}$. Consider $u_1,\ldots,u_m\in V(G)$
such that $f(u_i) = v_i$. Define $G'$ as $G[U]$ for $U = \{u_1,\ldots,u_m\}$.
By Lemma~2.3a and Lemma~2.4a, we have 
$\beta(H) \leq \beta(G') \leq \beta (G)$.

The following Theorem was initially proved in 2009 by D. Galvin~\cite{Galvin},
actually the proof was based on the result of V. Alekseev~\cite{Alekseev}.

{\bf Theorem 3.1} \cite{Alekseev,Davies,Galvin}.
Given a graph $G$ with $w = \omega(G)$, $n = |V(G)|$,
we have $C(G,x)\leq \big(1+\frac{nx}{w}\big)^w$ for all $x>0$
with equality if and only if $G$ is a complete multipartite graph with equal parts.

{\sc Proof}.
For $w = 1$, the statement is trivial. 
Let $w > 1$. Prove the statement by induction on $n = |V(G)|$.
The case $n = 1$ is also trivial.

Let $G$ equal $H_1+H_2+\ldots+H_k$ with connected $\bar{H}_i$, $i=1,\ldots,k$.
Define $n_i = |V(H_i)|$, $w_i = \omega(H_i)$.
It is clear that $\sum\limits_{i=1}^k n_i = n$, $\sum\limits_{i=1}^k w_i = w$.
Applying Lemma~1.5e, Jensen's inequality and the induction hypothesis, we get 
\begin{multline}\label{Jensen}
C(G,x) = \prod\limits_{i=1}^k C(H_i,x)
 \leq \prod\limits_{i=1}^k\left(1+\frac{n_i x}{w_i}\right)^{w_i} \\
 \leq \left({\sum\limits_{i=1}^kw_i\left(1+\frac{n_i x}{w_i}\right)}\big/{w}\right)^w
 = \left(1+\frac{nx}{w}\right)^w.
\end{multline}
We have equality in \eqref{Jensen} if and only if all $n_i$ are equal, 
it means all components of $G$ have the same order.

Let $G$ be a graph such that $\bar{G}$ is connected.
Choose a vertex $v\in V(G)$ of minimal degree 
$d(v) = n-\Delta-1$. By Lemma~1.5a, 
$C(G,x) = C(G \setminus v,x) + xD(G[N(v)],x)$.
Since $\omega(G \setminus v)\leq w$, $\omega(G[N(v)])\leq w-1$,
and the function $(1+nx/w)^w$ is increasing in $w>0$ for all $x > 0$, we have 
\begin{equation}\label{PCGlobalBound}
C(G,x) \leq \left(1+\frac{(n-1)x}{w}\right)^w
 + x\left(1+\frac{(n-\Delta-1)x}{w-1}\right)^{w-1}.
\end{equation}
It is easy to show the inequality $\Delta \geq \frac{n-1}{w}$ 
when $\bar{G}$ is either complete graph or a cycle of an odd length.
Otherwise, by Brooks' theorem 
$\Delta\geq \chi(\bar{G})\geq \frac{n}{\alpha(\bar{G})} = \frac{n}{w}$.
Inserting the inequality $\Delta \geq \frac{n-1}{w}$ into~\eqref{PCGlobalBound}, we obtain
$$
C(G,x) \leq \left(1+\frac{(n-1)x}{w}\right)^w
 + x\left(1+\frac{(n-1)x}{w}\right)^{w-1}
$$
and so 
$$
\frac{C(G,x)-\big(1+\frac{nx}{w}\big)^w}{\big(1+\frac{(n-1)x}{w}\big)^w}
 \leq 1+\frac{xw}{w+(n-1)x}-\left(1+\frac{x}{w+(n-1)x}\right)^w \leq 0
$$
with equality if and only if $w = 1$, i.e., $G$ is complete.
Theorem is proved.

{\bf Corollary 3.1} \cite{Alekseev,Erdos}.
Let $G$ be a graph with $w = \omega(G)$, $|V(G)| = n$
and $c(G)$ denotes the number of all cliques in $G$. Then 
$c(G)\leq \big(1+\frac{n}{w}\big)^w$.
We have equality if and only if $G$ is a complete multipartite 
graph with equal parts.

{\sc Proof}.
Note that $c(G) = c_0(G) + c_1(G) + \ldots + c_{w}(G) = D(G,-1)$.
It remains to apply Theorem~3.1.

{\bf Corollary 3.2}.
Given a graph $G$ with $w = \omega(G)$, $n = |V(G)|$,
for all $x<0$ we have $PC(G,x)\leq \big(x-\frac{n}{w}\big)^w$ for even~$w$
and $PC(G,x)\geq -\big(x-\frac{n}{w}\big)^w$ for odd~$w$.
We have equality if and only if $G$ is a complete multipartite 
graph with equal parts.

{\bf Remark 3.1}.
In \cite{Davies}, Theorem~3.1 was proved with the help of 
the Moon---Moser inequalities~\cite{MoonMoser}
on the numbers of cliques of sizes~$s-1$, $s$ and $s+1$ in a graph $G$:
\begin{equation}\label{MoonMoser}
c_{s+1}\geq \frac{s^2}{s^2-1}c_s\left(\frac{c_s}{c_{s-1}}-\frac{n}{s^2}\right),
\end{equation}
where $s\geq2$, $c_{s-1}\neq0$.

{\bf Remark 3.2}.
In \cite{Fisher1992-Clique}, Corollary~3.1 was stated directly 
with the help of the following inequalities~\cite{Fisher1992-Clique,Khadziivanov}
\begin{equation}\label{FisherClique}
\left(\frac{c_1}{\binom{w}{1}}\right)^{1/1}
 \geq \left(\frac{c_2}{\binom{w}{2}}\right)^{1/2}
 \geq \left(\frac{c_3}{\binom{w}{3}}\right)^{1/3}\geq \ldots
 \geq \left(\frac{c_w}{\binom{w}{2}}\right)^{1/w},\quad w = \omega(G).
\end{equation}
Indeed, from \eqref{FisherClique} it follows that $c_j\leq C_w^j\big(\frac{n}{w}\big)^j$
and $c(G) = 1 + c_1 + \ldots + c_w \leq \big(1+\frac{n}{w}\big)^w$.

{\bf Corollary 3.3 (Tur\'{a}n's Theorem)}.
Given a graph $G$ with $w = \omega(G)$, $|V(G)| = n$, 
we have $k = |E(G)|\leq \frac{n^2}{2}\big(1-\frac{1}{w}\big)$.

{\sc Proof}.
Let $w$ be even, then
\begin{multline}\label{TuranProof}
PC(G,x) = x^w -nx^{w-1}+kx^{w-2}-\ldots  \\
\leq
x^w -nx^{w-1}+\frac{n^2(w-1)}{2w}x^{w-2}-\ldots
= \left(x-\frac{n}{w}\right)^w.
\end{multline}
For $x\ll -1$, the inequality \eqref{TuranProof} is fulfilled only if 
$k\leq \frac{n^2}{2}\big(1-\frac{1}{w}\big)$.

The proof for odd $w$ is analogous.

{\bf Statement~3.4}.
Let $G$ be a graph with $w = \omega(G)$, $|V(G)| = n$. 
Tur\'{a}n's Theorem is equivalent to the fact that 
the polynomial $PC^{(w-2)}(G,x)$ has real roots.

{\sc Proof}. 
Let $k = |E(G)|$, then 
\begin{equation}\label{QuadTuran}
\frac{1}{(w-2)!}PC^{(w-2)}(G,x)
 = \frac{w(w-1)}{2}x^2 - n(w-1)x + k.
\end{equation}
The discriminant of the RHS from \eqref{QuadTuran} is nonnegative if and only if 
$k\leq \frac{n^2}{2}\big(1-\frac{1}{w}\big)$.

Initiated by an effort of  A. Granville to resolve the Cameron---Erd\H{o}s conjecture,
N.~Alon proposed in 1991 the following 

{\bf Conjecture 3.1} \cite{Alon}.
For any $n$-vertex $d$-regular graph~$G$, we have 
$i(G) \leq (2^{d+1} - 1)^{\frac{n}{2d}}$.

\noindent
Here $i(G) = c(\bar{G})$ equals the number of all independent sets in $G$.

Note that for $n$ divisible by $2d$,
we have the equality $i(G) = (2^{d+1} - 1)^{\frac{n}{2d}}$
when $G$~is a disjoint union of $n/(2d)$ complete bipartite graphs $K_{d,d}$.
In 2001, J. Kahn proved Conjecture 3.1 for any bipartite graph~$G$ \cite{Kahn}.
In 2009, Y. Zhao proved Conjecture 3.1 completely \cite{Zhao2}, we state here this proof
concerned with independence polynomial.

Let $\mathcal{I}(G)$ denote the set of all independent sets in $G$.

For $A, B \subset V(G)$, say that $A$ is independent from~$B$ 
if any $a \in A$ and any $b \in B$ are disconnected. 
Let $\mathcal{J}(G)$ denote the set of pairs $(A,B)$ of subsets of vertices of $G$ 
such that $A$ is independent from $B$ and $G[A \cup B]$ is bipartite. 
For a pair $(A, B)$ of subsets of $V(G)$, define its size as $|A| + |B|$.

{\bf Lemma 3.1} \cite{Zhao2}. 
For a graph $G$, there is a bijection between 
$\mathcal{I}(G) \times \mathcal{I}(G)$ and~$\mathcal{J}(G)$.

{\sc Proof}.
For every $W \subset V(G)$ such that $G[W]$ is bipartite, fix a bipartition 
$W = W_1 \cup W_2$ so that $W_1$ and $W_2$ are both independent in $G$. 
Let $\mathcal{K}(G)$ be the set of pairs $(A, B)$ of subsets of $V(G)$ such that 
$G[A\cup B]$ is bipartite. Note that 
$\mathcal{I}(G) \times \mathcal{I}(G) \subset \mathcal{K}(G)$ and 
$\mathcal{J}(G) \subset \mathcal{K}(G)$. 
Indeed, if $A, B \in \mathcal{I}(G)$ 
then $A \cup (B \setminus A)$ is already a bipartition.

Let us construct a bijection $\varphi\colon \mathcal{K}(G)\to\mathcal{K}(G)$ as follows. 
For any $(A, B) \in \mathcal{K}(G)$, let $W_1 \cup W_2$ be the chosen bipartition 
of $W = A \cup B$. We define 
$$
\varphi((A, B))
 = ((A \cap W_1) \cup (B \cap W_2),  (A \cap W_2) \cup (B \cap W_1)).
$$
It is easy to check that $\varphi$ is a bijection on $\mathcal{K}(G)$ that 
maps $\mathcal{I}(G) \times \mathcal{I}(G)$ to $\mathcal{J}(G)$ and vice versa. 

{\bf Theorem 3.2} \cite{Zhao2}.
Given a $d$-regular graph $G$ with $n = |V(G)|$, for all $x\geq 0$ we have 
$I(G,x) \leq (2(1+x)^d-1)^{n/(2d)}$.

{\sc Proof}.
Let $G \times K_2$ denote the bipartite graph with 
$$
V(G \times K_2) = \{(v, i)\mid v \in V(G), i =0,1\},\
E(G \times K_2) = \{((v,0),(w,1))\mid (v,w)\in E(G)\}.
$$
Thus, independent sets in $G \times K_2$ correspond to pairs $(A,B)$ 
of subsets of $V(G)$ such that $A$ is independent from $B$. 
Applying Lemma~3.1, for all $x\geq 0$ we have
\begin{multline} \label{ZhaoProof}
I(G \times K_2,x)
 = \sum_{I \in \mathcal{I}(G \times K_2)}  \hspace{-1em} x^{|I|}
 = \sum_{\substack{A, B \subset V(G) \\ A \text{ indep.~from } B}} x^{|A| + |B|} \\
 \geq \sum_{(A,B) \in \mathcal{J}(G)} x^{|A| + |B|}
 = \sum_{A, B \in \mathcal{I}(G)} x^{|A| + |B|}
 = I(G,x)^2.
\end{multline}

Note that $G \times K_2$ is also a $d$-regular graph. 
Applying the inequality of Statement holding for bipartite 
graphs~\cite{Galvin2}, by \eqref{ZhaoProof}
we have for all $x\leq 0$
$$
I(G,x)
 \leq I(G \times K_2,x)^{1/2} 
 \leq I(K_{d,d},x)^{n/(2d)}
 = (2(1+x)^d-1)^{n/(2d)}.
$$

{\bf Corollary 3.4} \cite{Zhao2}.
For any $n$-vertex $d$-regular graph $G$, 
$i(G)\leq (2^{d+1} - 1)^{\frac{n}{2d}}$.

{\bf Corollary 3.5} \cite{Zhao2}.
Given a $d$-regular graph $G$ with $n = |V(G)|$, for all $x\leq 0$ we have 
$C(G,x) \leq (2(1+x)^{n-d}-1)^{\frac{n}{2(n-d)}}$
and $c(G) = C(G,1)\leq (2^{n-d+1} - 1)^{\frac{n}{2(n-d)}}$.

Let us state without proofs some results devoted to the applications of 
(in)dependence polynomial in graph theory.

{\bf Theorem 3.3} \cite{Zhao2018}.
Let $G$ be a graph and $d_v$ denote the degree of vertex $v$ in $G$. 

a) Suppose that $G$ is without isolated vertices. Then for all $x>0$
$$
I(G,x) \leq \prod_{uv \in E(G)} I(K_{d_v,d_u},x)^{1/(d_ud_v)} 
 = \prod_{uv \in E(G)} ((1+x)^{d_u} + (1+x)^{d_v}-1)^{1/(d_ud_v)}.
$$
Equality holds if and only if $G$ is a disjoint union of complete bipartite graphs.

b) For all $x>0$, we have
$$
I(G,x) \geq \prod_{v \in V(G)} I(K_{d_v+1},x)^{1/(d_v+1)} 
 = \prod_{v \in V(G)} ((d_v+1)x + 1)^{1/(d_v+1)}.
$$
Equality holds if and only if $G$ is a disjoint union of cliques.

{\bf Corollary 3.6} \cite{CutlerRadcliffe,Zhao2018}.
Let $G$ be a graph. 

a) If $G$ is without isolated vertices, then
$i(G) \leq \prod_{uv \in E(G)} (2^{d_u} + 2^{d_v} - 1)^{1/(d_ud_v)}$.
Equality holds if and only if $G$ is a disjoint union of complete bipartite graphs.

b) We have $i(G)\geq \prod_{v\in V(G)}(d_v+2)^{1/(d_v+1)}$.
Equality holds if and only if $G$ is a disjoint union of cliques.

Let $H_{d,n}$ denote the $d$-regular, $n$-vertex graph that 
is the disjoint union of $n/(2d)$ copies of $K_{d,d}$.

{\bf Statement 3.5} \cite{Davies2017}.
For all $d$-regular graphs $G$ on $n$ vertices (where $2d$ divides $n$), 
$i_k(G) \le 2 \sqrt{n} \cdot i_k(H_{d,n})$. 

By the {\it occupancy fraction} we mean the expected fraction of vertices 
that appear in the random independent set
$$
\alpha(G,x) 
 = \frac{E(|I|)}{|V(G)|}
 = \frac{1}{|V(G)|}\sum\limits_{I\in\mathcal{I}(G)}|I|\cdot \mathrm{Pr}[I]
 = \frac{xI'(G,x)}{|V(G)|I(G,x)} = x\left(\frac{1}{|V(G)|}\ln I(G,x)\right)'.
$$
Here $\mathrm{Pr}[I] = \dfrac{x^{|I|}}{\sum\limits_{J\in\mathcal{I}(G)}x^{|J|}}$
is so called {\it hard-core distribution} which is simply the uniform
distribution over all independent sets of $G$ at fugacity~$x$.
The expression $\dfrac{\ln I(G,x)}{|V(G)|}$ is called the {\it free energy}.
The independence polynomial is interpreted as the {\it partition function} 
of the hard-core model on $G$ at fugacity $x$.

{\bf Example 3.1}. 
a) $\alpha(\bar{K}_n,x) = \dfrac{nx(1+x)^{n-1}}{n(1+x)^n} = \dfrac{x}{1+x}$;

b) $\alpha(K_n,x) = \dfrac{xn}{n(1+nx)} = \dfrac{x}{1+nx}$;

c) $\alpha(\bar{C}_n,x) = \dfrac{x(n+2nx)}{n(1+nx+nx^2)} = \dfrac{x(1+2x)}{1+nx+nx^2}$,
where $C_n$ denotes the $n$-vertex cycle.

{\bf Statement 3.6} \cite{Davies2017,Zhao2017}.
Let $G$ be a $d$-regular graph. For all $x\geq0$, 
we have the following inequality 
$\alpha(G,x)\leq \alpha(K_{d,d},x) = \dfrac{x(1+x)^{d-1}}{2(1+x)^d-1}$.

Let $o_d(1)$ denote a quantity that tends to zero as $d$ tends to infinity.

{\bf Statement 3.7}~\cite{Davies}.
a) For any graph $G$, $\alpha (G,x)$ is monotone increasing in~$x$.

b) Let $G$ be a triangle-free graph on $n$ vertices with maximum degree $d$,
we have $\alpha(G,x)\geq (1+o_d(1))\frac{\ln d}{d}$
for any $x\geq 1/\ln d$.

{\sc Proof}.
a) Show that the derivative of $n\alpha(G,x)$ is positive:
\begin{gather*}
(n\alpha(G,x))' 
 = \left(\frac{x I'(G,x)}{I(G,x)}\right)' = \frac{I'(G,x)}{I(G,x)} + \frac{xII''(G,x)-x(I'(G,x))^2}{I^2(G,x)}  \\
 = \frac{I'(G,x)}{I(G,x)} +\frac{1}{x} \left(\frac{x^2 I''(G,x)}{I(G,x)}-\left(\frac{x I'(G,x)}{I(G,x)}\right)^2 \right) \\
 = \frac{E(|I|) + E(|I|^2) - E(|I|) - (E(|I|))^2}{x}
 = \frac{D(|I|)}{x} \geq 0,
\end{gather*}
where $I$ is a random independent set drawn from the hard-core model at fugacity $x$. 

b) See in \cite{Davies} the proof based on a) and the lower bound 
on $\alpha(G,x)$ for triangle-free graphs via Lambert $W$-function.

{\bf Corollary 3.8}~\cite{Davies,Shearer2}.
For the Ramsey numbers $R(3,k)$, we have the upper bound
$R(3,k)\leq (1+o(1))\frac{k^2}{\ln k}$.

{\sc Proof}.
Suppose that $G$ is triangle-free graph with no independent set of size $k$. 
Then $G$ must have maximum degree less than $k$. 
By Statement~3.7b, on the one hand 
$n\alpha(G)\geq n\alpha(G,1)\geq (1+o(1))\frac{\ln k}{k}n$.
On the other hand, $\alpha(G)<k$, i.e.,
$k>(1+o(1))\frac{\ln k}{k}n$ as required.

The upper bound for $R(3,k)$ from Corollary~3.7 is still the best known one.
In \cite{RamseyA,RamseyB}, it was proved that $R(3,k)\geq (1/4 + o(1))k^2/\ln k$. 
Reducing the factor 4 gap between these bounds is an open hard problem.

In \cite{Davies}, the following conjecture about the ratio 
of the maximum and average indepen\-dent set sizes was formulated.

{\bf Conjecture 3.2}~\cite{Davies}.
a) For every $K_r$-free graph $G$, 
$\frac{\alpha(G)}{\alpha(G,1)}\geq 1+\frac{1}{r}$.

b) For every $K_r$-free graph $G$ of minimum degree $d$,
$\frac{\alpha(G)}{\alpha(G,1)}\geq 2-o_d(1)$
with $r$ fixed as $d\to \infty$.

{\bf Statement 3.8}.
Let $G$ be not empty graph with $n$ vertices and the spectral radius~$\rho$. 
Then $\alpha(\bar{G},x)\geq \frac{1}{n}$ for any $x\geq 1/\rho$.

{\sc Proof}.
Let us prove the inequality for $x = 1/\rho$. 
For greater~$x$ it will follow by Statement~3.7a.

Let $w = \omega(G)$, by the condition $w\geq2$.
In 2002, V. Nikiforov proved~\cite{NikiforovIneq} the inequality 
\begin{equation}\label{NikiforovRhoIneq}
\rho^w \leq c_2(G)\rho^{w-2}+2c_3(G)\rho^{w-2}+\ldots+(i-1)c_i(G)\rho^{w-i}+\ldots+(w-1)c_w(G),
\end{equation}
which is equivalent to $f'(\rho)\leq 0$ for the function
$$
f(x) = x+\sum\limits_{i=2}^w c_i x^{1-i}
 = x\left(1+\sum\limits_{i=2}^w \frac{c_i}{x^i}\right)
 = x(C(G,1/x) - n/x).
$$
Since 
$$
f'(x) = C(G,1/x) - n/x + x\left(-\frac{1}{x^2}C'(G,1/x) + \frac{n}{x^2}\right)
 = C(G,1/x) - (1/x)C'(G,1/x),
$$
we have $C(G,1/\rho)\leq (1/\rho)C'(G,1/\rho)$.
In terms of independence polynomial and occu\-pancy fraction,
we get the inequality
$$
\alpha(\bar{G},1/\rho)
 = \frac{I'(\bar{G},1/\rho)}{n\rho I(\bar{G},1/\rho)}\geq \frac{1}{n}
$$
and we are done.

See the exposition \cite{Zhao2017} for another applications 
of the stated results for counting colorings, graph homomorphisms 
and independent sets of graphs with different constraints.
Just note that for colorings, the Potts model 
plays the same role as (in)dependent polynomial for (co)cliques.

In the preprint \cite{Perkins}, the applications of 
the independence polynomial were also shown for sphere packings problems.

In the book \cite{Barvinok}, the different problems
concerned roots of matching and independence polynomials
were considered including applications for permanents.

\newpage
\section{Partially commutative Lie algebras}

It is well-known that any associative algebra $\langle A,\cdot\rangle$ under the commutator 
\begin{equation}\label{AsToLie}
[x,y] = x\cdot y - y\cdot x
\end{equation}
is a Lie algebra. Denote the Lie algebra obtained as $A^{(-)}$.
If a Lie algebra $g$ is a Lie subalgebra of $A^{(-)}$,
then $A$ is called an (associative) enveloping of~$g$.

Let $g$ be a Lie algebra and $X$ be a linear basis of~$g$.
The universal enveloping associative algebra $U(g)$ 
is defined as enveloping of~$g$ such that 
$U(g)$ is generated by $X$ and any enveloping of~$g$ generated by $X$
is a homomorphic image of~$U(g)$. 
The algebra $U(g)$ is unique up to isomorphism.

{\bf Poincar\'{e}---Birkhoff---Witt Theorem}.
Let $\{x_i\colon i \in I\}$ be a basis of a Lie algebra~$g$ totally ordered by a set~$I$.
Then words $x_{i_1}\ldots x_{i_n}$, $i_1 \leq \ldots\leq i_n$, form a basis of $U(g)$.

{\bf Remark 4.1}.
One can endow $U(g)$ with the operation $\circ$ such that
$\langle U(g),\circ\rangle \cong \Bbbk[x_i]_{i\in I}$,
the polynomial algebra on $x_i$.

{\bf Remark 4.2}.
Any Lie algebra~$g$ (injectively) embeds into $U^{(-)}(g)$.

Define a partially commutative Lie algebra $L(X,G)$ as 
$\Lie\langle X|[a,b] = 0,(a,b)\in E(G)\rangle$.
Denote the dimension of the homogeneous space of all products 
of length~$n$ in the alphabet~$X$ in $L(X,G)$ as $l_n$.
Consider in $\Ass^{(-)}(X,G)$ the Lie subalgebra $T$ generated by the set~$X$.

In 1992, G. Duchamp and D. Krob proved that 

{\bf Lemma 4.1} \cite{DuchampKrob1992}.
a) $U(L(X,G)) \cong\Ass(X,G)$,
b) $L(X,G)\cong T$.

{\sc Proof}.
In the algebra $U = U(L(X,G))$ the relations $ab = ba$ for $(a,b)\in E(G)$ hold. 
Since $U$ is generated by $X$, there exists 
a homomorphism $\varphi\colon \Ass(X,G)\to U$ such that $\varphi(\Ass(X,G)) = U$.
Moreover, $L(X,G)$ is a homomorphic image of $T$ under~$\varphi$.

On the other hand, in $\Ass^{(-)}(X,G)$ the relations $[a,b] = 0$ for $(a,b)\in E(G)$ are fulfilled.
Thus, $T$ is a homomorphic image of $L(X,G)$ and we have proved b).
Therefore, the algebra $\Ass(X,G)$ is an enveloping of~$L(X,G)$ generated by~$X$.
We have that $\Ass(X,G)$ is a homomorphic image of~$U$
and hence, $\Ass(X,G)\cong U$. Lemma is proved.

Given roots $x_1,\ldots,x_{t_0}$ of $PC(G,x)$, $t_0 = \omega(G)$, define the numbers 
\begin{equation}\label{preWitt}
p_n = \sum\limits_{j=1}^{t_0}x_j^n.
\end{equation}
The numbers $p_n$ could be expressed by the Newton's identities
via the coefficients of $PC(G,x)$, i.e., via the numbers $c_k(G)$.

In the same work~\cite{DuchampKrob1992}, 
G. Duchamp and D. Krob actually proved the following result 
(but not in the most comfortable form).

{\bf Theorem 4.1} \cite{DuchampKrob1992}.
We have
\begin{equation}\label{Witt}
l_n = \frac{1}{n}\sum\limits_{d|n}\mu(d)p_{n/d},
\end{equation}
where $\mu$ is the M\"{o}bius function and the numbers $p_i$ are defined by~\eqref{preWitt}.

{\sc Proof}.
Form in $L(X,G)$ a linear basis $F$ which is a union 
of linear bases of homo\-geneous spaces of products of the same length in the alphabet~$X$:
$$
f_{11},\ldots,f_{1l_1},\
f_{21},\ldots,f_{2l_2},\
\ldots,\
f_{n1},\ldots,f_{nl_n},\
\ldots.
$$

Given a word $w = x_{i_1}\ldots x_{i_s}$ from the basis of $U(L(X,G))$,
the sum of the lengths of letters $x_{i_j}$ is called a degree of $w$.
The degree of~$w$ equals the degree of the homogeneous expression 
$x_{i_1}\ldots x_{i_s}$ in the alphabet~$X$ after calculating it 
in $U(L(X,G))$ by \eqref{AsToLie}.
The generating function of the sequence formed 
by the numbers of words of the fixed degree in $U(L(X,G))$ 
consisting only of letters from $\{f_{n1},\ldots,f_{nl_n}\}$
equals $\left(\frac{1}{1-x^n}\right)^{l_n}$.
Since letters $x_{i_j}$ of the basic word $x_{i_1}\ldots x_{i_s}$ 
can be any Lie words of any length, 
the generating function of the sequence formed 
by the numbers of words of the fixed degree in $U(L(X,G))$ equals 
$$
\prod\limits_{n\geq1}\left(\frac{1}{1-x^n}\right)^{l_n}.
$$

Poincar\'{e}---Birkhoff---Witt Theorem implies the equality
$$
\prod\limits_{n\geq1}\left(\frac{1}{1-x^n}\right)^{l_n}
 = \sum\limits_{m\geq0}a_m x^m,
$$
where $a_m$ equals the dimension of the homogeneous space of all words of length~$m$
in $\Ass(X,G)$. Let $x_1,\ldots,x_{t_0}$ be roots (including complex) of $PC(G,x)$. By~\eqref{GF-PC}
$$
\prod\limits_{n\geq1}\left(\frac{1}{1-x^n}\right)^{l_n}
 = \frac{1}{D(G,x)}
 = \frac{1}{\prod\limits_{j=1}^{t_0}(1-x_j x)}.
$$
Calculating logarithm on both sides of the last equality, we get
$$
\ln (D(G,x)) = \sum\limits_{n\geq1}l_n\ln(1-x^n)
$$
or
$$
\sum\limits_{j=1}^{t_0}\sum\limits_{m\geq1}\frac{x_j^m x^m}{m}
 = \sum\limits_{n\geq1}l_n\sum\limits_{k\geq1}\frac{x^{nk}}{k}.
$$
Comparing the coefficients by the same degree $x^n$, find
$$
\sum\limits_{j=1}^{t_0}\frac{x_j^n}{n}
  = \sum\limits_{d|n}\frac{l_d}{n/d}
  = \frac{1}{n}\sum\limits_{d|n}dl_d
$$
which implies $\sum\limits_{j=1}^{t_0}x_j^n = \sum\limits_{d|n}dl_d$.
Applying the M\"{o}bius inversion formula, we have
$$
l_n = \frac{1}{n}\sum\limits_{d|n}\mu(d)\left(\sum\limits_{j=1}^{t_0}x_j^{n/d}\right)
 = \frac{1}{n}\sum\limits_{d|n}\mu(d)p_{n/d}.
$$

{\bf Corollary 4.1}.
If $G = K_n$, then the growth rate of partially commutative Lie algebra $L(X,G)$ equals~0.
Otherwise, it equals~$\beta(G)$.

{\sc Proof}.
For $G = K_n$, we calculate $\{l_n\} = \{1,|X|,0,0,\ldots,0,\ldots\}$,
thus, the growth rate equals~0. Otherwise, by Lemma~2.6 
$\beta(G)\geq2$ and due to \eqref{preWitt},~\eqref{Witt} 
we get the asymptotics 
$l_n\sim \frac{\beta^n(G)}{n}$ and
$\lim\limits_{n\to\infty}\sqrt[n]{l_n}
 = \lim\limits_{n\to\infty}\sqrt[n]{\frac{\beta^n(G)}{n}} = \beta(G)$.
Corollary is proved.

Note that the classical Witt's formula could be derived from~\eqref{Witt}
as a particular case. Indeed, the free Lie algebra generated by~$X$ 
is exactly partially commutative Lie algebra $L(X,\bar{K}_n)$.
Let $q = |X|$, then $PC(G,x) = x - q$ and
$$
l_n = \frac{1}{n}\sum\limits_{d|n}\mu(d)q^{n/d}.
$$

{\bf Remark 4.3}.
Applying Corollary~4.1, 
the part of Lemma~2.6 devoted to the upper bound $\beta(G)\leq n$
and when it is reached follows from~\cite{Bahturin}.
In \cite{Bahturin}, the growth rate of subalgebras, ideals and subideals
of a free finitely generated Lie algebra was studied.

{\bf Corollary 4.2}.
a) Let $G$ be a tree with $q$ vertices, then for $n\geq2$ 
$$
l_n(G) = \frac{1}{n}\sum\limits_{d|n}\mu(d)(q-1)^{n/d} = l_n(\bar{K}_{q-1}).
$$

b) For $G = K_{q_1,\ldots,q_s}$, we have for $n\geq2$ 
$$
l_n(G) = \frac{1}{n}\sum\limits_{j=1}^s\sum\limits_{d|n}\mu(d)q_j^{n/d}
 = l_n(\bar{K}_{q_1}) + \ldots + l_n(\bar{K}_{q_s}).
$$

{\sc Proof}.
It follows from Example~1.1, the formulas \eqref{preWitt}, \eqref{Witt} 
and the properties of the M\"{o}bius function. 

Let us write down some results found by computations.

{\bf Corollary 4.3}.
a) Let $G = K_2\cup K_1\cup K_1\cup K_1\cup K_1$.
For $n\geq3$ the sequence $l_n(G)$ coincides with A212443 \cite{Sloane}.

b) Let $G = K_2\cup K_1$.
For $n\geq3$ the sequence $l_n(G)$ coincides with A072337 \cite{Sloane}.

c) Let $G = K_2\cup K_1\cup K_1$.
For $n\geq3$ the sequence $l_n(G)$ coincides with A072279 \cite{Sloane}.

{\sc Proof}.
a) We have $PC(G,x) = x^2 - 6x + 1$, so $3\pm2\sqrt{2}$ are roots of $PC(G,x)$.
By~\eqref{Witt} 
$$
l_n(G) = \frac{1}{n}\sum\limits_{d|n}\mu(d)((3+2\sqrt{2})^{n/d}+(3-2\sqrt{2})^{n/d}).
$$
It is known that the sequence A212443 \cite{Sloane} equals
$$
a_n = \frac{1}{n}\sum\limits_{d|n}\mu(d)((1+\sqrt{2})^{n/d}+(1-\sqrt{2})^{n/d})^2,
$$
which for $n\geq3$ implies
$$
a_n = l_n(G) + \frac{2}{n}\sum\limits_{d|n}\mu(d)(-1)^{n/d}
    = l_n(G).
$$

b) By \cite{YM2}, A072337 equals a sequence of the dimensions of homogeneous spaces
of the Lie Yang---Mills algebra over an algebraically closed field of characteristic zero
$$
\mathfrak{nm}(k) = \Lie\bigg\langle x_1,\ldots,x_k\mid \sum\limits_{i=1}^{k}[x_i,[x_i,x_j]],\,j=1,\ldots,k\bigg\rangle
$$
for $k=3$.
Let $\mathrm{YM}(k) = U(\mathfrak{nm}(k))$ be its universal enveloping algebra. 

The generating function of the dimensions $b_n$ of $\mathrm{YM}(3)$
equals $\frac{1}{(1-x^2)(1-3x+x^2)}$.
For $G = K_2\cup K_1$, we have $PC(G,x) = x^2 - 3x + 1$
and $t_{1,2} = \frac{3\pm\sqrt{5}}{2}$ are roots of $PC(G,x)$.
By~\eqref{Witt} and the properties of the M\"{o}bius function,
for $n\geq3$ we have 
$$
b_n = \frac{1}{n}\sum\limits_{d|n}\mu(d)\big(t_1^{n/d}+t_2^{n/d}+1^{n/d}+(-1)^{n/d}\big)
 = \frac{1}{n}\sum\limits_{d|n}\mu(d)\big(t_1^{n/d}+t_2^{n/d}\big)
 = l_n(G).
$$

c) The statement follows from the properties of the sequence~A072279~\cite{Sloane} 
and the formula~\eqref{Witt}. The sequence appears as the sequence
of the dimensions of homogeneous spaces of the Lie algebra $\mathfrak{nm}(4)$. 
As in b), the generating function of the dimensions of $\mathrm{YM}(4)$ equals $1/((1-x^2)D(G,x))$.

{\bf Remark 4.4}.
By computational experiments, the following connection between 
par\-tially commutative Lie algebras and the Lie algebras $L^{(l)}$ 
of primitive elements in the connected cocommutative Hopf algebra 
$\mathrm{Sym}^{(l)}$~\cite{Thibon} was found.
Let $G = C_3\cup K_1$, $H = K_3\cup K_3$.
For $n\geq3$, the sequences $l_n(G)$ and $l_n(H)$ coincide with the 
sequences of the dimensions of $L^{(l)}$ for $l=2$ 
(A141312 \cite{Sloane}) and $l=3$ (A185162 \cite{Sloane})
respectively.

\newpage
\section{Random graphs}

\subsection{$PC$-polynomial of random graph}

Consider $PC$-polynomial and clique polynomial for 
the random graph $G_{n,p}$ with $n$ vertices and edge probability $p$ 
\begin{gather*}
PC(G_{n,p},x) = x^n {-} \binom{n}{1}x^{n-1} {+} \binom{n}{2} px^{n-2} 
 {-} \ldots {+} (-1)^k \binom{n}{k} p^{\frac{k(k-1)}{2}}x^k
 {+} \ldots {+} (-1)^n p^{\frac{n(n-1)}{2}}, \\
C(G_{n,p},x) = 1 + \binom{n}{1} x + \binom{n}{2}px^2
 + \ldots+\binom{n}{k} p^{\frac{k(k-1)}{2}}x^k
 + \ldots + p^{\frac{n(n-1)}{2}}x^n.
\end{gather*}

Denote by $\beta(G_{n,p})$ the largest real root of $PC(G_{n,p},x)$
(see Theorem~5.1 below for the correctness of the definition). 

For $p = 0$, we have $PC(G_{n,0},x) = x - n = 0$ and $\beta(G_{n,0}) = n$.
For $p = 1$, we have $PC(G_{n,1},x) = (x-1)^n = 0$ and $\beta(G_{n,1}) = 1$.
So, we are interested on the properties of $PC(G_{n,p})$ for $p\in(0,1)$.
Let us denote $G_{n,p}$ as $G_p$, when the number $n$ is known.

For $n = 2$, we have 
$C(G_p,x) = 1 + 2x + px^2$ and $PC(G_p,x) = x^2 - 2x + p$.
The largest root of $PC(G_p,x)$ equals
\begin{equation}\label{n=2}
\beta(G_{2,p}) = \frac{2+\sqrt{4-4p}}{2} = 1 +\sqrt{1-p}.
\end{equation}

Let us prove the analogue of Lemma~1.5 for random graph.

{\bf Lemma 5.1} \cite{Brown2}.
The following equalities hold:

a) $C(G_{n,p},x) = C(G_{n-1,p},x) + xC(G_{n-1,p},px)$,

b) $C'(G_{n,p},x) = nC(G_{n-1,p},px)$, 

c) $PC'(G_{n,p},x) = nPC(G_{n-1,p},x)$.

{\sc Proof}.
a) Fix a vertex~$v$. Divide all cliques in $G_{n,p}$ into 
two types, the first ones contain $v$ but not the second ones.
All cliques of the second type are counted in a~summand $C(G_{n-1,p},x)$.
All cliques of the first type are listed in the sum 
$$
x\sum\limits_{i=0}^{n-1} \binom{n-1}{i}p^{i}p^{\frac{i(i-1)}{2}}x^i
 = x\sum\limits_{i=0}^{n-1}\binom{n-1}{i} p^{\frac{i(i-1)}{2}}(px)^i
 = xC(G_{n-1,p},px),
$$
where the factor $p^i$ equals the probability of having 
all edges between~$v$ and all vertices of a clique of size~$i$ in $G_{n,p}\setminus v$.

The proof of b) and c) is analogous.

The following theorem was stated by J. Brown and R. Nowakowski in 2005~\cite{BrownNowak} for $p=1/2$
and by J. Brown et al in 2012~\cite{Brown2} for any $p$.

{\bf Theorem 5.1} \cite{Brown2}.
Let $p\in(0,1)$. 

a) All roots of $C(G_{n,p},x)$ are real and simple.

b) Write roots of $C(G_{n,p},x)$ in ascending order
$r_{n}<\ldots<r_1<0$. Then $pr_{i+1}<r_{i}$ for all $i=1,\ldots,n-1$.

{\sc Proof}.
We proceed by induction on $n$.
For $n = 1$, the polynomial $C(G_{1,p},x)$ has the unique root $-1$.
For $n = 2$, the polynomial $C(G_{2,p},x) = 1 + 2x + px^2$ has roots 
$r_{\pm} = \frac{-1\pm\sqrt{1-p}}{p}$
and we have the required inequality $pr_{-}<r_{+}$.

By the induction hypothesis, $C(G_{n-1,p},x)$ has $n-1$ real roots, 
$$
r_{n-1}<r_{n-2}<\ldots<r_1<0,
$$
and $pr_{i+1}<r_i$.
By Lemma~5.1a, $C(G_{n,p},x) = C(G_{n-1,p},x) + xC(G_{n-1,p},px)$.

The numbers $r_i/p$, $i=1,\ldots,n-1$, and 0 are roots of the polynomial $xC(G_{n-1,p},px)$.
Moreover, we know how the roots of $C(G_{n-1,p},x)$ and $xC(G_{n-1,p},px)$
lie on the real line:
$$
r_{n-1}/p<r_{n-1}<r_{n-2}/p<r_{n-2}<\ldots<r_2/p<r_2<r_1/p<r_1<0.
$$

It is easy to check that signs of $C(G_{n-1,p},x)$ 
and $xC(G_{n-1,p},px)$ coincide in all intervals
$$
d_1 = (r_{n-1}/p,r_{n-1}),\quad d_2=(r_{n-2}/p,r_{n-2}),\quad \ldots,\quad d_{n-1}=(r_1/p,r_1),\quad (0,+\infty),
$$
although $C(G_{n-1,p},x)$ and $xC(G_{n-1,p},px)$
have different signs in a pair of intervals $d_i,d_{i+1}$, $i=1,\ldots,n-2$.
Hence, $C(G_{n,p},x)$ has a real root in all intermediate intervals
$$
(r_{n-1},r_{n-2}/p),\quad (r_{n-2},r_{n-3}/p),\quad \ldots,\quad (r_2,r_1/p),\quad (r_1,0).
$$
The remaining $n$-th root of $C(G_{n,p},x)$ lies in $(-\infty,r_{n-1}/p)$. 
The required inequalities for neighbour roots of $C(G_{n,p},x)$ 
follow from the form of intervals in which the roots lie.

{\bf Remark 5.1}.
Real-rootedness of $C(G_{n,p},x)$ could be obtained from the result of E.~Laguerre.
By change of variables $q = \sqrt{p}$, $y = x/q$, we have
$C(G_{n,p},x)
 = \sum\limits_{k=0}^n \binom{n}{k} q^{k^2}y^k$.
Since all roots of the polynomial
$\sum\limits_{k=0}^n \binom{n}{k} y^k = (1+y)^n$ 
are real, we may apply~\cite{Laguerre} for $|q|\leq1$
to state that all roots of 
$\sum\limits_{k=0}^n \binom{n}{k} q^{k^2}y^k$
are real.

Let $q = \sqrt{p}$, $y = {q}^{n-1}x$. 
Define the polynomial
\begin{equation}\label{RandomChange}
\widetilde{C}(G_{n,q},y)
 = C(G_{n,p},x)
 = \sum\limits_{k=0}^n \binom{n}{k} \frac{y^k}{{q}^{k(n-k)}}.
\end{equation}

A polynomial $F(x) = \sum\limits_{i=0}^n a_i x^i$ of degree~$n$
is called symmetric if $a_i = a_{n-i}$ for all $i=0,\ldots,n$.
By \eqref{RandomChange}, we immediately get

{\bf Statement 5.1}.
a) The polynomial $\widetilde{C}(G_{n,q},y)$ is symmetric.

b) For odd $n$, $p^{\frac{n-1}{2}}$ is a middle root of $PC(G_{n,p},x)$.
All other roots of $PC(G_{n,p},x)$ for odd $n$
and all roots for even $n$ could be gathered in pairs
with the roots product equal~$p^{n-1}$.

{\bf Corollary 5.1}.
For $p\in(0,1)$, we have 
\begin{gather}
\beta(G_{3,p}) = 1+\frac{1-p}{2}+\frac{\sqrt{3(1-p)(3+p)}}{2},\label{n=3}\\
\beta(G_{4,p})
 = 1 + (1-p)\sqrt{\frac{2+p}{2}}
     + \sqrt{\frac{(1-p)(4+p+p^2+ 2\sqrt{2}\sqrt{2+p})}{2}}, \label{n=4}
\end{gather}
\begin{multline}\label{n=5}
\beta(G_{5,p})
 = 1+\frac{1-p^2}{4}+\frac{\sqrt5(1-p)\sqrt{5+2p+p^2}}{4} \\
 +\frac{\sqrt{1-p}\sqrt{5(5+p+p^2+p^3)+(5-p^2)\sqrt{5}\sqrt{5+2p+p^2}}}{2\sqrt2}.
\end{multline}

{\sc Proof}.
Let $q = \sqrt{p}$, $y = {q}^{n-1}x$.
For $n = 3$, we have
$$
\widetilde{C}(G_{q},y) = 1 + \frac{3y}{{q}^2} + \frac{3y^2}{{q}^2} + y^3
 = 1 + \frac{3y}{p} + \frac{3y^2}{p} + y^3
 = (1+y)\left(y^2+y\left(\frac{3}{p}-1\right)+1\right).
$$

The smallest root of $\widetilde{C}(G_p,x)$ equals 
$y_0 = \dfrac{p-3+\sqrt{3(3-2p-p^2)}}{2p}$, so,
$$
\beta(G_p) = -\frac{2p^2}{p-3+\sqrt{3(3-2p-p^2)}}
 = 1+\frac{1-p}{2}+\frac{\sqrt{3(1-p)(3+p)}}{2}.
$$

For $n = 4$, we have
$$
\widetilde{C}(G_{q},y) = 1 + \frac{4y}{q^3} + \frac{6y^2}{q^4} + \frac{4y^3}{q^3} + y^4
= y^2\left( \frac{1}{y^2} + \frac{4}{yq^3} + \frac{6}{q^4} + \frac{4y}{q^3} + y^2 \right) = 0.
$$
For a new variable $z = y + 1/y$, we get the equation
$z^2 + \frac{4z}{{q}^3} + \frac{6}{{q}^4}-2 = 0$
which roots are 
$$
z_{\pm} = \frac{-\frac{4}{q^3}\pm\sqrt{\frac{16}{q^6}
 - 4\left(\frac{6}{q^4}-2\right)}}{2}
 = -\frac{2}{q^3}\pm\frac{1}{q^3}\sqrt{4 - 6p + 2p^3}.
$$

To solve the equation $y_{\pm}^2-z_{\pm}y_{\pm}+1 = 0$, write down
$$
D_{\pm}
 = z_{\pm}^2 - 4
 = \frac{2(4-3p-p^3\mp2\sqrt{2}\sqrt{2-3p+p^3})}{p^3},\quad
(y_{\pm})_{\pm} = \frac{z_{\pm}\pm\sqrt{D_{\pm}}}{2}.
$$
Since all numbers $(y_{\pm})_{\pm}$ are negative,
the smallest root of $\widetilde{C}(G_{q},y)$
should be the one among
$(y_{\pm})_{+} = -\frac{1}{2}(|z_{\pm}|-\sqrt{D_{\pm}})$.
Since the function $x-\sqrt{x^2-4}$ is monotonically decreasing 
for $x\in[2,\infty)$, the $(y_-)_+$ is the smallest root of 
$\widetilde{C}(G_{q},y)$ and 
\begin{multline*}
\beta(G_p) = -\frac{2{q}^3}{z_- +\sqrt{D_-}}
 = \frac{{q}^3(-z_- +\sqrt{D_-})}{2} \\
 = 1 + (1-p)\sqrt{\frac{2+p}{2}}
     + \sqrt{\frac{(1-p)(4+p+p^2+2\sqrt{2}\sqrt{2+p})}{2}}.
\end{multline*}

For $n = 5$, express 
\begin{multline*}
\widetilde{C}(G_{q},y)
 = 1 + \frac{5y}{p^2} + \frac{10y^2}{p^3} + \frac{10y^3}{p^3} + \frac{5y^4}{p^2} + y^5 \\
 = (1+y)\left(y^4 + y^3\left(\frac{5}{p^2}-1\right) + y^2\left(\frac{10}{p^3}-\frac{5}{p^2}+1\right)
+ y\left(\frac{5}{p^2}-1\right) + 1\right).
\end{multline*}
Changing a variable $z = y + 1/y$, we get the equation 
$$
z^2 + z\left(\frac{5}{p^2}-1\right) + \frac{10}{p^3}-\frac{5}{p^2}-1 = 0
$$
which roots are 
$$
z_{\pm}
 = \frac{p^2-5\pm\sqrt{5}\sqrt{5-8p+2p^2+p^4} }{2p^2}.
$$

To solve the equation $y_{\pm}^2-z_{\pm}y_{\pm}+1 = 0$, we calculate its discriminant
$$
D_{\pm}
 = z_{\pm}^2 - 4
 = \frac{5(5-4p-p^4)\mp(p^2-5)\sqrt{5}\sqrt{5-8p+2p^2+p^4} }{2p^4},\quad
(y_{\pm})_{\pm} = \frac{z_{\pm}\pm\sqrt{D_{\pm}}}{2}.
$$
Analogously to the case $n = 4$ we find 
$\beta(G_p) = p^2(-z_- +\sqrt{D_-})/2$ 
which implies~\eqref{n=5}.

\subsection{Random algebra}

Let $X =\{x_1,\ldots,x_n\}$ be a finite set. Fix an order on $X$ such that $x_i>x_j$ if $i<j$.
Consider a word $w = w_1w_2\ldots w_m\in X^*$ of length $m = |w|$.
Let a letter $x_{i_j}$ occurs in $w$ exactly $m_i\geq1$ times, $i=1,\ldots,k$.
We suppose that $x_{i_1}>x_{i_2}>\ldots>x_{i_k}$.
Consider a new alphabet
$$
X' = X'(w) = \{x_{i_1}^1,\ldots,x_{i_1}^{m_1},x_{i_2}^1,\ldots,x_{i_2}^{m_2},\ldots,x_{i_k}^1,\ldots,x_{i_k}^{m_k}\}.
$$
Define an order on the set $X'$: 
$x_{i_a}^s>x_{i_b}^t$ 
if $a<b$ or $a=b$ and $s<t$.

Given a word $w\in X^*$, let us construct a word $w'\in (X')^*$ of the same length as follows.
If $w_j$ is the $t$-th occurrence (counting from the left) of a letter $x_s$ in~$w$,
then the $j$-th letter of $w'$ equals $x_s^t$.
For example, 
$(abccdaba)' = a^1 b^1 c^1 c^2 d^1 a^2 b^2 a^3$.

Denote the set of all multipartite graphs with parts 
$\{x_{i_1}^1,\ldots,x_{i_1}^{m_1}\}$, \ldots, $\{x_{i_k}^1,\ldots,x_{i_k}^{m_k}\}$
as $MP(w)$. Define $M = M(w)$ be equal to the product $m_1m_2\ldots m_k$.

Let $p\in[0;1]$. Define a weight $s_p(w)$ of a word $w$ as 
\begin{equation}\label{word-weight}
s_p(w) = \sum\limits_{G\in MP(w)}p^{|E(G)|}(1-p)^{M-|E(G)|}
 I(w'\ \mbox{is in n.f. in }\,M(X',G)),
\end{equation}
where n.f. means ``normal form'',
$M = \prod\limits_{i=1}^k m_i$,
$I(A) = \begin{cases}
1, & A\ \mbox{is true}, \\
0, & \mbox{otherwise}.
\end{cases}$

Actually $s_p(w)$ equals a probability of the event that 
$w'$ is in the normal form 
in hypothetical partially commutative monoid
with commutativity graph $G_{p}(w)$, 
the random multipartite graph with fixed parts 
with $m_1,\ldots,m_k$ vertices and edge probabi\-lity~$p$.

{\bf Example 5.1}.
Let $w = acb$ and $a>b>c$. Note that $MP(w) = G_{3,p}$.
To calculate $s_p(w)$, we should avoid an edge $(b,c)$ 
in a graph $G$ with 3 vertices. Thus,
$$
s_p(w) 
 = (1-p)^3 + p(1-p)^2 + p(1-p)^2 + p^2(1-p) = 1-p.
$$

Call a word $w =w_1w_2\ldots w_k\in X^*$ monotonic if 
$w_i\geq w_{i+1}$ for all $i=1,\ldots,k-1$.
For example, $abbcccd$ is monotonic but not $acb$ if $a>b>c>d$.

{\bf Lemma 5.2}.
a) For a monotonic word $w\in X^*$, $s_p(w) = 1$ for all $p\in [0,1]$.
If $w\in X^*$ is not monotonic, then $s_{p_1}(w) < s_{p_2}(w)$
for any $p_1,p_2\in[0;1]$, $p_1>p_2$.

b) For any word $w\in X^*$ of length $m$, we have 
$s_p(w)\geq (1-p)^{m(1-1/n)}$.

c) If there is precisely $t$ different 
pairs of not equal neighbour letters in $w\in X^*$, then $s_p(w)\geq (1-p)^t$.
(We do not distinguish the pairs $ab$ and $ba$ in the word $abba$.)

{\sc Proof}.
a) If a word $w\in X^*$ is monotonic, then 
$s_p(w) = 1$ by Lemma~1.2 and \eqref{word-weight}.

Let $w$ be not monotonic, $G\in MP(w)$.
By Lemma~1.2, $w'$ is not in normal form in $M(X,G)$
if and only if there exists a pair $a,b\in X'$
such that $a<b$, $(a,b)\in E(G)$, $w' = xaybz$ 
for $x,y,z\in (X')^*$ and $(t,b)\in E(G)$ for all $t\in y$.
The probability that a triple $(a,b,y)\in X\times X\times X^*$ 
breaks normality of $w'$ in $G_p(w)$ equals a probability that 
we have a star graph $K(a,b,y)$ connecting $a$ and all letters from~$y$ to $b$.
Interpret the probability $P(a,b,y)$ of appearance of $K(a,b,y)$
like a geometric probability in the space $[0;1]^{\binom{n}{2}}$, $n = |X|$, 
where we relate an edge $e$ in $K_n$ to a line segment $[0;p]$ 
if $e$ belongs to $K(a,b,y)$ or to a line segment $[0;1]$ otherwise. 
Thus, $P(a,b,y)$ equals a volume of the constructed body. 
Hence, the probability that $w'$ is not in normal form
equals a union of all such bodies constructed for every $K(a,b,y)$.
By this observation we conclude that the greater edge probability~$p$ 
implies the greater probability that $w'$ is not in normal form.

b) Let $G\in MP(w)$ and $w = w_1w_2\ldots w_m$. 
Draw in $G$ all anti-edges (i.e., we draw edges in $\bar{G}$)
between $w_i$ and $w_{i+1}$ if $w_i<w_{i+1}$.
Then $w'$ is in normal form in $M(X',G)$ by Lemma~1.2.

For $m\leq n$, we have equal or less than $m-1$ 
pairs of neighbour letters in $w'$
and the statement follows from the inequality $m-m/n \geq m-1$.

Let $m>n$. Since $|X| = n$, 
there exists a pair of letters $w_i$, $w_{i+1}$ 
in a subword $w_1w_2\ldots w_{n+1}$ such that $w_i \geq w_{i+1}$.
We can find a such pair in a subword $w_{n+1}\ldots w_{2n+1}$ and so on.
Therefore, we need equal or less $m-m/n$ anti-edges to provide 
that $w'$ is in normal form.

c) Let $G\in MP(w)$, draw anti-edges between 
all pairwise distinct neighbour letters in~$w$.
Then $w'$ is normal form in $M(X',G)$ and thus $s_p(w)\geq (1-p)^t$. Lemma is proved.

Define on the free associative algebra $\Ass\langle X\rangle$
as on the vector space a new product $\cdot$.
Let $X_n$ denote the set of all words of length $n$ in the alphabet~$X$.
Due to distributivity, it is enough to define the product $\cdot$ on elements 
$w_i\in X_{n_i}$, $i=1,2$,
\begin{multline}\label{RandomProduct}
w_1\cdot w_2
 = \frac{1}{s_p(w_1)s_p(w_2)}
\sum\limits_{G\in MP(w_1w_2)}p^{|E(G)|}(1-p)^{M-|E(G)|}
I(w_1',w_2'\ \mbox{in n.f.})[(w_1w_2)'] \\
 = \frac{1}{s_p(w_1)s_p(w_2)}\sum\limits_{u\in X_{n_1+n_2}}
P(w_1',w_2'\ \mbox{in n.f.},\, u = [(w_1w_2)'])u \\
 = \sum\limits_{u\in X_{n_1+n_2}}P(u = [(w_1w_2)'] \mid w_1',w_2'\ \mbox{in n.f.})u,
\end{multline}
where $[w]$ denotes the normal form of~$w$,
$P(A)$ denotes the probability of an event $A$ 
in the probability theory model constructed 
by the random multipartite graph~$G_p(w_1w_2)$
and $P(A\mid B)$ denotes the conditional probability of~$A$ given~$B$. 

To avoid the division by~0, we define $\cdot$ for $p = 1$
by the third line of \eqref{RandomProduct}.
For $p<1$, we have $s_p(w)>0$ for any $w\in X^*$ by Lemma 5.2.

Let us call the space $\Ass\langle X\rangle$ under the product $\cdot$
defined by~\eqref{RandomProduct}
as partially commutati\-ve algebra with random commutativity graph
or as random algebra for short denoting it by $\Ass(X,p)$.

Define a weight function $s_p$ on all elements from $\Ass(X,p)$ by linearity.

{\bf Lemma 5.3}.
a) For $0\leq p<1$, the algebra $\Ass(X,p)$ is isomorphic to 
the free asso\-ciative algebra $\Ass\langle X\rangle$.
For $p = 1$, $\Ass(X,p)$ is isomorphic to the polynomial algebra $\mathbb{R}[X]$.

b) The map $s_p\colon \Ass(X,p)\to \langle\mathbb{R},\cdot\rangle$
is a semigroup homomorphism.

{\sc Proof}.
a) Note that for $p = 1$,
we have $s_p(w) = 1$ if $w$ is monotonic
and $s_p(w) = 0$ otherwise. 
Hence, $\Ass(X,p)$ under the product $\cdot$ 
defined by the third line of~\eqref{RandomProduct}
coincides with the polynomial algebra $\mathbb{R}[X]$.

Let $p<1$. Define a linear map 
$\varphi\colon \Ass\langle X\rangle\to \Ass(X,p)$
acting on a word $w$ as follows
\begin{equation}\label{RandomIso}
\varphi(w)
 = \sum\limits_{G\in MP(w)} p^{|E(G)|}(1-p)^{M-|E(G)|}[w'].
\end{equation}

Prove that $\varphi$ is an isomorphism of the algebras.
Let us state that $\varphi$ is injective.
To the contrary, $\varphi(x) = 0$ for $x\neq0$.
Express $x$ via the basis:
$x = \sum \lambda_i w_i$, $\lambda_i\neq0$, $w_i\in X^*$.
Choose the smallest word $w_0$ among $\{w_i\}$
in the lexicographic ordering of words in $X^*$.
Expressing $\varphi(x)$ via the basis, 
$\varphi(x)$ has a nonzero coefficient on $w_0$,
a~contradiction.
By~\eqref{RandomIso}, spaces of words of given length and 
given set of letters (with the numbers of occurrences) 
are $\varphi$-invariant.
Since $\varphi$ is injective and all such spaces are finite-dimensional,
$\varphi$ is bijective.

Let us prove that $\varphi$ is a homomorphism. Write down
\begin{equation}\label{RandomHom1}
\varphi(w_1w_2)
 = \sum\limits_{G\in MP(w_1w_2)}p^{|E(G)|}(1-p)^{M-|E(G)|}[(w_1w_2)'],
\end{equation}
\vspace{-0.5cm}
\begin{multline}\label{RandomHom2}
\varphi(w_1)\cdot \varphi(w_2) \\
 = \sum\limits_{\substack{H_i\in MP(w_i),\\ G\in MP(w_1w_2)}}
 p^{|E(G)|+|E(H_1)|+|E(H_2)|}(1-p)^{M_G+M_{H_1}+M_{H_2}-|E(G)|-|E(H_1)|-|E(H_2)|}\\
\times
\frac{1}{s_p([w_1'])s_p([w_2'])}I_G([w_1'],[w_2'])[([w_1']_{H_1}[w_2']_{H_2})'],
\end{multline}
where $I_{G}(w)$ denotes $I(w\ \mbox{in n.f. in }M(X',G))$,
$M_G = M(w_1w_2)$, $M_{H_i} = M(w_i)$, $i=1,2$, 
$[w']_H$ equals the normal form of $w'$ in $M(X',H)$.

A choice of $G\in MP(w_1w_2)$ is equivalent 
to an independent choice of graphs $G_i\in MP(w_i)$, $i=1,2$,
and $G_{12}\in MP(w_1-w_2))$, where $MP(w_1-w_2)$
denotes the set of all subgraphs $MP(w_1w_2)$ 
in which all vertices of $w_1'$ as well as $w_2'$
are pairwise disconnected.

Rewrite~\eqref{RandomHom1},~\eqref{RandomHom2}:
\begin{multline}\label{RandomHom1-2}
\varphi(w_1w_2)
 = \sum\limits_{\substack{G_i\in MP(w_i),\\ G_{12}\in MP(w_1-w_2)}}
 p^{|E(G_1)|+|E(G_2)|+|E(G_{12})|} \\
(1-p)^{M_{G_1}+M_{G_2}+M_{G_{12}}-|E(G_1)|-|E(G_2)|-|E(G_{12})|}
[([w_1']_{G_1}[w_2']_{G_2})'],
\end{multline}
\begin{multline}\label{RandomHom2-2}\allowdisplaybreaks
\varphi(w_1)\cdot \varphi(w_2)
 = \sum\limits_{\substack{H_i\in MP(w_i),\\ G_{12}\in MP(w_1-w_2)}}
 p^{|E(H_1)|+|E(H_2)|+|E(G_{12})|} \\
 \times (1-p)^{M_{H_1}+M_{H_2}+M_{G_{12}}-|E(H_1)|-|E(H_2)|-|E(G_{12})|} \\
\times \prod\limits_{i=1}^2\bigg(
 \frac{1}{s_p([w_i'])}
 \sum\limits_{G_i\in MP(w_i)}p^{|E(G_i)|}(1-p)^{M_{G_i}-|E(G_i)|}
I_{G_i}([w_i'])\bigg)[([w_1']_{H_1}[w_2']_{H_2})'] \\
= \sum\limits_{\substack{H_i\in MP(w_i),\\ G_{12}\in MP(w_1-w_2)}}
 p^{|E(H_1)|+|E(H_2)|+|E(G_{12})|} \\
 (1-p)^{M_{H_1}+M_{H_2}+M_{G_{12}}-|E(H_1)|-|E(H_2)|-|E(G_{12})|} [([w_1']_{H_1}[w_2']_{H_2})'].
\end{multline}

By change of ``variables'' $G_i\leftrightarrow H_i$, $i=1,2$,
the RHS of~\eqref{RandomHom1-2}, \eqref{RandomHom2-2}
coincide. Thus, we get $\varphi(w_1w_2) = \varphi(w_1)\cdot\varphi(w_2)$.

b) It is enough to prove that 
$s_p(w_1\cdot w_2)/(s_p(w_1)s_p(w_2)) = 1$
for any $w_i\in X_{n_i}$, $i=1,2$.

By the third line from~\eqref{RandomProduct}, we calculate 
\begin{multline}\label{WeightCalcul}
s_p(w_1\cdot w_2)
 = \sum\limits_{u\in X_{n_1+n_2}}
  P(u = [(w_1w_2)']\mid w_1',w_2' \mbox{ in n.f.})s_p(u) \\
 = \sum\limits_{u\in X_{n_1+n_2}}
  P(u = [(w_1w_2)']\mid w_1',w_2' \mbox{ in n.f.})
  \sum\limits_{G\in MP(w_1w_2)}p^{|E(G)|}(1-p)^{M-|E(G)|}
  I(u\ \mbox{in n.f.}) \\
 {=} \sum\limits_{G\in MP(w_1w_2)}p^{|E(G)|}(1-p)^{M-|E(G)|}
 \bigg(
 \sum\limits_{u\in X_{n_1+n_2}}
  P(u = [(w_1w_2)']\mid w_1',w_2' \mbox{ in n.f.})I(u\ \mbox{in n.f.}) \bigg) \\
 =  \sum\limits_{G\in MP(w_1w_2)}p^{|E(G)|}(1-p)^{M-|E(G)|}
I(w_1',w_2'\ \mbox{in n.f.}).
\end{multline}

Applying \eqref{WeightCalcul}, we get 
\begin{multline*}
\frac{s_p(w_1\cdot w_2)}{s_p(w_1)s_p(w_2)} 
= \frac{1}{s_p(w_1)s_p(w_2)}
\sum\limits_{G\in MP(w_1w_2)}p^{|E(G)|}(1-p)^{M-|E(G)|}
I(w_1',w_2'\ \mbox{in n.f.}) \\
 = \frac{1}{s_p(w_1)s_p(w_2)}P(w_1',w_2'\ \mbox{in n.f. in }M(X',G(w_1w_2))) \\
 = \frac{1}{s_p(w_1)s_p(w_2)}P(w_1'\ \mbox{in n.f. in }M(X',G(w_1w_2)))
P(w_2'\ \mbox{in n.f. in }M(X',G(w_1w_2))) = 1,
\end{multline*}
as the events $\{w_1'\ \mbox{in n.f. in }M(X',G(w_1w_2))\}$
and $\{w_2'\ \mbox{in n.f. in }M(X',G(w_1w_2))\}$ are independent
and their probabilities do not interchange when we consider
the supergraph $G(w_1w_2)$ instead of $G(w_1)$ and $G(w_2)$ respectively.

{\bf Corollary 5.2}.
For $w_i\in X_{n_i}$, $i=1,2,3$, we have in $\Ass(X,p)$
$$
w_1\cdot w_2\cdot w_3
 = \sum\limits_{t\in X_{n_1+n_2+n_3}}
P(t = [(w_1w_2w_3)'] \mid w_1',w_2',w_3'\ \mbox{in n.f. in }M(X',G(w_1w_2w_3)))t.
$$

{\bf Remark 5.2}.
The algebra $\Ass(X,p)$ is a trivial deformation~\cite{Gerstenhaber}
of the free associative algebra $\Ass\langle X\rangle$,
as for $p = 0$ we have $w_1\cdot w_2 = w_1w_2$
and $\Ass(X,p) = \Ass\langle X\rangle$.

Given a subset $D =\{d_1,\ldots,d_k\}$ of $X^*$,
define a weight $s_p(D)$ of $D$ be equal to $s_p(d_1)+\ldots+s_p(d_k)$.
Let $m_t(p) = s_p(X_t)$.

{\bf Example 5.2}. 
Denote $n = |X|$. Let us calculate $m_t(p)$ for $t=1,2$: 
\begin{gather*}
m_1(p) = \sum\limits_{i=1}^n s_p(x_i) 
       = \sum\limits_{i=1}^n 1 	= n, \\  
m_2(p) = \sum\limits_{i=1}^n s_p(x_i^2)
       {+} \sum\limits_{x_i>x_j}s_p(x_i x_j)
       {+} \sum\limits_{x_i<x_j}s_p(x_i x_j)
       = n {+} \binom{n}{2} {+} (1-p)\binom{n}{2}
       = n^2 {-} \binom{n}{2}p.
\end{gather*}

{\bf Theorem 5.2}.
The numbers $m_t(p)$, $t\geq1$, satisfy the reccurence relation
\begin{multline}\label{rekurrent2}
m_t(p)
= \binom{n}{1} m_{t-1}(p) - \binom{n}{2}p m_{t-2}(p) + \ldots \\
 + (-1)^{k+1}\binom{n}{k} p^{\binom{k}{2}} 
 + \ldots + (-1)^{n+1}p^{\binom{n}{2}}m_{t-n}(p)
\end{multline}
with initial data $m_0(p) = 1$, $m_{-1}(p) = \ldots = m_{-n+1}(p) = 0$.

{\sc Proof}.
The proof will be very similar to the proof of Theorem~1.1.
Define the sets 
$$
N_i = \{p^{\binom{i}{2}}u\cdot(x_{j_1}x_{j_2}\ldots x_{j_i}) \mid u \in M_{t-i}, x_{j_s}\in X,
x_{j_1} > x_{j_2} > \ldots > x_{j_i}\},\quad i=1,\ldots,n,
$$
where $\cdot$ is the product in $\Ass(X,p)$.

Given a word $w\in X_t$, define
\begin{multline}\label{li-definition}
l_i(w) = \sum\limits_{\substack{u\in X_{t-i}, x_{j_s}\in X,\\ x_{j_1} > \ldots > x_{j_i}}}
P([(u x_{j_1}x_{j_2}\ldots x_{j_i})'] = w,\, u \mbox{ in n.f.}, x_{j_s}\sim x_{j_t},s\neq t) \\
= \sum\limits_{G\in MP(w)} p^{|E(G)|}(1-p)^{M_G-|E(G)|}I_G(w) \\
\times \bigg( 
\sum\limits_{\substack{u\in X_{t-i}, x_{j_s}\in X,\\ x_{j_1} > \ldots > x_{j_i}}}
I_G(w' = (u x_{j_1}x_{j_2}\ldots x_{j_i})')I_G(u)I_G(x_{j_s}\sim x_{j_t},s\neq t) \bigg),
\end{multline}
$$
l_0(w) = \sum\limits_{i=1}^n(-1)^{i+1}l_i(w),
$$
where $x\sim y$ denotes that $x$ and $y$ are connected.
Denote by $\tilde{k}_i(w)$ the expression from the last line of~\eqref{li-definition}.

Applying Lemma~5.3b, we get
\begin{multline}\label{rek2:preRightPart}\allowdisplaybreaks
\sum\limits_{w\in X_t}l_i(w)
 = \sum\limits_{w\in X_t}\sum\limits_{\substack{u\in X_{t-i}, x_{j_s}\in X,\\ x_{j_1} > \ldots > x_{j_i}}}
P([(u x_{j_1}x_{j_2}\ldots x_{j_i})'] = w,\, u \mbox{ in n.f.}, x_{j_s}\sim x_{j_t},s\neq t) \\
 = \sum\limits_{u\in X_{t-i}}\sum\limits_{x_{j_1} > \ldots > x_{j_i}}
   \sum\limits_{G\in MP(ux_{j_1}\ldots x_{j_i})} p^{|E(G)|}(1-p)^{M_G-|E(G)|}I_G(u)I_G(x_{j_s}\sim x_{j_t},s\neq t) \\
 = \sum\limits_{u\in X_{t-i}}\sum\limits_{H_1\in MP(u)} p^{|E(H_1)|}(1-p)^{M_{H_1}-|E(H_1)|}I_{H_1}(u) \\
   \times 
   \sum\limits_{x_{j_1} > \ldots > x_{j_i}}\sum\limits_{H_2\in MP(x_{j_1}x_{j_2}\ldots x_{j_i})} p^{|E(H_2)|}(1-p)^{M_{H_2}-|E(H_2)|}I_{H_2}(x_{j_s}\sim x_{j_t},s\neq t)
   \\
  = \bigg(\sum\limits_{u\in X_{t-i}}s_p(u)\bigg)\times p^{\binom{i}{2}}\sum\limits_{x_{j_1} > \ldots > x_{j_i}}1
  = \binom{n}{i}p^{\binom{i}{2}}m_{t-i}(p) = s_p(N_i).
\end{multline}

Therefore, by \eqref{rek2:preRightPart}
\begin{multline}\label{rek2:RightPart}
\binom{n}{1} m_{t-1}(p) - \binom{n}{2}p m_{t-2}(p) + \ldots + (-1)^{n+1}p^{\binom{n}{2}}m_{t-n}(p) \\
 = \sum\limits_{i=1}^n (-1)^{i+1}s_p(N_i)
 = \sum\limits_{w\in X_t}(-1)^{i+1}l_i(w) = \sum\limits_{w\in X_t}l_0(w).
\end{multline}

Suppose that $l_0(w) = s_p(w)$ for all $w \in X^*$,
then the RHS of~\eqref{rekurrent2} by~\eqref{rek2:RightPart}
equals 
$\sum\limits_{w\in X_t}l_0(w) = \sum\limits_{w\in X_t}s_p(w) = m_t(p)$,
and we are done.

Given a word $w\in X_t$, we prove that $l_0(w) = s_p(w)$
by the following equalities
\begin{multline*}
\sum\limits_{i=1}^n(-1)^{i+1}l_i(w)
 = \sum\limits_{i=1}^n(-1)^{i+1}
  \sum\limits_{G\in MP(w)}
   p^{|E(G)|}(1-p)^{M_G-|E(G)|}I_G(w)\tilde{k}_i(w) \\
 = \sum\limits_{G\in MP(w)}
 p^{|E(G)|}(1-p)^{M_G-|E(G)|}I_G(w)
 \sum\limits_{i=1}^n(-1)^{i+1}\tilde{k}_i(w) \\
 = \sum\limits_{G\in MP(w)}p^{|E(G)|}(1-p)^{M_G-|E(G)|}I_G(w) = s_p(w),
\end{multline*}
since for $I_G(w) = 1$ the coefficients $\tilde{k}_i(w)$ equal 
the numbers $k_i(w')$ for partially commuta\-tive monoid $M(X',G)$ 
from the proof of Theorem~1.1. By Lemma~1.3 we are done.

{\bf Corollary 5.3}.
The polynomial $PC(G_{n,p},x)$ is a characteristic polynomial
for the sequence $\{m_t(p)\}$ and $\beta(G_{n,p})$ equals its growth rate.

{\bf Lemma 5.4}.
a) The following inequalities hold
\begin{equation}\label{BoundsOnBetaNP}
1+(n-1)(1-p)\leq \beta(G_{n,p})\leq 1+(n-1)\sqrt{1-p}.
\end{equation}

b) The number $\beta(G_{n,p})$ for fixed~$n$
is strictly monotonic function on $p\in[0;1]$ decreasing from $n$ to~1.

{\sc Proof}.
a) The upper bound follows from Theorem~5.1a and Samuelson's Inequality. Prove the lower bound.

Let $X = V(G)$.
Consider the random graph $H$ with $n$ vertices, edge probability $1-p$
and loops in all vertices. Denote by $W_t = W_t(H)$ the expected value
of the numbers of walks in $H$ of length $t$.
Let us show the inequality $W_t\leq m_t(p)$ where $m_t(p) = s_p(X_t)$.
Given a word $w\in X_t$, denote by $t$ the number of different pairs
of not equal neighbour letters in~$w$.
Then the walk $w$ will be counted in a sum $W_s$ as $(1-p)^k$.
By Lemma~5.2c, $s_p(w)\geq (1-p)^k$.

By Fekete's Lemma, there exists the limit $\lim\limits_{t\to \infty}\sqrt[t]{W_t(G)}$.
By \cite{Cvetkovich,Harary}, 
$\lim\limits_{t\to \infty}\sqrt[t]{W_t(G)} = \rho(G)$
for any simple graph~$G$, where $\rho(G)$ equals the largest eigenvalue of the 
adjacency matrix~$A(G)$ (so called spectral radius of~$G$).
This result could be easily proved for~$H$,
since it is the corollary of the Perron---Frobenius theory.

The adjacency matrix of $H$
$$
A = A(H) = (a_{ij}),\quad
a_{ij} = \begin{cases}
1,& i=j, \\
1-p, & i\neq j,\end{cases}
$$
is symmetric. It is easy to state that 
$\rho(H) = 1 + (n-1)(1-p)$.
Finally, we have 
$$
\beta(G_{n,p}) = \lim\limits_{t\to \infty}\sqrt[t]{m_t(p)}
\geq \lim\limits_{t\to \infty}\sqrt[t]{W_t(H)}
= \rho(H) = 1 + (n-1)(1-p).
$$

b) By Lemma~5.2a and Corollary~5.3, we have that 
$\beta(G_{n,p})$ is not strictly monotonic function.
Suppose that there exist $p_1<p_2$, $p_1,p_2\in[0;1]$, such that
$\beta(G_{n,p_1}) = \beta(G_{n,p_2})$. Then we have 
$\beta(G_{n,p_1}) = \beta(G_{n,p})$ for all $p\in[p_1;p_2]$.
Thus, polynomials $PC(G_{n,p_1},x)$ and $PC(G_{n,p},x)$
are not coprime for all $p\in[p_1;p_2]$.
It means that the resultant of this pair of polynomials
as a polynomial on $p$ has infinite number of zeros, a contradiction. 

{\bf Remark 5.3}.
Let us show that two natural lower bounds are worse than~\eqref{BoundsOnBetaNP}.
At first, the lower bound on $\beta(G_{n,p})$ arisen from Lemma~5.2b
is always not better than the lower bound from Lemma~5.4a, as 
\begin{multline*}
n(1-p)^{1-1/n}
= n(1-p)\left(1+\frac{p}{1-p}\right)^{1/n} \\
\leq n(1-p)\left(1+\frac{p}{n(1-p)}\right)
= n(1-p)+p
= 1+(n-1)(1-p).
\end{multline*}
By Theorem~5.1b, we can get the bound 
$\beta(G_{n,p})>n(1-p)$ with the help of geometric progression.
But this bound is still worse than one we have considered above.

{\bf Corollary 5.4}.
a) For $p = 1/2$, we have 
\begin{equation}\label{RateRandomHalf}
\frac{n}{2}<\beta(G_{n,1/2})<\frac{n}{\sqrt{2}}+1.
\end{equation}

b) \cite{BrownNowak}
The polynomial $C(G_{n,1/2})$
equals the arithmetic mean of clique polynomials 
computed for all graphs with $n$ vertices.

{\sc Proof}. a) It follows from Lemma~5.4.

b) Any clique $C$ of size $k$ is calculated in $C(G_{n,1/2})$ 
with coefficient $2^{-\binom{k}{2}}$. The clique $C$
is calculated in the sum 
$\frac{1}{2^{\binom{n}{2}}}\sum\limits_{G\colon|V(G)|=n} C(G,x)$
with the same coefficient
$\frac{2^{\binom{n}{2}-\binom{k}{2}}}{2^{\binom{n}{2}}} = 2^{-\binom{k}{2}}$.

In 2012, J. Brown et al~\cite{Brown2} got the following bounds 
on the smallest root $\alpha(G_{n,p})$ of $PC(G_{n,p},x)$:
$$
\frac{p^{n-1}}{n}\leq \alpha(G_{n,p})\leq \frac{p^{n-1}}{1+\sqrt{1-p}}.
$$

Below we improve both bounds.

{\bf Corollary 5.5}.
For any $p\in(0,1)$, we have the following bounds
$$
\frac{p^{n-1}}{1+(n-1)\sqrt{1-p}}
\leq \alpha(G_{n,p})\leq \frac{p^{n-1}}{1+(n-1)(1-p)}.
$$

{\sc Proof}.
It follows from Statement~5.1b and Lemma~5.4.

\newpage
\section{Planar graphs}

Denote the set of all planar graphs with $n$ vertices and $k$ edges as $Pl(n,k)$.

Let find minimal and maximal values of $\beta(G)$ for $G\in Pl(n,k)$.
We also show graphs for these extremal values.

For any planar graph, the degree of $PC(G,x)$ is not higher than~4.
Nore the following upper bounds on the clique numbers~\cite{Wood2007}
for any planar graph with $n$ vertices:
\begin{equation}\label{Wood}
c_2\leq 3n-6,\quad c_3\leq 3n-8,\quad c_4\leq n-3.
\end{equation}

{\sc Case 1}. 
For $n<3$ and any~$k$ or $n\geq3$ and $0\leq k\leq 2$
either $PC(G,x) = x-n$ or $PC(G,x) = x^2 - nx + k = 0$.
In both variants, we have 
\begin{equation}\label{k0-2}
\lambda_-(n,k) = \lambda_+(n,k)
 = \frac{n+\sqrt{n^2-4k}}{2}.
\end{equation}

{\sc Case 2}.
Consider the four particular variants of values of $n,k$.
For $n = k = 3$, we have $PC(G,x) = (x-1)^3$ and
\begin{equation}\label{kn3}
\lambda_-(3,3) = \lambda_+(3,3) = 1.
\end{equation}
Analogously, 
\begin{equation}\label{k6n4}
\lambda_-(4,6) = \lambda_+(4,6) = 1.
\end{equation}

For $n = 4$, $k = 5$, we have the planar graph $K_{1,1,2}$.
Thus, by Example~2.1c, 
\begin{equation}\label{k5n4}
\lambda_-(4,5) = \lambda_+(4,5) = 2.
\end{equation}

For $n = 5$, $k = 9$, we have the planar graph $K_{1,1,1,2}$.
By Example~2.1c, 
\begin{equation}\label{k9n5}
\lambda_-(5,9) = \lambda_+(5,9) = 2.
\end{equation}

{\sc Case 3}.
Let $n\geq4$, $3\leq k\leq 2n-4$.
We can construct a bipartite graph $G_-$ consisting of two parts 
with 2 and $n-2$ vertices. Indeed, the complete bipartite 
graph $K_{2,n-2}$ contains $2(n-2)\geq k$ edges.
Then $PC(G_-,x) = x^2 - nx + k = 0$ and
\begin{equation}\label{2min}
\lambda_-(n,k) = \frac{n+\sqrt{n^2-4k}}{2}.
\end{equation}
Note that $\lambda_-$ as the function on $k$ for fixed~$n$
is strictly monotonic decreasing.
The formula~\eqref{2min} gives $\lambda_- = n-1$ for $k = n-1$ and 
$\lambda_- = n-2$ for $k = 2n-4$.

{\bf Lemma 6.1}.
Let $G\in Pl(n,k)$,
$n\geq4$, $k\leq 2n-4$, $c_3(G)>0$, then 
$\beta(G)>\lambda_-(n,k)$.

{\sc Proof}.
Let $\lambda = \lambda_-(n,k)$, then 
$$
PC(G,\lambda) = - c_3(G)\lambda+c_4(G) < 0,
$$
as $c_4\leq n-3$ by~\eqref{Wood} and $\lambda\geq n-2$ by Lemma~2.4.

{\sc Case 4}.
Let $n\geq5$, $2n-4 < k < 3n-6$.
Construct the graph $G_-$ as a supergraph of $K_{2,n-2}$
with parts $V_1 = \{e_1,e_2\}$ and $V_2 = \{s_1,s_2,\ldots,s_{n-2}\}$ 
in which $k-2n+4 < n-2$ edges in $V_2$ are drawn in such way that
these edges form a tree with $k-2n+5$ vertices inside $V_2$
(see Picture~1). 

Then 
\begin{equation}\label{3a-min}
PC(G_-,x) = x^3 - nx^2 + kx - 2(k-2n+4) = 0.
\end{equation}
One of the roots of the equation~\eqref{3a-min} equals~2, so,
\begin{equation}\label{3a-min'}
\lambda_-(n,k) = \frac{n-2+\sqrt{(n-2)^2-4(k-2n+4)}}{2}
  = -1+\frac{n+\sqrt{n^2+4n-4k-12}}{2}.
\end{equation}
In particular, $\lambda_-(n,3n-7) = n-3$.

{\bf Lemma 6.2}.
Let $G\in Pl(n,k)$, $n\geq5$, $k = 2n-4+s$, $1\leq s\leq n-3$. Then 

a) $c_3(G)\geq 2s$, 

b) $c_3(G)>2s$ if $c_4(G)>0$.

{\sc Proof}.
a) To the contrary, suppose that $c_3 = 2s-t$, $t>0$. 
By the Euler's formula, $G$ has $f \geq 2+k-n=n-2+s$ faces. 
Since there are not greater than $2s-t$ triangle faces in $G$,
remaining $n-2-s+t$ (or more) faces contain at least 4 vertices on their border.
Then the doubled number of all edges in $G$ equals, 
from the one hand, $2k = 4n-8+2s$.
From the other hand, $2k\geq N$, where $N$ equals the number of all edges
of which consist borders of all faces. By the following
$$
N\geq 3(2s-t)+4(n-2-s+t) 
 = 4n-8+2s+t>4n-8+2s
$$
we arrive at a contradiction.

b) Let $c_4(G)>0$, find a subgraph $K_4$ in $G$.
In every of four faces on which the subgraph $K_4$ divides the plane,
there is $k_i$ vertices and $l_i$ edges of $G$, $i=1,2,3,4$, 
different from the ones of the chosen $K_4$. 
We have that every face of $K_4$ together with the border 
forms a planar graph with $k_i+3$ vertices and $l_i+3$ edges.

By a), the joint number of all triangles in all four faces 
is not less than the following number doubled:
$$
\sum\limits_{i=1}^4 (l_i+3-2(k_i+3)+4)
 = \sum\limits_{i=1}^4(l_i -2k_i +1)
 = (2n-4+s-6) -2(n-4) + 4
 = s + 2.
$$
Thus, $c_3(G)\geq 2s+4$.

{\bf Corollary 6.1}.
Let $G\in Pl(n,k)$, $n\geq5$, $2n-4<k<3n-6$, then 
$\beta(G)\geq\lambda_-(n,k)$.

{\sc Proof}.
We apply Lemma~6.2. The proof is analogous to the 
proof of Lemma~6.1 with only one refinement: 
$\lambda_-(n,k)\geq n-3$, so we have a not strict upper bound
$\beta(G)\geq\lambda_-(n,k)$.

{\sc Case 5}. 
Let $n\geq6$, $k = 3n-6$.
Construct the graph $G_-$ as a supergraph of $K_{2,n-2}$
with parts $V_1 = \{e_1,e_2\}$, $V_2 = \{s_1,s_2,\ldots,s_{n-2}\}$
in which the part $V_2$ is a cycle with $n-2$ vertices.
To draw this graph on plane, we put the vertex $e_1$ 
inside the cycle and we put $e_2$ outside (see Picture~2).

\begin{figure}[h]
\centering
\includegraphics[height = 4cm]{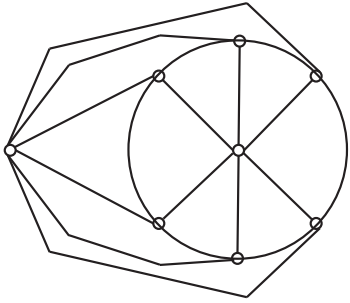}
\caption{Graph $G_-$ for the case $n\geq5$, $k = 3n-6$.} \label{Pic2}
\end{figure}

We have
\begin{equation}\label{3b-min}
PC(G_-,x) = x^3 - nx^2 + (3n-6)x - 2n+4 = 0.
\end{equation}
We find that 2 is a root of the equation~\eqref{3b-min}, so
\begin{equation}\label{3b-min'}
\lambda_-(n,3n-6) = \frac{n-2+\sqrt{(n-2)^2-4(n-2)}}{2}
= \frac{n+\sqrt{n^2-8n+12}}{2} -1.
\end{equation}
Note that 
$\beta(G_-) = \lambda_-(n,3n-6) = \beta(C_{n-2})$
for the cycle $C_{n-2}$ with $n-2$ vertices.

{\bf Lemma 6.3}.
Let $G\in Pl(n,3n-6)$, $n\geq6$, then $\beta(G)\geq\lambda_-(n,3n-6)$.

{\sc Proof}.
The graph~$G$ is planar, otherwise we can draw at least one additional edge
preserving planarity, a contradiction to \eqref{Wood}.
By the Euler's formula, $G$ has $2n-4$ faces.
Since $G$ is maximal planar graph, all faces are triangle.
Thus, $c_3(G)\geq 2n-4 = c_3(G_-)$.

Let $\lambda = \lambda_-(n,3n-6)$.
If $c_4(G) = 0$, then $PC(G,\lambda)\leq 0$
and $\beta(G)\geq \lambda$.
Otherwise,
\begin{equation}\label{min3n-6}
PC(G,\lambda) = -(c_3(G)-2n+4)\lambda +c_4(G).
\end{equation}
Since $G$ contains $K_4$, hence,
$G$ contains a separating triangle (it means which is a border of any face).
So, $c_3(G)-2n+4\geq1$.
For $n\geq6$, we have $\lambda\geq n-4$.
If $c_4(G)\leq n-4$, then $PC(G,\lambda)\leq 0$ by~\eqref{min3n-6}.
If $c_4(G) = n-3$, then $c_4(G)\geq3$ for all $n\geq6$.
Therefore, we can find at least two separating triangles in $G$ and 
$$
PC(G,\lambda) = -(c_3(G)-2n+4)\lambda +c_4(G)
 \leq -2(n-4)+n-3 = 5 - n < 0.
$$

{\sc Case 6}.
Let $n\geq4$, $k\geq3$.
Construct a graph $G_+$ as follows.
Initially we draw $K_3$.
On each step, we add one new vertex inside of some 
(triangle) face of the graph and connect it with each vertex of the face.
We proceed on while we have edges (see Picture~3). 
Sometimes, the constructed graph is called an Apollonian network. 
For $k = 3n - 6$, this graph maximizes the number 
of all cliques for planar graphs~\cite{Wood2007}.

\begin{figure}[h]
\centering
\includegraphics[height = 3.8cm]{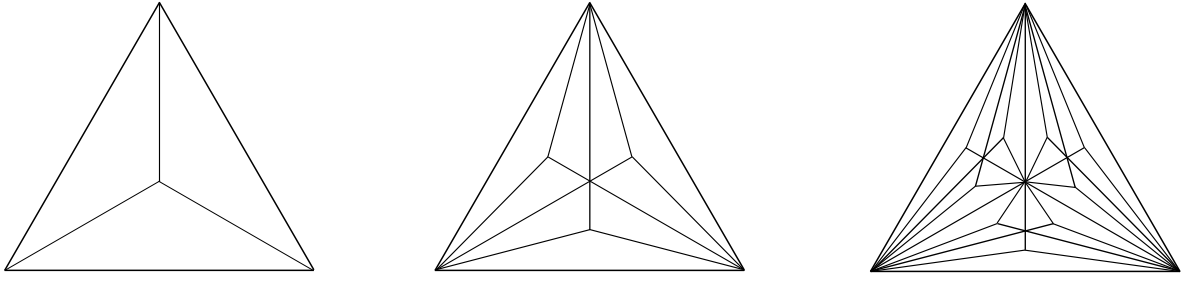}
\caption{Graph $G_+$ for the case $n\geq4$, $k\geq3$.} \label{Pic3}
\end{figure}

For $3\leq k<6$, $G_+$ contains no cliques of size~4, so
\begin{equation}\label{2a-max}
PC(G_+,x) = x^3 - nx^2 + kx - (1+[(k-3)/2].
\end{equation}

If $6\leq k\leq 3n-6$, we have
\begin{equation}\label{2b-max}
PC(G_+,x) = x^4 - nx^3 + kx^2
- \left(1+\bigg[\frac{k-3}{3}\bigg]+\bigg[\frac{2(k-3)}{3}\bigg]\right)x
+ \bigg[\frac{k}{3}\bigg] - 1.
\end{equation}
In particular, for $k = 3s$:
\begin{equation}\label{2b-max0}
PC(G_+,x) = x^4 - nx^3 + 3s x^2 - (3s-2)x + s - 1,
\end{equation}
for $k = 3s+1$:
\begin{equation}\label{2b-max1}
PC(G_+,x) = x^4 - nx^3 + (3s+1) x^2 - (3s-2)x + s - 1,
\end{equation}
for $k = 3s+2$:
\begin{equation}\label{2b-max2}
PC(G_+,x) = x^4 - nx^3 + (3s+2) x^2 - (3s-1)x + s - 1,
\end{equation}

Denote the largest root of the equations \eqref{2a-max} and \eqref{2b-max} 
(i.e., of $PC(G_+,x)$) as $\lambda_+(n,k)$.

For $k = 3n-6$, by~\eqref{2b-max0} we write down
\begin{equation}\label{3-max}
PC(G_+,x) = x^4 - nx^3 + 3(n-2)x^2 - (3n-8)x + n-3 = 0.
\end{equation}
The value $x = 1$ is a root of $PC(G_+,x)$.
Thus, we find that $\lambda_{+}(n,3n-6) = n-3$.

{\bf Lemma 6.4}.
Let $G\in Pl(n,k)$, $n\geq4$, $k\geq3$,
then $c_3(G)\leq c_3(G_+)$.

{\sc Proof}.
We prove the statement by induction on~$n$.
For $n = 4$, it is easy to show that the statement is true.

If $k\leq 2n-1$, then there exists a vertex~$v$
of degree $d(v)\leq 3$ in a planar graph~$G$.
Hence, we can reorder edges in $G$
in such way that the number of triangles will be not less and 
$v$ will be isolated. So, we may apply the induction hypothesis.

If $G$ is disconnected, then we are done
by the induction hypothesis and by~\eqref{2b-max}.

For $k\geq 2n$, we have $c_3(G_+)\geq2n-3$.
We may assume that $G$ is connected. 
If $c_3(G)>c_3(G_+)$, then 
by the Euler's formula $G$ has a separating triangle~$T$.
Let $G_1$ ($G_2$) be a subgraph of $G$ induced by 
$V(T)$ and all vertices lying inside (outside)~$T$.
Applying the induction hypothesis for $G_1$ and $G_2$, 
we prove the inductive step.

{\bf Lemma 6.5}.
Let $G\in Pl(n,k)$, $n\geq4$, $k\geq3$.
If $c_4(G) < c_4(G_+)$, then $c_3(G) < c_3(G_+)$.

{\sc Proof}.
Prove the statement by induction on~$n$. 
For $n = 4$, the statement if checked directly.
We may assume that $G$ is connected, otherwise
we apply the induction hypothesis.
If $G$ contains no separating triangles, then 
by the Euler's formula $n + c_3(G)\leq k-2$ and so $c_3(G)<c_3(G_+)$.
Let $T$ be a separating triangle in~$G$.
Let $G_1$ ($G_2$) be a subgraph of $G$ induced by 
$V(T)$ and all vertices lying inside (outside)~$T$.
Let $G_i$ contain $n_i$ vertices and $k_i$ edges,
$G_{+,i} = G_+(n_i,k_i)$, $i=1,2$.

There are $k = k_1+k_2+3$ edges and 
$[\frac{k_1+k_2}{3}]$ cliques of size~4 in~$G_+$.
Let $c_4(G) < c_4(G_+)$, then 
we have either $c_4(G_i)<c_4(H_i)$ for some $i\in\{1,2\}$
and we may apply the induction hypothesis and Lemma~6.4
or $c_4(G_i) = c_4(H_i)$, $i=1,2$,
\begin{equation}\label{C4forG+}
c_4(G_1)+c_4(G_2)
 = \bigg[\frac{k_1}{3}\bigg] + \bigg[\frac{k_2}{3}\bigg]
 < \bigg[\frac{k_1+k_2}{3}\bigg].
\end{equation}
The last inequality in~\eqref{C4forG+} holds only if 
by division on~3 the numbers $k_1,k_2$ have reminders 
either 1,2, or 2,1, or 2,2 respectively.
In all cases, we can rearrange one edge in $G$ 
to get a new graph $G'$ such that
the numbers of edges inside (outside) $T$ 
will give reminders 0,0 or 0,1 by division on~3.
Applying Lemma~6.4, we have inequalities $c_3(G)<c_3(G')\leq c_3(G_+)$.

{\bf Corollary 6.2}.
Let $G\in Pl(n,k)$, $n\geq4$, $k\geq3$,
then $\beta(G)\leq \lambda_{+}(n,k)$.

{\sc Proof}.
Note that by Lemma~2.4, 
$\lambda = \lambda_{+}(n,k)\geq \lambda_{+}(n,3n-6) = n-3$.
If $c_3(G) = c_3(G_+)$, then by Lemma~6.5
$c_4(G)\geq c_4(G_+)$. So, 
$$
PC(G_+,\beta(G))
 = c_4(G_+) - c_4(G)\leq0
$$
and $\lambda\geq \beta(G)$.

Let $c_3(G)<c_3(G_+)$. 
Supposing that $\beta(G)>\lambda\geq n-3$, we have by~\eqref{Wood}
$$
PC(G_+,\beta(G))
 = -(c_3(G)-c_3(G_+))\beta(G) + (c_4(G_+)-c_4(G))
 < -(n-3) + n-3 = 0,
$$
a contradiction.

{\bf Remark 6.1}.
The following equality holds
$$
\lambda_-(n,3n-7) = \lambda_+(n,3n-6),
$$
i.e., moving edges in a planar graph we can ``eat'' the whole edge.

{\bf Remark 6.2}. 
In 1991 N.B. Boots with G.F. Royle~\cite{Geography}
and in 1993 D. Cao with A.~Vince~\cite{PlanarSpectrum} independently posed 
the conjecture that the graph 
$G_+(n,3n-6) = K_2\cup P_{n-2}$ 
has the maximal spectral radius among $n$-vertex planar graphs.
Although the conjecture fails for some small~$n$, in 2017 
M. Tait and J. Tobin solved this conjecture~\cite{ThreeConj}
showing that the graph $G_+(n,3n-6)$ is the unique planar graph on $n$~vertices with 
maximum spectral radius for all sufficiently large $n$.
We believe that this result could be proved (with concrete lower bound on $n$)
with the help of the Kelmans transformation and beyond (see~\S8.1).

Another conjecture of~D. Cao with A.~Vince~\cite{PlanarSpectrum}
says that the graph 
$G_-(n,3n-6) = \bar{K}_2\cup C_{n-2}$ 
has the maximal spectral radius among $n$-vertex planar graphs
with the minimal degree~4.

{\bf Theorem 6.1}.
Let $G\in Pl(n,k)$, then
$$
\lambda_-(n,k) = \min\limits_{G\in Pl(n,k)}\beta(G),\quad
\lambda_+(n,k) = \max\limits_{G\in Pl(n,k)}\beta(G),
$$
where the numbers $\lambda_{\pm}(n,k)$ are defined by \eqref{k0-2}--\eqref{2min},
\eqref{3a-min'}, \eqref{3b-min'}, \eqref{2a-max} and \eqref{2b-max}.
All these bounds are reached.

{\sc Proof}.
If $n < 3$ or $k\leq2$, see Case~1.
If $n = 3$ and $k = 3$, see Case~2.

Thus, we may consider only $n\geq4$ and $k\geq3$.
For the maximal value of $\beta(G)$ for $G\in Pl(n,k)$,
see Corollary~6.2.

For the minimal value of $\beta(G)$ for $G\in Pl(n,k)$,
for $n = 4$ see Case~2 when $k=5,6$ and Lemma~6.1 when $k\leq 4$.
For $n = 5$ see Case~2 when $k = 9$ and Corollary~6.2 when $k\leq 8$
(the case $n = 5$, $k = 10$ gives $K_5$ which is not planar).
For $n\geq6$ see Lemma~6.3.

Let $Pl(n)$ denote the set of all planar graphs with $n$ vertices.
In 2007, O. Gim\'{e}nez and M. Noy stated~\cite{Noy} that
$|Pl(n)|\sim C_0n^{-7/2}\gamma^n n!$
for $\gamma \approx 27.227$ and 
a planar graph in average contains $\kappa \approx 2.213 n$ edges. 
Moreover, the expected values and dispersions of the values 
$c_3$ and $c_4$ in a planar graph with $n$ vertices equal respectively
\begin{equation}\label{RandomPl}
\mu_n(K_3)\sim\frac{n}{6\gamma^3},\quad \mu_n(K_4)\sim\frac{n}{24\gamma^4},\quad
\sigma_n^2(K_3),\sigma_n^2(K_4)\sim \gamma n.
\end{equation}

Introduce the $PC$-polynomial of random planar graph as follows
$$
PC(Pl(n),x)
 = x^4-nx^3+\kappa n x^2-\frac{nx}{6\gamma^3}+\frac{n}{24\gamma^4}.
$$

It is easy to verify that for $n\gg1$, $PC(Pl(n),x)$ has a root 
$$
n-\kappa-\left(\frac{6\gamma^3\kappa^2-1}{6\gamma^3}\right)\frac{1}{n}
 + O\bigg(\frac{1}{n^2}\bigg)\approx n-\kappa.
$$

Since $PC(Pl(n),\kappa)<0$ for $n\gg1$ and $PC(Pl(n),0)>0$,
another root lies in the interval $(0,\kappa)$.

{\bf Lemma 6.6}. Let $n\gg1$.

a) The polynomial $PC(Pl(n),x)$ has two real and two complex non-real roots.
Complex roots lie in the right half-plane, real roots are simple.

b) Let $\rho$ be a root of $PC(Pl(n),x)$, then 
$|\rho|\geq 1/(6\gamma^3\kappa)$.

{\sc Proof}.
a) By Hurwitz stability criterion~\cite[Theorem 11.4.5]{Hurwitz}, 
to show that all roots of $PC(Pl(n),x)$ lie in the right half-plane,
it is enough to check the following inequalities:
$$
n>0,\quad
\kappa n^2>\frac{n}{6\gamma^3},\quad
\frac{\kappa n^3}{6\gamma^3} > \frac{n^3}{24\gamma^4}+\frac{n^2}{6\gamma^3},\quad
\frac{n}{24\gamma^4}>0,
$$
which are clearly fulfilled.

The Sturm series calculated for $PC(Pl(n),x)$ via its derivative
gives two real roots and two complex non-real roots.

b) It follows from the Enestr\"{o}m---Kakeya Theorem applied for 
$PC(Pl(n),-1/x)$.

{\bf Theorem 6.2}.
a) $PC(G,x)$ of almost all planar graphs is a polynomial of 4-th degree,
has two complex and real roots. Complex roots lie in the right half-plane,
real roots are simple.

b) Let $\varepsilon>0$, then for almost all planar graphs
we have $|\rho|>\frac{1-\varepsilon}{6\gamma^3\kappa}$
for any root $\rho$ of $PC$-polynomial.

{\sc Proof}.
By~\eqref{RandomPl} and Chebyshev's inequality,
we nay assume that for almost all planar graphs
with $n$ vertices 
$$
c_k(G)=\varepsilon_k c_k(PC(Pl(n),x)),\quad k=2,3,4,
$$
where $c_k(PC(Pl(n),x))$ are the coefficients 
of $PC(Pl(n),x)$ by the degree $x^{4-k}$
and $\varepsilon_k$ lie in some neighbourhood of~1.
It remains to repeat the proof of Lemma~6.6.

{\bf Remark 6.3}.
Let us consider the example showing that roots 
of $PC$-polynomial of planar graph could lie in the left half-plane.
For $G = K_4\cup \bar{K}_{96}$, the polynomial 
$PC(G,x) = x^4-100x^3+6x^2-4x+1$ has the 
following approximate values of the roots due to~\cite{wolfram}
$$
0,17028,\quad 99,940,\quad -0,05532\pm 0,23601i.
$$

{\bf Statement 6.1}.
The average value of the growth rate of partially commutative monoid
with planar commutativity graph equals
$$
\beta_{\mathrm{ev},Pl}(n)
 =\frac{1}{|Pl(n)|}\sum\limits_{G\in Pl(n)}\beta(G) = n-\kappa + O(1/n).
$$

{\sc Proof}.
By the construction of $\lambda_{\pm}(n,k)$, we have
\begin{equation}\label{lambda-pm}
\begin{gathered}
\lambda_{-} = n-\frac{k+4}{n} + \frac{2k}{n^2}-\frac{k^2}{n^3} + O(1/n^2),\quad k>2n-4, \\
\lambda_{+} = n-\frac{k}{n} + \frac{k}{n^2}-\frac{k^2}{n^3} + O(1/n^2),
\end{gathered}
\end{equation}
from which the statement follows.

{\bf Remark 6.4}.
By~\eqref{lambda-pm}, we get that the difference
$\lambda_{+}(n,k)-\lambda_{-}(n,k)$ in average asymptotically equals 
$\dfrac{4}{n}-\dfrac{\kappa}{n}\approx \dfrac{1.787}{n}$.

{\bf Statement 6.2}.
Let $G$ be a planar graph, $|V(G)| = n$, $|E(G)|\leq\sqrt{\frac{16n}{13}}$,
$w(G) = 3$. Then $PC(G,x)$ has one real root and two complex non-real roots.

{\sc Proof}.
It is clear that $n\geq 8$, otherwise $k = |E(G)| < 3$ and $c_3 = 0$.
Applying the inequality $1\leq c_3 < k$,
compute the discriminant of the cubic equation
$-PC(G,-x) = x^3+nx^2+kx+c_3 = 0$:
$$
\Delta
 = n^2k^2-4k^3-4n^3c_3-27c_3^2+18nkc_3
 < \frac{16}{13}n^3 -4n^3+18\cdot\frac{16}{13}n^2
 = \frac{36}{13}n^2(8-n)\leq0.
$$
Hence, $PC(G,-x)$ as well as $PC(G,x)$ has only one real root.

\newpage
\section{Lower bounds on $\beta(G)$}

By analogy with the well-known Nordhaus---Gadddum inequalities~\cite{Zykov,NordhausGaddum}
\begin{equation}\label{NordhausGadddum}
\begin{gathered}
2 \sqrt{n} \leq \chi(G) + \chi(\bar{G}) \leq n+1,\\
n \leq \chi(G) \chi(\bar{G}) \leq \left(\frac{n+1}{2}\right)^2
\end{gathered}
\end{equation}
for chromatic number of a graph~$G$ and its complement, we are interested on 

{\bf Problem~1}. To find the tight bounds for 
the expressions $\beta(G)+\beta(\bar{G})$ and $\beta(G)\beta(\bar{G})$.

Denote by $G(n,k)$ the set of all graphs with $n$ vertices and $k$ edges. Introduce 
$$
\beta_-(n,k) = \min\limits_{G\in G(n,k)}\beta(G),\quad
\beta_+(n,k) = \max\limits_{G\in G(n,k)}\beta(G).
$$

{\bf Problem~2}.
To find the values $\beta_{\pm}(n,k)$ and the graphs on which they are reached.

Given $G\in G(n,k)$, we denote by $\bar{k}$ the number $|E(\bar{G})| = \binom{n}{2}-k$.
For fixed $n$ and~$k$, the value $\beta(G)$ could essentially exchange depending on~$G$.
For example, let $\bar{G}_1 = K_{2,2,\ldots,2}\in G(n,k)$ with $\bar{k} = n/2$
and $\bar{G}_2 = K_{s,1,1,\ldots,1}\in G(n,k)$ with $n/2 = \binom{s}{2}$.
Then the ratio $\frac{\beta(G_2)}{\beta(G_1)} = \frac{s}{2}$
could be any arbitrarily large natural number when $n\gg1$.

For a graph $G\in G(n,k)$, define $e(G)$, the edge $PC$-density of $G$, as follows:
$$
e(G) = \frac{n-\beta(G)}{k}.
$$
We will see that this parameter could not essentially exchange and has excellent bounds.
For $k = 0$, define $e(\bar{K}_n) = \frac{1}{n}$.

{\bf Example 7.1}.
We have the following values of $e(G)$ 

a) for the complete multipartite graph with $p$ equal parts, 
$$
e(K_{n/p,n/p,\ldots,n/p}) = \frac{n-\frac{n}{p}}{\frac{n(n-n/p)}{2}} = \frac{2}{n};
$$

b) for a tree with $n$ vertices,
$$
e(T_n) = \frac{n-(n-1)}{n-1} = \frac{1}{n-1};
$$

c) for the planar graph $G_+(n,3n-6)$ (see Case~6 from \S6),
$$
e(G_+) = \frac{n-(n-3)}{3n-6} = \frac{1}{n-2};
$$

d) for the random planar graph $PC(Pl(n),x)$,
$$
e(PC(Pl(n),x)) 
 \sim \frac{n-(n-\kappa)}{\kappa n} = \frac{1}{n};
$$

d) for the graph $\overline{K_{s,1,1,\ldots,1}}$ with $n$ vertices,
$$
e(\overline{K_{s,1,1,\ldots,1}})
 = \frac{n-s}{\binom{n}{2}-\binom{s}{2}}
 = \frac{2(n-s)}{(n-s)(n+s-1)}
 = \frac{2}{n+s-1};
$$

e) for random graph with $n$ vertices, by Lemma~5.4, 
$$
\frac{1}{n}\leq e(G_{n,p})\leq \frac{2}{n}.
$$

Given a simple graph~$G$, define the graph $\bar{G}^{loop}$ with loops (without multiedges):
$$
V(\bar{G}^{loop}) = V(G),\quad
E(\bar{G}^{loop}) = \{(u,v)\mid (u,v)\not\in E(G),u\neq v\}\cup \{(u,u)\mid u\in V(G)\}.
$$

The following lower bound for $\beta(G)$ (part a)) was proved by D.C. Fisher in 1989.

{\bf Theorem 7.1}.
a) \cite{Fisher1989} Given a graph $G$ with $n$ vertices and $k$ edges, 
$\beta(G)\geq n-\frac{2k}{n}$. 

b) The bound from a) is reached if and only if $G$ is an empty graph or 
$G$ is a complete multipartite graph with equal parts.

{\sc Proof}.
a) The proof will be very similar to the proof of Lemma~5.4a. 
Let $H = \bar{G}^{loop}$. 
Denote by $W_s = W_s(H)$ the number of walks of length~$s$ in $H$.
The inequality $W_s\leq m_s(G)$ holds, since 
a word from partially commutative monoid $M(V(G),G)$
corresponding to a walk from $W_s$ is always in normal form by Lemma~1.2.
As it was noticed in the proof of Lemma~5.4 , we have
$\lim\limits_{s\to \infty}\sqrt[s]{W_s(G)} = \rho(G)$,
where $\rho(G)$ is the spectral radius of~$G$.
It is easy to show \cite[Theorem~1.1]{StevanovicMonograph}
that $W_s$ equals the sum of all entries of the matrix $A^s$,
where $A = A(H)$ is the adjacency matrix of~$H$.

By the construction of~$H$, we have 
$\rho(H) = 1 + \rho(\bar{G})$.
Applying the known lower bound 
$\rho(\bar{G})\geq \frac{2\bar{k}}{n}$~\cite{Collatz} for simple graphs,
we conclude
$$
\beta(G) = \lim\limits_{s\to \infty}\sqrt[s]{m_s(G)}
 \geq \lim\limits_{s\to \infty}\sqrt[s]{W_s(G)}
 = \rho(H) = 1 + \rho(\bar{G})\geq 1 + \frac{2\bar{k}}{n}
 = n - \frac{2k}{n}.
$$

b) By Example~1.1, the lower bound is reached 
on empty graph or on a complete multipartite graph with equal parts.
Let $G$ be a graph with $n$ vertices, $k$ edges, $k>0$,
and $\beta(G) = n-\frac{2k}{n} = 1 + \frac{2\bar{k}}{n} = 1 +\rho(\bar{G})$.
It is known that \cite{Collatz} the bound $\rho(\bar{G}) = \frac{2\bar{k}}{n}$
is attained only if $\bar{G}$ is a regular graph.
If $G$ is a complete multipartite graph (with not necessary equal parts),
then regularity implies that all parts of~$G$ are equal. 

Let $G$ be not complete multipartite graph, so 
$G$ is a supergraph of a complete multipartite graph with parts 
$G_1,G_2,\ldots,G_r$, and $\bar{G}_i$ is connected for all $i=1,\ldots,r$.
By Lemma~1.5e, $PC(G,x) = \prod\limits_{i=1}^r PC(G_i,x)$.
Define $n_i = |V(G_i)|$, $\bar{k}_i = |E(\bar{G}_i)|$, $i=1,\ldots,r$.
From the following relations
$$
1 +\frac{2\bar{k}}{n}
 = \beta(G) \geq 1+\max\limits_{i}\frac{2\bar{k}_i}{n_i}
 \geq 1 + \frac{2(\bar{k}_1+\ldots +\bar{k}_r)}{n_1+\ldots+n_r}
 = 1 +\frac{2\bar{k}}{n},
$$
we conclude that 
$\frac{2\bar{k}_1}{n_1} = \frac{2\bar{k}_2}{n_2} = \ldots = \frac{2\bar{k}_r}{n_r}$.
Thus, $\beta(G) = \beta(G_1) = \ldots = \beta(G_r)$ and all $G_i$ are regular
and not empty. 

Consider one of parts $G_i$. We may assume that $i=1$
Let us show that for any $v\in V(G_1)$
there exist $u,w\in V(G_1)$ such that 
\begin{equation}\label{K-1,2}
(v,u),(v,w)\in E(\bar{G}_1),\quad (u,w)\not \in E(\bar{G}_1).
\end{equation}
Consider $d = \frac{2\bar{k}_1}{n_1}$ vertices from the neighbourhood $N(v)$ of~$v$.
If any two vertices from $N(v)$ are connected, then 
by regularity and connectedness of $\bar{G}_1$, we get $\bar{G}_1 = K_{d+1}$,
a contradiction. 

Let $H_1 = \bar{G_1}^{loop}$, $\rho(H_1) = 1 +\rho(\bar{G}_1)$.
To get a contradiction, it is enough to prove that 
$\rho(H_1)<\beta(G_1)$.

Denote by~$B$ the adjacency matrix of~$H_1$,
$$
\mathrm{Spec}(B)
 = \{\rho(H_1) = \lambda_1>\lambda_2\geq \lambda_3\geq\ldots\geq\lambda_{n_1}\},
$$
the inequality $\lambda_1>\lambda_2$ holds due to 
the Perron---Frobenius theory~\cite[p.~8]{StevanovicMonograph},
as the graph $H_1$ is connected.
For $D = \mathrm{diag}\{\lambda_1,\lambda_2,\ldots,\lambda_{n_1}\}$,
there exists an orthogonal matrix~$S$ such that $B = S^{-1}DS$.

A walk $w_1w_2\ldots w_t$ is called closed if $w_1 = w_t$.
The number of closed walks of length~$t$ equals 
$\tr(B^t) = \sum\limits_{i=1}^{n_1}\lambda_i^t$~\cite{Harary}.
Since 
\begin{multline*}
B^t = S^{-1}D^tS
 = S^{-1}\mathrm{diag}\{\lambda_1^t,\lambda_2^t,\ldots,\lambda_{n_1}^t\}S \\
 = \lambda_1^t S^{-1}\mathrm{diag}\left\{1,\frac{\lambda_2^t}{\lambda_1^t},\ldots,\frac{\lambda_{n_1}^t}{\lambda_1^t}\right\}S
 \sim \lambda_1^tS^{-1}\mathrm{diag}\{1,0,\ldots,0\}S
\end{multline*}
for $t\gg1$, we can find a vertex $v\in V(H_1)$ such that 
the number $W_t(v)$ of closed walks of length $t$ starting with $v$ 
is not less than $\lambda_1^t/(n+1)$ for $t\gg1$.
Thus, there exists a positive constant $z$ such that 
$W_t(v)\geq \lambda_1^t/z$ for all $t\geq1$.

Consider $t\gg1$. For $v$, find $u,w\in V(H_1)$ 
such that~\eqref{K-1,2} holds. Therefore, $(u,w)\in E(G_1)$
and $(v,u),(v,w)\not\in E(G_1)$. Define on the set $V(G_1)$ an order with $u > w$.
By Lemma~1.2, the word $Y = vuwv$ is in normal form.
Construct a set of normal words from $M_t(V(G_1),G_1)$ as follows.
For any $0\leq s\leq t/4$, we consider the word $Y^s$ of length $4s$.
In all $s+1$ positions between the blocks of~$Y$ 
we put $t-4s$ letters in such way that a~word~$x$ written between neighbour $Y$-blocks
safisfies the condition that the word $vxv$ is a walk in $H_1$
(a word $x$ could be empty).
All obtained words of length~$t$ are in normal form and pairwise distinct, 
so we may estimate 
\begin{multline*}
m_t(V(G_1),G_1)
 \geq \sum\limits_{s=0}^{t/4}
 \binom{s+t-4s}{s} \frac{\lambda_1^{t-4s} }{z^{s+1}}
 = \frac{\lambda_1^t}{z}\sum\limits_{s=0}^{t/4}
 \binom{t-3s}{s} \frac{1}{\big(\lambda_1^4z\big)^s} \\
 \geq \frac{\lambda_1^t}{z}\sum\limits_{s=0}^{t/4}
 \binom{t/4}{s} \frac{1}{\big(\lambda_1^4z\big)^s}
 = \frac{\lambda_1^t}{z}\left(1+\frac{1}{\lambda_1^4z}\right)^{t/4}.
\end{multline*}
Hence, 
$$
\beta(G_1)
 \geq \lambda_1\left(1+\frac{1}{\lambda_1^4z}\right)^{1/4} > \lambda_1 = \rho(H_1).
$$

{\bf Corollary 7.1}.
For any graph $G$ with $n$ vertices, 
$e(G)\leq\frac{2}{n}$.

{\bf Corollary 7.2}.
Let $G$ be a graph with $n$ vertices and $k$ edges. 

a) If $k\leq \frac{n^2}{4}$, then $\beta(G)\geq\frac{n}{2}$.

b) If $k\leq \frac{n(n-1)}{4}$, then $\beta(G)\geq\frac{n+1}{2}$.

{\bf Corollary 7.3}.
For any graph $G$ with $n$ vertices, 

a) $n+1\leq \beta(G)+\beta(\bar{G})$,

b) $n\leq \beta(G)\beta(\bar{G})$.

\noindent Moreover, the bounds are reached only if and only if 
$\{G,\bar{G}\} = \{K_n,\bar{K}_n\}$.

{\sc Proof}.
a) By Theorem~7.1a, 
$$
\beta(G)+\beta(\bar{G})\geq n-\frac{2k}{n} + n-\frac{2\bar{k}}{n}
= 2n-\frac{2 \binom{n}{2} }{n} = 2n - (n-1) = n+1.
$$
We have equality if and only if both $G$ and $\bar{G}$
are empty or complete multipartite grfaph with equal parts.
It coud happen only if $\{G,\bar{G}\} = \{K_n,\bar{K}_n\}$.

b) Analogously, by Theorem~7.1 we have
$$
\beta(G)\beta(\bar{G})\geq \bigg(n-\frac{2k}{n}\bigg)\bigg(n-\frac{2\bar{k}}{n}\bigg)
= n^2 - n(n-1)+\frac{4k\bar{k}}{n^2} \geq n.
$$
The equality holds only if $k\bar{k} = 0$, i.e., $\{G,\bar{G}\} = \{K_n,\bar{K}_n\}$.

{\bf Remark  7.1}.
Corollary~7.3b could be derived from Corollary~7.2b and Lemma~2.6b.

{\bf Corollary 7.4}.
For any graph $G$ with $n$ vertices and $k$ edges,

a) $n+\frac{4k\bar{k}}{n^2}\leq\beta(G)\beta(\bar{G})$,
in particular for $k=\frac{n(n-1)}{4}$ we have 
$\big(\frac{n+1}{2}\big)^2\leq\beta(G)\beta(\bar{G})$,

b) $2(n-1)\leq\beta(G)\beta(\bar{G})$ if $\{G,\bar{G}\} \neq \{K_n,\bar{K}_n\}$.

{\sc Proof}.
a) The statement was actually proved in Corollary~7.3b.

b) It follows from Lemma~2.6b and Corollary~7.3a.

{\bf Remark 7.2}.
Let $G$ be a graph with $n$ vertices and $k$ edges.
In light of Theorem~7.1 and Lemma~2.6a, 
Tur\'{a}n's Theorem could be interpreted as follows:
the lower bound $\frac{n}{\omega(G)}$ on $\beta(G)$
is always not better than the one $n-\frac{2k}{n}$.

{\bf Corollary 7.5} \cite{Fisher1989}.
For any graph $G$, we have $\beta(G)\geq 1 +\rho(\bar{G})$.

{\bf Statement 7.1}.
Let $k\leq\frac{n^2}{4}$, then 
$\beta_-(n,k) = \frac{n+\sqrt{n^2-4k}}{2}$.
This value is reached for a~graph $G\in G(n,k)$
if and only if $G$ is a triangle-free graph.

{\sc Proof}.
If $G$ is a triangle-free graph with $n$ vertices and $k$ edges,
then $PC(G,x) = x^2-nx+k$ and 
$\beta(G) = A = \frac{n+\sqrt{n^2-4k}}{2}\geq\frac{n}{2}$.
Let $H$ be a graph with $n$ vertices and $k\leq \frac{n^2}{4}$ edges.
By Lemma~2.7 we have $\beta^2-n\beta+k\geq0$
which by Corollary~7.2 implies $\beta\geq A$. 
Suppose that $c_3 = c_3(H) > 0$.
If $c_4 = c_4(H) = 0$, then $PC(H,A) = -c_3 A < 0$ and hence $\beta>A$.
For $c_4 > 0$, by Lemma~2.7, we have
$\beta^4-n\beta^3+k\beta^2-c_3\beta+c_4\geq0$.
If $\beta = A$, then $c_4\geq c_3A\geq \frac{nc_3}{2}$,
a contradiction to the simple inequality $c_4\leq \frac{n-3}{4}c_3$.

{\bf Statement 7.2}.
Let $k = \frac{n^2}{2}\big(1-\frac{1}{w}\big)$, $w\in \mathbb{N}$, $w\geq3$,
then $\beta_-(n,k) = n-\frac{2k}{n}$.
This value is reached for a graph $G\in G(n,k)$
if and only if $G$ is the complete multipartite graph with equal $w$ parts.

{\sc Proof}.
It follows from Theorem~7.1 and Example 7.1 .

{\bf Corollary 7.6} \cite{Fisher1989,Razborov}.
Let $G$ be a graph with $n$ vertices and 
$\frac{n^2}2\leq k\leq \frac{n^2}3$ edges.
Then $c_3\geq (9nk-2n^3-2(n^2-3k)^{3/2})/27$.

{\sc Proof}.
Consider the function $g(x) = x^3 - nx^2 + kx - c_3$. 
By Lemma~2.7,
$g(\beta)\leq 0$, where $\beta = \beta(G)$.
By the conditions and Theorem~7.1, we have
$n/3\leq s = n-2k/n\leq \beta$, so,
$\min\limits_{x\geq s} g(x)\leq 0$.
Since $g''(x) = 3x-n \geq0$ for $x\geq s$,
the unique critical point of $g(x)$ for $x\geq s$
is the point $x_0 = (n+\sqrt{n^2-3k})/3$. Thus, 
$$
0\geq \min\limits_{x\geq s} g(x) = g(x_0)
 = \frac{9nk-2n^3-2(n^2-3k)^{3/2}}{27} - c_3.
$$

Let us finish the paragraph with the lower bound on $\beta(G)$ 
proven by D.C. Fisher and J.M. Nonis 
which is stronger than the obtained one in Theorem~7.1.

{\bf Lemma 7.1}~\cite{FisherUpper}.
Let $A_1,A_2,\ldots,A_n$ be a sequence of finite sets. Then for all $s\in\mathbb{N}$,
\begin{equation}\label{NonisPreIneq}
\binom{s}{2}\bigg|\bigcup\limits_{1\leq i\leq n}A_i \bigg|  
 - s\sum\limits_{1\leq i\leq n}|A_i| + \sum\limits_{1\leq i<j\leq n}|A_i\cap A_j|\geq0.
\end{equation}

{\sc Proof}.
Let an element $a$ be in exactly $m$ sets of the $A_i$'s. Then $a$ is counted 
$\binom{s+1}{2} - sm + \binom{m}{2} = \binom{s-m+1}{2}\geq 0$ 
times on the left side of~\eqref{NonisPreIneq}.

Given a partially commutative monoid $M(X,G)$, let us denote by $M_j(w)$
the set of words in $M(X,G)$ of the length $j$ that can end with $w$.
So, $M_j = M_j(\emptyset)$ and remind that $m_j = |M_j|$.
Oce can prove (in a manner similiar to Lemma~1.1) that
\begin{equation}\label{MjX}
M_j(x)\cap M_j(y) = \begin{cases}
M_j(xy), & (x,y) \in E, \\
\emptyset, & \mbox{otherwise}.
\end{cases}
\end{equation}

{\bf Lemma 7.2}~\cite{FisherUpper}. 
Given a graph~$G$ with $n$ vertices and $k$ edges, we have for all $s\in\mathbb{N}$,

a) $\binom{s+1}{2}m_j(G)-sn m_{j-1}(G) + km_{j-2}(G)\geq 0$, 

b) $\binom{s+1}{2}\beta^2(G)-sn \beta(G) + k\geq 0$.

{\sc Proof}.
a)  With the help of~\eqref{MjX}, we state that 
\begin{multline*}
\binom{s+1}{2} m_j = \binom{s+1}{2} \bigg|\bigcup\limits_{x\in X}M_j(x)\bigg|
 \geq s\sum\limits_{x\in X}|M_j(x)| - \sum\limits_{x,y\in X,\,x\neq y}|M_j(x)\cap M_j(y)| \\
  = s \sum\limits_{x\in X}|M_j(x)| - \sum\limits_{(x,y)\in E}|M_j(xy)|
  = sn m_{j-1} - km_{j-2}.
\end{multline*}

b) It is enough to divide the inequalities from a) by $m_{j-2}$ 
and consider the limit $j\to\infty$.

{\bf Theorem 7.2}~\cite{FisherUpper}.
Given a graph $G$ with $n$ vertices and $k$ edges, let $w$ be such a~natural number that
$\big(1-\frac{1}{w-1}\big)\frac{n^2}{2}<k\leq \big(1-\frac{1}{w}\big)\frac{n^2}{2}$. 
Then 
\begin{equation}\label{LowerBoundBest}
\beta(G)\geq \frac{n}{w}\left(1+\sqrt{1 -\frac{2kw}{n^2(w-1)}}\right).
\end{equation}

{\sc Proof}.
By Theorem~7.1,
$$
\beta(G) \geq n-\frac{2k}{n} \geq n - n\left(1-\frac{1}{w}\right) = \frac{n}{w}.
$$
Thus, the statement follows from Lemma~7.2b for $s = w-1$.

\newpage
\section{Upper bound on $\beta(G)$}

\subsection{Maximum value $\beta_+(n,k)$}

Let us define a relation $\geq$ on simple graphs as follows.
We write that $G\geq H$ if $D(G,x)\leq D(H,x)$ 
on the line segment $[0;1/\beta(G)]$. 
From $G\geq H$, we have $\beta(G)\geq \beta(H)$.

In the proofs of Lemmas~2.3 and 2.4 (via Lemma~1.5),
we actually have stated that

{\bf Lemma 8.1} \cite{CsikvariKelman}.
a) Given an induced subgraph $H$ of $G$, we have $G\geq H$; \\
b) Given a spanning subgraph $H$ of $G$, we have $H\geq G$.

Let $G$ be a graph, $u,v\in V(G)$, $u\neq v$.
In 1981, A. Kelmans defined \cite{Kelmans} 
so called Kelmans transformation
which transfers a graph $G$ into a graph $G' = KT(G,u,v)$.
To get $G'$, we erase all edges between $v$ and $N(v)\setminus (N(u)\cup\{u\})$
and add all edges between $u$ and $N(v)\setminus (N(u)\cup\{u\})$
(see Picture~4).
The obtained graph has the same number of edges as $G$.
In~2011, P. Csikv\'{a}ri stated~\cite{CsikvariKelman} 
some important properties of Kelmans transformations.

\begin{figure}[h]
\centering
\includegraphics[height = 5cm]{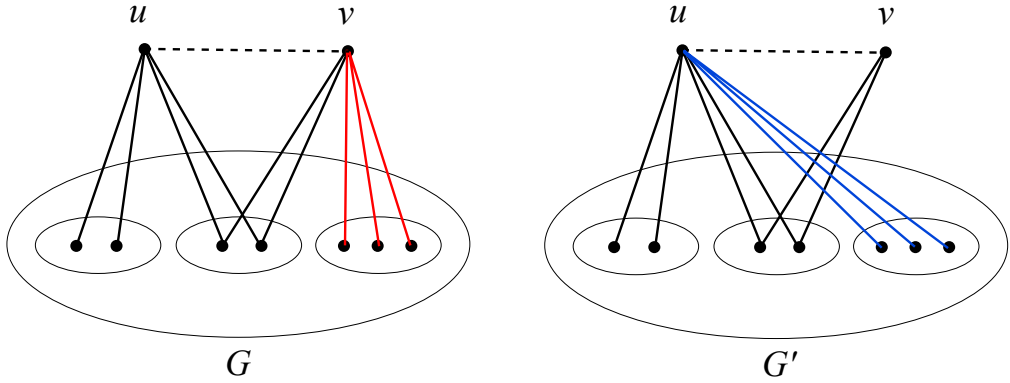}
\caption{The Kelmans transformation $G' = KT(G,u,v)$.} \label{Pic4}
\end{figure}

{\bf Lemma 8.2} \cite{CsikvariKelman,CsikvariThesis}.
Let $G$ be a graph and $G'$ be a graph obtained from~$G$ by a Kelmans transformation.
Then 

a) $G'\geq G$ and so, $\beta(G')\geq\beta(G)$,

b) $c_k(G')\geq c_k(G)$ for all $k$.

{\sc Proof}.
Let $G' = KT(G,u,v)$ and we may assume that 
$N(u)\setminus N(v)\neq\emptyset$,
otherwise $G$ and $G'$ are isomorphic.

a) Prove the statement by induction on $n = |V(G)|$.
For $n=1,2$, the statement is true as $G' = G$.
For $w\in N(u)\setminus N(v)$, write down by Lemma~1.5 ,
\begin{gather}\label{KT:Dep}
D(G,x)  = D(G \setminus w,x) - xD(G[N(w)],x),\\
D(G',x) = D(G'\setminus w,x) - xD(G'[N(w)],x).
\end{gather}

Since $G'\setminus w = KT(G\setminus w,u,v)$,
we have $G\setminus w\leq G'\setminus w$
by the induction hypothesis.
Note that $G[N(w)]$ is isomorphic to the spanning subgraph of $G'[N(w)]$.
Hence, by Lemma~8.1b and~\eqref{KT:Dep} we derive 
$$
D(G,x)\geq D(G',x),\quad x\in[0;1/\max\{\beta(G'\setminus w),\beta(G[N(w)])\}].
$$

We have $G'\geq G'\setminus w$ by Lemma~8.1a, as $G'\setminus w$ is an induced subgraph of $G'$.
Since $G'[N(w)\cup\{v\}]\setminus \{u\}$ is a spanning subgraph of $G[N(w)]$,
by Lemma~8.1 we have 
$$
G'\geq G'[N(w)\cup\{v\}]\setminus \{u\}\geq G[N(w)].
$$
Altogether, we get 
$\beta(G')\geq \max\{\beta(G'\setminus w),\beta(G[N(w)])\}$.
It implies $D(G,x)\geq D(G',x)$ for $x\in[0;1/\beta(G')]$, i.e., $G'\geq G$.

b) Let $A_k$ and $A_k'$ denote the set of all cliques of size $k$ in $G$
and $G'$ respectively. It is easy to see that $S\in A_k\setminus A_k'$
has to contain $v$ and at least one vertex from $N_G(v)\setminus N_G(u)$.
So, $u$ is not a vertex of~$S$.
For any such clique~$S$, we have a new clique $S'\in A_k'\setminus A_k$
obtained by the replacement of $v$ on $u$. 
Thus, $|A_k'|\geq|A_k|$. Moreover, we can say that 
$|A_k'|-|A_k|$ equals the number of cliques of size~$k$ in $G'$
such that contain $u$, at least one vertex from $N_G(v)\setminus N_G(u)$
and at least one vertex from $N_G(u)\setminus N_G(v)$.
Lemma is proved.

A threshold graph is a graph that can be constructed from 
a single vertex by repeated applications of the following operations:
a) addition of a isolated vertex,
b) addition of a vertex connected to all previous vertices.
By the definition, we may identify a threshold graph $G$ with $n$ vertices 
with a vector from $\thr(G)\in\mathbb{Z}_2^{n-1}$ as follows.
The $i$-th coordinate of the vector equals~0 if 
the $i$-th vertex appeared in~$G$ is isolated and~1, otherwise. 

It is easy to show that (see \cite{CsikvariNikif})
any graph by a series of Kelmans transformations
can be transformed to a threshold graph.

Let $G$ be a graph such that $V(G) = V(G_1)\cdot \hspace{-9.5pt}\cup V(G_2)$.
Moreover, let any vertex of $G_2$ be either connected or
disconnected with all vertices of $G_1$ and
$u\in V(G_1)$ be such hanging vertex in $G_1$ that $G_1\setminus u$ is not complete.
Define the isolating transformation which transforms~$G$ to a graph 
$G'=I(G,u)$ as follows.
We obtain a $G'$ by arising the only edge in $G_1$
incident to~$u$ and adding an edge in $G_1\setminus u$.

{\bf Lemma 8.3}.
Let $G$ be a graph such that $V(G) = V(G_1)\cdot \hspace{-8pt}\cup V(G_2)$.
Moreover, let any vertex of $G_2$ be either connected or
disconnected with all vertices of $G_1$ and
$u\in V(G_1)$ be such hanging vertex in $G_1$ that $G_1\setminus u$ is not complete.
There exists an isolating transformation $G' = I(G,u)$
such that $\beta(G')\geq\beta(G)$.

{\sc Proof}.
Let $V = V(G)$, $V(G_2) = V_2'\cup V_2''$,
where all vertices from $V_2'$ are connected with all vertices from $V_1 = V(G_1)$
and all vertices from $V_2''$ are disconnected with all vertices from $V_1$.
Denote by $e = (u,w)$ the only edge in $G_1$ incident to~$u$. 

If there exists a vertex $t\in V_1$ such that $(w,t)\not\in E(G_1)$,
we consider the isolating transformation which maps $e$ to $(w,t)$.
Otherwise, we construct the isolating transforma\-tion which maps $e$ 
to an edge $(v,t)$ with $v,t$ distinct from $u,w$.
Let us consider the second case (see Picture~5), the proof in the first one is analogous.
Note that $E(G_1)\ni (w,v),(w,t)$.

Let us define an order on $V$: all vertices from $V_2$ 
are greater than all vertices from $V_1$.
Further, $w > u > v > t$ and $t$ is greater 
than all vertices from $V_1\setminus\{u,v,w,t\}$.

\begin{figure}[h]
\centering
\includegraphics[height = 5.5cm]{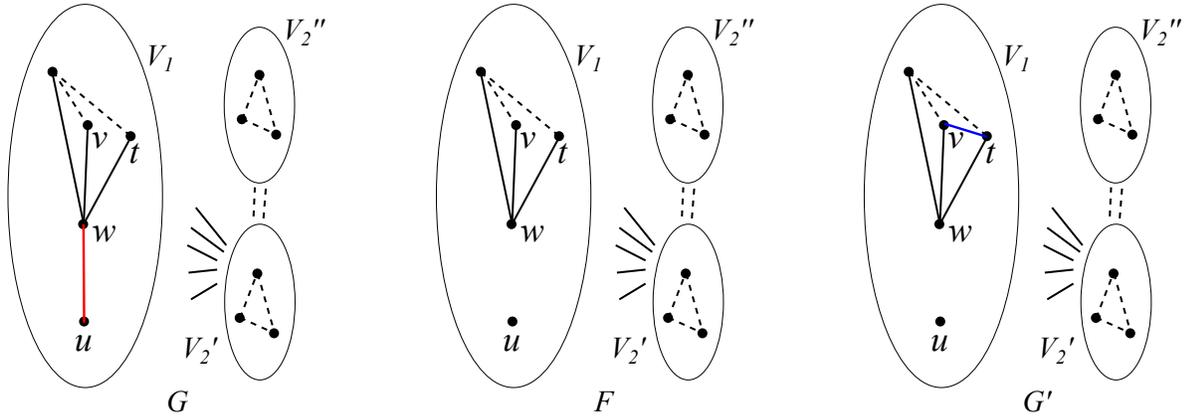}
\caption{The isolating transformation $G' = I(G,u)$.} \label{Pic5}
\end{figure}

Fix a natural number~$n$.
Consider the graph $F = G\setminus e$ and $M_n(V,F)$, 
the set of normal words of length~$n$ in 
the partially commutative monoid~$M(V,F)$.
The set $M_n(V,G)$ could be obtained from $M_n(V,F)$ by removing 
the subset $S_n$ of all words which are not in normal form if adding the edge~$e$.
By Lemma~1.2, $S_n$ consists of words from $M_n(V,F)$ 
containing $uw$ as a subword.
Analogously, $M_n(V,G')$ could be obtained from $M_n(V,F)$ 
by removing the set $T_n$ of all words $M_n(V,F)$ 
containing a subword $tv$ or a subword $vt$ in 
such a position that the entire word is not in normal form iff adding the edge $e'$
(call such a subword $vt$ in words from $M_n(V,F)$ as special one).
Prove that $|S_n|\geq |T_n|$ for all $n$,
it means that $m_n(V,G)\leq m_n(V,G')$
and so, $\beta(G)\leq \beta(G')$.

Introduce $P_n = S_n\cap T_n$.
We want to show that 
$|S_n\setminus P_n|\geq |T_n\setminus P_n|$.
The set $S_n\setminus P_n$ consists of words from $M_n(V,F)$ of the form
$$
\chi_0 uw \chi_1 uw\chi_2 \ldots \chi_{k-1} uw\chi_k, \quad k\geq1,
$$
where $\chi_i$, $i = 0,1,\ldots,k$, are any words in normal form from $M(V,F)$
not containing subwords $uw$, $tv$ and special $vt$.
Some of words $\chi_i$ could be empty.

The set $T_n\setminus P_n$ consists of words from $M_n(V,F)$ of the form
$$
\chi_0 (tv)^{\varepsilon_1} \chi_1 (tv)^{\varepsilon_2}\chi_2 
 \ldots  \chi_{k-1}(tv)^{\varepsilon_k}\chi_k,
\quad \varepsilon_i\in\{\pm1\},\ i=1,\ldots,k,\ k\geq1,
$$
where $(tv)^1 := tv$, $(tv)^{-1} := vt$;
$\chi_i$, $i = 0,1,\ldots,k$, are any words in normal form from $M(V,F)$
not containing subwords $uw$, $tv$ and special $vt$.
Moreover, $\chi_{i-1} (tv)^{\varepsilon_i} \chi_i$, $i=1,\ldots,k$,
should be in normal from in $M(V,F)$.

Consider the set $A$ of all words from $S_n\setminus P_n$
with fixed positions of occurrences of subwords $uw$
and the set $B$ of words from $T_n\setminus P_n$
with the same positions of occurrences of subwords $tv$ (or special $vt$).
Define the map $\varphi\colon B\to A$:
$$
\varphi(\chi_0 (tv)^{\varepsilon_1} \chi_1 (tv)^{\varepsilon_2}\chi_2 
 \ldots  \chi_{k-1}(tv)^{\varepsilon_k}\chi_k)
 = \chi_0 uw \chi_1 uw\chi_2 \ldots \chi_{k-1} uw\chi_k.
$$
The map $\varphi$ is well-defined because of order on $V(G)$.
Let us state that words $\chi_{i-1}tv\chi_i$, $\chi_{i-1}vt\chi_i\in T_n\setminus P_n$ 
could not be in normal form simultaneously in $M(V,F)$,
it will imply that signs $\varepsilon_i$ are uniquely determined.
Indeed, if a word $\chi_{i-1}vt\chi_i\in T_n\setminus P_n$ 
is in normal form in $M(V,F)$ but not in normal form by adding the edge $e'$,
then by Lemma~1.2, $\chi_{i-1} = cyd$ for $y\in V$, $c,d,\in V^*$,
$y<t$, $t$ is connected with $y$ and all letters from $d$. 
Thus, the word $\chi_{i-1}tv\chi_i$ is not in normal form in $M_n(V,F)$.

Hence, $\varphi(B)\subset A$ and $|A|\geq|B|$.
So, $|S_n\setminus P_n|\geq |T_n\setminus P_n|$ for all $n\geq1$.
Lemma is proved.

Now we are ready to prove the generalization of Conjecture~1 from~\cite{FisherUpper}. 

{\bf Theorem 8.1}.
Let $n,k$ be natural numbers, $k = \binom{d}{2} + e\leq \binom{n}{2}$ for $0\leq e<d$.
Construct a graph $G$ with $n$ vertices and $k$ edges as follows.
We add a vertex of degree $e$ to the complete graph $K_d$
and leave all other vertices to be isolated. Then $\beta_+(n,k) = \beta(G)$.

{\sc Proof}.
Note that the constructed graph $G$ is the threshold graph
with the corres\-ponding vector record 
$\thr(G) = (\underbrace{1,\ldots,1}_{d-e-1},0,\underbrace{1,\ldots,1}_{e},0,\ldots,0)$.

Let $H$ be a graph with $n$ vertices, $k$ edges and $\beta(H) = \beta_+(n,k)$.
By Lemma~8.2, we may assume that $H$ is a threshold graph.
Consider the maximal right position of zero coordinate in~$\thr(H)$
such that there is a~1 to the right of the zero.
Denote by~$s$ the position of this zero
and by $u$ the vertex corresponding to the $s$-th coordinate.
Let $V_1$ be a subset of $V(H)$ formed by the initial vertex and all vertices
corresponding to the coordinates of~$\thr(H)$ from the 1st to $(s+1)$-th
and define $V_2$ as $V(H)\setminus V_1$. If $H$ is not isomorphic to~$G$,
we can apply the isolating transformation $H' = I(H,u)$ with 
the respect to the parts $V_1,V_2$.
By Lemma~8.3, $\beta(H') = \beta(H) = \beta_+(n,k)$.
Repeating such procedure, on some step we get~$G$.
On any step, $\beta$ is the same. Theorem is proved. 

Now we can deduce some sort of Kruskal---Katona theorem.

{\bf Corollary 8.1} \cite{CutlerRadcliffe,Wood2007}.
The constructed graph $G$ from Theorem~8.1 maximizes
all numbers $c_k$ among graphs from $G(n,k)$.
In particular, $C(H,x)\leq C(G,x)$ for any graph $H\in G(n,k)$
and for all $x\geq 0$.

{\sc Proof}. 
It follows from Lemma~8.2b, the proof of Theorem~8.1 and 
the fact that any isolating transformation does not decrease
$c_k$ for any $k$.

{\bf Remark 8.1}.
By the same strategy of the proof what we did in Theorem~8.1, 
one can reprove the analogous result for the spectral radius 
(the problem initially posed by Brualdi and Hoffman in 1976~\cite{BrualdiHoffman}
and solved by P. Rowlinson in 1988 \cite{Rowlinson}\footnote{
More precisely, P. Rowlinson stated that the graph $G$ is a unique graph
with the maximal spectral radius among graphs from $G(n,k)$.
We will prove the same refinement for $\beta(G)$ in Corollary~9.1.}).
Indeed, P.~Csikv\'{a}ri in \cite{CsikvariKelman} showed that 
the spectral radius of a graph is not decreasing by a~Kelmans transformation. 
Let us show the same property for an isolating transformation $G' = I(G,u)$ in terms of Lemma~8.3.
Consider the set $W_n(V,H)$ of all walks of length~$n$ in a graph 
$H = G\cup\{(w,t)\}$. To get $W_n(V,G')$, we remove all walks from $W_n(V,H)$ 
containing the edge $e$ and to get $W_n(V,G)$, 
we remove all walks containing $(w,t)$.
By the definition, for $W_n(V,G')$ we avoid all walks having subwords 
$auwb$, $awub$ and $wuw$, where $a,b\in V_2'$.
But for $W_n(V,G)$ we avoid not less walks, since 
we forbid at least all subwords $atwb$, $awtb$ and $wtw$, $a,b\in V_2'$.
So, we have 
$\rho(G) = \lim\limits_{s\to\infty}\sqrt[s]{W_s(V,G)}
\leq \lim\limits_{s\to\infty}\sqrt[s]{W_s(V,G)} = \rho(G')$
and it remains to apply the algorithm of the proof of Theorem~8.1.

\subsection{Upper bound on $\beta(G)$ and Nordhaus---Gadddum inequalities}

{\bf Statement 8.1}.
Let $G$ be a graph with $n$ vertices such that 
all roots of $PC(G,x)$ are real. Then $\beta\leq n-\frac{k}{n}$
and $e(G)\geq \frac{1}{n}$.

{\sc Proof}.
Let $w = \omega(G)$, $\beta = \beta(G)$. By Samuelson's Inequality,
$$
\beta \leq \frac{n}{w}+\frac{w-1}{w}\sqrt{n^2-\frac{2wk}{w-1}}
 \leq \frac{n}{w}+\frac{n(w-1)}{w}\left(1-\frac{wk}{n^2(w-1)}\right)
 = n-\frac{k}{n},
$$
hence, $e(G) = \frac{n-\beta}{k}\geq \frac{1}{n}$.

{\bf Corollary 8.2}.
Given a graph $G$ with $n$ vertices and $e(G)<\frac{1}{n}$,
the polynomial $PC(G,x)$ has complex non-real roots.

Now we want to derive the upper bound on $\beta(G)$
of the form $\beta(G)\leq n-\frac{\alpha k}{n}$
applying Theorem~8.1. Before the proof we need

{\bf Lemma 8.4}.
Let $G$ be a graph with $n\gg1$ vertices and 
$k\geq \frac{n^2}{2}\big(1-\frac{1}{pe^{p+2}}\big)^2$ edges,
then $\beta(G)<\frac{n}{p}$ and $e(G)>\frac{2}{n}\big(1-\frac{1}{p}\big)$.

{\sc Proof}.
Since, $k\leq \binom{n}{2}$,
we get the inequality
$pe^{p+2} \leq n(2+o(n))<en$.
It implies that there exists $\varepsilon>1/2$ such that 
$s = \big[n\big(1-\frac{1}{pe^{p+2}}\big)\big] 
 = n\big(1-\frac{1}{pe^{p+\varepsilon}}\big)$
for the floor function.
By Lemma~2.4 and Theorem~8.1, 
$\beta(G)\leq \beta(H_s)$ 
for $H_s = K_s \cup \bar{K}_{n-s}$.
Since $PC(H_s,x) = (x-1)^s - (n-s)x^{s-1}$, 
we bound for $x\geq \frac{n}{p}$ the following expression
\begin{multline*}
\frac{PC(H_s,x)}{(x-1)^s}
 = 1 - (n-s)\frac{x^{s-1}}{x^s}
 = 1 - \frac{n}{xpe^{p+\varepsilon}}\left(\frac{x}{x-1}\right)^s \\
 > 1 - \frac{1}{e^{p+1/2}}\left(\left(1+\frac{1}{x-1}\right)^{x-1}\right)^{n/(x-1)} 
 > 1 - \frac{1}{e^{p+1/2}}e^{n/(x-1)} \\
 \geq 1 - \frac{e^{\frac{np}{n-p}-p}}{e^{1/2}}
 = 1 - \frac{e^{p^2/(n-p)}}{e^{1/2}} > 0.
\end{multline*}

To prove the second part of the statement, we bound
$$
e(G) = \frac{n-\beta}{k} 
 > \frac{n}{k}\left(1-\frac{1}{p}\right) 
 > \frac{2n}{n^2}\left(1-\frac{1}{p}\right) 
 = \frac{2}{n}\left(1-\frac{1}{p}\right).
$$

{\bf Theorem 8.2}.
Let $n\gg1$.
For any graph $G$ with $n$ vertices and $k$ edges, we have

a) $\beta(G)\leq n-\frac{\alpha k}{n}$,

b) $e(G)\geq \frac{\alpha}{n}$, \\
where $\alpha\approx 0.9408008$.

{\sc Proof}.
For $k=0$ or $k=\binom{n}{2}$, the statement is true.
Suppose that $0<k<\binom{n}{2}$.
By Theorem~8.1, we know a graph on which the maximal value $\beta(n,k)$ is reached.
At first, let $k = \binom{s}{2}$ for some $1<s<n$.
Consider a graph $H_s = K_s \cup \bar{K}_{n-s}$, then 
$PC(H_s,x) = (x-1)^s - (n-s)x^{s-1}$. For $\beta = \beta(H_s)$, we have
$\big(\frac{\beta-1}{\beta}\big)^s = \frac{n-s}{\beta}$.
Thus,
$$
\frac{e^{s/\beta}(n-s)}{\beta} = \theta, \quad
0 < \theta = (e^{1/\beta}(1-1/\beta))^s<1.
$$

Let $n = st$, $t>1$. Then 
$
e^{s/\beta}\frac{s}{\beta}(t-1) = \theta,
$
which implies 
\begin{equation}\label{UpperBound:beta}
\beta = \frac{s}{W(\frac{\theta}{t-1})}
 = \frac{n}{tW(\frac{\theta}{t-1})},
\end{equation}
where $W(x)$ is the Lambert $W$-function, the inverse function to $f(y) = ye^y$.

{\sc Case 1}: $s\geq n(1-\frac{1}{2e^4})$.
By Lemma~8.4, we have $e(G)>1/n$ and we are done. 

{\sc Case 2}: $s<n(1-\frac{1}{2e^4})$.
Thus, $n - s > \frac{n}{110} = O(n)$. By Theorem~7.1, we have $\beta = O(n)$. Compute 
$$
ne(G) = n\frac{n-\beta}{k}
 = \frac{n^2}{\binom{s}{2}}\bigg(1-\frac{1}{tW(\frac{\theta}{t-1})}\bigg)
 = \frac{2ts}{s-1}\bigg(t-\frac{1}{W(\frac{\theta}{t-1})}\bigg)
 \sim 2t\bigg(t-\frac{1}{W(\frac{1}{t-1})}\bigg).
$$
We apply $\theta\sim1$ as 
$$
\theta = \left(e^{1/\beta}\left(1-\frac{1}{\beta}\right)\right)^s
\geq \left(\left(1+\frac{1}{\beta}\right)\left(1-\frac{1}{\beta}\right)\right)^s
 = \left(1-\frac{1}{\beta^2}\right)^s 
 \geq 1 - \frac{s}{\beta^2}.
$$

See the plot of $f(t)$ in Picture~6.

\begin{figure}[h]
\centering
\includegraphics[height = 6.2cm]{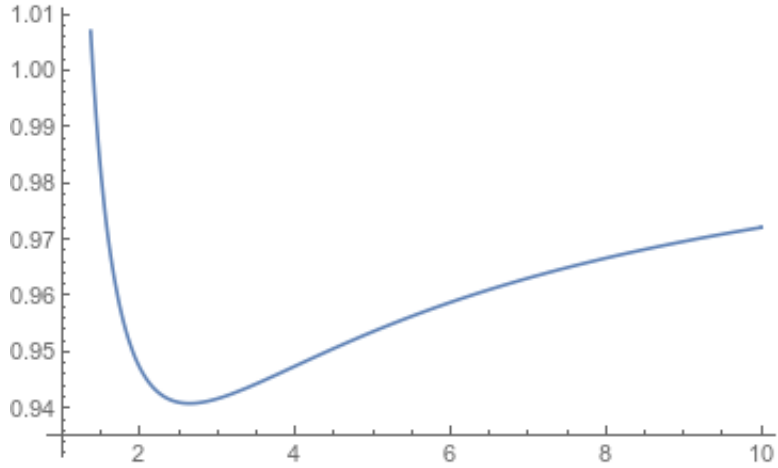}
\caption{The plot of $f(t) = 2t\big(t-\frac{1}{W(\frac{1}{t-1})}\big)$ for $t\in (1;10)$~\cite{wolfram}.}
\label{Pic6}
\end{figure}

Let us proceed on the proof of Theorem~8.2 with the following result.

{\bf Lemma 8.5}.
The function $f(t) = 2t\big(t-\frac{1}{W(\frac{1}{t-1})}\big)$ 
satisfies the following properties:

a) $\lim\limits_{t\to+1}f(t) = 2$,

b) $\lim\limits_{t\to\infty}f(t) = 1$,

c) $f(t)$ is monotonic when $t\to+1$ and $t\to\infty$,

d) $f(t)$ has the minimum value 
$f(t_0) \approx 0.9408008$ 
in the point $t_0\approx 2.6390005$.

{\sc Proof of Lemma~8.5}.
The Taylor series of $W(x)$ around~$+0$ is given by
\begin{equation}\label{Lambert:Taylor}
W(x) 
 = \sum\limits_{i=1}^{\infty}\frac{(-i)^{i-1}}{i!}x^i
 = x - x^2 + \frac{3}{2}x^3 - \frac{8}{3}x^4+\ldots\!.
\end{equation}

a) It follows from the fact that $W(x)\to\infty$ when $x\to\infty$.

b) Applying~\eqref{Lambert:Taylor} and the formula $1/(1-x) = 1 + x + x^2 +\ldots$, 
we prove the statement by the following calculations
\begin{multline}\label{f(t):infty}\allowdisplaybreaks
f(t) = 2t\left(t-\frac{1}{W(\frac{1}{t-1})}\right)
 = 2t\left(t-\frac{1}{\frac{1}{t-1}(1-(\frac{1}{t-1}-\frac{3}{2(t-1)^2} + \frac{8}{3(t-1)^3} + O(\frac{1}{t^4})))}\right) \\
 = 2t\bigg(t-(t-1)\\
 \times\left(1+\frac{1}{t-1}-\frac{3}{2(t-1)^2}+\frac{8}{3(t-1)^3}+\frac{1}{(t-1)^2}-\frac{3}{(t-1)^3}+\frac{1}{(t-1)^3}
  + O\left(\frac{1}{t^4}\right) \right)\bigg) \\
  = 2t\left(t- (t-1)-1+\frac{1}{2(t-1)}-\frac{2}{3(t-1)^2}+O\left(\frac{1}{t^3}\right)\right) \\
  = \frac{t}{t-1}-\frac{4t}{3(t-1)^2} + O\left(\frac{1}{t^2}\right)
  = 1 - \frac{1}{3t} + O\left(\frac{1}{t^2}\right).
\end{multline}

c) The function $f(t)$ monotonically increases 
when $t\gg1$ by \eqref{f(t):infty}. 
Let $a = 1/(t-1)$, then 
$$
f(t) = g(a) = \left(1+\frac{1}{a}\right)^2 - \left(1+\frac{1}{a}\right)\frac{1}{W(a)}.
$$
To prove that $f(t)$ monotonically decreases 
when $t$ is in a right neighbourhood of~$1$,
we will show that $f'(t)<0$ in such neighbourhood.
Equivalently, let us show that $g'(a)>0$ when $a\gg1$.
The derivative of~$W(x)$ equals $W'(x) = \frac{W(x)}{x(1+W(x))}$~\cite{Lambert}.
We compute
$$
g'(a) = -2\left(1+\frac{1}{a}\right)\frac{1}{a^2} + \frac{1}{a^2W(a)} 
 + \left(1+\frac{1}{a}\right)\frac{1}{W^2(a)}\cdot\frac{W(a)}{a(1+W(a))}
$$
or 
\begin{multline*}
a^2g'(a) 
 = -2\left(1+\frac{1}{a}\right) + \frac{1}{W(a)} + \frac{a+1}{W^2(a)+W(a)} \\
 \sim -2 + \frac{a+1}{W^2(a)+W(a)} = -2 + \frac{W(a)e^{W(a)}+1}{W^2(a)+W(a)} > 0
\end{multline*}
for $a\gg1$.

d) It follows due to~\cite{wolfram}. Lemma is proved.

By Lemma~8.5, we are done in Case~2.

At second, consider the case $G\in G(n,k)$ for 
$k = \binom{s}{2} + e$, $0 < e < s$.
For $s = O(n)$, the claim follows from the above stated results.
Indeed, by Lemmas~2.3 and~2.4, 
$\beta(H_{s+1})\leq \beta(G)\leq(H_s)$
and also $|E(G)|,|E(H_s)|,|E(H_{s+1})|\sim \frac{s^2}{2}$.
For $s = o(n)$, we consider two variants.
If $s\geq 100$, then by previous part, we may bound
$$
e(G)
 = \frac{n-\beta}{\binom{s}{2}+e}
 \geq \frac{n-\beta(H_s)}{\binom{s}{2}+e}
 = \frac{n-\beta(H_s)}{\binom{s}{2}}\frac{\binom{s}{2}}{\binom{s}{2}+e}
 \geq \frac{0.99}{n}\cdot\frac{1}{1.02}>\frac{0.97}{n}.
$$
If $s<100$, consider $PC(G,x) = (x-1)^s - (n-s-1)x^{s-1} - x^{s-e-1}(x-1)^e$.
We have $\beta = n(1-\varepsilon)$, where by Theorem~7.1,
$\varepsilon \sim \frac{\alpha}{n^2}$. So, 
\begin{multline}\label{small-k}\allowdisplaybreaks
0 = \frac{PC(G,\beta)}{(\beta-1)^s}
  = 1 - \frac{n-s-1}{n(1-\varepsilon)}\left(1+\frac{1}{n(1-\varepsilon)-1}\right)^s
  - \frac{1}{n(1-\varepsilon)}\left(1+\frac{1}{n(1-\varepsilon)-1}\right)^{s-e} \\
  = 1 - \frac{n-s-1}{n(1-\varepsilon)}
  \left(1+s\left(\frac{1}{n(1-\varepsilon)}+\frac{1}{n^2(1-\varepsilon)^2}\right)
  +\frac{s(s-1)}{2n^2(1-\varepsilon)^2}\right) \\
  - \frac{1}{n(1-\varepsilon)}\left( 1+\frac{s-e}{n(1-\varepsilon)}\right) + O\left(\frac{1}{n^3}\right) \\
  = -\frac{\varepsilon}{1-\varepsilon}-\frac{\varepsilon s}{n(1-\varepsilon)^2}
   +\frac{s^2-s+2e}{2n^2(1-\varepsilon)^2}+ O\left(\frac{1}{n^3}\right).
\end{multline}
Hence, $\alpha = \frac{s^2-s}{2}+e = k$.
Thus, $\beta \sim n\big(1-\frac{\alpha}{n^2}\big) = n - \frac{k}{n}$ and $e(G)\sim \frac{1}{n}$. 
Theorem is proved.

Statement~6.1 could be derived from the following 

{\bf Corollary 8.3}.
Let $G$ be a graph with $n\gg1$ vertices and $k = O(n^{2\theta})$, $\theta<1$, edges.
Then $\beta(G)\sim n -\frac{k}{n}$. 

{\sc Proof}.
By~Statement~7.1, 
$\beta(G)\geq \beta_-(n,k) = \frac{n+\sqrt{n^2-4k}}{2}\sim n -\frac{k}{n}$.
For $s = O(n^\theta)$, we can analogously to \eqref{small-k} 
state that $\beta(G)\lesssim n -\frac{k}{n}$.

{\bf Remark 8.2}.
One can get sufficiently good upper bounds on $\beta(G)$
of the form $\beta(G)\leq n-\frac{\alpha k}{n}$
with the help of Moon---Moser inequalities~\eqref{MoonMoser}
or Fisher---Khadziivanov inequalities~\eqref{FisherClique}~\cite{Fisher1992-Clique,Khadziivanov}
(with $\alpha\in[0.7;0.8]$).
Let us show two upper bounds proven in 1990 by D.C. Fisher and J.~Nonis.

{\bf Statement 8.2}~\cite{FisherUpper}. Given a graph $G$ with $n$ vertices and $k$ edges, we have

a) $\beta^3 - n\beta^2 + k\beta - \sqrt{2}k^{3/2}/3\leq 0$,

b) $\beta\leq \sqrt{n^2-1.5k}$.

{\sc Proof}.
a) By~\eqref{FisherClique}, we have $c_3(G)\leq \sqrt{2}k^{3/2}/3$.
Substitute this bound on $c_3(G)$ in the inequality 
$\beta^3 - n\beta^2 + k\beta - c_3(G)\leq0$ holding by Lemma~2.7a
and we are done.

b) Given a word $w\in X^*$, define $s(w)$ as the number 
of pairs of neighbour letters in $w$ which are connected. 
For example, we have $s(w) = 4$ for the word $w = abccdaba$
and the graph $G$ with $V(G) = \{a,b,c,d\}$ and $E(G) = \{(a,b),(c,d)\}$. 

We also define $h_j = \sum\limits_{w\in X,\,|w|=j}{1}/{2^{s(w)}}$.
Let $A$ be the adjacency matrix of $G$ and we introduce the matrix 
$B = (b_{ij})\in M_n(\mathbb{Z})$ as follows: $b_{ij} = 2 - a_{ij}$.
By the definition of~$h_j$, we have
$$
h_j = \frac{1}{2^{j-1}}(1,1,\ldots,1) B (1,1,\ldots,1)^T,
$$
so in the ``sum-of-squares'' norm, we get
$$
h_j\leq \frac{\|(1,1,\ldots,1)\|^2 \|B\|^{j-1}}{2^{j-1}}
 = n(\sqrt{n^2-1.5k})^{j-1}.
$$
The statement follows from the inequality
$\beta(G) = \lim\limits_{j\to\infty}\sqrt[j]{m_j}\leq \lim\limits_{j\to\infty}\sqrt[j]{h_j}$.

Note that Statement~8.2b implies the asymptotic upper bound 
$\beta(G)\lesssim n-\frac{3k}{4n}$ for $k = o(n^2)$.

{\bf Corollary 8.4}.
Let $G$ be a graph with $n\gg1$ vertices and $k$ edges.

a) For $k\geq 0.256736n^2$, we have 
$e(G)\geq\frac{1}{n}$ and $\beta(G)\leq n-\frac{k}{n} < \frac{3n}{4}$.

b) For $k\geq n^2/4$, we have 
$e(G)\geq\frac{0.996}{n}$
and $\beta(G)\leq n-\frac{0.996k}{n} < 0.751n$.

{\sc Proof}.
a) It follows from the proof of Theorem~8.1 
and the equality $f(\delta) = 1$ for $\delta\approx 1.3955366$.
Hence, $e(G)\geq \frac{1}{n}$ and $\beta(G)\leq n-\frac{k}{n}$
are fulfilled for $\sqrt{2k}\approx s\geq n/\delta$
or equivalently for 
$k\geq \frac{n^2}{2\delta^2}\approx 0.256736n^2$.

b) Note that $f(\sqrt{2})\approx 0.99607$. 
The rest proof is analogous to a). 

{\bf Corollary 8.5}.
For any graph $G$ with $n\gg1$ vertices the following inequalities hold

a) $\beta(G)+\beta(\bar{G})<1.50197n$,

b) $\beta(G)\beta(\bar{G})<0.56398n^2$.

{\sc Proof}.
Let us consider the case $k,\bar{k}\in O(n^2)$.
We may suppose that $k = \binom{s}{2}$ and 
$\bar{k} = \binom{s'}{2}$. Define $t = n/s\in(1;\infty)$
and $t' = n/s' \approx t/\sqrt{t^2-1}$.
In the same way as we did in the proof of Theorem~8.2 (see \eqref{UpperBound:beta}), 
we have 
$$
\frac{\beta(G)+\beta(\bar{G})}{n} 
 \lesssim \frac{1}{tW(\frac{1}{t-1})} + \frac{\sqrt{t^2-1}}{tW(t^2-1+t\sqrt{t^2-1})}.
$$
The function in the RHS has the maximum value $1.50197$ when $t = \sqrt{2}$~\cite{wolfram}.

If one of $k,\bar{k}$ is $o(n^2)$, then by~Lemma~8.4, $\beta(G)+\beta(\bar{G}) = n(1+o(1))$.

b) Analogously to a), we find the maximum of the RHS 
$$
\frac{\beta(G)\beta(\bar{G})}{n^2} 
 \lesssim \frac{1}{tW(\frac{1}{t-1})}\cdot\frac{\sqrt{t^2-1}}{tW(t^2-1+t\sqrt{t^2-1})},
$$
which is equal to~$0.56398$ when $t = \sqrt{2}$~\cite{wolfram}.

{\bf Example 8.1}.
The condition $n\gg1$ in Corollary~8.5 is essential. Indeed, 
\begin{gather*}
\beta(K_{1,2})+\beta(\overline{K_{1,2}})
 = 2+\frac{3+\sqrt{5}}{2} \approx 4.618 \approx 1.539n , \\
\beta(K_{1,2})\beta(\overline{K_{1,2}})
 = 3+\sqrt{5} \approx 0.5818n^2, \\
\beta(K_{1,3})+\beta(\overline{K_{1,3}})
 \approx 3 + 3.1479 = 6.1479 \approx 1.537n,\\
\beta(K_{1,3})\beta(\overline{K_{1,3}})
 \approx 9.4437 \approx 0.59n^2.
\end{gather*}

Lemma~8.2 says that the Kelmans transformation increases (in weak sense)
the value of $\beta$. Note that if $G' = KT(G,u,v)$, then
$\bar{G}' = KT(\bar{G},u,v)$. Thus, the maximal values
of Nordhaus---Gadddum expressions $\beta(G)+\beta(\bar{G})$
and $\beta(G)\beta(\bar{G})$ are attained when both $G,\bar{G}$
are threshold. We conjecture that the following example gives
their maximal values. 

{\bf Example 8.2}.
Let $G = K_s \cup \bar{K}_{n-s}$, $s\geq 2$, $st = n \gg1$,
then by~\eqref{UpperBound:beta} compute
$$
\frac{\beta(G)+\beta(\bar{G})}{n}
 \sim \frac{1}{tW(\frac{1}{t-1})}+\frac{1}{t}.
$$
provided that $n-s = O(n)$. 
Due to~\cite{wolfram} the function 
$h(t) = \frac{1}{t}\big(1+\frac{1}{W(\frac{1}{t-1})}\big)$
has the maximal value $h(t_1)\approx 1.46594$
in the point $t_1\approx 1.2773044$.

We deal with $\beta(G)\beta(\bar{G})$ analogously.
Due to~\cite{wolfram} the function 
$l(t) = \frac{1}{t^2W((\frac{1}{t-1}))}$
has the maximal value $l(t_2)\approx 0.535919$
in the point $t_2\approx 1.31745$.
Thus, the following values can be reached
\begin{gather*}
\beta(G)+\beta(\bar{G})\approx 1.46594n, \\
\beta(G)\beta(\bar{G})\approx 0.535919n^2.
\end{gather*}

{\bf Example 8.3}.
Let $n\gg1$.
Define a graph $G$ such that $V(G) = V_1\cup V_2\cup V_3\cup V_4$,
$V_i\cap V_j = \emptyset$, $i\neq j$, $|V_i| = n/4$, $i,j=1,2,3,4$.
We have the following edges:
$$
E(V) = \{(u,v)\mid u,v\in V_3\cup V_4\}\cup 
\{(u,v)\mid u\in V_1,v\in V_3\}\cup
\{(u,v)\mid u\in V_2,v\in V_4\}.
$$
Then
$$
PC(G,x) = (x-1)^{n/2} - \frac{n}{2}(x-1)^{n/4}x^{n/4-1}
 = (x-1)^{n/4}\bigg((x-1)^{n/4} - \frac{n}{2}x^{n/4-1}\bigg).
$$
For $\beta = \beta(G)$, we have 
$\big(\frac{\beta-1}{\beta}\big)^{n/4} = \frac{n}{2\beta}$.
Since $\beta>n/2$, we have 
$e^{\frac{n}{4\beta}}\frac{n}{4\beta} \approx \frac{1}{2}$.
Hence, $\beta = \frac{n}{4W(0.5)}\approx 0.710765n$.

Since $\bar{G}$ is isomorphic to $G$, we compute
\begin{gather*}
\beta(G)+\beta(\bar{G})\approx 1.42153n, \\
\beta(G)\beta(\bar{G})\approx 0.50519n^2.
\end{gather*}

{\bf Remark 8.3}.
Let $n\gg1$. By Corollary~7.5 and Corollary~8.5 we have 
the upper bound $\rho(G)+\rho(\bar{G})< 1.5158n$
for any graph $G$ with $n$ vertices.

Of course, this bound is too bad. 
It was a series of works devoted to Nordhaus---Gadddum
problem for the spectral radius.
Analogously to the proof of Theorem~7.1a,
one can state that $n-1\leq\rho(G)+\rho(\bar{G})$.
In 1970, E.~Nosal proved that
$\rho(G)+\rho(\bar{G})\leq \sqrt{2}n$~\cite{Nosal}.
In 2007, V.~Nikiforov improved this one to $(\sqrt{2}-8\cdot10^{-7})n$
and conjectured that $\rho(G)+\rho(\bar{G})\leq 4n/3$~\cite{Nikiforov}.
In 2009, P.~Csikv\'{a}ri showed that the Kelmans transformation
increases (in weak sense) the spectral radius of~$G$ as well as of~$\bar{G}$
and dealing with threshold graphs proved 
$\rho(G)+\rho(\bar{G})\leq \frac{1+\sqrt{3}}{2}n$~\cite{CsikvariNikif}.
Finally, in 2011 T. Terpai proved Nikiforov's conjecture~\cite{Terpai}
with the help of analytic methods.
The common example for the asymptotically tight upper bound for 
both $\rho(G)+\rho(\bar{G})$ and $\rho(G)\rho(\bar{G})$ is the following:
it is a clique of size $2n/3$ and its complement.
It is interesting to compare with the situation with $\beta(G)$.
If Example~8.2 gives the asymptotically tight upper bound for 
the Nordhaus---Gadddum expressions for~$\beta$,
they are attained in essentially different graphs.

\newpage
\section{Minimum value $\beta_{-}(n,k)$}

After Statement~7.1, we are interested to find $\beta_{-}(n,k)$ for $k>n^2/4$.

Given a graph $G$ and two distinct vertices $u,v\in V(G)$, let us call 
a Kelmans transformation $KT(G,u,v)$ as a nontrivial one, if $c(G')-c(G)>0$,
where $c(H)$ denotes the number of all cliques in a graph $H$.
Nontriviality of a Kelmans transformation $KT(G,u,v)$
is equivalent to the following condition: 
there exist $x\in N_G(u)\setminus N_G(v)$ 
and $y\in N_G(v)\setminus N_G(u)$
such that $(x,y)\in E(G)$.

{\bf Lemma 9.1}.
Let $G$ be a graph with connected complement.
If a graph $G'$ is a result of a nontrivial Kelmans transformation of~$G$,
then $\beta(G')>\beta(G)$.

{\sc Proof}.
Let $G' = KT(G,u,v)$. Denote 
$A = N_G(u)\setminus N_G(v) =\{a_1,\ldots,a_p\}$,
$B = N_G(u)\cap N_G(v)$
and $C = N_G(v)\setminus N_G(u) = \{c_1,\ldots,c_r\}$.
Due to the definition of the Kelmans transformation,
\begin{multline}\label{KelmansDiff}
D(G',x) - D(G,x) \\
 = D(G[A\cup B\cup C],x) - D(G[A\cup B],x) - D(G[B\cup C],x) + D(G[B],x).
\end{multline}

Applying Lemma~1.5a successively, we get
\begin{multline}\label{KelmansDiffProof}
D(G[A\cup B\cup C],x)
 = D(G[A\cup B\cup C]\setminus c_1,x) - xD(G[N(c_1)\cap (A\cup B\cup C)],x) \\
 = D(G[A\cup B\cup C]\setminus \{c_1,c_2\},x) \\
 - xD(G[N(c_1)\cap (A\cup B\cup C)],x) 
 - xD(G[N(c_2)\cap (A\cup B\cup C_2)],x) \\
 = \ldots = D(G[A\cup B],x) - x\left( 
  \sum\limits_{i=1}^{r} D(G[N(c_i)\cap (A\cup B\cup C_i)],x)  \right), 
\end{multline}
where $C_i = C\setminus\{c_1,c_2,\ldots,c_i\}$.

Dealing with $D(G[B\cup C],x)$ in a similar way as we did in \eqref{KelmansDiffProof}
and substituting~\eqref{KelmansDiffProof} and the obtained formula for $D(G[B\cup C],x)$
in \eqref{KelmansDiff}, we have
\begin{multline}\label{KelmansDiffProof2}
D(G[A\cup B\cup C],x) - D(G[A\cup B],x) - D(G[B\cup C],x) + D(G[B],x) \\
 = -x\sum\limits_{i=1}^{r} ( D(G[N(c_i)\cap (A\cup B\cup C_i)],x) 
 - D(G[N(c_i)\cap (B\cup C_i)],x) ).
\end{multline}
Now we apply Lemma~1.5a for every 
$D(G[N(c_i)\cap (A\cup B\cup C_i)],x)$ 
considering successively vertices from $A$:
\begin{multline*}
D(G[N(c_i)\cap (A\cup B\cup C_i)],x) - D(G[N(c_i)\cap (B\cup C_i)],x)  \\
 = -x\sum\limits_{j=1}^{p} D(G[N(a_j)\cap N(c_i)\cap (A_j\cup B\cup C_i)],x),
\end{multline*}
where $A_j = A\setminus\{a_1,a_2,\ldots,a_j\}$. Finally, we write down 
\begin{equation}\label{KelmansDiffProofSubst}
D(G',x) - D(G,x) 
 = x^2\sum\limits_{i=1}^{r}\sum\limits_{j=1}^{p} 
  D(G[N(a_j)\cap N(c_i)\cap (A_j\cup B\cup C_i)],x). 
\end{equation}

By Lemma~8.2a, $\beta(G')\geq\beta(G)$.
Substituting $x = 1/\beta(G')$ in \eqref{KelmansDiffProofSubst},
on the one hand we get a nonpositive number in the LHS.
On the other hand, by Lemma~2.3a ($\bar{G}$ is connected!)
we have a positive number in the RHS of~\eqref{KelmansDiffProofSubst}
provided that at least one of sets 
$M_{ij} = N(a_j)\cap N(c_i)\cap (A_j\cup B\cup C_i)$ is not empty.
Suppose that $M_{ij} = \emptyset$ for all $i$ and~$j$. 
Thus, $c_k(G') = c_k(G)$ for all $k\geq4$. 
Since we consider nontrivial Kelmans transformation,
$c_3(G')>c_3(G)$ and hence $\beta(G')>\beta(G)$. Lemma is proved.

{\bf Lemma 9.2}.
Let $G$ be a graph with connected complement.
Let $e\in E(G)$ be an edge lying in a clique of size~$t\geq3$
and $a,b\in V(G)$ such vertices that $(a,b)\not\in E(G)$
and $N_G(a)\cap N_G(b) = \emptyset$.
Consider the graph $G'$ obtained by removing an edge $e$
and adding an edge $(a,b)$. Then $G'\leq G$ and $\beta(G')<\beta(G)$.

{\sc Proof}.
Let $e = (u,v)$, denote $A = N_G(u)\cap N_G(v)$. 
It is easy to get the following equality (see, e.g., Lemma~1.5b)
$$
D(G',x) - D(G,x) = x^2 - x^2D(G[A],x).
$$

Denote $\beta = \beta(G)$.
Let us show that the function $f(x) = 1 - D(G[A],x)$ 
increases strictly monotonically on $x\in [0;1/\beta]$,
it will imply that $G\geq G'$ and also
$$
\beta^2 D(G',1/\beta) = 1 - D(G[A],1/\beta) > 1 - D(G[A],0) = 0,
$$
hence $\beta(G')<\beta(G)$.
Indeed, by Lemma~1.5c, 
$$
f'(x) = - D'(G[A],x) = \sum\limits_{w\in A}D(G[A\cap N(w)],x) > 0, 
\quad x\in [0;1/\beta],
$$
since the edge $e$ lies in a clique of size~$t\geq3$
and $D(H,x)>0$, $x\in [0;1/\beta]$, 
for every proper induced subgraph $H$ of $G$ by Lemma~2.3b.

{\bf Corollary 9.1}.
The graph $G$ constructed in Theorem~8.1 is a unique graph with the maximal~$\beta(G)$
among all graphs with $n$ vertices and $k$ edges with one exception:
when $k = \binom{d}{2}+1$ for some~$d$. In this case, 
the set $\{H\in G(n,k)\mid \beta(H) = \beta_+(n,k)\}$
consists of all graphs obtained from $K_d\cup \bar{K}_{n-d}$
by adding one edge.

{\sc Proof}.
The statement is trivial for $n<3$ or $k<3$. 
If $\bar{k}<n-1$, then we have
$PC(G,x) = (x-1)^{n-1-\bar{k}}PC(G',x)$, where $G'$ obtained from $G$
by removing $n-1-\bar{k}$ vertices of degree $n-1$. 
So, we may assume that $3\leq k\leq \binom{n}{2}-n+1$.

Let $H$ be a graph with $\beta(H) = \beta_+(n,k)$.
Suppose that $\bar{H}$ is disconnected, i.e.,
$G = H_1+\ldots+H_s$ with connected $\bar{H}_i$.
Let $H_1$ be a maximal part of $H$.
If at least one of $|V(H_2)|,\ldots,|V(H_s)|$
is greater than~1 (say $H_2$), then 
we apply the Kelmans transformation $G' = KT(G,u,v)$
for any $v\in V(H_1)$ and $u\in V(H_2)$. 
Since the maximal part of $G'$ contains $H_1$ 
as the proper induced subgraph we get 
$\beta(G') > \beta(G)$ by Lemma~2.3, a contradiction.

Let $|V(H_2)| = \ldots = |V(H_s)| = 1$.
Let us show that this case is contradictary for $s = 2$,
the proof for $s>2$ is analogous.
By Theorem~8.1, $\beta(H) = \beta(H_1) = \beta(F_1)$,
where the graph $F_1$ consists of the complete graph $K_{d'}$,
a vertex $v$ of degree $e'$ and remaning 
$n-1-d'-e'\geq1$ isolated vertices.
Here $\binom{d'}{2}+e' = k-n+1$ (we have fixed $s =2$), $0\leq e'<d'$. 
Let $V(H_2) = \{u\}$. Introduce the graph $F = F_1 + \{u\}$
and the graph $J$ obtained from $F$ by removing 
an edge $(u,k)$ for a vertex $k$ from the clique $K_{d'}$ 
and adding an edge $(v,i)$ for some remaining isolated vertex~$i$.
Note that both $F$ and $J$ have connected complements.
Therefore, we can apply the proof of Lemma~9.2
for the pair of graphs $F$ (instead of $G'$) and~$J$
(instead of $G$). If $d'>1$, we get $\beta(F)<\beta(J)$, a contradiction.
If $d'=1$, it means that $F_1$ ($=H_1$) is an empty graph.
Since $\bar{k}\geq n-1\geq 2$, we have $|V(H_1)|\geq3$.
Thus, $F$ is a triangle-free graph and we have by Statement~7.1
$\beta(F)<\beta_+(n,k)$, a contradiction.

We have proved that $H$ has to be connected. 
If $\omega(H)<d$, then some of the Kelmans or isolating transformations 
(see proof of Theorem~8.1) imply $\beta(H)<\beta_+(n,k)$ by Lemma~9.1 or by Lemma~9.2.
So, $\omega(H) = d$.

For $k = \binom{d}{2} + e$, if $e=0$ or $e=1$, we are done. Suppose that $e\geq2$. 
We get the extremal graph $G$ either from at least one Kelmans transformation
or after at least one isolating transformation (and, maybe Kelmans transformations before;
see the proof of Theorem~8.1).
In the first case, consider the last Kelmans transformation $G = KT(H,u,v)$,
where $\beta(H) = \beta(G)$. Because of the the structure of~$G$ and $\omega(H) = d$,
this Kelmans transformation is nontrivial and $\beta(H)<\beta(G)$, a contradiction.
In the second case, consider the last isolating transformation $G = I(H,u)$.
In the vector notation of threshold graphs, we transformed the graph 
$\thr(H) = (\underbrace{1,\ldots,1}_{d-e},0,\underbrace{1,\ldots,1}_{e-1},0,1,0,\ldots,0)$
to the graph 
$\thr(G) = (\underbrace{1,\ldots,1}_{d-e-1},0,\underbrace{1,\ldots,1}_{e},0,\ldots,0)$.
As $d\geq 3$, we have $\beta(G)>\beta(H)$ by Lemma~9.2.

{\bf Corollary 9.2}.
Let $k>[n^2/4]$ and $G$ be a graph such that $\beta(G) = \beta_-(n,k)$.
Then $G$ is connected graph having diameter~2.

{\sc Proof}.
If $\bar{G}$ is disconnected, then $G$ is connected and has diameter~2, 
we are done. Let $\bar{G}$ be connected. 

If $G$ has more than one connected components, by Mantel's Theorem 
(Statement~3.1), at least one component contains a triangle. 
Applying Lemma~9.2, we get a graph $G'$ such that $\beta(G')<\beta(G)$,
a contradiction.

Suppose $G$ has diameter~3 or greater. 
It means that one can find a pair of vertices $u,v\in V(G)$
such that $(u,v)\not\in E(G)$ and there no a vertex $w\in V(G)$ with the property 
$(u,w),(w,v)\in E(G)$. It remains to apply Lemma~9.2 for a graph $G'$ 
obtained by removing an edge $e$ from any triangle of $G$
and adding an edge $(u,v)$. We get $\beta(G')<\beta(G)$, a contradiction.
Corollary is proved.

Now we formulate the following conjecture modulo which we will find
the exact value of $\beta_-(n,k)$.

{\bf Conjecture 9.1}.
Let $k>[n^2/4]$ and $G$ be a graph such that $\beta(G) = \beta_-(n,k)$.
Then $\bar{G}$ is disconnected.

Consider the case when $n^2/4 < k < n^2/3$.
At first, let $n = 2l+1$ and $k<l(l+2)$. 
In this case we construct the graph $G_-$ as follows. 
It is a supergraph of $K_{l+1,l}$
in which the $k' = k-l(l+1)$ edges form 
a tree inside the part with $l+1$ vertices. Then 
\begin{equation}\label{BetaMinValue3Spec}
\beta(G_-) = \frac{l+1}{2}\left(1+\sqrt{1 - \frac{4k'}{(l+1)^2}}\right).
\end{equation}

At second, let $n$ be either even or odd with $k\geq n^2/4+n/2$.
Find a natural number $n_1$ such that
\begin{equation}\label{BetaMinPartition3}
n_1(2n-3n_1)\leq k< (n_1-1)(2n-3n_1+3).
\end{equation}
Note that $n/3<n_1\leq n/2$.
The inequalities \eqref{BetaMinPartition3} are equivalent to the following ones:
$$
|E(K_{n_1,n_1,n-2n_1})|\leq k<|E(K_{n_1-1,n_1-1,n-2n_1+2})|.
$$

Denote $k' = k-n_1(2n-3n_1)$. From~\eqref{BetaMinPartition3} we conclude that
$0\leq k'<6n_1-2n-3$. 
Let us construct the graph $G_-$ as follows. 
It is a supergraph of $K_{n_1,n_1,n-2n_1}$
in which the $[k'/2]$ edges form a tree 
in the first part and remaining $k' - [k'/2]$ edges 
form a tree in the second part. So, $\beta(G_-)$ is the largest root of the equation
$x^2 - n_1x + [k'/2] = 0$, i.e., 
\begin{equation}\label{BetaMinValue3}
\beta(G_-) = \frac{n_1}{2}\left(1 + \sqrt{1-\frac{2k'}{n_1^2}+\frac{2\varepsilon}{n_1^2}}\right),
\end{equation}
where $\varepsilon = k'/2 - [k'/2]\in\{0,1\}$.
Note that we allow for $n-2n_1$ be equal to zero.

{\bf Theorem 9.1}.
Let $n^2/4<k<n^2/3$. If Conjecture~9.1 holds, then $\beta_-(n,k) = \beta(G_-)$
for the constructed graph $G_-$. Thus, $\beta_-(n,k)$
is defined by \eqref{BetaMinValue3Spec} if $n = 2l+1$ and $k<l(l+2)$
or by \eqref{BetaMinValue3} otherwise.

{\sc Proof}.
By Conjecture~9.1, a graph $H\in G(n,k)$ with $\beta(H) = \beta_-(n,k)$
has disconnec\-ted complement. Assume that $H$ is a supergraph of 
the complete multipartite graph with parts $H_1,\ldots,H_s$
such that $\bar{H}_i$ are connected. 

If $n = 2l+1$ and $k<l(l+2)$, then $s = 2$, since we have 
too less edges to form either more parts or to form bipartite graph
with another capacities of parts. It remains to apply Statement~7.1.

Let $n$ be even and $k<\frac{n^2}{4}+n-3$, then $s \leq 3$.
We have two possibilities: to construct a supergraph of $K_{n/2,n/2}$
(with connected complements) 
or to construct a supergraph of $K_{n/2,n/2-1,1}$ 
(provided $k\geq \frac{n^2}{4}+\frac{n}{2}-1$). 
It is easy to see that the smallest value of $\beta$
is reached on the graph $G$ obtained from 
inserting $[k'/2]$ (where $k' = k - n^2/4$) 
edges in one part of $K_{n/2,n/2}$
and inserting $k' - [k'/2]$ edges in the another.
Actually, it is \eqref{BetaMinValue3} for $n_1 = n/2$
and the same value $k'$.

Now we consider the case when $k$ is enough big to form 
a tripartite graph with parts $a,a,b$ with $a\geq b$.
Let us show that there exists a graph $H'\in G(n,k)$
with $\beta(H') = \beta_-(n,k)$ and $\bar{H'}$ 
has exactly three connected components.
Suppose there exists a graph $H\in G(n,k)$ with the minimal $\beta$
and $s = 4$. Denote $a_i = |H_i|$, $1\leq i\leq 4$,
and order the numbers in such way that $a_1\geq a_2\geq a_3\geq a_4$. 
If $a_1>n_1$, then $\beta(H)\geq\beta(G_-)$.
Indeed, if $H$ has more than $a_1-1$ edges over 
the edges of the graph $K_{a_1,a_2,a_3,a_4}$,
then we can construct a supergraph of 
$K_{a_1-1,a_2,a_3,a_4+1}$ with smaller $a_1$.
If $|E(H)|-|E(K_{a_1,a_2,a_3,a_4})|<a_1-1$,
by Statement~7.1, $\beta(H)>a_1-1\geq \beta(G_-)$.
If $a_1<n_1$, then 
\begin{multline*}
k = |E(H)| = a_1(a_2+a_3+a_4) + a_2(a_3 + a_4) + a_3a_4 \\
 > (n_1-1)(n-n_1+1) + (n_1-1)(n-2n_1+2)
 = 2nn_1-3n_1^2 + 6n_1 - 2n-3,
\end{multline*}
a contradiction to \eqref{BetaMinPartition3}. 
If $a_1 = n_1$, in a similar way we can show that $a_2 = n_1$
and the unique case when $H$ has the minimal value of $\beta_-(n,k)$
is the following: $a_3 = a_4 = 1$ and $k' = k - n_1(2n-3n_1)$ is odd.
In this case we can move the only edge between parts $a_3,a_4$
inside of one parts of $a_1,a_2$.

For $H$ which complement has three connected components,
the minimality of $\beta(H)$ leads to form the small as possible 
maximal parts of $H$, it is exactly the condition \eqref{BetaMinPartition3}.
After that, we may apply Statement~7.1.

{\bf Corollary 9.3}.
Let $n^2/4<k<n^2/3$. 

a) If Conjecture~9.1 holds, then the graph 
$H\in G(n,k)$ such that $\beta(H) = \beta_-(n,k)$
is unique up to the choice of triangle-free graphs
with $[k'/2]$ and $k' - [k'/2]$ vertices inside the two largest parts.
If $k'$ is odd, we also may form triangle-free graphs with $(k'-1)/2$
edges inside the two largest parts and draw the last edge inside the third part. 

b) We have 
$\frac{n}{3}+\frac{1}{3}\sqrt{n^2-3k}\leq \beta_{-}(n,k)<\frac{n}{3}+\frac{1}{3}\sqrt{n^2-3k}+1$.

{\sc Proof}.
a) It follows from the proof of Theorem~9.1.

b) The lower bound follows from Theorem~7.2, the upper one 
is true by the construction of the graph $G_-$.

{\bf Remark 9.1}.
For the constructed graph $G_-$, we have
\begin{multline*}
c_3(G_-) 
 \approx n_1^2(n-2n_1) 
 \approx \frac{1}{27}(n+\sqrt{n^2-3k})^2(n-2\sqrt{n^2-3k}) \\
 = \frac{9nk-2n^3-2(n^2\sqrt{n^2-3k})^{3/2}}{27},
\end{multline*}
which is equal to the lower bound from Corollary~7.6.
Thus, the graph $G_-$ has asymptoti\-cally the smallest number of triangles.

Now let us consider the general case for $n^2/4<k$.
Let $\big(1-\frac{1}{w-1}\big)\frac{n^2}{2}<k<\big(1-\frac{1}{w}\big)\frac{n^2}{2}$
for a natural number $w$.

Define $l = \big[\frac{n}{w-1}\big]$, $p = n-l(w-1)$ and $q = w-1-p$.
Suppose that $k < \binom{w-1}{2}l^2+pl(w-1)$, i.e.,
$k$ is not enough large to construct a supergraph of $K_{a,a,\ldots,a,b}$
(with $w-1$ parts of $a$ and $a\geq b\geq0$).
Thus, the graph $H$ with $\beta(H) = \beta_-(n,k)$
is a supergraph of $K_{l+1,l+1,\ldots,l+1,l,\ldots,l}$ 
with $p$ parts with $l+1$ edges and $q$ parts with $l$ edges.
Introduce 
$$
k' = k - \left(\binom{p}{2}(l+1)^2+\binom{q}{2}l^2+pql(l+1)\right).
$$
We construct the graph $G_-$ as follows.
We draw a triangle-free graphs with $[k'/p]$ vertices inside all $p$
parts with $l+1$ vertices and draw remaining edges anywhere.
So, 
\begin{equation}\label{BetaMinValueNSpec}
\beta(G_-) = \frac{l+1}{2}\left(1+\sqrt{1 - \frac{4[k'/p]}{(l+1)^2}}\right).
\end{equation}

Let $k$ be enough large to construct a supergraph of $K_{a,a,\ldots,a,b}$ 
with prescribed conditions on parts.
Find a natural number $n_1$ such that
\begin{equation}\label{BetaMinPartitionN}
(w-1)n_1\left(n-\frac{wn_1}{2}\right)\leq k 
< (w-1)(n_1-1)\left(n-\frac{w(n_1-1)}{2}\right).
\end{equation}
We also have $\frac{n}{w}<n_1\leq \frac{n}{w-1}$ and 
the inequalities \eqref{BetaMinPartitionN} are equivalent to the following ones:
$$
|E(K_{n_1,\ldots,n_1,n-(w-1)n_1})|\leq k<|E(K_{n_1-1,\ldots,n_1-1,n-(w-1)(n_1-1)})|.
$$

Denote $k' = k - \big((w-1)n_1n-\binom{w}{2}n_1^2\big)$ and 
$0\leq k'< (w-1)(wn_1-n-w/2)$. 
We construct the graph $G_-$ as a supergraph of $K_{n_1,\ldots,n_1,n-(w-1)n_1}$
in which the $[k'/(w-1)]$ edges form a triangle-free graph
in each part with $n_1$ vertices and we put remaining 
$k' - (w-1)[k'/(w-1)]$ edges anywhere. 
Hence, $\beta(G_-)$ as the largest root of the equation
$x^2 - n_1x + \left[\frac{k'}{w-1}\right] = 0$ equals
\begin{equation}\label{BetaMinValueN}
\beta(G_-) = \frac{n_1}{2}\left(1 + \sqrt{1-\frac{4k'}{(w-1)n_1^2}+\frac{4\varepsilon}{(w-1)n_1^2}}\right),
\end{equation}
where $\varepsilon = k'/(w-1) - [k'/(w-1)]\in\{0,1\}$.

{\bf Theorem 9.2}.
Let $\big(1-\frac{1}{w-1}\big)\frac{n^2}{2}<k<\big(1-\frac{1}{w}\big)\frac{n^2}{2}$. 
If Conjecture~9.1 holds, then $\beta_-(n,k) = \beta(G_-)$
for the constructed graph $G_-$. Thus, $\beta_-(n,k)$
is defined by \eqref{BetaMinValueNSpec}, \eqref{BetaMinValueN}.

{\sc Proof}.
The proof is analogous to the proof of Theorem~9.1.

{\bf Corollary 9.4}.
Let 
$\big(1-\frac{1}{w-1}\big)\frac{n^2}{2}<k\leq\big(1-\frac{1}{w}\big)\frac{n^2}{2}$, $w\geq2$, then 
\begin{equation}\label{LowerBoundGen}
\frac{n}{w}+\frac{1}{w}\sqrt{n^2-\frac{2kw}{w-1}}\leq 
 \beta_{-}(n,k) < \frac{n}{w}+\frac{1}{w}\sqrt{n^2-\frac{2kw}{w-1}}+1.
\end{equation}

{\sc Proof}. 
The lower bound follows from Theorem~7.2. 
For $k = \big(1-\frac{1}{w}\big)\frac{n^2}{2}$, the inequalities follow from Statement~7.1 and Statement 7.2. 
For $k < \big(1-\frac{1}{w}\big)\frac{n^2}{2}$,
the upper bound follows from the construction of $G_-$.
Indeed, if $k \geq \binom{w-1}{2}l^2+pl(w-1)$, then
we have 
$$
\beta_{-}(n,k)\leq \beta(G_-)\leq n_1 = \left\lceil\frac{n}{w}+\frac{1}{w}\sqrt{n^2-2kw/(w-1)}\right\rceil
$$
by~\eqref{BetaMinPartitionN}. Otherwise, it is easy to show that 
$l < \frac{n}{w}+\frac{1}{w}\sqrt{n^2-\frac{2kw}{w-1}}$ and thus, 
$$
\beta_{-}(n,k)\leq \beta(G_-)< l+1 < \frac{n}{w}+\frac{1}{w}\sqrt{n^2-\frac{2kw}{w-1}}+1.
$$

{\bf Remark 9.2}.
Note that the construction of $G_-$ is almost the same as the constructions
from the works of C. Reiher \cite{Reiher} and A.A. Razborov~\cite{Razborov} 
devoted to the problem of finding the minimal values of the numbers of cliques in a graph from $G(n,k)$ 
(see also the article \cite{TrianglesBig} about the exact evaluation 
of the number of triangles).

See in Picture~\ref{Pic8} the asymptotics of the borders of possible values of $\beta(G)/n$,
the upper bound is derived from the proof of Theorem~8.2 and 
the lower bound follows from Corollary~9.4.

\newpage
\section{Average value of $\beta(G)$}

\subsection{Average value}

We are interested on the asymptotics of the expression
$$
\beta_{\evv}(n)
 =\frac{1}{2^{\binom{n}{2}}}\sum\limits_{G\colon |V(G)| = n}\beta(G),
$$
i.e., we want to find the average value of growth rate 
of partially commutative monoid with $n$-vertex commutativity graph. 

{\bf Lemma 10.1}.
For all $p\in[0;1]$, there exists the limit 
$\lim\limits_{n\to\infty}\frac{\beta(G_{n,p})}{n}$.

{\sc Proof}.
For $p = 0$, the limit equals~1, for $p = 1$ the limit equals~0.
Let $p\in(0;1)$ and $b = 1/p$.

Consider the polynomials
$$
P_t(x) = \sum\limits_{i=0}^t \frac{(-1)^i}{i!b^{\binom{i}{2}}}x^{t-i},\quad t\geq1.
$$
Define $x_t$ as the largest positive root of $P_t(x)$.
For even $t$, if $P_t(x)$ has no such roots, define $x_t = 0$.
It is clear that $x_k\in [0;1]$ and $x_{2k} < x_{2k+1}$, $x_{2k} < x_{2k-1}$ 
for all $k\in\mathbb{N}$.

For polynomials
\begin{equation}\label{Lem:BetaAve-def:D}
D_t(x) = \sum\limits_{i=0}^t \frac{(-1)^i}{i!b^{\binom{i}{2}}}x^i 
 = x^t P_t(1/x),\quad t\geq1,
\end{equation}
we have the following relation 
\begin{equation}\label{Lem:BetaAve-D}
D_{n+1}'(x) = -D_n(px).
\end{equation}

From the definition of $D_t(x)$ and~\eqref{Lem:BetaAve-D}, we conclude that
\begin{equation}\label{Lem:BetaAve-P}
P_{n}'(x) = \frac{1}{xb^{n-1}} P_{n-1}(bx) + \frac{n}{x}P_{n}(x).
\end{equation}

Suppose that $x_{2k+1}$ is not a simple root of $P_{2k+1}(x)$,
then by \eqref{Lem:BetaAve-P} we have that $bx_{2k+1}$ is a root of 
the polynomial $P_{2k}(x)$, which contradicts the inequality $x_{2k}<x_{2k+1}$.
Hence, $x_{2k+1}$ is a simple root of $P_{2k+1}(x)$.

Introduce polynomials
$$
R_t(n,x) = \sum\limits_{i=0}^t \frac{(-1)^i\binom{n}{i}}{b^{\binom{i}{2}}} x^{t-i},\quad t\geq1,
$$
and $x_t(n)$ as the largest positive root of $R_t(n,x)$.
For even $t$, if $R_t(n,x)$ has no such roots, define $x_t(n) = 0$.

Let us show that $x_{2k}(n)\leq \beta(G_{n,p})\leq x_{2k+1}(n)$ for $k\gg1$.
Indeed,
\begin{multline}\label{Lem:BetaAve-R}
PC(G_{n,p},x)
 = x^{n-2k}R_{2k}(n,x)
 - \left(\frac{\binom{n}{2k+1}x^{n-2k-1}}{b^{\binom{2k+1}{2}}}
- \frac{\binom{n}{2k+2}x^{n-2k-2}}{b^{\binom{2k+2}{2}}}\right) \\
 - \left(\frac{\binom{n}{2k+3}x^{n-2k-3}}{b^{\binom{2k+3}{2}}}
- \frac{\binom{n}{2k+4}x^{n-2k-4}}{b^{\binom{2k+4}{2}}}\right)
 - \ldots\!.
\end{multline}
By Lemma~5.4a, $\beta(G_{n,p}) > n(1-p)$.
The inequality $x_{2k}(n)\leq \beta(G_{n,p})$ 
follows from \eqref{Lem:BetaAve-R} considered for $x\geq \beta(G_{n,p})$.
Analogously we prove that $\beta(G_{n,p})\leq x_{2k+1}(n)$.

Consider 
$$
R_t(n,x/n)
 = \sum\limits_{i=0}^t \frac{(-1)^i(1-\frac{1}{n})\ldots
 (1-\frac{i-1}{n})}{i! b^{\binom{i}{2}}}x^{t-i},\quad t\geq1.
$$
Since $P_{2k+1}(x)$ is strict monotonic for $x\geq x_{2k+1}$ by \eqref{Lem:BetaAve-P}
and $R_{2k+1}(n,x/n)$ as the function on $n$ is continuous, we have
$$
1-p\leq \varlimsup\limits_{n\to\infty}\frac{\beta(G_{n,p})}{n}
 \leq \varlimsup\limits_{n\to\infty}\frac{x_{2k+1}(n)}{n} \leq x_{2k+1}.
$$
From the inequality $1-p\leq x_{2k+1}$ and the relation
$$
P_{2k+1}(x) = x^2P_{2k-1}(x) + \frac{1}{ (2k)!b^{\binom{2k}{2}}}\left(x-\frac{1}{(2k+1)b^{2k}}\right),
$$
we state that the sequence $x_{2k+1}$ is monotonically decreasing starting with some~$k$.
This sequence is bounded, thus, there exists the limit
$\lim\limits_{k\to \infty} x_{2k+1} = x_0 \geq 1-p$.

Suppose that $P_{2k+1}(x)$ has another positive real roots 
distinct from $x_{2k+1}$. Define $z_{2k+1}$ as 
the second largest positive real root of $P_{2k+1}(x)$.
Then, $D_{2k+1}'(x)$ has a root $a$ in the interval $(1/x_{2k+1};1/z_{2k+1})$.
Note that all roots of $P_t(x)$ and $D_t(x)$
are in one-to-one correspondence by the rule $x\leftrightarrow 1/x$.
So, by \eqref{Lem:BetaAve-D}, we have $x_{2k}\geq 1/(ap)$. It implies
$$
z_{2k+1}<\frac{1}{a}\leq px_{2k}<p x_{2k+1},
$$
thus, $z_{2k+1} < px_{2k+1}$.
Therefore, all other real roots of $P_{2k+1}(x)$
could be separated from the root $x_{2k+1}$ for all $k$.

By Ostrowski's Theorem applied for the polynomials $P_{2k+1}(x)$ and $xP_{2k}(x)$, we have 
\begin{equation}\label{OstrowskiAppl}
|x_{2k+1}-\xi|
 < 4\left(\frac{1}{(2k+1)!b^{\binom{2k+1}{2}}}\right)^{1/(2k+1)}
 \sim \frac{4e}{(2k+1)b^k} = \varepsilon
\end{equation}
for some root $\xi$ of $xP_{2k}(x)$.
Let us prove that there exists a real root $x'$ of $xP_{2k}(x)$
such that $x_{2k+1}-x'<4e/b^k$.
Indeed, suppose that $\xi = \xi_1\in\mathbb{C}\setminus\mathbb{R}$,
then by \eqref{OstrowskiAppl},
there exists a~non-real complex root $\nu_1$ of $P_{2k+1}(x)$
in a $\varepsilon$-neighbourhood of $\bar{\xi}_1$, the complex conjugate of $\xi_1$.
So, by \eqref{OstrowskiAppl}, there exists a root $\xi_2$
in a $\varepsilon$-neighbourhood of $\bar{\xi}_1$.
Since $x_{2k+1}$ is a separated simple root of $P_{2k+1}(x)$,
we can find the required real root $x'$ of $xP_{2k}(x)$.
By the definition of $x_{2k}$ and the equality $x_{2k}<x_{2k+1}$,
we have that $x_{2k+1}-x_{2k}<4e/b^k$.
Hence, the sequence $x_{2k}$ for $k\gg1$ also has the limit equal to~$x_0$.
By~\eqref{Lem:BetaAve-P} we can prove that $x_{2k}$ is a simple root of $P_{2k}(x)$ for $k\gg1$.
We can state that all other real roots of $P_{2k}(x)$ are separated from the root $x_{2k}$ for $k\gg1$
(in a similar way as we did for $x_{2k+1}$).
By Ostrowski's Theorem applied for the polynomials
$P_{2k}(x)$ and $R_{2k}(n,x)$, $n,k\gg1$, we get 
$$
x_{2k} \leq \varliminf\limits_{n\to\infty}\frac{x_{2k}(n)}{n}
 \leq \varliminf\limits_{n\to\infty}\frac{\beta(G_{n,p})}{n}
 \leq \varlimsup\limits_{n\to\infty}\frac{\beta(G_{n,p})}{n}
 \leq \varlimsup\limits_{n\to\infty}\frac{x_{2k+1}(n)}{n} \leq x_{2k+1}.
$$
By the squeeze theorem, there exists the limit 
$\lim\limits_{n\to\infty}\frac{\beta(G_{n,p})}{n}$.

{\bf Theorem 10.1}.
The average value of $\beta(G)$ on graphs with $n\gg1$ vertices
asymptoti\-cally equals
$$
\beta_{\mathrm{ev}}(n)
 \sim n\lim\limits_{n\to\infty}\frac{\beta(G_{n,1/2})}{n} = \beta_0 n \approx 0.672008n.
$$

{\sc Proof}.
The sequence $x_n$ in Lemma~10.1 was defined.
For $p = 1/2$, we can compute that
$\lim\limits_{n\to\infty}x_n = \beta_0\approx 0.6720076538$~\cite{magma,wolfram}.

In \cite{BolobasErdos} it was proved that with probability tending to~1,
the clique number of a graph $G$ with $n$ vertices equals
$\omega(G) = 2\log_2(n)+O(\log_2\log_2(n))$.
Fix a natural number $t$ and put $n > 2^t$.
For almost all graphs with $n$ vertices the number of cliques 
of size $2\leq i\leq t$ equals 
$\frac{(1+\varepsilon_i)\binom{n}{i}}{2^{ \binom{i}{2}}}$,
where $\varepsilon_i\in\big(-\frac{\ln n}{n};\frac{\ln n}{n}\big)$ \cite{Bolobas}.
Denote $\varepsilon = \{\varepsilon_i\}_{i\geq2}$
and define $z_{t}(n,\varepsilon)$ as the largest positive real root of the polynomial
$$
Q_t(n,\varepsilon,x)
 = \sum\limits_{i=0}^t\frac{(-1)^i (1+\varepsilon_i) \binom{n}{i}}{2^{\binom{i}{2}}}x^{t-i},\quad t\geq1.
$$
If $Q_t(n,\varepsilon,x)$ has no such roots, define $z_{t}(n,\varepsilon) = 0$.
Below while writting $Q_t(n,\varepsilon,x)$ as well as $z_{t}(n,\varepsilon)$,
we assume that $n>2^t$.

By Corollary~7.3a, we have $\beta_{\evv}(n)\geq \frac{n+1}{2}$.
Analogously to \eqref{Lem:BetaAve-R},
we get the inequalities
$z_{2k}(n,\varepsilon)\leq \beta_{\evv}(n)\leq z_{2k+1}(n,\varepsilon)$
for all sufficiently large~$n$.
Let $t\gg1$. Since the coefficients of $Q_t(n,\varepsilon,x/n)$ tend to 
the corresponding coefficients of $P_t(x)$ when $n\to\infty$
and $P_t(x)>0$ for $x>x_t$, we have 
$\varlimsup\limits_{n\to\infty}\frac{z_t(n,\varepsilon)}{n} \leq x_t$.
On the other hand, $x_t$ is a~simple separated root of $P_t(x)$,
so Ostrowski's Theorem considered for the polynomials
$P_t(x)$ and $Q_t(n,\varepsilon,x/n)$ gives
$x_t\leq \varliminf\limits_{n\to\infty}\frac{z_t(n,\varepsilon)}{n}$.
Hence, 
$\lim\limits_{n\to\infty}\frac{z_t(n,\varepsilon)}{n} = x_t$
holds independently of the choice of a tuple 
$\varepsilon\in \big(-\frac{\ln n}{n};\frac{\ln n}{n}\big)^{t-1}$
and we have 
$\lim\limits_{n\to\infty}\frac{\beta_{\evv}(n)}{n} = \beta_0$.

{\bf Corollary 10.1}.
The average value of $e(G)$ for graphs with $n\gg1$ vertices asymptoti\-cally equals
$\frac{4(1-\beta_0)}{n} \approx \frac{1.312}{n}$.

{\bf Remark 10.1}.
The constant $\beta_0 = 0.6720076538\ldots$ firstly appeared 
in the article of R. Stanley~\cite{Stanley} in 1973.
In \cite{Stanley}, the number of all acyclic orientations of a digraph was counted as 
$A\cdot n!2^{\binom{n}{2}}\beta_0^n$ for $A\approx 1.741$.
The number $1/\beta_0$ was obtained in \cite{Stanley} as the unique root of the function 
$F(x) = \sum\limits_{n=0}^{\infty}\dfrac{x^n}{n!2^{\binom{n}{2}}}$
which modulus is not greater than~2.

{\bf Statement 10.1}.
For almost all graphs with $n$ vertices,
$\beta(G)$ lies in a neighbourhood of $\beta_0 n$
and is the unique root $PC(G,x)$ which modulus is greater than $n/2$.

{\sc Proof}.
The first part of the statement follows from Theorem~9.1. 
Let us prove the second one. Denote $p = 1/2$, $b = 2$. 
Consider a polynomial $D_t(x)$, $t\geq3$,
defined by \eqref{Lem:BetaAve-def:D},
and $H(x) = 1-x+\frac{x^2}{4}-\frac{x^3}{48}$.
It is easy to show that in the circle~$x = 2e^{i\varphi}$, we have 
$$
\min\limits_{|x|=2}\{|H(x)|\}\geq \left(\frac{\cos \varphi}{6}\right)^2
+ \left(\frac{\sin \varphi}{2}\right)^2
\geq \frac{1}{36} 
 > \frac{1}{48} = 2\cdot\frac{2^4}{4!2^6}
  > \max\limits_{|x|=2}\{|D_t(x)-H(x)|\}.
$$
By Rouch\'{e}'s theorem, the polynomial $D_t(x)$, $t\geq3$, 
has a unique root in a circle $|x| = 2$.
We can apply Rouch\'{e}'s theorem in the same way 
for the polynomials $H(x)$ and
$$
\widetilde{Q}_t(n,\varepsilon,x)
 = \sum\limits_{i=0}^t\frac{(-1)^i (1+\varepsilon_i)\binom{n}{i}}{2^{\binom{i}{2}}}x^i
$$
for $n,t\gg1$, where $\varepsilon$ is chosen 
as we did in the proof of Theorem~10.1. We are done.

{\bf Remark 10.2}.
Statement~10.1  implies that that Theorem~2.1 for almost all graphs could be
derived from Rouch\'{e}'s theorem. 

Define the number $\beta_{\evv}(n,k)$
as the average value of $\beta(G)$ for the set $G(n,k)$
of all graphs with $n$ vertices and $k$ edges:
$$
\beta_{\evv}(n,k)
 =\frac{1}{\binom{\binom{n}{2} }{k}}\sum\limits_{G\in G(n,k)}\beta(G).
$$

{\bf Statement 10.2}.
Let $k(n)$ be such integer-valued function that $0\leq k(n)\leq \frac{n(n-1)}{2}$
and there exists $\lim\limits_{n\to\infty}\frac{2k(n)}{n^2} = k_0<1$. Then
$\beta_{\evv}(n,k) \sim n\lim\limits_{n\to\infty}\frac{\beta(G_{n,k_0})}{n}$.

{\sc Proof}.
In the random graph model $G_{n,k}$, 
(initial Erd\H{o}s---R\'{e}nyi model, slightly diffe\-rent from $G_{n,p}$)
a graph is chosen uniformly at random from the $G(n,k)$.
In $G_{n,k}$, the expected value of cliques of size~$t$ equals
$$
\frac{\binom{n}{t} \binom{\binom{n}{2}-\binom{t}{2} }{k-\binom{t}{2} }  }{\binom{\binom{n}{2} }{k} }
 = \binom{n}{t} \frac{k}{\binom{n}{2}}\cdot\frac{k-1}{\binom{n}{2}-1}\cdot\ldots \cdot\frac{k-\binom{t}{2}+1}{\binom{n}{2}-\binom{t}{2}+1}
 = \binom{n}{t} p^{\binom{t}{2}}\theta(n,k,t)
$$
for $k\geq \binom{t}{2}$ and~0, otherwise. 
Here $\theta(n,k,t)\to 1$ when $t$ is fixed and $n\to\infty$.
The condition $k_0<1$ and Theorem~7.1 imply that $\frac{\beta_{\evv}(n,k)}{n}>\frac{1-k_0}{2}>0$.
The rest of proof is analogous to the proof of Theorem~10.1.
The correctness of transfer from one random graph model to another one
see in \cite{Bolobas,FriezeKaronski}.

{\bf Remark 10.3}.
If $k_0 = 1$, Statement~10.2 is not true. 
For example, for $k(n) = \binom{n}{2}$ we have $\beta_{\evv}(n,k) = 1$,
for $k(n) = \binom{n}{2}-1$ we have 
$\beta_{\evv}(n,k) = 2$.
For $k(n)\ll \binom{n}{2}-n$ we have 
$\beta_{\evv}(n,k) \geq 3$, since 
$p = \frac{1}{n}$ is a threshold for $G_{n,p}$
to contain at least one triangle~\cite{ErdosRenyi}.
Thus, the complement graph $\bar{G}$ to $G\in G(n,k)$ 
contains a triangle and $\beta(G)\geq 3$ by Statement~3.2a.

\subsection{Roots of $PC$-polynomial of random graph}

Let $p>0$ and $\beta_r = \beta_r(G_{n,p})$ denote the $r$-th
largest root of $PC(G_{n,p},x)$.

{\bf Statement~10.3}. 
For any $p\in(0;1]$, we have 
$\frac{(1-p)}{2}p(n-1)\leq \beta_2(G_{n,p})\leq p(n-1)$.

{\sc Proof}.
By~\eqref{BoundsOnBetaNP}, we have the upper bound on 
$\beta_2 = \beta_2(G_{n,p})$:
$$
\beta_2 \leq p n -\beta_1 \leq n - 1 - (n-1)(1-p) = p(n-1).
$$

By Theorem~5.1b, we can bound 
$\beta_2$ as follows:
$$
n - \beta_1 = (\beta_1+\ldots+\beta_n) - \beta_1
 \leq \beta_2\left(1+\frac{1}{b}+\ldots+\frac{1}{b^{n-2}}\right)
 \leq \frac{\beta_2}{1-p}.
$$
By the upper bound on $\beta_1$, we have
\begin{equation}\label{LowerBoundBeta2}
\beta_2 \geq (1-p)(n-1)(1-\sqrt{1-p})
 = \frac{p(1-p)(n-1)}{1+\sqrt{1-p}}\geq \frac{p(1-p)(n-1)}{2}.
\end{equation}

{\bf Lemma 10.2}.
For all $p\in(0;1]$, there exists the limit 
$\lim\limits_{n\to\infty}\frac{\beta_2(G_{n,p})}{n}$.

{\sc Proof}.
For $p = 1$, the limit equals~0.
Let $p\in(0;1)$ and $b = 1/p$.
The idea of the proof is the same as of the proof of Lemma~10.1.
For $t\gg1$, define $y_t$ as the second largest real root of $P_t(x)$.
For $k\gg1$, the $y_{2k}$ is well-defined as the largest positive real 
root of the polynomial $P_{2k}(x)/(x-x_{2k})$.
If $y_{2k+1}$ is not defined, put $y_{2k+1} = 0$.

It is easy to show that $y_{2k}>y_{2k-1}$ and $y_{2k}>y_{2k+1}$ for all $k$.
Suppose that $y_{2k}$ is not a simple root of $P_{2k}(x)$.
Then by \eqref{Lem:BetaAve-P}, we have $by_{2k} = x_{2k-1}$
as $by_{2k}$ is a real root of $P_{2k-1}(x)$ larger than $y_{2k-1}$.
Let us find the root $a$ of $D_{2k-1}(x)$ in the interval $(1/x_{2k};1/y_{2k})$.
By \eqref{Lem:BetaAve-D}, $1/(pa)$ is a root of $P_{2k-1}(x)$
which has to coincide with $x_{2k-1}$.
Thus, $x_{2k-1} = \frac{y_{2k}}{p}<\frac{1}{pa}$, a~contradiction.
So, $y_{2k}$ is a simple root of $P_{2k}(x)$.

Define $y_t(n)$ as the second largest positive real root of the polynomial
$R_t(n,x)$. We apply the lower bound 
$\frac{(1-p)}{2}p(n-1)\leq \beta_2(G_{n,p})$ 
and the inequality $P_{2k}(x)<0$ for $x\in(y_{2k},x_{2k})$ 
to prove that
$y_{2k+1}(n)\leq \beta_2(G_{n,p})\leq y_{2k}(n)$ for $k\gg1$.
From the relation
$$
P_{2k+2}(x) = x^2P_{2k}(x) 
 - \frac{1}{ (2k+1)!b^{\binom{2k+1}{2}}}\left(x-\frac{1}{(2k+2)b^{2k+1}}\right),
$$
we state that the sequence $y_{2k}$ is monotonically decreasing starting with some~$k$.
This bounded sequence has the limit
$\lim\limits_{k\to \infty} y_{2k} = y_0 \geq p(1-p)/2$.

Introduce $z_{2k}$ as the third largest real root of $P_{2k}(x)$.
So, $D'_{2k}(x)$ has a~root $a$ in $(1/y_{2k};1/z_{2k})$.
By \eqref{Lem:BetaAve-D}, $D_{2k-1}(pa) = 0$.
Hence, $1/(pa)\leq y_{2k-1}<y_{2k}$.
We conclude that $z_{2k}<\frac{1}{a}<p y_{2k}$
and all other real roots of $P_{2k}(x)$ are separated from $y_{2k}$.
The rest of the proof is the same as of Lemma~10.1.

{\bf Theorem 10.2}.
For all $p\in(0;1]$ and $r\geq1$, there exists the limit 
$\lim\limits_{n\to\infty}\frac{\beta_r(G_{n,p})}{n}$.

{\sc Proof}.
For $p = 1$, the limit equals~0. Let $p\in(0;1)$ and $b = 1/p$.
Prove the statement by induction on $r$. For $r=1,2$
it follows from Lemmas~10.1 and~10.2. 
Let $r>2$. We define the $r$-th largest positive real 
roots $w_r$ of $P_t(x)$ and $w_r(n)$ of $R_t(n,x)$
respectively.
For odd~$r$ we repeat the proof of Lemma~10.1,
for even~$r$ we apply the proof of Lemma~10.2.

For completeness, we should show only one thing:
the lower bound $q_r n<\beta_r(G_{n,p})$ for $n\gg1$ 
with a positive constant $q_r$. For $r = 1,2$,
we had $q_1 = 1-p$, $q_2 = p(1-p)/3$.
To the contrary, suppose that $\beta_{r}(G_{n,p}) = o(n)$. 
By Theorem~5.1b and Vieta's formula, 
\begin{multline}\label{LowerBoundBetaR}
\binom{n}{r}p^{\binom{r}{2}}
 = \sum\limits_{1\leq i_1<\ldots<i_r\leq n}\beta_{i_1}\ldots \beta_{i_r} 
 < \sum\limits_{j=1}^{r-1}(\beta_1+\ldots+\beta_{r-1})^j(\beta_r+\ldots+\beta_n)^{r-j} \\
 < \sum\limits_{j=1}^{r-1}n^j\left(\frac{\beta_r}{1-p}\right)^{r-j}
 = o(n^r),
\end{multline}
a contradiction. Theorem is proved.

{\bf Corollary 10.2}.
Let $r>0$. For almost all graphs with $n\gg1$ vertices, 
the real roots of $PC$-polynomial which moduli is not less than $n/r$
lie in neighbourhoods of the roots of $PC(G_{n,1/2},x)$.

{\sc Proof}. The proof is analogous to the proof of Theorem~10.1.

With the help of engines~\cite{wolfram,magma},
one can approximately compute six largest roots 
of $PC(G_{n,1/2},x)$: 
\begin{gather*}
0.672008n,\quad
0.204871n,\quad
0.073744n,\quad
0.028756n,\quad
0.011768n,\quad
0.004975n.
\end{gather*}

Let us show how all coefficients of the series expansion of 
$\lim\limits_{n\to\infty}\frac{\beta_r(G_{n,p})}{n}$ on $p$
could be calculated.
With the help of expressions of power sums in terms of 
elementary symmetric polynomials and Vieta's formulas 
we have
\begin{multline}\label{FindSeriesBeta1}
\beta_1^5+\beta_2^5+\ldots+\beta_n^5 \\
 = n^5 - 5n^3\binom{n}{2}p + 5n\left(\binom{n}{2}p\right)^2
 + 5n^2\binom{n}{3}p^3 
 - 5\binom{n}{2}\binom{n}{3}p^4
 -5n\binom{n}{4}p^6 + 5\binom{n}{5}p^{10} \\
 \sim n^5\left(1 - \frac{5p}{2} + \frac{5p^2}{4} + \frac{5p^3}{6} 
 - \frac{5p^4}{12} - \frac{5p^6}{24} + \frac{p^{10}}{24}\right).
\end{multline}
At first, $\beta_1^5$ is less than the RHS of \eqref{FindSeriesBeta1}.
At second, by Theorem~5.1b, $\beta_1^5$
is greater than the RHS of \eqref{FindSeriesBeta1}
multiplied by $1-p^5$. So, we can find the first five coefficients
of the series expansion of $\beta_1$ on $p$:
$\beta_1(G_{n,p})\sim n(1-p/2-p^2/4-p^3/12-p^4/16)$.
Continuing on, we get the following expression
\begin{multline}\label{Beta1Series}
\frac{\beta_1(G_{n,p})}{n}\sim 
1-\frac{p}{2}-\frac{p^2}{4}-\frac{p^3}{12}-\frac{p^4}{16}
 -\frac{p^5}{48}-\frac{7p^6}{288}-\frac{p^7}{96} -\frac{7p^8}{768} \\
 - \frac{49p^9}{6912}-\frac{113p^{10}}{23040}-\frac{17p^{11}}{4608}
 -\frac{293p^{12}}{92160}-\frac{737p^{13}}{276480}-\frac{3107p^{14}}{1658880}+O(p^{15}).
\end{multline}
For $p = 1/2$, the equality \eqref{Beta1Series}
gives the average value $\beta_0$ with mistake less than $3\cdot10^{-8}$.

\begin{figure}[h]
\centering
\includegraphics[height = 6.5cm]{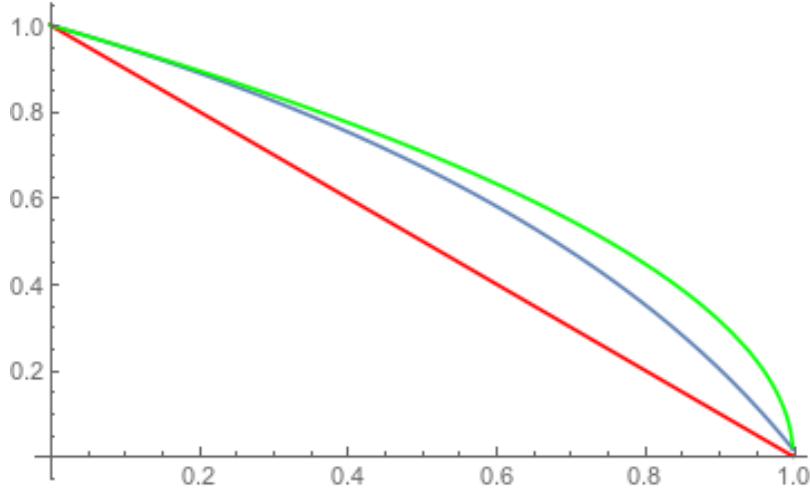}
\caption{The plots of the RHS of \eqref{Beta1Series} (blue), 
 $1-p$ (red) and $\sqrt{1-p}$ (green) for $p\in(0,1)$.}
\label{Pic7}
\end{figure}

In Picture~\ref{Pic7} the plots of the RHS of \eqref{Beta1Series} 
and asymptotic bounds on $\beta_1(G_{n,p})/n$ from Corollary~5.4 are drawn.

In terms of $q = p/2$, we get 
\begin{multline}\label{SeriesBeta1InQ}
\frac{\beta_1(G_{n,p})}{n}\sim 
1-q-q^2-\frac{2q^3}{3}-q^4-\frac{2q^5}{3}-\frac{14q^6}{9}-\frac{4q^7}{3} -\frac{7q^8}{3} \\
 - \frac{98q^9}{27}-\frac{226q^{10}}{45}-\frac{68q^{11}}{9}-\frac{586q^{12}}{45}-\frac{2948q^{13}}{135} 
 - \frac{12428 q^{14}}{405} + O(q^{15}).
\end{multline}

To calculate $\beta_2(G_{n,p})$,
we may consider the inequalities
$$
\big(\big(\beta_1^7+\beta_2^7+\ldots+\beta_n^7\big) - \beta_1^7\big)
(1-p^7) < \beta_2^7<\big(\beta_1^7+\beta_2^7+\ldots+\beta_n^7\big) - \beta_1^7,
$$
which imply
\begin{equation}\label{Beta2Series}
\beta_2^7 = \big(\big(\beta_1^7+\beta_2^7+\ldots+\beta_n^7\big) - \beta_1^7\big)(1+O(p^7)).
\end{equation}
Inserting in \eqref{Beta2Series} 
the RHS of the 7th power analogue of~\eqref{FindSeriesBeta1} and 
the series \eqref{Beta1Series}, we find
\begin{equation}\label{Beta2Formula}
\frac{\beta_2(G_{n,p})}{n}\sim 
\frac{p}{2}\left(1-\frac{p}{6}-\frac{5p^2}{18} -\frac{29p^3}{216} 
 - \frac{85p^4}{648}-\frac{163p^5}{3888}-\frac{1387p^6}{19440}\right) + O(p^8),
\end{equation}
which in terms of $q = p/6$ gives
$$
\beta_2/n\sim 3q(1-q-10q^2-29q^3-170q^4-326q^5-16644q^6/5) + O(q^8).
$$
Note that the lower bound from Statement~10.3 is closer to $\beta_2(G_{n,p})$ 
than the upper one.

Further, we calculate
$$
\frac{\beta_3(G_{n,p})}{n}\sim \frac{p^2}{3}\left(1-\frac{p}{12}-\frac{25p^2}{144}+O(p^3)\right),\quad
\frac{\beta_4(G_{n,p})}{n}\sim \frac{p^3}{4}(1+O(p)).
$$

Thus, by such algorithm we can compute all coefficients of the series expansion of 
$\lim\limits_{n\to\infty}\frac{\beta_r(G_{n,p})}{n}$ on $p$ for any~$r$. 
Note that with the help of Statement~5.1b, we automatically get the expressions
for the smallest roots of $PC(G_{n,p},x)$.

By the expressions of four largest roots of $PC(G_{n,p})$
and in light of Statement~5.1b, we formulate

{\bf Conjecture 10.1}.
For any~$r\geq1$, we have 
$$
\lim\limits_{n\to\infty}\frac{\beta_r(G_{n,p})}{n}
 = \frac{p^{r-1}}{r}\left(1-\frac{p}{r(r+1)}+O(p^2)\right).
$$ 

Now let us consider the following symmetric polynomial 
on roots of $PC(G_{n,p})$:
$$
(1+\beta_1)(1+\beta_2)\ldots(1+\beta_n)
 = \sum\limits_{k=0}^n\binom{n}{k}p^{\binom{k}{2}}
 = C(G_{n,p},1).
$$
Since $\beta_i = O(np^{i-1})$ by Theorem~5.1b 
and $np^{i-1} = O(1)$ already for 
$i = \frac{\ln n}{\ln(1/p)}$,
we may assume that $\ln C(G_{n,p},1) = O(\frac{\ln^2 n}{\ln(1/p)})$.
In \cite{BrownAsy}, it was proven that
$\ln C(G_{n,p},1) \sim \frac{\ln^2 n}{2\ln(1/p)}$.

In 2014, W. Gawronski and T. Neuschel stated the following result:
(see also \cite{Paris})

{\bf Theorem 10.3}~\cite{Gawronski}.
For a fixed $p\in(0;1)$, as $n\to\infty$, we have 
$$
C(G_{n,p},1)
 = \frac{1}{\sqrt{r(n)}}
 \left(\theta_3 \left(\frac{\pi r(n)}{\ln(1/p)}, e^{-2\pi^2/{\ln(1/p)}}\right)+o(1)\right)
 \exp\left(\frac{r(n)^2 + 2r(n)}{2\ln(1/p)}\right),
$$
where $\theta_3(z,q)$ is the Jacobi's third Theta function,
$r(n)$ is defined as the positive solution $t$ of the equation
$t\left(e^t+\sqrt{1/p}\right) = n\sqrt{1/p}\ln(1/p)$.

In 2012--13 \cite{BrownAsy,BrownFn}, 
J. Brown, K. Dilcher and V. Manna initiated to study the poly\-nomials 
$f_n(z) = \sum\limits_{k=0}^n\binom{n}{k}z^{\binom{k}{2}}$
for complex variable $z$.
For real $z\in[0;1]$, $f_n(z) = C(G_{n,z},1)$.

{\bf Statement 10.4} \cite{BrownFn}.
a) For each $n\geq 0$, $f_{2n+1}(z)$ is divisible by $z^n + 1$.

b) For each $n\geq 3$ the roots of $f_n(z)$ lie inside a circle of radius $1 + \frac{3\ln n}{n}$.

c) For each $n\geq 3$ the roots of $f_n(z)$ lie outside a circle of radius $2/n$.

d) For each $n\geq 4$ there is a negative real root of $f_n(z)$ in the interval
$-\frac{2+\frac{2}{n}}{n}< z < -\frac{2}{n}$.

{\sc Proof}. a) It follows from Statement~5.1b, 
as the proof works for any (not necessary real) $p$.

b), c), d) See in \cite{BrownFn}.

In \cite{BrownFn}, there were formulated a lot of conjectures
about the distributions of roots, irreducibility of $f_n(z)$ etc. 

{\bf Remark 10.5}.
We may suppose that just for complex variable~$p$,
we may have 
\begin{equation}\label{BrownFunctViaRoots}
f_n(z) = \prod\limits_{i=1}^n(1+\beta_i(G_{n,z})),
\end{equation}
where $\beta_i(G_{n,z})$ are defined by the described above series expansions.
But we should be very careful. Indeed, if we try to insert 
the expression~\eqref{Beta2Formula} in \eqref{BrownFunctViaRoots} to find a root of $f_n(z)$,
we will make a mistake.
From the equality $1 + \frac{zn}{2}\big(1-\frac{z}{6}+O(z^2)\big) = 0$,
we find the root $z = -\frac{2-\frac{2}{3n}}{n} + O\big(\frac{1}{n^3}\big)$,
a contradiction to Statement~10.4c. The reason of the mistake is the following:
$\beta_2(n,z) = \frac{zn}{2}\big(1-\frac{z}{6}+O(z^2)\big) + o(n)$,
where we do not control the tail $o(n)$.
If our assumption about the correctness of the decomposition \eqref{BrownFunctViaRoots}
is true, we conclude that the tail $o(n)$ is at least $1/(3n)$.

\begin{figure}[h]
\centering
\includegraphics[height = 6.5cm]{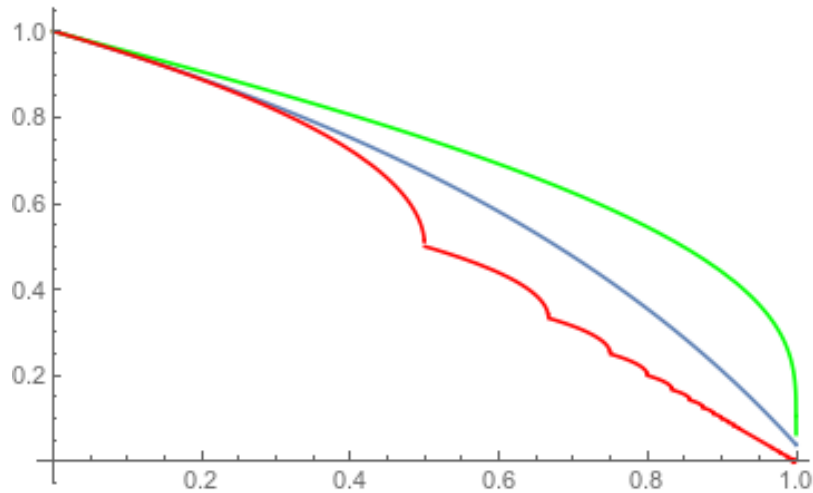}
\caption{The plots of $f_1(x)$ (green), $f_2(x)$ (red) and $f_3(x)$ (blue) for $x\in(0,1)$.}
\label{Pic8}
\end{figure}

Let us finish the section with Picture~\ref{Pic8} 
on which the borders of the possible values of $\beta(G)/n$ are drawn (asymptotically), 
compare it with Fig.~3 from~\cite{FisherUpper}. 
Define 
\begin{gather*}
f_1(x) = \frac{\sqrt{x}}{W(\frac{\sqrt{x}}{1-\sqrt{x}})},\\
f_2(x) = \frac{1}{\lceil\frac{1}{1-x}\rceil}\left(1+\sqrt{1-\frac{\lceil\frac{1}{1-x}\rceil x}{\lceil\frac{1}{1-x}\rceil-1}}\right),
\end{gather*}
\vspace{-0.5cm}
\begin{multline*}
f_3(x) =  1-\frac{x}{2}-\frac{x^2}{4}-\frac{x^3}{12}-\frac{x^4}{16}
 -\frac{x^5}{48}-\frac{7x^6}{288}-\frac{x^7}{96} -\frac{7x^8}{768}\\
 - \frac{49x^9}{6912}-\frac{113x^{10}}{23040}-\frac{17x^{11}}{4608}
 -\frac{293x^{12}}{92160}-\frac{737x^{13}}{276480}-\frac{3107x^{14}}{1658880}.
\end{multline*}

Let $n = |V(G)|\gg1$ and $x = \frac{2k}{n^2}\in(0,1)$.
The function $f_1(x)$ is asymptotically equal to the upper bound of $\beta(G)/n$ 
arising from the proof of Theorem~8.2, 
the function $f_2(x)$ asymptotically equals to the lower bound of $\beta(G)/n$ from Corollary~9.4.
And $f_3(x)$ is the approximation of the asymptotics of the average value of $\beta(G)/n$~\eqref{Beta1Series}. 

To compare $f_1(x)$ with $f_3(x)$, we have 
$$
f_1(x) = 1-\frac{x}{2}+\frac{x^{3/2}}{6}-\frac{7x^2}{24}+\frac{31x^{5/2}}{120} + O(x^3).
$$

Note that the maximum value of the difference $f_1(x) - f_3(x)$ is attained in 
a point~$x$ close to~1,
the maximum of the difference $f_3(x) - f_2(x)$ is attained either in $x = 1/2$ or in $x = 2/3$.

\newpage
\section{Weighted case. Lov\'{a}sz local lemma}

Let $G$ be a simple graph and $Cl(G)$ denote the set of all cliques of $G$.
Let us define the weighted dependence polynomial of $G$ as 
$$
D_w(G,x) = \sum\limits_{B\in Cl(G)}(-1)^{|B|}w(B),\quad 
w(B) = \prod\limits_{v\in B}\alpha_v x^{d_v},\
\alpha_v,d_v\in \mathbb{R}_{+}.
$$
If $\alpha_v = d_v = 1$ for all $v\in V(G)$, 
then the weighted dependence polynomial of $G$
coincides with the dependence polynomial $D(G,x)$.

Define the weighted $PC$-polynomial of $G$ as 
$PC_w(G,x) = x^{w_0}D_w(G,1/x)$ for $w_0 = \deg(D_w(G,x))$.

S. Lavall\'{e}e and C. Reutenauer in 2009 and 
S. Lavall\'{e}e in 2010 proved the following results
which we gathered in one statement.

{\bf Theorem 11.1}.
a) \cite{Reutenauer} 
The set of all weighted dependence polynomials 
with $\alpha_v,d_v\in \mathbb{N}_{>0}$
coincide with the set of all polynomials of the 
form $\det(E-xM)$, where $M$ is a square matrix with natural entries.

b) \cite{Lavallee}
The set of all weighted dependence polynomials 
with $d_v\in \mathbb{N}_{>0}$
coincide with the set of all polynomials of the 
form $\det(E-xM)$, where $M$ is a square matrix with nonnegative entries.

{\sc Proof}.
Prove both parts of the statement simultaneously.
When the difference of the conditions is important,
we will consider both cases.

Let $f(x) = \det(E-xM)$, where $M$ is a square matrix of an order $n$ 
with nonnegative entries.
Consider a weighted directed graph $D$ (maybe, with loops)
such that $M = (m_{uv})$ is its adjacency matrix;
every edge $(u,v)\in E(D)$ has a weight $m_{uv}$.
Construct a simple graph $G$ as follows.
The vertices of $G$ are simple (no repeated edges) 
oriented cycles of $D$ (including loops).
We put that $u,v\in V(G)$ are connected if 
their corresponding cycles in $D$ have no common vertices.
Let $c = (x_1,x_2,\ldots,x_k)$ be a simple cycle in $D$
and $v = v_c\in V(G)$ be a corresponding vertex.
We define the $\alpha_v$ to equal 
$m_{x_1 x_2}m_{x_2 x_3}\ldots m_{x_{k-1} x_k}m_{x_k x_1}$
and the $d_v$ to equal $k$, a length of the cycle.
For more details, see Example 8 from \cite{Lavallee}.

By expanding the determinant $\det(E-xM)$ in terms of permutations, 
and then decomposing the latter into not intersecting cycles, we get
$$
\det(E-xM)
 = \sum\limits_{k\geq 0}(-1)^k\sum\limits_{(c_1,\ldots,c_k)
 \in S_n(k)}m(c_1,\ldots,c_k) x^{|c_1|+\ldots+|c_k|},
$$
where $S_n(k)$ denotes the set of all $k$ mutually disjoint cycles 
from the symmetric group $S_n$,
$|c_i|$ is a length of a cycle $c_i$ and 
$m(c_1,\ldots,c_k)$ is equal to the product
$\prod\limits_{i=1}^k \alpha_{v_{c_i}}$ by the construction.
Thus, $f(x) = D_w(G,x)$.

Let $g(x) = D_w(G,x)$.
At first, prove that 
for any graph $G$ such that $\bar{G}$ is connected
and a set of natural numbers $d_v$, $v\in V(G)$, is coprime,
there exists a square matrix~$M$ such that $g(x) = \det(E-xM)$, 
and the entries of~$M$ satisfy the conditions of Theorem.
To state this we will apply the Boyle---Handelman Theorem,
but before this we need more definitions and results.

Given $X = V(G)$ with a total order $<$, 
a word $w = w_1\ldots w_t\in M_t(X,G)$,
define the two functions 
$d(w) = d_{w_1}+\ldots+d_{w_t}$,
$\alpha(w) = \alpha_{w_1}\alpha_{w_2}\ldots \alpha_{w_t}$
and the numbers 
\begin{equation}\label{m_t-weighted}
m_t(X,G) = \sum\limits_{w\in M(X,G),\,d(w) = t}\alpha(w),\quad t\in\mathbb{N}.
\end{equation}
Given a clique $B = \{b_1,\ldots,b_k\}$ in $G$, define 
$d(B) = d(b_1b_2\ldots b_k)$, 
$\alpha(B) = \alpha(b_1b_2\ldots b_k)$.
One can prove the analogue of Theorem~1.1 for the numbers $m_t(X,G)$:
$$
m_t = \sum\limits_{B\in Cl(G)}(-1)^{|B|+1}\alpha(B)m_{t-d(B)}.
$$

Also, it is easy to state the analogues 
of Lemmas~2.1--2.3 for weighted $PC$-polynomials (dependence polynomials).
Denote by $\beta_w(G)$ the largest real root of $PC_w(G,x)$.
In weighted case, the direct analogue of Theorem~2.1 is not true.
Thus, we need some additional restrictions on a graph $G$ and weights.

{\bf Lemma 11.1}.
Let $G$ be a graph such that $\bar{G}$ is connected.
Let $D_w(G,x)$ be a weighted dependence polynomial
such that the set $\{d_v\mid v\in V(G)\}$ is coprime.
Then the number $1/\beta_w(G)$ is the only complex root 
of $D_w(G,x)$ with modulus less or equal to $1/\beta_w(G)$.

{\sc Proof}.
Let $u$ be a vertex of $G$. 
Since $\bar{G}$ is connected, there exists a vertex $v\in V(G)$ 
such that $(v,u)\not\in E(G)$.
By the same reason, $G[N(v)]$ is a proper induced subgraph of $G\setminus v$. 
Consider the function 
$$
h(x)
 = \frac{D_w(G\setminus v,x)}{D_w(G,x)}
 = \frac{1}{1 - \alpha_v x^{d_v}\frac{D_w(G[N(v)],x)}{D_w(G\setminus v,x)}},\quad
f(x) = \alpha_v x^{d_v}\frac{D_w(G[N(v)],x)}{D_w(G\setminus v,x)}.
$$

By the weighted analogue of Lemma~2.3 , 
the convergence radius of the series expansion of $f(x)$ 
is greater or equal to $1/\beta_w(G\setminus v)>1/\beta_w(G)$.
By the weighted analogue of Lemma~2.1,
any complex root of $D_w(G,x)$ has a modulus not less than $1/\beta_w(G)$.
Let $\rho$ be a complex non-real root of $D_w(G,x)$ with modulus $1/\beta_w(G)$.
By the analogue of Lemma~2.3b, $\rho$ is a pole of $h(x)$.
Hence, $f(\rho) = 1$ and $|f(\rho)| = f(\rho)$.
By Daffodil Lemma, $\rho = \beta_w^{-1}(G)e^{2\pi ir/p}$,
where $p>1$ is a span of $f(x)$. 
At the same time, by the weighted analogue of Lemma~2.2,
all monomials $x^k$ of the series expansion of $f(x)$ with $d_u|k$
have nonzero coefficients. Thus, $p|d_u$.
Since the vertex $u$ was chosen arbitrarily, 
$p$ is a common divisor of the set $\{d_r\mid r\in V(G)\}$.
We arrive at a contradiction. Lemma is proved.

Let $G$ be a graph with connected complement 
and the set $\{d_v\mid v\in V(G)\}$ is coprime.
Now check the conditions of Boyle---Handelman Theorem
for the polynomial $D_w(G,x)$.
The condition (1) holds by the construction, (2) --- by Lemma~10.1, 
(3) --- by the weighted analogue of Theorem~4.1 for 
partially commutative Lie algebra graded by the weight $d(w)$.
By the M\"{o}bius inversion formula,
$s(\Lambda,n) = \sum\limits_{k|n}t(\Lambda,k)$,
so the condition ($3'$) is fulfilled
and we proved the statement with these restrictions.

At second, consider the general case of $G$ and weights $d_v$.
Let $\bar{G}$ consist of the connected components 
$\bar{H}_1,\ldots,\bar{H}_s$, 
then 
$$
D_w(G,x) = \prod\limits_{i=1}^s D_w(H_i,x)
 = \prod\limits_{i=1}^s \det(E-xM_i).
$$
For a block diagonal matrix~$M$ formed by the blocks $M_i$, $i=1,\ldots,s$,
we have $D_w(G,x) = \det(E-xM)$.

So, we may assume that $\bar{G}$ is connected.
Let $p$ be the greatest common divisor of the set 
$\{d_v\mid v\in V(G)\}$ and suppose that $p>1$.
Consider the weighted dependence polynomial $D_{w'}(G,x)$
such that $\alpha'_v = \alpha_v$, $d'_v = d_v/p$ for all $v\in V(G)$.
We have already proved that 
$D_{w'}(G,x) = \det(1-xM')$ 
for a square matrix $M'$ with required entries.
Moreover, $D_w(G,x) = \det(1-x^p M')$. 
Let $r$ be an order of the matrix $M'$.
Construct a weighted directed graph $G$ as follows.
The set $V(G)$ will include a set of $r$ vertices 
(say $1,2\ldots,r$) which index the rows and the columns of $M'$. 
For any nonzero entry $a_{ij}x^p$ of $x^p M'$,
we construct $p$ simple paths from the vertex~$i$ to~$j$ 
with the lengths $l=1,\ldots,p$ adding to $V(G)$ all required intermediate vertices.
To the first edge of any of the constructed paths from $i$ to $j$
we assign a weight $a_{ij}$, to all other edges we assign~1.
Denote by $M$ the adjacency matrix of the directed graph $G$.
It is not hard to verify~\cite{Boyle91} that $\det(1-x^pM') = \det(1-xM)$
and all other conditions on entries of $M$ are also fulfilled.

{\bf Remark 11.1}.
Note that Theorem~2.1 and its weighted analogue Lemma~11.1
are equivalent to the Perron---Frobenius theorem 
for nonnegative square matrices via Boyle---Handelman Theorem, if-part.
Indeed, in \cite{Lavallee} the second part of the proof Theorem~11.1 
was derived from the Perron---Frobenius theorem.
On the other hand, the first part of the proof of Theorem~11.1
and Lemma~11.1 provide the proof of the Perron---Frobenius theorem.

{\bf Corollary 11.1}.
a) \cite{Lavallee}
The set of all weighted $PC$-polynomials with $d_v\in \mathbb{N}_{>0}$
multiplied by $x^k$, $k\in\mathbb{N}$,
coincide with the set of all characteristic polynomials
of square matrices with nonnegative entries.

b) \cite{Reutenauer} 
The set of all weighted $PC$-polyno\-mials with $\alpha_v,d_v\in \mathbb{N}_{>0}$
multiplied by $x^k$, $k\in\mathbb{N}$,
coincide with the set of all characteristic polynomials of
square matrices with natural entries.

In 2015, C. McMullen considered~\cite{McMullen} weighted dependence polynomial 
$D_w(G,x)$ with $\alpha_v=1$ for all $v\in V(G)$.
Let us list some results from \cite{McMullen} in this case.

The growth rate of $(G,w)$ is defined by 
$\lambda(G,w) = \lim\limits_{t\to\infty} n_t^{1/t}$,
where 
$n_t = \sum\limits_{i=1}^t m_i$ 
and numbers $m_i$ are defined by~\eqref{m_t-weighted}.
If $d_v = 1$ for all $v\in V(G)$, then 
$\lambda(G,w) = \beta(G)$.

{\bf Theorem 11.2} \cite{McMullen}.
The largest positive root of $PC_w(G,x)$ equals $\lambda(G,w)$,
and the function 
$h(w) = \ln \lambda(G,w)$ is convex.
Provided $\bar{G}$ is connected and $|V(G)|>1$, 
the function $h(w)$ is strictly convex, real-analytic, 
and $h(w)\to\infty$ at the boundary of the cone of positive weights.

A weight $w$ is called {\it admissible} if $d(K)\leq 1$ 
for all cliques $K$ of $G$. The number 
$\lambda(G) = \inf \{\lambda(G,w)\mid w\ \mbox{is admissible}\}$
is called the {\it minimal growth rate} of~$G$.
An admissible weight $w$ is called {\it optimal} if $\lambda(G) = \lambda(G,w)$.
An admissible weight $w$ is {\it symmetric} if it is invariant under $\Aut(G)$
and {\it maximal} if there is no other admissible weight $w'$ 
with $d_v\leq d'_v$ for all $v$.

For any induced subgraph $H$ of $G$,
we have $\lambda(G,w)\geq \lambda(H,w|_{H})$ (analogue of Lemma~2.3a)
and hence $\lambda(G)\geq \lambda(H)$~\cite{McMullen}.

{\bf Statement 11.1} \cite{McMullen}.
Given $M > 0$, there are only finitely many graphs $G$ 
with connected complement and $\lambda(G)\leq M$.

{\sc Proof}.
It suffices to show that $\lambda(G)$ is large whenever $|V(G)|$ is large.
By Ramsey theory, if $|V(G)|\geq \binom{2n-2}{n-1}$, 
either $G$ or $\bar{G}$ contains a clique $K$ of size~$n$.
If $K$ is a clique of $\bar{G}$ then $G$~contains 
an empty induced subgraph $\bar{K}$ with $n$ vertices and 
$\lambda(G)\geq n$ (see Statement~3.2a).

Suppose that $K$ is a clique of $G$.
The optimal weight~$w$ for $G$ satisfies $d(K)\leq1$,
so $d_u\leq 1/n$ for some $u\in K$.
Since $\bar{G}$ is connected, there is a vertex 
$v\in V(G)$ such that $(u,v)\not\in E(G)$ and $d_v\leq1$.
Let $H$ be an induced subgraph of $G$ with $V(H) = \{u,v\}$. 
Then 
$\lambda(G) = \lambda(G,w) \geq \lambda(H,w|_H) \geq \lambda_n$, where 
$t=1/\lambda_n$ is the smallest root of 
$1 - t^{1/n} - t$.
It~is easy to show that $\lambda_n = O(n/\ln n)$ and we are done.

{\bf Statement 11.2} \cite{McMullen}.
If $\bar{G}$ has no isolated vertices, 
then $G$ carries a~unique optimal weight.
This weight is maximal and symmetric.

{\sc Proof}.
Suppose that $\bar{G}$ is connected and $|V(G)|>1$.
By Theorem~11.2, the function 
$h(w) = \ln\lambda(G,w)$
is strictly convex and tends to infinity at the boundary 
of the cone of positive weights.
The admissibility conditions are convex and
bound $h(w)$ from below, so there is a unique optimal weight.

In general, $G$ is a supergraph of complete multipartite graph
with components $H_1$, \ldots, $H_k$ such that $\bar{H}_i$ is connected for all~$i$.
We have $\lambda(G,w) = \max\limits_{i=1,\ldots,k}\lambda(H_i,w|_{H_i})$.
Since $\bar{G}$ has no isolated vertices, 
$|V(H_i)|>1$ for all $i=1,\ldots,k$.
So each $H_i$ has a unique optimal weight $w_i$.
Then the optimal weight for $G$ is the unique convex combination
$w = \sum\limits_{i=1}^k \alpha_i w_i$ such that 
$\lambda(H_i,w|_{H_i}) = \lambda(H_i)^{1/\alpha_i}=\lambda(H_j)^{1/\alpha_j}$,
$i,j=1,\ldots,k$.
By convexity and uniqueness, the optimal weight $w$ must be maximal
and symmetric.

{\bf Statement 11.3} \cite{McMullen}.
For any graph $G$, the growth rate $\lambda(G)$ is algebraic over $\mathbb{Q}$.

{\bf Example 11.1} \cite{McMullen}.
a) We have $\lambda(K_2\cup K_2) = 6 + 4\sqrt{2}$. 
By symmetry the optimal weights are all $1/2$. 
Thus $D_w(t) = 1 - 4t^{1/2} + 2t$.

b) We have $\lambda(K_n\cup K_1) = 2^n$. 
The optimal weights are $1/n$ on the vertices of $K_n$ and 1 on the remaining vertex. 
Thus, $D_w(t) = (1-t^{1/n})^n - t$
and its smallest real root is~$1/2^n$.

c) We have $\lambda(P_4) = 8$ for the path graph with 4~vertices. 
In this case by symmetry and maximality there is an 
$a \in (0;1)$ such that the optimal weights on the 
consecutive vertices of $P_4$ are $(a,1-a,1-a,a)$. 
Hence, $D_w(t) = 1 - 2t^a - 2t^{1-a} + t^{2-2a} + 2t$.
To determine the value of $a$, we observe that both $D_w$ and 
$\dfrac{dD_w}{da}$ must vanish at $(a,1/\lambda)$; hence 
$a = 2/3$ and $1/\lambda = 1/8$.

d) For $G = P_3\cup K_1$ we have 
$\lambda = \lambda(G) = 3 + 2\sqrt2$.
By symmetry and maximality, there is an 
$a\in (0;1)$ such that $w$ takes the values
$(a,1-a,a)$ along the vertices of $P_3$ 
and $w(v) = 1$ at the remaining vertex of $G$. 
Thus, $D_w(t) = 1-2t^a - t^{1-a} + t$.
By optimality we must have 
$\dfrac{dD_w}{da} = 0$ at $t = 1/\lambda$. 
This gives a second relation 
$2(1/\lambda)^a = (1/\lambda)^{1-a}$; 
when combined with $D_w(1/\lambda) = 0$,
it implies $\lambda = 3 + 2\sqrt2$ and $a = \ln(2\lambda)/\ln(\lambda^2)$.

C. McMullen used the stated results to apply 
for spectral radius of called reciprocal positive matrices. 
A square matrix $A = (a_{ij})$ is called reciprocal if $a_{ij} = 1/a_{ji}$
for all $i,j$.

{\bf Theorem 11.3} \cite{McMullen}. 
The minimum value of the spectral radius $\rho(A)$ 
over all recipro\-cal matrices $A \in M_{2g}(\mathbb{N})$, $g\geq 2$, 
is given by the largest root of the polynomial
$t^{2g} - t^g(1 + t + t^{-1}) + 1$.
Consequently $\rho(A)^g \geq (3 + \sqrt{5})/2$ for all such $A$.

\newpage
In 2005, A. Scott and A. Sokal found~\cite{Scott} the deep connection
between weighted (in)dependence polynomials and Lov\'{a}sz local lemma.
Below we stated the particular case of their results.

For a graph~$G$, we define the weighted independence polynomial 
$I_w(G,x)$ as $D_w(\bar{G},x)$\footnote{It is the only place 
where we use the definition of $I(G,x)$ as $D(\bar{G},x)$
instead of $C(\bar{G},x)$.}\!.

{\bf Theorem 11.4} \cite{CsikvariThesis,Scott}.
Assume that given a graph $G$ and there is an event $A_i$
assigned to each vertex~$i$. 
Assume that $A_i$ is totally independent of the events
$\{A_k\mid (i,k)\not\in E(G)\}$.
Set $P(A_i) = p_i$ and define the weighted independence polynomial
$I_w(G,x)$ with $\alpha_i = p_i$ and $d_i = 1$ for all $i=1,\ldots,n$.

a) Assume that $I_w(G,t)>0$ for $t\in[0;1]$, then 
$P\big(\bigcap\limits_{i=1}^n \overline{A_i}\big)\geq I_w(G,1)>0$.

b) Assume that $I_w(G,t) = 0$ for some $t\in [0;1]$.
Then there exist a probability space and a~family of events $B_i$, $i=1,\ldots,n$,
with $P(B_i)\leq p_i$ and with dependency graph~$G$ such that
$P\big(\bigcap\limits_{i=1}^n \overline{B_i}\big) = 0$.

{\sc Proof}.
a) Let us define the events $B_i$, $i=1,\ldots,n$, 
on a new probability space as follows
$$
P\bigg(\bigcap\limits_{i\in S}B_i\bigg)
= \begin{cases}
\prod\limits_{i\in S}p_i, & S\ \mbox{is independent in}\ G,\\
0, & \mbox{otherwise}.
\end{cases}
$$
The expression
$P((\bigcap\limits_{i\in S}B_i)\cap(\bigcap\limits_{i\not\in S}B_i))$
equals zero if $S$ is not an independent set.
So assume that $S$ is an independent set, then 
\begin{multline}\label{Scott-Socal:ex1}
P\bigg(\bigg(\bigcap\limits_{i\in S}B_i\bigg)
 \cap \bigg(\bigcap\limits_{i\not\in S}\overline{B_i}\bigg)\bigg)
 = \sum\limits_{S\subseteq I}(-1)^{|I|-|S|}P\bigg(\bigcap\limits_{i\in I}B_i\bigg) \\
 = \sum\limits_{S\subseteq I,\,I\in \mathcal{I}}(-1)^{|I|-|S|}\prod\limits_{i\in I}p_i
 = \bigg(\prod\limits_{i\in S}p_i\bigg)I_w(G\setminus N[S],1),
\end{multline}
where $\mathcal{I}$ is the set of all independent sets of $G$
and $N[S]$ denotes the set $S$ together with all their neighbors.
Let $\beta = \beta_w(\bar{G})$.
By the conditions, $\beta<1$.
By the weighted analogue of Lemma~2.3 \cite{CsikvariThesis}, we have
$\beta_w(\overline{G\setminus N[S]})\leq \beta < 1$.
This means that the RHS of \eqref{Scott-Socal:ex1} is nonnegative for all $S\subseteq\{1,\ldots,n\}$.
Hence, we have defined a probability measure on the generated
$\sigma$-algebra $\sigma(B_i\mid i=1,\ldots,n)$.

Let us state that $(B_i)_{i=1}^n$ minimizes the expression
$P\big(\bigcap\limits_{i=1}^n B_i\big)$ 
among the families of events with dependency graph~$G$.
For $S\subseteq\{1,\ldots,n\}$, we set 
$P_S = P\big(\bigcap\limits_{i\in S}\overline{A_i}\big)$ and
$Q_S = P\big(\bigcap\limits_{i\in S}\overline{B_i}\big)$.
Prove by induction on $|S|$ that $P_S/Q_S$ is monotonically increasing in $S$.
For $Q_S$, we have
\begin{gather*}
Q_S
 = \sum\limits_{I\subseteq S}(-1)^{|I|}P\bigg(\bigcap\limits_{i\in I}B_i\bigg)
 = \sum\limits_{I\subseteq S,\,I\in \mathcal{I}}(-1)^{|I|}\prod\limits_{i\in I}p_i
 = I_w(G[S],1)>0; \\
Q_{S\cup \{j\}}
 = I_w(G[S\cup \{j\}],1) = I_w(G[S],1) - p_j I_w((G[S\setminus N[j]],1)
 = Q_S-p_j Q_{S\setminus N[j]},
\end{gather*}
where $j\not\in S$. On the other hand, 
$$
P_{S\cup \{j\}}
 = P_S-P\bigg(A_j\cap\bigg(\bigcap\limits_{i\in S}\overline{A_i}\bigg)\bigg)
 \geq P_S-P\bigg(A_j\cap\bigg(\bigcap\limits_{i\in S\setminus N[j]}\overline{A_i}\bigg)\bigg)
 \geq P_S-p_j P_{S\setminus N[j]}.
$$

Now show that 
$P_{S\cup \{j}\}/Q_{S\cup \{j\}}\geq P_S/Q_S$ 
for all $j\not\in S$.
By the induction hypothesis, we compute
\begin{multline*}
P_{S\cup \{j}\}Q_S-Q_{S\cup \{j\}}P_S
 \geq (P_S-p_j P_{S\setminus N[j]})Q_S
 -(Q_S-p_j Q_{S\setminus N[j]})P_S \\
 = p_j(P_S Q_{S\setminus N[j]}-Q_S P_{S\setminus N[j]})\geq 0.
\end{multline*}

Since $P_S/Q_S$ is monotonically increasing in $S$, we have
$P_{V(G)}/Q_{V(G)}\geq P_{\emptyset}/Q_{\emptyset} = 1$. We have proved a).

b) Let us use the construction of the events $B_i$ with probability
$\beta^{-1} p_i$, $i=1,\ldots,n$.
Then this will define a probability measure again the same way (as in a). So, we have
$P\big(\bigcap\limits_{i=1}^n\overline{B_i}\big) = I_w(G,1/\beta) = 0$.

{\bf Corollary 11.2} \cite{CsikvariThesis,Scott}.
Let $A_i$, $i=1,\ldots,n$, be events with dependence graph $G$
such that $P(A_i)\leq t$ for all $i=1,\ldots,n$. Then 
$P\big(\bigcap\limits_{i=1}^n \overline{A_i}\big)>0$
if and only if $t\leq 1/\beta(\bar{G})$.

{\bf Corollary 11.3}.
a) \cite{Scott}
For any graph $G$, 
$\beta(G)\leq \frac{d^{d}}{(d-1)^{d-1}}<ed$,
where $d = \Delta(\bar{G}){\geq}2$.
For $d = 1$, we have $\beta(G) = 2$.

b) Let $G$ be a $(n-d-1)$-regular graph, then 
$1+d\leq\beta(G) < ed$.

{\sc Proof}.
a) The statement follows from the Lov\'{a}sz local lemma in the Shearer's ver\-sion~\cite{Shearer}, 
Corollary~11.2 and Lemma~2.6b.

b) The upper bound follows from a), the lower bound follows from Theorem~7.1:
$\beta(G)\geq 1+\frac{2\bar{k}}{n} = 1 + d$.

{\bf Remark 11.2}.
About more precise bounds on $\beta(G)$ in terms of the maximal degree of a graph 
see in \cite{BencsZeroFree,PetersRegts}.

{\bf Corollary 11.4}.
Let $G$ be a graph with $n\gg1$ vertices and $k$ edges,
$A_i$, $i=1,\ldots,n$, are events with dependence graph~$G$.
Then $P\big(\bigcap\limits_{i=1}^n \overline{A_i}\big)>0$, if 

a) $k\leq 0.24326n^2$ and $P(A_i)\leq 4/(3n)$, $i=1,\ldots,n$, or

b) $k\leq n^2/4$ and $P(A_i)\leq 1.3316/n$, $i=1,\ldots,n$.

{\sc Proof}.
It follows from Corollary~8.4 and Corollary~11.2.

Note that in analogous way one can produce another global versions
of Lov\'{a}sz local lemma in case $k\leq \alpha n^2$.
It is an open question: 
Wether such version of Lov\'{a}sz local lemma like in Corollary~11.4 is useful?

\newpage
\section{Claw-free graphs. Matching polynomial}

\subsection{Chudnovsky---Seymour theorem}

A graph $G$ is called claw-free if $G$ contains no 
induced subgraphs isomorphic to $K_{1,3}$.

Remind that the independence polynomial $I(G,x)$ is defined as follows:
$I(G,x) = 1+\sum\limits_{k=1}^{\alpha(G)}s_k(G)x^k$,
where $\alpha(G)$ is the size of maximal coclique in $G$
and $s_k$ denotes the number of all cocliques of size~$k$ in $G$.
By Lemma~1.5, we get immediately the formulas
\begin{gather}
I(G,x) = I(G \setminus v,x) + xI(G\setminus [N(v)],x),\quad v\in V(G);
\label{IndepMainFormula} \\
I'(G,x) = \sum\limits_{v\in V(G)}I(G\setminus [N(v)],x);
\label{IndepDiffFormula} \\
I(G_1\cup G_2,x) = I(G_1,x)I(G_2,x).
\label{IndepUnionFormula}
\end{gather}

As earlier, given a subset $H\subset V(G)$,
we denote the set $\big(\bigcup\limits_{v\in H}N[v]\big)\cup H$ by $N[H]$.
The analogues of the Christoffel---Darboux identities 
were successfully applied in \cite{Heilman72} for matching polynomial. 
Let us state them in more general case.

{\bf Lemma 12.1} \cite{Bencs}.
Let $G$ be a graph, $v\in V(G)$,
$B_v$ be the set of induced connected, 
bipartite subgraphs of $G$ containing the vertex~$v$.
For $H\in B_v$, let $A(H)$ be a part containing $v$ 
and $B(H)$ be another one.
Let $a(H) = |A(H)|$, $b(H) = |B(H)|$, then 
\begin{multline}\label{Christoffel-Darboux}
I(G, x)I(G \setminus v, y) - I(G \setminus v, x)I(G,y) \\
 = \sum\limits_{H\in B_v}
(x^{a(H)}y^{b(H)} - x^{b(H)}y^{a(H)})I(G \setminus N[H],x)I(G \setminus N[H],y).
\end{multline}

{\sc Proof}.
Denote by $\mathcal{I} = \mathcal{I}(G)$ the set of all independent sets in~$G$,
$\mathcal{I}_1 = \mathcal{I}(G)\times \mathcal{I}(G\setminus v)$,
$\mathcal{I}_2 = \mathcal{I}(G\setminus v)\times \mathcal{I}(G)$.
The LHS of \eqref{Christoffel-Darboux} equals
\begin{equation}\label{Christoffel-Darboux:eq1}
\sum\limits_{(A,B)\in \mathcal{I}_1}x^{|A|}y^{|B|}
- \sum\limits_{(A,B)\in \mathcal{I}_2}x^{|A|}y^{|B|}
= \sum\limits_{(A,B)\in \mathcal{I}_0}x^{|A|}y^{|B|}
- \sum\limits_{(A,B)\in \mathcal{I}_0}x^{|B|}y^{|A|},
\end{equation}
where $\mathcal{I}_0 = \mathcal{I}_1\setminus
(\mathcal{I}(G\setminus v)\times \mathcal{I}(G\setminus v))$.

Note that for any pair $(A, B) \in \mathcal{I}_0$ we have $v\in A$, $v\not\in B$.
Let $P(A,B)$ be the connected component of the induced subgraph $G[A\cup B]$
which contains~$v$,
$A' = A\cap V(P(A,B))$, $B' = B\cap V(P(A,B))$. 
Since $A$ and $B$ are independent, we get $A'\cap B' = \emptyset$,
and $P(A,B)$ is a connected bipartite subgraph containing~$v$.
Hence, 
\begin{multline*}
\sum\limits_{(A,B)\in\mathcal{I}_0}x^{|A|}y^{|B|}
 = \sum\limits_{(A,B)\in\mathcal{I}_0}x^{a(P(A,B))}y^{b(P(A,B))}
 x^{|A|-a(P(A,B))}y^{|B|-b(P(A,B))} \\
 = \sum\limits_{H\in B_v}x^{a(H)}y^{b(H)}
\sum\limits_{K,L\in\mathcal{I}(G\setminus N[H])}x^{|K|}y^{|L|} \\
= \sum\limits_{H\in B_v}x^{a(H)}y^{b(H)}I(G\setminus N[H],x)I(G\setminus N[H],y).
\end{multline*}

Analogously, we express
$$
\sum\limits_{(A,B)\in\mathcal{I}_0}x^{|B|}y^{|A|}
 = \sum\limits_{H\in B_v}x^{b(H)}y^{a(H)}
I(G\setminus N[H],x)I(G\setminus N[H],y).
$$
It remains to insert the obtained relations in~\eqref{Christoffel-Darboux:eq1}.

{\bf Lemma 12.2} \cite{Bencs}.
Let $G$ be a graph, $B$ be the set of induced connected, 
bipartite subgraphs of $G$.
For any $H\in B$, let $P(H)$ denote one part of $H$
and $R(H)$ denote the other.
Let $p(H) = |P(H)|$, $r(H) = |R(H)|$, then
\begin{multline}\label{Christoffel-Darboux2}
yI(G, x)I'(G, y) - xI'(G, x)I(G, y)  \\
= \sum\limits_{H\in B}
(p(H) - r(H))(x^{p(H)}y^{r(H)} - x^{r(H)}y^{p(H)})
I(G\setminus N[H], x)I(G \setminus N[H], y).
\end{multline}

{\sc Proof}.
Let $n = |V (G)|$. 
From the formulas \eqref{IndepMainFormula}, \eqref{IndepDiffFormula},
we deduce 
\begin{equation}\label{IndepPolyConseq}
\sum\limits_{v\in V(G)}
I(G\setminus v,x) = nI(G, x) - xI'(G,x).
\end{equation}
Let us sum the identity \eqref{Christoffel-Darboux} for all $v\in V(G)$
and apply \eqref{IndepPolyConseq}:
\begin{multline*}
\sum\limits_{v\in V(G)}(I(G,x)I(G\setminus v,y) - I(G\setminus v,x)I(G,y)) \\
= I(G,x)\sum\limits_{v\in V(G)}I(G\setminus v, y)
 - I(G,y)\sum\limits_{v\in V(G)}I(G\setminus v,x) \\
= I(G,x)(nI(G,y) - yI'(G,y)) - I(G,y)(nI(G,x) - xI'(G,x)) \\
= xI'(G,x)I(G,y) - yI(G,x)I'(G,y) \\
= \sum\limits_{v\in V(G)}\sum\limits_{H\in B_v}
(x^{a(H)}y^{b(H)} - y^{a(H)}x^{b(H)})I(G\setminus N[H],x)I(G\setminus N[H],y) \\
= \sum\limits_{H\in B}(p(H) - r(H))(x^{p(H)}y^{r(H)} - y^{p(H)}x^{r(H)})
I(G\setminus N[H],x)I(G\setminus N[H],y).
\end{multline*}

Y. Hamidoune in 1990 \cite{Hamidoune} and R. Stanley in 1998 \cite{Stanley98}
conjectured and M. Chudnov\-sky and P. Seymour in 2004 proved that

{\bf Theorem 12.1} \cite{Bencs,Chudnovsk}.
If $G$ is a claw-free graph then all roots of $I(G,x)$ are real.

{\sc Proof}.
Note that every induced connected bipartite subgraph $H$ 
of a claw-free graph $G$ is a path or a cycle. 
Indeed, claw-freeness implies that degree of every vertex of $H$
is at most 2, and connectedness of $H$ implies that it is a path or a cycle. 
Let $H$ be an even path or a cycle, then $p(H) - r(H) = 0$.
Suppose that $H$ is an odd path, then choose 
$P(H)$ and $R(H)$ such that $p(H) - r(H) = 1$. 
Denote by $\mathcal{P}$ the set of all odd paths in $G$, 
then by~\eqref{Christoffel-Darboux2} we have
\begin{equation}\label{ChudnovskyProof}
\frac{yI(G, x)I'(G, y) - xI'(G, x)I(G, y)}{x - y}
 = \sum\limits_{H\in \mathcal{P}}x^{r(H)}y^{r(H)}I(G\setminus N[H], x)I(G\setminus N[H], y).
\end{equation}

Suppose that $G$ is a counterexample with 
the smallest number of vertices and $z$ is a complex non-real root of $I(G,x)$.
By~\eqref{IndepUnionFormula}, $G$ is connected.
Putting in~\eqref{ChudnovskyProof} the values $x = z$, $y = \bar{z}$,
we get~0 on the LHS and a positive number on the RHS, 
since $G \setminus N[H]$ is an induced subgraph 
without claws and with less vertices than~$G$.
Let us show that the set $\mathcal{P}$ is not empty.
Indeed, if $G$ has no paths of length~3, 
then $G$ is a complete graph~\cite{ChudnovskInduced,Lozin}.
But $I(K_n,x) = 1 + nx$ has only one real root.

{\bf Corollary 12.1}.
If $\bar{G}$ is a claw-free graph then all roots of $PC(G,x)$ are real.

{\bf Corollary 12.2}.
Let $G$ be a graph with $n$ vertices and $k$ edges.

a) If $\bar{G}$ is claw-free then 
$n-\frac{2k}{n}\leq \beta(G)\leq n-\frac{k}{n}$
and $\frac{1}{n}\leq e(G)\leq \frac{2}{n}$.

b) If $\beta(G)>n-\frac{k}{n}$ then $\bar{G}$ has a claw.

c) Let $G$, $\bar{G}$ be claw-free~\cite{Pouzet}, then
$\beta(G)+\beta(\bar{G})\leq \frac{3n+1}{2}$,
$\beta(G)\beta(\bar{G})\leq \big(\frac{3n+1}{4}\big)^2$.

{\sc Proof}.
a) The lower bound for $\beta(G)$ was stated in Theorem~7.1,
the upper one follows from Theorem~11.1 and Statement~8.1.

b), c) It follows fom a). 

\subsection{Matching polynomial}

For a graph~$G$, $m_k$ denotes the number of matchings with $k$ edges in $G$
($m_0 = 1$, $m_1 = |E(G)|$). 
Define the matching polynomial of~$G$ as 
$\mu(G,x)=\sum\limits_{k=0}^{\nu(G)}(-1)^k m_k x^{n-2k}$,
where $n = |V(G)|$ and $\nu(G)$ is the maximum of sizes of matchings in $G$. 
A perfect matching exists if and only if $\nu(G) = n/2$.

Another very close polynomial, so called matching-generating polynomial,
was defined: 
$M(G,x)=\sum\limits_{k=0}^{\nu(G)}m_kx^k$.
Both polynomials are connected in the following way:
\begin{equation}\label{MatchingPoly}
\mu(G,x) = x^n M(G,-x^{-2}),\quad M(G,x) = I(L(G),x),
\end{equation}
where $L(G)$ denotes the line graph of $G$: 
$V(L(G)) = E(G)$,
and two vertices in $L(G)$ are connected 
if and only if corresponding edges have a common vertex in~$G$.

Let $T_n$, $U_n$ be the Chebyshev functions of the first and second kind,
let $He_n$, $H_n$ be the two standard forms of the Hermit\'{e} polynomials,
let $L_n$ be the Laguerre polynomial. 

{\bf Example 12.1} \cite{Heilman72,Gutman1979}.
a) $\mu(K_n,x) = He_n(x) = 2^{-n/2}H_n(x/\sqrt{2})$,

b) $\mu(K_{n,n},x) = (-1)^n L_n(x^2)$,

c) $\mu(C_n,x) = 2T_n(x/2)$, where $C_n$ is a cycle of length $n$,

d) $\mu(P_n,x) = \frac{2}{\sqrt{4-x^2}}U_{n+1}(x/2)$,
where $P_n$ is a path of length~$n$.

From Theorem~12.1 follows the result of O. Heilmann and E. Lieb (1972):

{\bf Corollary 12.3} \cite{Heilman72}.
All roots of $\mu(G,x)$ and $M(G,x)$ are real.

{\sc Proof}.
By the construction, the line graph $L(G)$ is a claw-free graph.
Hence, by Theorem~12.1, the polynomial $M(G,x) = I(L(G),x)$
has only real (negative) roots. By the formula~\eqref{MatchingPoly},
all roots of $\mu(G,x)$ are real. Corollary is proved.

Denote the largest root of $\mu(G,x)$ as $t(G)$. 
By the definition, $t^2(G) = \beta(\overline{L(G)})$. 
The next bound for $t(G)$ was stated by D.C. Fisher and J. Ryan in 1992:

{\bf Statement 12.1} \cite{Fisher1992-Match}.
For a graph~$G$ with $n$ vertices and $k$ edges, 
$t^2(G)\geq \frac{4k}{n}-1$.

{\sc Proof}.
By \eqref{MatchingPoly}, $t^2(G) = \beta(\overline{L(G)})$.
Applying Corollary~7.5, the lower bound on the spectral radius~\cite{Collatz}
of the form $\rho(G)\geq 2|E(G)|/|V(G)|$ and the inequality between means, we get
\begin{multline*}
t^2(G)\geq 1+\rho(L(G))
 \geq 1+\frac{2}{k}\bigg(\sum\limits_{v\in V(G)}\frac{d(v)(d(v)-1)}{2}\bigg) \\
 = 1+\frac{2}{k}\bigg(\sum\limits_{v\in V(G)} d^2(v) - k\bigg)
 = \frac{1}{k}\sum\limits_{v\in V(G)} d^2(v) - 1
 \geq \frac{4k}{n}-1,
\end{multline*}
where $d(v)$ denotes the degree of a vertex $v\in V(G)$.
Statement is proved.

In the famous article of O. Heilmann and E. Lieb written in 1972~\cite{Heilman72},
the bounds on $t(G)$ in terms of the maximal degree were obtained.

{\bf Statement 12.2} \cite{Fisher1992-Match,Heilman72}.
Let $\Delta>1$ denote the maximal degree of graph~$G$, then 
the following inequalities hold
$\sqrt{\Delta}\leq t(G)\leq 2\sqrt{\Delta-1}$.

{\sc Proof}.
For any (not necessary induced) subgraph~$H$ of~$G$,
the induced subgraph $\overline{L(H)}$ of $\overline{L(G)}$
corresponds. By Lemma~2.3a, we have $t(H)\leq t(G)$
for any subgraph~$H$ of~$G$. 
Consider the subgraph $H$ consisting on the vertex 
of the maximal degree~$\Delta$ and all its outgoing edges. 
Thus, we have the lower bound. 

The detailed proof of the upper bound $t(G)\leq 2\sqrt{\Delta-1}$ 
see in \cite{LovaszPlummer}, it contained three steps:
a) matching polynomial of a forest coincides with its characteristic polynomial, 
b) matching polynomial of any connected graph $G$ is a factor 
of the matching polynomial of the path-tree associated with~$G$,
c) $\rho(G)\leq 2\sqrt{\Delta-1}$
for a forest $G$ with the maximal degree~$\Delta$.

We will clarify only the part a). 
Let $A$ denote the adjacency matrix of a forest~$T$
and $E$ denote the unit matrix.
Expanding the determinant $|\lambda E-A|$ in terms of permutations, 
we get a nonzero summand only if 
the corresponding permutation consists of $s$ cycles of length~2 
and of $n-2s$ cycles of length~1, $0\leq s\leq \big[\frac{n}{2}\big]$.
Cycles of length~1 give the factor $\lambda^{n-2s}$
and cycles of length~2 correspond to a matching of size $s$.
This observation implies~a).

{\bf Remark 12.1}.
By Corollary~11.3a, we immediately get the upper bound
$t(G)\leq \sqrt{2e(\Delta-1)}$, since 
the degree of all vertices of $L(G)$ is not greater than $2(\Delta-1)$.

In a similar to the proof of Statement~12.2 way,
with exchange of associated path-tree on simplicial clique,
J. Leake and N. Ryder in 2016 stated the bounds for $\beta(G)$ 
in a claw-free case:

{\bf Statement 12.3} \cite{LeakeRyder}.
If $\bar{G}$ is a claw-free graph then 
$\alpha(G)\leq \beta(G)\leq 4\max\{\alpha(G)-1,\Delta(G)\}$.

The lower bound of Statement~12.3 follows from Statement~3.2a.
We leave the proof of the upper bound.

Let us list another interesting results about roots of matching polynomial.

{\bf Statement 12.4} \cite{Heilman72}.
If a graph $G$ contains a Hamiltonian path, then all roots of $\mu(G,x)$ are simple.

{\bf Statement 12.5} \cite{Godsil93}.
The maximum multiplicity of a root of $\mu(G,x)$ is at most 
equal to the number of vertex-disjoint paths required to cover~$G$.
The number of distinct roots of $\mu(G,x)$ is at least equal to the length of the longest path in~$G$.

In 2014, D. Bevan interpreted the number $t^2(G)$ 
as the growth rate of a class of geometric grids
corresponding to~$G$~\cite{Bevan}.

Let us briefly state some results about matchings 
which are very close to the results from \S3 about the bounds 
on the number of all (co)cliques in a graph. 
Denote by $m(G)$ the number of matchings in a graph~$G$
and by $pm(G) = m_{|V(G)|/2}(G)$ the number of perfect matchings. 

{\bf Statement 12.6} \cite{Davies2017}.
For all $d$-regular graphs $G$ on $n$ vertices (where $2d$ divides $n$), 
$m_k(G) \le 2 \sqrt{n} \cdot m_k(H_{d,n})$,
where $H_{d,n}$ denotes the $d$-regular, $n$-vertex graph that 
is the disjoint union of $n/(2d)$ copies of $K_{d,d}$.

By the edge occupancy fraction \cite{Davies2017} we mean 
$\alpha^M(G,x) = \dfrac{xM'(G,x)}{|E(G)|M(G,x)}$.

{\bf Statement 12.7} \cite{Davies2017}.
For any $d$-regular graph $G$, the following inequalities hold 

a) $\alpha^M(G,x)\leq\alpha^M(K_{d,d},x)$,

b) $M(G,x)\leq M(K_{d,d},x)^{n/(2d)}$.

Moreover, the maximum is achieved only by unions of copies of $K_{d,d}$.

{\bf Corollary 12.4} \cite{Davies2017}.
For any $d$-regular graph~$G$, 
$m(G)\leq m(K_{d,d})^{n/(2d)}$.

{\sc Proof}. 
It is enough to consider Statement~12.6 with $x = 1$.

{\bf Corollary 12.5} \cite{Bregman,Davies2017}.
For any $d$-regular graph~$G$, 
$pm(G)\leq (d!)^{n/(2d)}$.

{\sc Proof}. 
For $x\to\infty$, the inequality 
$\alpha^M(G,x)\leq\alpha^M(K_{d,d},x)$
implies the inequality between leading terms
$pm(G)\leq pm(K_{d,d})^{n/(2d)} = (d!)^{n/(2d)}$.

{\bf Theorem 12.2}~\cite{CsikvariMatch}.
Let $G$ be a $d$-regular bipartite graph on $2n$ vertices,
let $p = \frac{k}{n}$, and $p_{\mu}$ be the probability 
that a random variable with distribution Binomial$(n,p)$ 
takes its mean value $\mu=k$. Then 
$$
m_k(G)\geq p_{\mu}{n \choose k}^2\left(\frac{d-p}{d}\right)^{n(d-p)}(dp)^{np}.
$$

{\bf Corollary 12.6} \cite{CsikvariMatch,Schrijver}. 
Let $G$ be a $d$--regular bipartite graph on $2n$ vertices, then
$pm(G)\geq \left(\frac{(d-1)^{d-1}}{d^{d-2}}\right)^n$.

In 2008, L. Gurvits \cite{Gurvitz} proved the general result
from which both van der Waerden conjecture about the permanents
of double stochastic matrices \cite{Egorychev,Waerden}
and Corollary~12.6 follow (see exposition \cite{Laurent}).

\newpage
\section{Survey. Open problems}

Below we consider only some (not all) questions concerned
with clique-type polynomi\-als.
Note the recent works devoted to independence
and matching polynomials defined for hypergraphs~\cite{Barvinok,Guo,Makowsky,HyperSupertrees,Hypergraph}. 

\subsection{Real-rootedness of $PC$-polynomial}

Let us once again note the nice Theorem~11.1 stating that 
all roots of $PC$-polynomial of complement to a claw-free graph are real. 
This result was proved by M. Chudnovsky and P. Seymour in 2004 (published in 2007).
After that, a lot of other proofs of Theorem~11.1 appear~\cite{Bencs,Lass,LeakeRyder}.

In 2014~\cite{Balhs2014}, P. Bahls, E. Bailey, and M. Olsen showed that given a graph $G$ 
with some constraints there exists a construction to get a graph $G'$ 
which $PC$-polynomial has only real roots. 
In 2017~\cite{Bencs17}, F. Bencs stated real-rootedness 
of $PC$-polynomials complement to some families of trees.

In the paper of J. Brown and R. Nowakowski of 2005~\cite{BrownNowak},
it was stated that for almost all graphs
$PC$-polynomial has a complex non-real root.
It was claimed that the Sturm sequence (starting with $PC(G,x)$ and $PC'(G,x)$)
has less sign changes in $(-\infty,+\infty)$ than it should be 
for a real-rooted polynomial. 
For this claim, the clique numbers $c_k(G)$ lie in some neighbourhoods 
of $\binom{n}{k}2^{-\binom{k}{2}}$, the expected values of the numbers of cliques of size~$k$.
Unfortunately, it looks like the proof is not complete right.
The problem is the following: the coefficients of $PC(G_{n,p})$,
the $PC$-polynomial of random graph~$G_{n,p}$, lie in these neighbourhoods.
But all roots of $PC(G_{n,p})$ are real (Theorem~5.1a).

By this reason, we state

{\bf Problem 13.1}.
To state if for almost all graphs $PC$-polynomial has a non-real root.

Maybe, there exists the constant $a\in(0,1]$ such that 
the part of real roots of $PC$-polynomial of almost all graphs equals~$a$. 
If yes, then Problem~13.1 could be reformulated as follows:
Is it true that $a = 1$?

The results about part of non-real roots of $PC$-polynomials of all graphs with $n$~verti\-ces
for $1\leq n\leq6$ are gathered below. 

\begin{center}
Table: Part of non-real roots of $PC$-polynomials of small graphs
\end{center}
\renewcommand{\arraystretch}{1.5}
\begin{center}
\begin{tabular}{|c|c|c|}
\hline
$n$ & \quad part of $PC$-polynomials with non-real roots \quad & \quad part of non-real roots \quad \\
\hline
1 & $0/2^0 = 0\%$ & $0/1 = 0\%$ \\
2 & $0/2^1 = 0\%$ & $0/3 = 0\%$ \\
3 & $0/2^3 = 0\%$ & $0/16 = 0\%$ \\
4 & $4/2^6 \approx 6.25\%$ & $8/151\approx 5.3\%$ \\
5 & $135/2^{10}\approx 13.18\%$ & $270/2750\approx 9.81\%$ \\
6 & $4666/2^{15}\approx 14.24\%$ & $9344/97839\approx 9.55\%$ \\
\hline
\end{tabular}
\end{center}

\subsection{Unimodality and log-concavity of clique polynomial}

A polynomial $F(x) = \sum\limits_{i=0}^{n}a_i x^i$, $a_i\in\mathbb{R}$, is called 

\begin{itemize}
 \item
unimodal if 
$a_0\leq \ldots \leq a_{k-1}\leq a_k\geq a_{k+1}\geq \ldots \geq a_n$
for some $k\in\{0,1,\ldots,n\}$,

 \item
logarithmically concave (log-concave) if $a_i^2\geq a_{i-1}a_{i+1}$ for all $i=1,\ldots,n-1$,

 \item
symmetric if $a_i = a_{n-i}$ for all $i=0,1,\ldots,n$.
\end{itemize}

{\bf Theorem 13.1} \cite{Bona}.
a) If a polynomial $F(x)$ with positive coefficients 
has only real roots, then $F(x)$ is log-concave.

b) Any log-concave polynomial is unimodal.

Thus, the condition of log-concativity of polynomial $C(G,x)$ 
or especially condition of unimodality are weaker
than the real-rootedness.
For example, for the graph $G = K_{7,7,7}\cup \overline{K_{43}}$
the clique polynomial $C(G,x) = 1 + 64x + 147x^2 + 343x^3$ \cite{Levit2005}
is unimodal but not log-concave.

By Theorem~13.1 and Corollary~12.3, the matching-generating polynomial $M(G,x)$
is log-concave and unimodal~\cite{Heilman72}.
If $\bar{G}$ is a claw-free graph, then polynomial $C(G,x)$ 
is log-concave and unimodal~\cite{Hamidoune} by Theorem~12.1.

In 1987, Y. Alavi, J. Malde, A. Schwenk, P. Erd\H{o}s proved that

{\bf Theorem 13.2} \cite{Schwenk}.
For every permutation $\pi\in S_{n}$ there exists a graph~$G$
with $\omega(G) = n$ such that
$c_{\pi(1)}<c_{\pi(2)}<\ldots <c_{\pi(n)}$.

{\bf Conjecture 13.1} \cite{Schwenk}.
Let $\bar{G}$ be a tree. Then $C(G,x)$ is unimodal.

In 2013, A. Bhattacharyya and J. Kahn~\cite{Bhattacharyya} found the example 
of the bipartite graph $G$ such that $C(G,x)$ is no unimodal.

A graph $G$ is said to be {\it well-covered} if every maximal independent set of $G$ 
is also a maximum independent set.
In 2003, T. Michael and W. Traves stated that 

{\bf Theorem 13.3} \cite{Michael}.
Let $\bar{G}$ be a well-covered graph and $\omega(G) = 3$. 
Then $C(G,x)$ is unimodal.

{\bf Theorem 13.4} \cite{Michael}.
Let $\bar{G}$ be a well-covered graph and $w = \omega(G)$. Then 
$c_1(G)\leq c_2(G)\leq \ldots\leq c_{\lceil w/2\rceil}(G)$.

In 2014, J. Cutler and L. Pebody proved an analogue of Theorem~13.2 
for well-covered graphs.

{\bf Theorem 13.5} \cite{Cutler}.
For every permutation $\pi$ of the set 
$\{\lceil w/2\rceil,\lceil w/2\rceil+1,\ldots,w\}$
there exists a well-covered graph~$G$ such that
$c_{\pi(\lceil w/2\rceil)}(G)
 <c_{\pi(\lceil w/2\rceil+1)}(G)
 {<}\ldots <c_{\pi(w)}(G)$
and $w = \omega(G)$.

{\bf Corollary 13.1} \cite{Cutler,Levit2006}.
Let $k\geq4$. There exists a well-covered 
graph $\bar{G}$ with $\omega(G) = k$ such that $C(G,x)$ is not unimodal.

A well-covered graph~$G$ is called a {\it very well-covered graph}~\cite{Favaron},
if $G$ contains no isolated vertices and $|V(G)| = 2\alpha(G)$.

In 2006, V. Levit and E. Mandrescu stated the following theorem.

{\bf Theorem 13.6} \cite{Levit2006very}.
Let $\bar{G}$ be a very well-covered graph,
$|V(G)|\geq2$, $w = \omega(G)$. Then 
a) $c_1(G)\leq c_2(G)\leq \ldots\leq c_{\lceil w/2\rceil}(G)$ and
$c_{\lceil (2w-1)/3\rceil}(G)\geq \ldots\geq c_{w-1}(G)\geq c_w(G)$;

b) $C(G,x)$ is unimodal for $w\leq 9$ and log-concave for $w\leq 5$.

Given graphs $G$ and $H$, a {\it corona} $G\circ H$~\cite{corona} 
of them is defined as the graph obtained by taking $|V(G)|$ copies of $H$
and for each $i$ inserting edges between the $i$th vertex of $G$ 
and each vertex of the $i$th copy of $H$.

In 2016, V. Levit and E. Mandrescu proved that 

{\bf Theorem 13.7} \cite{Levit2016}.
Let $H = K_r\setminus e$, $r\geq2$.
Then the polynomial $C(\overline{G\circ H},x)$ 
is unimodal and symmetric for any graph~$G$.

{\bf Corollary 13.2} \cite{Stevanovic98,Mandrescu12}.
The polynomial $C(\overline{G\circ 2K_1},x)$ 
is unimodal and symmetric for any graph~$G$.

Let us note the works of B.-X. Zhu~\cite{Zhu13}, B.-X. Zhu, Q. Lu~\cite{Zhu17}, 
D. Galvin~\cite{Galvin2011}, D. Galvin, J. Hilyard~\cite{Galvin2018}
and P. Bahls, B. Ethridge, L. Szabo~\cite{Balhs2018}.

\subsection{Recognizability by $PC$-polynomial}

If two graphs $G$ and $H$ are isomorphic, then $PC(G,x) = PC(H,x)$.
A class of graphs $\mathcal{K}$ is called $PC$-recognizable
if for any two graphs $G,H\in\mathcal{K}$
the equality $PC(G,x) = PC(H,x)$ implies $G\cong H$.
In 1994, C. Hoede and X. Li formulated~\cite{Clique1994}
the question: Which classes of graphs are $PC$-recognizable?

It is easy to show~\cite{Dohmen} that 
classes of trees and complements to trees are not $PC$-recognizable.
However, D. Stevanovic in 1997 proved that 

{\bf Theorem 13.8} \cite{Stevanovic97}.
The class of threshold graphs is $PC$-recognizable.

{\bf Remark 13.1}. 
Maybe, the second case of the last paragraph of the proof of Corol\-lary~9.1 
could be derived from Theorem~13.8.

In 2017, J.A. Makowsky and V. Rakita stated that for almost all graphs 
there exists a non-isomorphic graph with the same $PC$-polynomial~\cite{MakowskyRakita}.
Later J.A. Makowsky and R.X. Zhang extended this result for 
independence polynomial of hypergraphs.

A {\it spider}~\cite{Spiderman} is a tree having at most one vertex of degree greater than~3.
In 2008, V.~Levit and E. Mandrescu stated 

{\bf Theorem 13.9} \cite{Levit2008}.
Let $\overline{G\circ K_1}$ be connected,
$PC(\overline{G\circ K_1},x) = PC(\bar{T},x)$ and 
$\bar{T}$ is well-covered spider.
Then $G\circ K_1$ is isomorphic to $T$.

{\bf Conjecture 13.2} \cite{Levit2008}.
Let $\bar{G}$ be a connected graph, $\bar{T}$ be a well-covered tree. If 
$PC(\bar{G},x) = PC(\bar{T},x)$, then $\bar{G}$ is a well-covered tree.

About another results on the problem see in the survey~\cite{Levit2005}.

\subsection{Bounds on roots of $PC$-polynomial}

Given a graph~$G$, by the Enestr\"{o}m---Kakeya theorem ($\S0$),
modulus of any root of $PC(G,x)$ is not less than 
$\frac{(w-1)!}{n^{w-1}}$, where $n = |V(G|$ and $w = \omega(G)$.
J. Brown and R.~Nowakowski in 2001 stated the asymptotically best lower bound.

{\bf Theorem 13.10} \cite{BrownNowak01}.
For any graph $G$ with $n$ vertices and clique number 
$w\geq2$, modulus of any root of $PC(G,x)$ is not less than 
$\big(\frac{w-1}{n}\big)^{w-1} + O(n^{-w})$.
This bound is tight. 

{\bf Corollary 13.3}.
For almost all graphs, moduli of all roots of $PC$-polynomial are greater than the smallest
root $\alpha(G_{n,1/2},x) = \beta_{n}(G_{n,1/2})$ of $PC(G_{n,1/2},x)$.

{\sc Proof}.
In \cite{BolobasErdos}, it was proven that with probability tending to~1,
the clique number of a $n$-vertex graph equals 
$\omega(G) = 2\log_2(n) + O(\log_2\log_2(n))$.
Hence, modulus of any root of $PC$-polynomial of almost all graphs
is not less than $n^{-2(1+o(1))\log_2(n)}$.
By Corollary~5.5, we have $\alpha(G_{n,1/2},x) < 2^{1-n}$.
Calculating logarithms of both expressions, we prove the statement.

In 2000, J. Brown, K. Dilcher and R. Nowakowski stated that 

{\bf Theorem 13.11} \cite{BrownNowak00}.
Given a well-covered graph~$G$, modulus of any root of $PC(G,x)$ 
is not less than $1/\omega(G)$.

In 2018, J. Brown and B. Cameron studied stability of independence polynomial.

{\bf Theorem 13.12} \cite{BrownStable}.
a) All roots of $PC(K_n\cup \ldots \cup K_n,x)$ lie in the right half-plane.

b) Not all roots of $PC(K_1\cup K_2\cup \ldots \cup K_n,x)$ lie in the right half-plane for $n\geq15$.

Let us formulate some problems close to considered ones in~$\S9$.

{\bf Problem 13.2}.
To find the average value $\beta(n,k,w)$ over all graphs with $n$ vertices,
$k$ edges and clique number not greater than~$w$.

{\bf Problem 13.3}.
To find the average value of the maximal root $t(G)$ 
of matching polynomial over all graphs with $n$ vertices.

To confirm Conjecture~10.1, it is enough to solve the following problem.

{\bf Problem 13.4}.
To find the series expansion for 
$\lim\limits_{n\to\infty}\frac{\beta_k(n,p)}{n}$
for any $k$.

{\bf Problem 13.5}.
Let $r\gg 1$. Is it true that for almost all graphs with $n$ vertices, 
the roots of $PC$-polynomial which moduli is not less than $n/r$
lie in neighbourhoods of the roots of $PC(G_{n,1/2},x)$?

It is interesting to mention the figures~1 and~2 from the recent work~\cite{BrownStable}
of J. Brown and B. Cameron.

We can suggest the following approach to solve Problem~13.4.
At first, we should clarify the maximum natural number~$k$ such that 
$1/r<\lim\limits_{n\to\infty}\frac{\beta_k(n,1/2)}{n}$. 
Further, we need some intermediate number $q$ between 
$\frac{\beta_{k+1}(n,1/2)}{n}$ and $\frac{\beta_{k}(n,1/2)}{n}$.
Finally, we can prove Problem~13.5 in affirma\-tive way 
applying the Rouch\'{e}'s theorem in the circle of the radius $1/q$ 
for appropriate polynomials $D_t(x)$ and $H(x)$ (see the proof of Statement~11.1).

Is it true that $\beta_r(n,p_1)>\beta_r(n,p_2)$ provided that $1\geq p_1>p_2\geq0$?

\subsection{Adjoint polynomial}

In 1987, R.-Y. Liu defined~\cite{Adjont} for a graph $G$ with $n$ vertices 
the {\it adjoint polynomial} $h(G,x)=\sum\limits_{k=1}^{n}(-1)^{n-k}a_k(G)x^k$,
where $a_k(G)$ denotes the number of ways one can cover 
all vertices of~$G$ by exactly $k$ disjoint cliques of $G$.

The chromatic polynomial of a graph~$G$ is connected with adjoint polynomial of~$G$
in the following was follows
$ch(\bar{G}, x) = \sum\limits_{k=1}^n a_k(G)x(x-1)\ldots (x-k+1)$.
It implies that a graph $G$ is recognizable by its chromatic polynomial
if and only if $G$ is recognizable by its adjoint polynomial.
Thie remark was used by R.-Y. Liu and L.-C. Zhao in 1997~\cite{Adjont97}
to construct new families of graphs recognizable by chromatic polynomial.

In 2017, F. Bencs stated the following result. 

{\bf Theorem 13.13} \cite{BencsAdjont}.
For any graph $G$ there exists a graph $\hat{G}$ such that
$h(G,x) = x^n I(\hat{G},1/x) = x^{n-w}PC(\overline{\hat{G}},x)$,
where $n = |V(G)|$ and $w = \omega(G)$.

Let $V (G) = \{u_1,\ldots,u_n\}$. 
The graph $\hat{G}$ satisfying the conclusion of Theorem~12.12, 
could be constructed as follows. Put
$V(\hat{G}) = \{(u_i,u_j)\in E(G)\mid 1\leq i<j\leq n\}$.
Let $(u_i,u_j)$ and $(u_k,u_l)$ be two distinct vertices of $\hat{G}$.
Assume that $j\geq l$. Define that they are connected in $\hat{G}$ if 
$k\in\{i,j\}$ or $j=l$ and $(u_i,u_k)\not\in E(G)$.
It is clear that $\hat{G}$ is a spanning subgraph of the lie graph of $G$.

{\bf Corollary 13.4} \cite{BencsAdjont,CsikvariAdjont,Zhao}.
For any graph~$G$, there exists a positive real root $\gamma(G)$ 
of~$h(G,x)$ such that $\gamma(G)$ is larger than moduli of all other roots of $h(G,x)$.

For $\gamma(G)$, we can prove the analogues of Lemmas~2.3--2.5 
in the similar way~\cite{BencsAdjont,CsikvariAdjont,Zhao}.
Since $\hat{G}$ is a spanning subgraph of~$L(G)$, from Lemma~2.4 we deduce

{\bf Corollary 13.5} \cite{BencsAdjont,CsikvariAdjont}.
For any graph~$G$, $\gamma(G)\leq t^2(G)$, where $t(G)$ 
is the largest root of the matching polynomial of~$G$.

By Corollary~13.5 and Statement~12.2, we can derive bounds on $\gamma(G)$.

Although $\hat{G}$ is a spanning subgraph of $L(G)$
and by Theorem~11.1 all roots of $PC(L(G),x)$ are real,
we may not conclude that the adjoint polynomial is real-rooted. 
For triangle-free graphs and comparability graphs, all roots 
of adjoint polynomial are indeed real~\cite{Brenti,Brenti94}.

\subsection{The value of independence polynomial in $x = -1$}

The sum 
$I(G,-1) = 1 - s_1(G)+s_2(G)-\ldots +(-1)^{\alpha(G)}s_{\alpha(G)}(G)$ 
was called in~\cite{Bousquet} as the {\it alternating number of independent sets}.

Given a graph~$G$, the {\it decycling number} $\varphi(G)$~\cite{Beineke}
is the minimum number of vertices that need to be removed in order 
to eliminate all cycles in~$G$.

The following result was initially proved by A. Engstr\"{o}m in 2009,
V.~Levit and E.~Mandrescu found the elementary proof of it in 2011.

{\bf Theorem 13.13} \cite{Engstroem,Levit2011}.
For any graph~$G$, $|I(G,-1)|\leq 2^{\varphi(G)}$.

{\sc Proof}.
Let us prove the statement by induction on $\varphi(G)$.

If $\varphi(G) = 0$, then $G$ is a forest.
Let us state the inequality $|I(G,-1)|\leq 1$ in this case by induction on $n = |V(G)|$.
For $n = 1$, $I(G,x) = 1+x$ and thus $|I(G,-1)| = 0 \leq 1$.
Let $G$ be a forest with $n\geq2$. 
If $G$ has no leaves, then $I(G,x) = (1+x)^n$ and $I(G,-1) = 0$.
Otherwise, let $v$ be a leaf of~$G$, $(v,u)\in E(G)$.
By~\eqref{IndepMainFormula}, we have
$$
I(G,x) = I(G \setminus u,x) + xI(G\setminus [N(u)],x)
 = (1+x)I(G \setminus \{u,v\},x) + xI(G\setminus [N(u)],x),
$$
which implies 
$|I(G,-1)| = |I(G\setminus [N(u)],-1)|\leq1$
by the induction hypothesis.

Suppose that the statement is proven for all graphs satisfying 
$\varphi(G)\leq k$. Let $G$ be a graph with $\varphi(G) = k+1$.
It is clear that there exists a vertex $v\in V(G)$ such that 
$\varphi(G\setminus v)<\varphi(G)$.
So, the statement follows from~\eqref{IndepMainFormula}
applying for $G$ and $v$ and the induction hypothesis.

J. Cutler and N. Kahl in 2016 proved that 

{\bf Theorem 13.14} \cite{CutlerKahl}.
Given a positive integer $k$ and an integer $q$ with $|q| \leq 2^k$,
there is a connected graph~$G$ with $\varphi(G) = k$ and $I(G,-1) = q$.

\subsection{Extremal values of $\beta_{\pm}(n,k)$}

Let is list some works devoted to the problems close to
the one about the values $\beta_{\pm}(n,k)$.
In \cite{Optimiz}, unicyclic $n$-vertex graphs were studied
which independence polynomials have minimal coefficients.
The analogue of the problem of $\beta_+(n,k)$
for the spectral radius of a graph was solved by P.~Rowlinson in 1988~\cite{Rowlinson}.
See the monograph~\cite{Stevanovic} of D.~Steva\-novic
about the extremal values of the spectral radius of a graph 
with fixed number of vertices and fixed graph invariants 
(including clique number, independence number, matching number etc). 
In 2014, J. Cutler and A.J. Radcliffe formulated a problem to find 
a graph with minimum number of cliques among graphs from $G(n,k)$.
They solved it for $\bar{k}\leq n$. For this, they used some 
graph transformations which are similar to the transformations from 
$\S8$ and $\S9$.

The problem of finding the values $\beta_{\pm}(n,k)$ 
could be formulated for graphs which are embedded into different surfaces
(not necessary in plane as it was studied in $\S6$).
For the graph~$G$ embeddable in torus
with maximal number of cliques~\cite{Wood2011}, we have 
\begin{multline*}
PC(G,x) 
 = x^7-nx^6+3nx^5-(3n+14)x^4+(n+28)x^3-21x^2+7x-1 \\
 = (x-1)^3(x^4-(n-3)x^3+6x^2-4x+1),
\end{multline*}
since $c_8 = 0$ for any graph embeddable in torus. 

For the graph~$H$ embeddable in the projective plane
with maximal number of cli\-ques~\cite{Wood2011}, we have 
\begin{multline*}
PC(H,x) 
 = x^6-nx^5+3(n-1)x^4-(3n+2)x^3+(n+9)x^2-6x+1 \\
 = (x-1)^3(x^3-(n-3)x^2+3x-1),
\end{multline*}
since $c_7 = 0$ for any graph embeddable in the projective plane.  
We suppose that the graphs $G$ and $H$ are analogues of the graphs 
$G_+$ (see Pic.~3). So, to construct the graph for the value 
$\beta_+(n,k)$ in the class of all graph embeddable into
torus (the projective plane),
we should start with the clique $K_7$ ($K_6$)
and then repeat so called splitting of triangles \cite{Wood2011}.

Let us formulate the following problem. 

{\sc Clique game}. 
There are two players: Min and Max.
Before the start, the numbers of vertices ($n$) and edges ($k$) 
of the future graph are fixed, and we have initially empty graph. 
Each player takes turn adding one edge. 
The goal of Min (Max) to min(max)imize $\beta(G)$ of the final graph $G\in G(n,k)$.
Fix that Min makes the first turn.

{\bf Problem 13.6}.
a) What strategies should the players follow?

b) What is the final value $\beta(G)\in G(n,k)$ 
if both players follow the best strategies?

The clique game could be considered in some class of graphs, for example
for planar graphs. Then we have an additional rule: the graph should be planar after each turn. 

Finally, let us list the problems which were considered in this paper but not solved.

{\bf Problem 13.7}.
To prove Conjecture~9.1 to confirm the results about $\beta_-(n,k)$ for $k>n^2/4$.

By Corollary~9.2 we know that a graph $G\in G(n,k)$ with minimal~$\beta$
has to be connected graph of diameter~2.
It seems that some kind of Zykov's symmetrization could be helpful to solve this Problem.
Let us introduce the Zykov's EP (edge-preserving)-transformation as follows.
Let $u,v$ be two disconnected vertices of a graph~$G$.
Suppose that $|N_G(u)\setminus N_G(v)|\geq |N_G(v)\setminus N_G(u)|$,
then the graph~$G'$ is obtained from $G$ by removing 
edges $(v,w)$ for $w\in |N_G(v)\setminus N_G(u)|$ and adding the edges 
$(u,w)$, $w\in |N_G(v)\setminus N_G(u)|$. 
By some conditions, the Zykov's EP-transformation does not increase $\beta(G)$. 

Another possible strategy is to prove that if $G\in G(n,k)$ with $\beta(G) = \beta_-(n,k)$,
then there exists a graph $H\in G(n,k+1)$ such that $\beta(H) = \beta_-(n,k+1)$
and $G = H\setminus e$ for some edge $e$. 
Hence, Conjecture~9.1 could be derived from Statement~7.2.

The third idea is following. Let us introduce the graph invariants for a graph $G$
\begin{gather*}
t_1 = 1,\quad t_2 = 2n,\quad t_3 = 3n^2-2k,\\
t_4 = 4n^3-6nk+c_3,\quad t_5 = 5n^4-12n^2k+3k^2+3nc_3-c_4,\quad \ldots,
\end{gather*}
where $n = |V(G)|$, $k = |E(G)|$, $c_k = c_k(G)$.

They appeared as the coefficients for the numbers $m_i(G)$ 
expressed by the linear recurrence~\eqref{rekurrent} via $n,k,c_3,c_4,\ldots$:
\begin{gather*}
m_3 = nm_2 - km_1 + c_3 = (n^3-2nk) + t_1 c_3,\\
m_4 = nm_3 - km_2 + c_3m_1 - c_4 = (n^4-3kn^2 + k^2) + t_2 c_3 - t_1c_4,\quad \ldots, \\
m_7 = (n^7-6kn^5+10k^2n^3-4k^3n) + t_5 c_3 - t_4 c_4 + t_3 c_5 - t_2 c_6 + t_1 c_7,\quad \ldots, \\
m_s = f_s(n,k) + \sum\limits_{i=1}^{s-2}t_{s-i-1}c_{i+2}
\end{gather*}
for some homogeneous polynomials $f_s(n,k)$ on two variables $n,k$ of the degree~$s$ 
(we prescribe to $k$ a degree~2). Something close was done in~\cite{Beezer} for matching polynomial. 

Suppose that we transform $G$ to a new graph $G'\in G(n,k)$ 
(by a Kelmans or another transformation), then
we want to clarify if $m_i(G)\geq m_i(G')$ holds for all $i$.
For this, it can be useful to know the properties of the coefficients $t_i$.
 
{\bf Problem 13.7}.
To find the exact upper bounds for the expressions 
$\beta(G)+\beta(\bar{G})$ and 
$\beta(G)\beta(\bar{G})$.

We conjecture that the upper bounds obtained in Example~8.2 are asymptotically the best ones.

\section*{Acknowledgments}
I would like to express my sincere gratitude to Mikhail Novikov,
in 2016 a pupil of Novosibirsk's lyceum~9 
(now a student of Saint-Petersburg State University),
with whom we began to study the growth rate of partially commutative Lie algebras.  
The work would not have been appeared without this initial common interest.
The author is very grateful to Elena Konstantinova and Peter Csikv\'{a}ri
for the help of various kinds. 

The author acknowledges support by the Austrian Science Foun\-da\-tion FWF, 
grant P28079.

\noindent Vsevolod Gubarev \\
University of Vienna \\
Oskar-Morgenstern-Platz 1, 1090, Vienna, Austria \\
Sobolev Institute of mathematics \\
Acad. Koptyug ave., 4, 630090 Novosibirsk, Russia \\
e-mail: vsevolod.gubarev@univie.ac.at

\end{document}